\documentclass[moor,nonblindrev]{informs1} 

\usepackage{natbib}
 \NatBibNumeric
 \def\bibfont{\small}%
 \def\BIBand{and}%
 \bibpunct[, ]{[}{]}{,}{n}{}{,}%

\usepackage[colorlinks=true,breaklinks=true,bookmarks=true,urlcolor=blue,
     citecolor=blue,linkcolor=blue,bookmarksopen=false,draft=false]{hyperref}

\def\EMAIL#1{\href{mailto:#1}{#1}}  

\TheoremsNumberedThrough     

\EquationsNumberedThrough    

\usepackage{booktabs} 
\usepackage{multicol, multirow} 
\usepackage{float} 
\usepackage{caption} 
\usepackage{amsmath,amssymb} 
\usepackage{url,hyperref}
\usepackage[usenames,dvipsnames]{xcolor}
\usepackage{mathtools}
\usepackage{multirow}
\usepackage{graphicx}
\usepackage{color}
\usepackage{xfrac}
\usepackage{url,hyperref}
\usepackage{algorithm}
\usepackage{pifont}
\usepackage{booktabs,paralist}
\usepackage{algorithmicx,algpseudocode,algcompatible}
\usepackage{wrapfig}
\usepackage{makecell}

\newcommand{\Id}{\mathbb{I}}
\newcommand{\R}{\mathbb{R}}

\newcommand{\set}[1]{\left\{#1\right\}}
\newcommand{\sets}[1]{\{#1\}}
\newcommand{\norm}[1]{\left\Vert#1\right\Vert}
\newcommand{\norms}[1]{\Vert#1\Vert}

\newcommand{\tnorms}[1]{\vert\!\Vert#1\vert\!\Vert}
\newcommand{\Eproof}{\hfill $\square$}
\newcommand{\prox}{\mathrm{prox}}

\newcommand{\dom}[1]{\mathrm{dom}(#1)}
\newcommand{\gra}[1]{\mathrm{gra}(#1)}

\newcommand{\iprods}[1]{\langle #1\rangle}
\newcommand{\Exp}[1]{\mathbb{E}\left[#1\right]}
\newcommand{\Exps}[2]{\mathbb{E}_{#1}\left[#2\right]}
\newcommand{\Expsn}[2]{\mathbb{E}_{#1}\big[#2\big]}
\newcommand{\Expn}[1]{\mathbb{E}\big[#1\big]}

\newcommand{\Bc}{\mathcal{B}}
\newcommand{\Xc}{\mathcal{X}}

\newcommand{\Zc}{\mathcal{Z}}
\newcommand{\Sc}{\mathcal{S}}

\newcommand{\Dc}{\mathcal{D}}
\newcommand{\Lc}{\mathcal{L}}

\newcommand{\Tc}{\mathcal{T}}
\newcommand{\Rc}{\mathcal{R}}

\newcommand{\Fc}{\mathcal{F}}
\newcommand{\Uc}{\mathcal{U}}

\newcommand{\Nc}{\mathcal{N}}
\newcommand{\Wc}{\mathcal{W}}

\newcommand{\Pc}{\mathcal{P}}

\newcommand{\BigO}[1]{\mathcal{O}\left(#1\right)}
\newcommand{\BigOs}[1]{\mathcal{O}\big(#1\big)}
\newcommand{\SmallO}[1]{o\left(#1\right)}
\newcommand{\SmallOs}[1]{o\big(#1\big)}

\newcommand{\zer}[1]{\mathrm{zer}(#1)}

\newcommand{\myeq}[2]{\vspace{-0.5ex}
\begin{equation}\label{#1}
{#2}
\vspace{-0.5ex}
\end{equation}
}
\newcommand{\myeqn}[1]{\vspace{-0.3ex}
\begin{equation*}
{#1}
\vspace{-0.5ex}
\end{equation*}
}

\newcommand{\mytb}[1]{\textbf{#1}}

\newcommand{\mbf}[1]{\mathbf{#1}}
\newcommand{\mcal}[1]{\mathcal{#1}}
\newcommand{\mred}[1]{\textcolor{black}{#1}}
\newcommand{\mblue}[1]{\textcolor{blue}{#1}}

\newcommand{\revise}[1]{#1}

\newcommand{\beforesec}{\vspace{-0ex}}
\newcommand{\aftersec}{\vspace{-1ex}}
\newcommand{\beforesubsec}{\vspace{-1ex}}
\newcommand{\aftersubsec}{\vspace{-1ex}}
\newcommand{\beforesubsubsec}{\vspace{-1ex}}
\newcommand{\aftersubsubsec}{\vspace{-1ex}}

\begin{document}

\RUNAUTHOR{Q. Tran-Dinh \textit{and} N. Nguyen-Trung}

\RUNTITLE{
Accelerated Past-Extragradient Methods with Variance-Reduction
}

\TITLE{
New Accelerated Past-Extragradient Methods with Variance Reduction  for Generalized Equations
}

\ARTICLEAUTHORS{%
\AUTHOR{Quoc Tran-Dinh \textit{and} Nghia Nguyen-Trung}
\AFF{Department of Statistics and Operations Research \\ The University of North Carolina at Chapel Hill (UNC), 318 Hanes Hall, Chapel Hill, NC 27599-3260, USA.\\ 
\EMAIL{quoctd@email.unc.edu; nghiant@unc.edu}}
} 

\ABSTRACT{
We develop a novel past-extragradient-type algorithmic framework, combining both Nesterov's \textit{acceleration} and \textit{variance-reduction} techniques, to solve a class of generalized equations involving possibly \textit{nonmonotone operators} in data-driven applications. 
Our framework covers a wide class of stochastic variance-reduced schemes, including mini-batching and both unbiased and biased control-variate estimators.
We establish that our method achieves $\BigOs{1/k^2}$ convergence rates in expectation for the squared norm of the residual under Lipschitz continuity and  a ``co-hypomonotonicity-type'' assumption,  significantly  improving upon non-accelerated counterparts by a factor of $1/k$. 
We also prove faster $\SmallOs{1/k^2}$ convergence rates, both in expectation and almost surely. 
In addition, we show that the sequence of iterates generated by our method almost surely converges  to a solution of the underlying problem.
We demonstrate the applicability of our method using general error approximation criteria, covering mini-batch stochastic estimators as well as three well-known control variate estimators: Loopless SVRG, SAGA, and Loopless SARAH.
The resulting three variants attain significantly better oracle complexities than existing methods.
We validate our framework and theoretical results through three numerical examples. 
The numerical results illustrate promising performance of our accelerated method over its non-accelerated counterparts.
}

\MSCCLASS{90C25, 90C06, 90-08}
\ORMSCLASS{Mathematical programming; operations research; game theory}

\maketitle

\beforesec
\section{Introduction.}\label{sec:intro}
\aftersec
The goal of this paper is to develop a novel \textit{\textbf{accelerated past-extragradient-type framework} with both \textbf{unbiased and biased variance reduction} techniques} to solve the following possibly nonmonotone \textbf{generalized equation}:
\myeq{eq:GE}{
\textrm{Find $x^{\star}\in \R^p$ such that:~ $0 \in \Phi{x}^{\star}$ with $\Phi{x} := G{x} + Tx$},
\tag{GE}
}
where $G : \R^p \to \R^p$ is a single-valued mapping and $T : \R^p \rightrightarrows \R^p$ is a possibly multivalued mapping.
Throughout this paper, we assume that $\zer{\Phi} :=  \sets{x^{\star} \in \R^p : 0 \in \Phi{x}^{\star} } \neq\emptyset$, i.e., \eqref{eq:GE} has a solution.

Following  \cite{Rockafellar2004}, we refer to \eqref{eq:GE} as a composite \textit{generalized equation}, but it is also known as a composite \textit{[non]linear inclusion}, see \cite{Bauschke2011}.
In particular, if $T = 0$, then \eqref{eq:GE} reduces to the following [non]linear equation (also known as a \textit{root-finding problem}):
\myeq{eq:NE}{
\textrm{Find $x^{\star}\in \R^p$ such that:~ $Gx^{\star}  = 0$}.
\tag{NE}
}
This equation is equivalent to the following fixed-point problem: 
\myeq{eq:fixed_point_prob}{
\textrm{Find $x^{\star} \in \R^p$ such that $x^{\star} = Fx^{\star}$, where $Fx := x - \lambda Gx$ for any $\lambda > 0$}.
\tag{FP}
}
However, under appropriate conditions (e.g., maximal [co-hypo]monotonicity of $T$), we can reformulate \eqref{eq:GE} equivalently as \eqref{eq:NE} using, e.g., a forward-backward splitting, backward-forward splitting, or Tseng's forward-backward-forward splitting operator, see, e.g., \cite{Bauschke2011,TranDinh2025a}. 
Under such conditions, we have the following relations: \eqref{eq:GE} $\Leftrightarrow$ \eqref{eq:NE} $\Leftrightarrow$ \eqref{eq:fixed_point_prob}.

Problem  \eqref{eq:GE} looks simple, but it is well-known that \eqref{eq:GE} covers the optimality conditions of  unconstrained and constrained minimization, minimax, and [mixed] variational inequality problems.
It also comprises Nash equilibrium problems  in economics \cite{Facchinei2003}.

\vspace{0.75ex}
\noindent\textbf{$\mathrm{1.1.}$~Our problem classes and assumptions.}
Before specifying our assumptions, we highlight that we will focus on two classes of \eqref{eq:GE} as follows:
\begin{compactenum}
\item\textbf{Class 1.} \revise{$G$} is \textit{$L$-Lipschitz continuous} and \revise{$G + T$} is \textit{$\rho$-co-hypomonotone}, i.e., \revise{$G$} satisfies Assumption~\ref{as:A1} with $\alpha = 1$ and \revise{$G+T$} satisfies Assumption~\ref{as:A3} below with $\rho_n \geq 0$ and $\rho_c = 0$.
\item\textbf{Class 2.} \revise{$G$} satisfies \textit{Assumption~\ref{as:A1} with any $\alpha \in [0, 1)$} and \revise{$G+T$} satisfies both \textit{Assumption~\ref{as:A2} and Assumption~\ref{as:A3}  below with $\rho_n \geq 0$ and $\rho_c > 0$}.
\end{compactenum}
In both classes, we assume that $G$ is not directly accessible and that we have access only to an unbiased stochastic oracle $\mbf{G}$ of $G$.
More specifically, we impose the following assumptions.

\begin{assumption}[Unbiased oracle]{\!\!\!\!\!}\label{as:A0}
$G$ is equipped with an unbiased stochastic oracle $\mbf{G}$ as
\myeq{eq:unbiased_setting}{
Gx = \Expsn{\xi}{ \mbf{G}(x, \xi) }, \quad   \forall x \in \R^p,
}
where $\xi$ is a random variable on a given probability space $(\Omega, \Sigma,\mathbb{P})$ and $\mathbb{E}_{\xi}[ \mbf{G}(\cdot, \xi) ] = \int_{\Omega}\mbf{G}(\cdot, \xi)\mathbb{P}(d\xi)$. 
Moreover, there exists $\sigma \geq 0$ such that $\Expsn{\xi}{\norms{\mbf{G}(x, \xi) - Gx}^2 } \leq \sigma^2$.
\end{assumption}
Assumption~\ref{as:A0} is widely used in the literature, especially in data-driven optimization applications, see, e.g., \cite{lan2020first,Nemirovski2009}.
In particular, if we choose $\Omega = [n] := \sets{1, \cdots, n}$, $\mathbb{P}(\xi = i) = \frac{1}{n}$ for all $i \in \Omega$, and $\mbf{G}(x, \xi) := G_{\xi}x$ with $G_i : \R^p \to \R^p$ for $i \in [n]$ are single-valued mappings, then $G$ becomes
\myeq{eq:finite_sum}{
Gx = \Exps{\xi}{\mbf{G}(x, \xi)} = \sum_{i=1}^n\mathbb{P}(\xi = i)\mbf{G}(x, i) = \frac{1}{n}\sum_{i=1}^nG_ix.
}
This is called a finite-sum representation.
It can also be obtained by sample average approximation techniques.
In many applications, $n$ is very large so that evaluating $G$ is expensive or even infeasible.

\begin{assumption}[\textbf{Lipschitz continuity-type}]\label{as:A1}
There exist two constants $L \in (0, +\infty)$ and $\alpha \in [0, 1]$ such that for all $x, y \in \R^p$, we have
\myeq{eq:G_Lipschitz}{
\arraycolsep=0.2em
\begin{array}{lcl}
\alpha\norms{Gx - Gy}^2 + (1-\alpha)  \Expsn{\xi}{\norms{\mbf{G}(x, \xi) - \mbf{G}(y, \xi) }^2 }   \leq L^2\norms{x - y}^2.
\end{array}
}
\end{assumption}
Note that we allow $\alpha = 0$ and $\alpha = 1$ in Assumption~\ref{as:A1}.
If $\alpha=1$, then \eqref{eq:G_Lipschitz} reduces to $\norms{Gx - Gy} \leq L\norms{x - y}$, showing that $G$ is $L$-\revise{Lipschitz} continuous.
If $\alpha= 0$, then \eqref{eq:G_Lipschitz} becomes $\Expsn{\xi}{\norms{\mbf{G}(x, \xi)  - \mbf{G}(y, \xi) }^2} \leq L^2\norms{x - y}^2$, showing that $G$ is $L$-Lipschitz continuous in expectation.
While the former is standard in classical stochastic approximation methods, the latter condition is often used in variance-reduction algorithms with control variate techniques, especially in optimization, see, e.g.,  \cite{driggs2019accelerating,nguyen2017sarah,Pham2019,Tran-Dinh2019a}.
In fact, \eqref{eq:G_Lipschitz} is a convex combination of these two extreme cases.
Note that $\norms{Gx - Gy}^2 = \norms{ \Expsn{\xi}{\mbf{G}(x, \xi)  - \mbf{G}(y, \xi)} }^2 \leq \Expsn{\xi}{\norms{\mbf{G}(x, \xi) - \mbf{G}(y, \xi) }^2 }$ by Jensen's inequality, \eqref{eq:G_Lipschitz} implies $\norms{Gx - Gy} \leq L\norms{x - y}$ for all $x, y \in \R^p$, i.e., $G$ is $L$-Lipschitz continuous.

In the finite-sum setting \eqref{eq:finite_sum}, our condition \eqref{eq:G_Lipschitz} reduces to 
\myeq{eq:G_Lipschitz_finite_sum}{
\arraycolsep=0.2em
\begin{array}{lcl}
\alpha\norms{Gx - Gy}^2 + \frac{1-\alpha}{n} \sum_{i=1}^n \norms{G_ix - G_iy}^2   \leq L^2\norms{x - y}^2.
\end{array}
}
If $\alpha = 0$, then $G$ is said to be $L$-average Lipschitz continuous.

\begin{assumption}[\textbf{Weak Minty solution}]\label{as:A2}
There exists a solution $x^{\star} \in \zer{\Phi}$ of \eqref{eq:GE} and a constant $\rho_{*} \geq 0$ such that 
\myeq{eq:weak_Minty}{
\arraycolsep=0.2em
\begin{array}{lcl}
\iprods{w, x - x^{\star}} \geq -\rho_{*}\norms{w}^2, \quad \forall(x, w) \in \gra{\Phi}.
\end{array}
}
\end{assumption}
This assumption also covers a class of nonmonotone operators $\Phi$ and has been widely used in the literature in recent years, see, e.g., \cite{diakonikolas2021efficient,tran2024revisiting} for examples.
It is weaker than the so-called ``star-monotonicity'' of $\Phi$ (i.e., $\rho_{*} = 0$) used in some recent works, e.g., \cite{gorbunov2022extragradient,kotsalis2022simple}.

\begin{assumption}[\textbf{Co-hypomonotonicity-type}]\label{as:A3}
There exist $\rho_n \geq 0$ and $\rho_c \geq 0$ such that $\rho_n \geq \rho_c$ and  for all $(x, w_x), (y, w_y) \in\gra{\Phi}$, we have 
\myeq{eq:G_cohypomonotonicity}{
\arraycolsep=0.2em
\begin{array}{lcl}
\iprods{w_x - w_y, x - y} \geq - \rho_n\norms{w_x - w_y}^2 + \rho_c \Expsn{\xi}{  \norms{\mbf{G}(x, \xi) - \mbf{G}(y, \xi) }^2 }.
\end{array}
}
\end{assumption}
The condition \eqref{eq:G_cohypomonotonicity} in Assumption~\ref{as:A3} looks technical and nonstandard.
It deserves further discussion.
Note that since we allow $\rho_c = 0$, \eqref{eq:G_cohypomonotonicity} can reduce to  
\myeqn{
\arraycolsep=0.2em
\begin{array}{lcl}
\iprods{w_x - w_y, x - y} \geq - \rho_n\norms{w_x - w_y}^2.
\end{array}
}
This relation shows that $\Phi$ is \textbf{$\rho_n$-co-hypomonotone}.
Clearly, Assumption~\ref{as:A3} includes ordinary co-hypomonotonicity as a special case and, when $T = 0$, also encompasses average co-coercive operators.
To develop variance-reduction methods via control-variate techniques in Section~\ref{sec:VFOG}, we will require $\rho_c > 0$.
In the finite-sum setting \eqref{eq:finite_sum}, the condition \eqref{eq:G_cohypomonotonicity} reduces to $\iprods{w_x - w_y, x - y} \geq - \rho_n\norms{w_x - w_y}^2 + \frac{\rho_c}{ n} \sum_{i=1}^n  \norms{G_ix - G_iy}^2$.

In particular, let us discuss the case when $T = 0$ and $\rho_c > 0$.
Then, \eqref{eq:G_cohypomonotonicity} reduces to
\myeq{eq:G_cohypomonotonicity2}{
\arraycolsep=0.2em
\begin{array}{lcl}
\iprods{Gx - Gy, x - y} \geq - \rho_n\norms{Gx - Gy}^2 + \rho_c \Expsn{\xi}{  \norms{\mbf{G}(x, \xi) - \mbf{G}(y, \xi) }^2 }, \quad \forall x, y \in \dom{G}.
\end{array}
}
We have the following facts corresponding to $G$ in this case.
\begin{compactitem}
\item First, if $\mbf{G}$ is deterministic, i.e., $Gx \equiv \mbf{G}(x, \xi)$ in \eqref{eq:unbiased_setting} or $n=1$ in \eqref{eq:finite_sum}, then \eqref{eq:G_cohypomonotonicity2} reduces to $\iprods{Gx - Gy, x - y} \geq -\rho\norms{Gx - Gy}^2$ for $\rho := \rho_n - \rho_c$.
Clearly, if $\rho_n < \rho_c$, then $G$ becomes $|\rho|$-co-coercive, a case excluded by Assumption~\ref{as:A3} since we require $\rho_n \ge \rho_c$.
If  $\rho_c = \rho_n$, then $G$ is monotone.
If  $\rho_n > \rho_c$, then $G$ is $\rho$-co-hypomonotone (see Subsection~\ref{subsec:background} for the definition of these concepts).

\item Second, since $\Expsn{\xi}{\norms{\mbf{G}(x, \xi) - \mbf{G}(y, \xi) }^2} \geq \norms{\Expsn{\xi}{\mbf{G}(x, \xi) - \mbf{G}(y, \xi) }}^2$ by Jensen's inequality, \eqref{eq:G_cohypomonotonicity2} implies that $\iprods{Gx - Gy, x - y} \geq -\rho\norms{Gx - Gy}^2$ for $\rho := \rho_n - \rho_c$.
Therefore, if $\rho_n < \rho_c$, then \eqref{eq:G_cohypomonotonicity2} is slightly stronger than the $\vert\rho\vert$-co-coercivity of $G$.
If $\rho_n > \rho_c$, then \eqref{eq:G_cohypomonotonicity2} is slightly stronger than the $\rho$-co-hypomonotonicity of $G$.

\item Third, since $\rho_n \geq 0$, if $\rho_c > 0$, then \eqref{eq:G_cohypomonotonicity2} is considerably weaker than the $\rho_c$-average co-coercivity of $G$, i.e., $\iprods{Gx - Gy, x - y} \geq \rho_c \Expsn{\xi}{ \norms{ \mbf{G}(x, \xi) - \mbf{G}(y, \xi)}^2 }$. 
This $\rho_c$\mred{-average}-co-coercivity assumption has been widely used in many recent works to develop variance reduction methods for \eqref{eq:GE} and its special cases, see, e.g., \cite{cai2023variance,davis2022variance,tran2024accelerated,TranDinh2025a}.
\end{compactitem}

\begin{wrapfigure}{r}{0.38\textwidth}
\vspace{-4ex}
\centering
\includegraphics[width=0.38\textwidth]{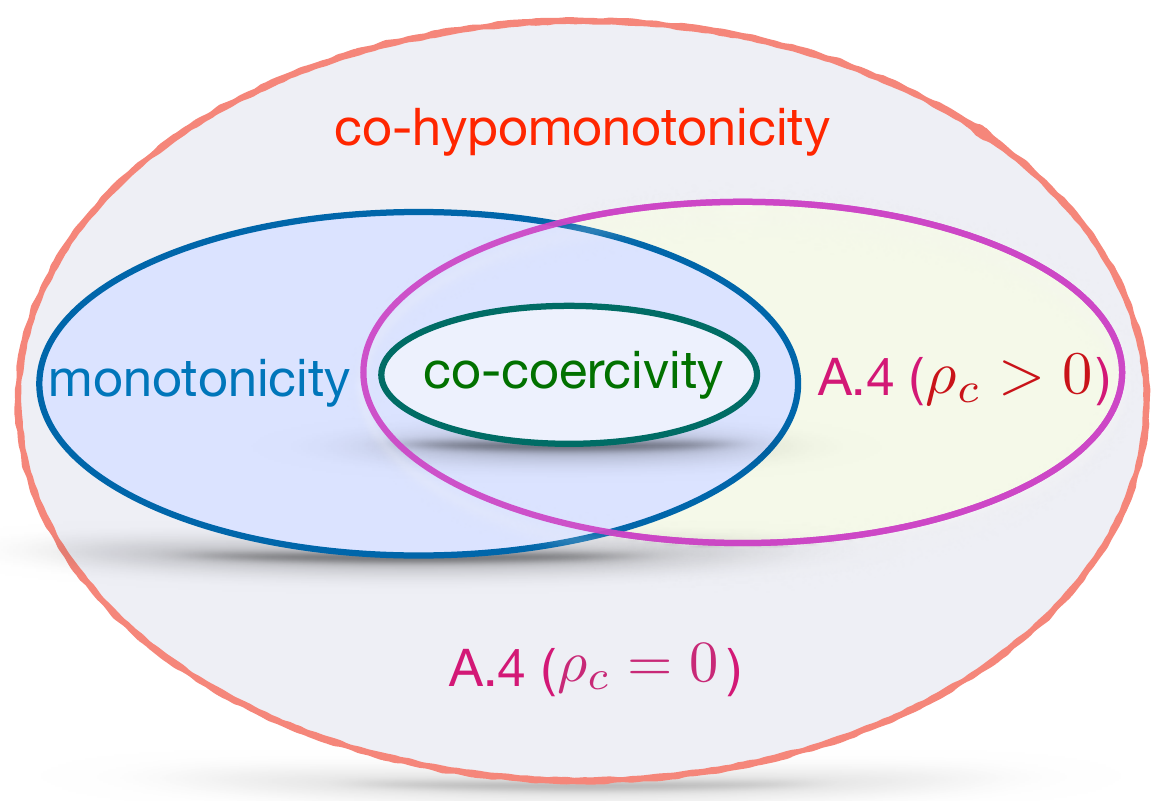}
\vspace{-4ex}
\caption{\footnotesize
Different classes of operators $G$.
The class ``\textbf{A.4. $(\rho_c > 0$)}'' represents the case $G$ satisfying \eqref{eq:G_cohypomonotonicity2} with $\rho_c > 0$.
The class ``\textbf{A.4. $(\rho_c = 0$)}'' covers all \textbf{co-hypomonotone operators}.
}
\label{fig:nonmono_ops}
\vspace{-3ex}
\end{wrapfigure}

In this paper, \textbf{we do not focus on the case $\rho_n < \rho_c$ in \eqref{eq:G_cohypomonotonicity2}} because $G$ is $(\rho_c-\rho_n)$-co-coercive, which has been studied in previous works, including our papers \cite{tran2024accelerated,TranDinh2025a}.
\revise{\textbf{We only consider the case $\rho_n \geq \rho_c \geq 0$}, which covers a class of monotone and co-hypomonotone operators.
Figure \ref{fig:nonmono_ops} illustrates the class of $\rho$-co-hypomonotone operators $G$ as well as the class of $G$ satisfying condition \eqref{eq:G_cohypomonotonicity} when $\rho_c > 0$.
The latter intersects with both monotone and nonmonotone operators.}
Subsections~\ref{subsec:co-hypomonotone_example} and \ref{subsec:nonmonotone_example} provide many examples of $\Phi$ and $G$ satisfying \eqref{eq:G_cohypomonotonicity} and \eqref{eq:G_cohypomonotonicity2}.

Note that if $T=0$, then Assumption~\ref{as:A3} implies the $\rho$-co-hypomonotonicity of $G$ with $\rho := \rho_n-\rho_c$. 
Thus, Assumption~\ref{as:A2} automatically holds with $\rho_{*} = \rho$. 
Alternatively, if $\rho_c = 0$, then Assumption~\ref{as:A3} also implies the $\rho_n$-co-hypomonotonicity of $\Phi$, and hence, Assumption~\ref{as:A2} also automatically holds with $\rho_{*} = \rho_n$. 
In these two cases, Assumption~\ref{as:A2} is redundant and can be removed.

\vspace{0.75ex}
\noindent\textbf{$\mathrm{1.2.}$~Motivation and challenges.}
Both minimax optimization and variational inequality (VIP) problems  are fundamental and powerful tools for modeling two-person games, Nash equilibria, robust optimization \cite{Bauschke2011,reginaset2008,Facchinei2003,phelps2009convex,Rockafellar2004,Rockafellar1976b,ryu2016primer}, and more recently, generative adversarial learning, online learning, reinforcement learning, and distributionally robust optimization \cite{arjovsky2017wasserstein,azar2017minimax,Ben-Tal2009,bhatia2020online,goodfellow2014generative,jabbar2021survey,levy2020large,lin2022distributionally,madry2018towards,rahimian2019distributionally,wei2021last}.
Their optimality conditions can be expressed as a generalized equation \eqref{eq:GE}. 
Yet, unlike classical models, modern applications are data-driven, large-scale, and often involve nonconvex, nonmonotone, and/or nonsmooth structures, limiting the use of deterministic and  primal-dual approaches. 
Mathematically, these problems are typically formulated as large finite-sums or as stochastic expectations induced by sampled data. 
\revise{
Consequently, contemporary methods face significant challenges, including restrictive assumptions, high computational operation cost, low efficiency, and slow convergence rates.
}

Hitherto, most methods developed for \eqref{eq:GE} and its special cases are deterministic and primarily address monotone problems, see, e.g., \cite{Bauschke2011,Facchinei2003,Konnov2001}. 
They are inefficient for the finite-sum setting \eqref{eq:finite_sum} when $n$ is large, and inapplicable to the expectation setting \eqref{eq:unbiased_setting}, where $G$ is not directly accessible. 
Randomized and stochastic methods exist for \eqref{eq:GE}, but under restrictive assumptions \cite{boct2021minibatch,combettes2015stochastic,cui2021analysis,kotsalis2022simple,pethick2023solving}. 
Many of these still rely on the classical Robbins-Monro approach without acceleration, resulting in slow convergence rates and oracle complexities often ranging from $\BigOs{\epsilon^{-6}}$ to $\BigOs{\epsilon^{-4}}$, depending on the method and assumptions, where $\epsilon$ is a desired accuracy level.

Recently, several works have applied variance-reduction techniques using control variates, commonly developed in optimization, to \eqref{eq:GE}, improving oracle complexity over the classical stochastic approximation. 
However, these algorithms are non-accelerated with slow convergence rates, typically $\BigOs{1/k}$, see, e.g., \cite{alacaoglu2021stochastic,alacaoglu2021forward,davis2022variance}. 
\revise{
While \cite{alacaoglu2021stochastic,alacaoglu2021forward} rely on Lipschitz continuity of $G$ and monotonicity of $\Phi$, others such as \cite{cai2023variance,davis2022variance,tran2024accelerated,TranDinh2025a} assume stronger conditions: co-coercivity and/or strong quasi-monotonicity. 
Co-coercivity is strong and hard to verify in practice; it is stronger than monotonicity plus Lipschitz continuity and fails for a bilinear mapping $Fx = [\mbf{L}v, -\mbf{L}^Tu]$, where $x = (u, v)$ and $\mbf{L}$ is given. 
Moreover, checking co-coercivity of the saddle mapping in minimax optimization applications is non-trivial, even in the convex-concave case.
}

\vspace{0.75ex}
\noindent\textbf{$\mathrm{1.3.}$~Our goal and approach.}
The above shortcomings raise the following research question:
\textit{Can we develop a \textbf{single-loop algorithmic framework} capable of achieving \textbf{faster $\BigOs{1/k^2}$} and $\SmallOs{1/k^2}$ convergence rates via \textbf{acceleration}, covering a wide range of variants $($using both unbiased and biased estimators$)$, attaining \textbf{improved oracle complexity} through \textbf{variance reduction}, and requiring \textbf{weaker assumptions} $($including a subclass of \textbf{nonmonotone} problems$)$?}

This paper aims to address this question.
However, we face two main technical challenges. 
\revise{
First, classical Krasnosel'ki\v{i}-Mann (KM) iteration is not directly applicable under our assumptions, requiring a shift to extragradient-type (EG) schemes. 
Nevertheless, EG involves the composite mapping $G(x - \eta Gx)$, which complicates the design and analysis of stochastic algorithms. 
We instead exploit Popov's \textbf{past-extragradient method} \cite{popov1980modification} (also closely related to the optimistic gradient scheme) to overcome this limitation. 
Second, developing variance-reduction methods with control variate techniques requires a structural assumption on $G$, often called ``average co-coercivity’’ \cite{cai2023variance,davis2022variance,tran2024accelerated,TranDinh2025a}, which excludes certain problem classes, particularly nonmonotone ones. 
We identify a new class of problems that includes both monotone and nonmonotone operators. 
This class typically satisfies Assumptions \ref{as:A1}, \ref{as:A2}, and \ref{as:A3}.
}

Unlike previous approaches, our algorithmic development simultaneously combines three main ideas: the past-extragradient method, Nesterov's acceleration, and unbiased and biased variance-reduction techniques. 
While each has been studied extensively, their unified combination to solve \eqref{eq:GE} under Assumptions~\ref{as:A1}-\ref{as:A3} remains largely unexplored, especially for nonmonotone problems.

\vspace{0.75ex}
\noindent\textbf{$\mathrm{1.4.}$~Our contributions.}
To this end, our contributions consist of the following:
\begin{compactitem}
\item[$\mathrm{(a)}$] 
First, we introduce Assumption~\ref{as:A3} to identify a new class of problems in~\eqref{eq:GE} that goes beyond standard monotonicity and co-coercivity assumptions.
Concrete examples in Subsection~\ref{subsec:nonmonotone_example} illustrate the scope and relevance of this assumption, showing that it covers subclasses of both monotone and nonmonotone operators.

\item[$\mathrm{(b)}$] 
Second, we develop a novel stochastic Nesterov-type accelerated past-extragradient framework to approximate solutions of \eqref{eq:GE}. 
Our framework only requires one stochastic approximation $\widetilde{G}y^k$ of $Gy^k$ per iteration, evaluated at the iterate $y^k$.
We introduce a general \textbf{\textit{error approximation condition}} connecting $Gy^k$ to its (unbiased or biased) stochastic variance-reduced estimator $\widetilde{G}y^k$, constructed via the stochastic oracle $\mathbf{G}(\cdot,\xi)$, which forms the basis of our convergence analysis.
Alternatively, we also propose a class of control-variate-based estimators of $G$ for our framework, covering both unbiased and biased instances, including SVRG \cite{SVRG}, SAGA \cite{Defazio2014}, and SARAH \cite{nguyen2017sarah}.

\item[$\mathrm{(c)}$]
Third, we establish $\BigOs{1/k^2}$ convergence rates in expectation for $\Expn{\norms{Gx^k + v^k}^2}$ using either the proposed error approximation criterion or our class of variance-reduced estimators, where $x^k$ is the iterate and $v^k \in Tx^k$. 
Our rate is significantly better than non-accelerated counterparts, e.g., \cite{alacaoglu2021stochastic,alacaoglu2021forward,TranDinh2024}, by a factor of $1/k$.
We further prove $\SmallOs{1/k^2}$ convergence of $\norms{Gx^k + v^k}^2$ both in expectation and almost surely, and show that the iterate sequences $\sets{x^k}$, $\sets{y^k}$, and $\sets{z^k}$ converge almost surely to a $\zer{\Phi}$-valued random variable $x^{\star}$, a solution of \eqref{eq:GE}.

\item[$\mathrm{(d)}$] 
Finally, we demonstrate our framework with three common estimators: Loopless SVRG, SAGA, and Loopless SARAH, and analyze their oracle complexities for the finite-sum case \eqref{eq:finite_sum}. 
For Loopless SVRG and SAGA, we obtain $\BigOs{n\log(n) + n^{2/3}\epsilon^{-1}}$, while for Loopless SARAH we achieve $\BigOs{n\log(n) + n^{1/2}\epsilon^{-1}}$, where $\epsilon > 0$ is the target accuracy.
\end{compactitem}

\noindent
We believe that our framework is simple to implement, as it is single-loop and uses elementary operations.
Moreover, it includes a broad class of variance-reduced and accelerated algorithms, covering mini-batch, unbiased, and biased estimators, a feature that is rarely seen in optimization.
Unlike most existing methods limited to the monotone setting or special cases of \eqref{eq:GE}, our approach can handle a class of nonmonotone generalized equations. 
Its convergence guarantees are notable, covering $\BigOs{1/k^2}$ and $\SmallOs{1/k^2}$ convergence rates, and almost sure convergence of iterates.
The $\SmallOs{1/k^2}$ rates and almost-sure convergence of the iterates have received limited attention, even in variance-reduction methods for convex optimization.

\vspace{0.75ex}
\noindent\textbf{$\mathrm{1.5.}$~Related work and comparison.}
Problem \eqref{eq:GE} is classical, with its theory and numerical methods extensively studied over decades, as presented in monographs such as \cite{Bauschke2011,Facchinei2003,minty1962monotone,Rockafellar2004} and references therein.
In this paper, we focus on a different class of \eqref{eq:GE} and develop novel stochastic single-loop methods with acceleration and variance reduction.
Accordingly, we restrict our review and comparison to works along this direction.

\begin{table}[ht!]
\vspace{-0ex}
\newcommand{\cell}[1]{{\!\!}{#1}{\!\!}}
\begin{center}
\caption{Comparison of existing variance-reduction single-loop methods and our algorithms}
\label{tbl:existing_vs_ours}
\vspace{0.25ex}
\begin{small}
\resizebox{\textwidth}{!}{  
\begin{tabular}{|c|c|c|c|c|c|} \hline
\textbf{Paper} & \cell{\textbf{Assumptions}} & \textbf{Estimators} & {\!\!\!} \textbf{Metric} {\!\!\!} &  \cell{\textbf{Convergence rates}} & \textbf{Complexity} \\ \hline
\cite{davis2022variance}                          & $^{*}$co-coercive and SQM     &  {\!\!\!\!} SVRG \& SAGA {\!\!\!\!} & Residual & linear rate & $\BigOs{(L/\mu)\log(\epsilon^{-1})}$  \\ \hline
\cite{TranDinh2024}                                 & weak-Minty+Lipschitz     & \cell{{\!\!\!}unbiased class }{\!\!}  &  Residual  & $\BigOs{1/k}$ & $\mcal{O}\big( n + n^{2/3}\epsilon^{-2} \big)$ \\ \hline
\cite{alacaoglu2021stochastic}                & monotone+Lipschitz       & SVRG & Gap & $\BigOs{1/k}$  & $\mcal{O}\big( n + n^{1/2}\epsilon^{-2} \big)$  \\ \hline
\cite{alacaoglu2021forward}                    & monotone+Lipschitz       & SVRG & Gap & $\BigOs{1/k}$  & $\mcal{O}\big( n\epsilon^{-2} \big)$  \\ \hline
\cite{yu2022fast}                                      & monotone+Lipschitz       & SVRG & Gap & $\BigOs{1/k}$  & $\mcal{O}\big( n\epsilon^{-2} \big)$  \\ \hline
\cite{tran2024accelerated}                      & co-coercive     & \cell{unbiased class} &  Residual & $\BigOs{1/k^2}$  & $\mcal{O}( n + n^{2/3}\epsilon^{-1} )$  \\ \hline
\cite{cai2023variance}                            & co-coercive     & SARAH  &  Residual  & $\BigOs{1/k^2}$ & $\widetilde{\mcal{O}}\big( n + n^{1/2}\epsilon^{-1} \big)$ \\ \hline
\cite{TranDinh2025a}                              & co-coercive       &   {\!\!\!} \makecell{unbiased and \\ biased class} &  Residual &{\!\!\!\!\!}  \makecell{ $\BigOs{1/k^2}$, $\SmallOs{1/k^2}$ \\ $x^k \to x^{\star}$ a.s.  } {\!\!\!\!\!} &{\!\!\!}  $\widetilde{\mcal{O}}\big(n + (n^{2/3} \vee n^{1/2})\epsilon^{-1}\big)$ {\!\!\!}  \\ \hline
\mblue{\textbf{Ours}}                               & \makecell{{\!\!\!}Assumptions \ref{as:A0}, \ref{as:A1}, and \ref{as:A3} {\!\!\!}\\ with $\alpha = 1$ and $\rho_c = 0$} & {\!\!\!} \makecell{\mblue{unbiased and} \\ \mblue{biased class}} {\!\!\!} & Residual &  \makecell{ \mblue{ $\BigOs{1/k^2}$, $\SmallOs{1/k^2}$} \\ \mblue{$x^k \to x^{\star}$ a.s. } } {\!\!\!\!\!} &  \makecell{$\widetilde{\mcal{O}}\big(\epsilon^{-4}\big)$} {\!\!\!} \\ \hline
\mblue{\textbf{Ours}}                               &{\!\!\!\!} \makecell{Assumptions \ref{as:A0} to \ref{as:A3} \\ with $\alpha \in [0, 1)$ and $\rho_c > 0$} {\!\!\!} & {\!\!\!} \makecell{ \mblue{unbiased and} \\ {\!\!\!\!}  \mblue{biased class} {\!\!\!\!} } & Residual &  \makecell{ \mblue{ $\BigOs{1/k^2}$, $\SmallOs{1/k^2}$} \\ \mblue{$x^k \to x^{\star}$ a.s. } } {\!\!\!\!\!} &  \makecell{$\widetilde{\mcal{O}}\big(n + n^{2/3}\epsilon^{-1}\big)$ \\ $\to \widetilde{\mcal{O}}\big(n + n^{1/2}\epsilon^{-1}\big)$ } {\!\!\!} \\ \hline
\end{tabular}}
\end{small}
\end{center}
{\footnotesize
\vspace{1.5ex}
\textbf{Notes:} 
\textbf{SQM} = strong quasi-monotonicity; 
\textbf{Residual} is $\Expn{\norms{Gx^k + v^k}^2}$ for $v^k \in Tx^k$; 
\textbf{Gap} is the metric $\mathcal{G}(x) := \max_{y \in \Bc}\iprods{Gy, x - y}$;
and \textbf{unbiased [and biased] class} = a class of unbiased [and biased] variance-reduced estimators.
}
\vspace{-4ex}
\end{table}

\noindent\mytb{$\mathrm{(a)}$~From non-accelerated to accelerated methods.}
The key advantage of accelerated methods over non-accelerated ones is their faster convergence rates: $\BigOs{1/k^2}$ vs. $\BigOs{1/k}$.
Nesterov first introduced acceleration for smooth convex optimization \cite{Nesterov1983}, later extended to other problems and monotone inclusions, including \eqref{eq:GE} \cite{attouch2020convergence,attouch2019convergence}. This line of research has recently gained strong attention, especially in deterministic settings \cite{bot2022fast,bot2022bfast,kim2021accelerated,mainge2021accelerated,tran2021halpern,tran2025accelerated}.

As an alternative, Halpern’s fixed-point iteration \cite{halpern1967fixed} also achieves accelerated $\BigOs{1/k^2}$ rates on the squared norm of residual for \eqref{eq:GE} under co-coercivity \cite{diakonikolas2020halpern,lieder2021convergence,sabach2017first}, improving on the classical $\BigOs{1/k}$ rate of the Krasnosel'ki\v{i}-Mann (KM) method and its variants.

These advances have driven extensive work in operator theory, mostly in deterministic settings \cite{cai2022accelerated,cai2022baccelerated,lee2021fast,tran2021halpern,yoon2021accelerated}. Though Nesterov- and Halpern-type accelerations were developed independently, their connection was later established in \cite{tran2022connection} and also in \cite{partkryu2022}.
In this paper, we build on Nesterov’s acceleration to design stochastic variance-reduction algorithms for \eqref{eq:GE} under Assumptions~\ref{as:A1}–\ref{as:A3} and establish their rigorous convergence guarantees.

\vspace{0.5ex}
\noindent\mytb{$\mathrm{(b)}$~From convex optimization to generalized equations.}
In convex optimization, the $L$-smoothness of the objective function is equivalent to its gradient co-coercivity \cite{Bauschke2011,Nesterov2004}, making smooth convex minimization problems special cases of \eqref{eq:NE} with $G$ being co-coercive. 
By contrast, smooth convex-concave minimax problems often lack this property (see, e.g., bilinear game models \cite{Facchinei2003}), highlighting a key distinction between optimization and generalized equations \eqref{eq:GE}. This shift complicates algorithm design and convergence analysis, as many standard analysis tools in convex optimization no longer apply to \eqref{eq:GE} \cite{attouch2019convergence,kim2021accelerated,tran2022connection}.
One technical challenge is the absence of an objective function in \eqref{eq:GE}, complicating the construction of a suitable Lyapunov function that underpins convergence guarantees in accelerated methods. 
In addition, unlike optimization, developing accelerated randomized and variance-reduced methods for data-driven \eqref{eq:GE} remains challenging but crucial.

Beyond co-coercivity, extragradient (EG) and its variants \cite{Korpelevic1976,popov1980modification,daskalakis2018training} are prominent, but their stochastic extensions face difficulties: dependence on intermediate updates \cite{yu2022fast} and the breakdown of gap function-based analysis under nonmonotonicity \cite{alacaoglu2021stochastic,alacaoglu2021forward}. 
Dual averaging and mirror descent \cite{Nemirovskii2004,kotsalis2022simple} are similarly restricted to monotone problems. 
In this work, we extend Popov’s past-extragradient variant to the stochastic setting, addressing cases with biased estimators.

\vspace{0.5ex}
\noindent\mytb{$\mathrm{(c)}$~Comparison with existing stochastic methods.}
Stochastic methods for \eqref{eq:GE} often follow the Robbins-Monro stochastic approximation framework, but converge slowly, typically at an $\BigOs{1/k}$ rate  with inefficient oracle complexity, see, e.g.,  \cite{boob2023first,censor2012extensions,cui2021analysis,jiang2008stochastic,juditsky2011solving,kannan2019optimal,koshal2012regularized}. 
Recent work improves this limitation via variance reduction techniques, e.g., adaptive or growing mini-batches \cite{boct2021minibatch,pethick2023solving}.
Alternative approaches such as dual averaging and mirror descent \cite{chen2017accelerated,kotsalis2022simple} extend to VIPs but require monotonicity, limiting applicability. 
Other advances include stochastic splitting and mirror-prox methods \cite{beznosikov2023stochastic,gorbunov2022stochastic,loizou2021stochastic}.
In contrast, our work employs accelerated methods with a broad class of unbiased and biased variance-reduced estimators, including several well-known instances.

\vspace{0.5ex}
\noindent\mytb{$\mathrm{(d)}$~Comparison with variance-reduction methods using control variate techniques.}
In the co-coercive setting, several variance-reduced methods using control variates have been proposed \cite{cai2023variance,davis2022variance,tran2024accelerated,TranDinh2025a}, exploiting co-coercivity for faster rates and improved oracle complexity. Table~\ref{tbl:existing_vs_ours} summarizes recent works by assumptions, estimators, and guarantees.
Beyond co-coercivity, EG variants, often based on Popov’s scheme \cite{popov1980modification}, have been studied for monotone VIPs. 
Examples include SVRG-based methods in \cite{alacaoglu2021stochastic,alacaoglu2021forward,yu2022fast}, which use gap functions for convergence analysis and rely on bounded-variance or unbiased estimators.

Our approach retains a simple single-loop structure, similar to the methods in \cite{alacaoglu2021stochastic,alacaoglu2021forward,cai2023variance,davis2022variance}, but it further (i) incorporates acceleration, (ii) applies to ``co-hypomonotone-type'' operators, (iii) accommodates both unbiased and biased estimators under general error-bound conditions, (iv) establishes stronger convergence guarantees, including $\SmallOs{1/k^2}$ rates in expectation and almost surely, as well as almost sure convergence of the iterates to a solution, and (v) achieves improved oracle complexity under weaker assumptions.
Unlike \cite{alacaoglu2021stochastic,yu2022fast}, our method also provides residual convergence guarantees with better oracle complexity in the finite-sum setting.

\vspace{0.75ex}
\noindent\textbf{$\mathrm{1.6.}$~Paper organization.}
The remainder of this paper is organized as follows.
Section~\ref{sec:background} reviews the necessary background and concepts, and presents examples of $\Phi$ and $G$ satisfying Assumption~\ref{as:A3}.
Section~\ref{sec:SFOG} develops a unified accelerated past-extragradient-type framework with variance reduction to solve \eqref{eq:GE} for the expectation setting  \eqref{eq:unbiased_setting} under Assumption~\ref{as:A0}, Assumption~\ref{as:A1} with $\alpha = 1$, and Assumption~\ref{as:A3} with $\rho_c = 0$, establishing convergence under a general error bound condition.
Section~\ref{sec:VFOG} introduces an alternative variant using control-variate-based variance reduction to handle both \eqref{eq:unbiased_setting} and \eqref{eq:finite_sum}, with similar convergence guarantees.
Section~\ref{sec:convergence_analysis} provides a complete convergence analysis of the proposed methods.
Section~\ref{sec:numerical_experiments} presents three numerical examples to validate our theoretical results and benchmark against state-of-the-art methods.
Preliminary results and technical proofs are deferred to the appendix.

\beforesec
\section{Background, Preliminaries, and Examples.}\label{sec:background}
\aftersec
We recall necessary notation and concepts used in this paper, and  briefly describe Popov's method.
Finally, we provide different examples to illustrate Assumption~\ref{as:A3}.

\beforesubsec
\subsection{\textbf{Notation and concepts.}}
\label{subsec:background}
\aftersubsec
First, we recall necessary notation and concepts.

\noindent\textbf{Basic concepts.}
We work with a finite-dimensional space $\R^p$ equipped with the standard inner product $\iprods{\cdot, \cdot}$ and norm $\norms{\cdot}$.
For a function $f : \R^p\to\R$, $\nabla{f}$ denotes the gradient of $f$, $\partial{f}$ denotes the subdifferential of $f$, and $\prox_{f}$ is the proximal operator of $f$.
For a given random vector $\mbf{u}_{\xi}$ in $\R^p$ and $v \in \R^p$, we define $u := \mathbb{E}_{\xi}[\mbf{u}_{\xi}]$ and $\tnorms{u + v}^2 := \mathbb{E}_{\xi}[ \norms{ \mbf{u}_{\xi} + v}^2 ]$.
By Jensen's inequality, we also have $\tnorms{u + v}^2 \geq \norms{u + v}^2$.
For a given symmetric matrix $\mbf{Q}\in\R^{p\times p}$, $\mbf{Q} \succeq 0$ (resp., $\mbf{Q} \succ 0$) means that $\mbf{Q}$ is positive semidefinite (resp., positive definite).
For two symmetric matrices $\mbf{P}$ and $\mbf{Q}$, $\mbf{P} \succeq \mbf{Q}$ means $\mbf{P} - \mbf{Q} \succeq 0$.
For given functions $g(t)$ and $h(t)$, we say that $g(t) = \BigOs{h(t)}$ if there exist $M > 0$ and $t_0 \geq 0$ such that $g(t) \leq Mh(t)$ for $t \geq t_0$.
We use $\widetilde{\mcal{O}}(g(t))$ to hide a poly-log factor of $g(t)$.

\vspace{0.5ex}
\noindent\textbf{Filtration and expectations.}
Let $\Fc_k$ denote the $\sigma$-algebra generated by all the randomness of our algorithm up to iteration $k$, including the iterates $x^0, \cdots, x^k$.
We also denote by $\Expsn{k}{\cdot} := \Expn{\cdot \mid \Fc_k}$ the conditional expectation conditioned on $\Fc_k$, and by $\Expn{\cdot}$ the full expectation.
For a random variable $\xi$, we denote by $\Expsn{\xi}{\cdot}$ the expectation w.r.t. $\xi$.

\vspace{0.5ex}
\noindent\textbf{Lipschitz continuity.}
For a single-valued mapping $G : \R^p\to\R^p$, we say that $G$ is $L$-Lipschitz continuous if there exists $L \geq 0$ such that $\norms{Gx - Gy} \leq L\norms{x - y}$ for all $x, y \in \dom{G}$.
If $L < 1$, then $G$ is contractive, and if $L=1$, then $G$ is nonexpansive.

\vspace{0.75ex}
\noindent\textbf{Co-monotonicity, monotonicity, and co-coercivity.}
For a single-valued or multivalued mapping $G : \R^p \rightrightarrows \R^p$, we denote by $\dom{G} := \sets{x \in \R^p : Gx \neq\emptyset}$ the domain of $G$ and by $\gra{G} := \sets{(x, u) \in \R^p\times\R^p : u \in Gx}$ the graph of $G$.
If $G$ is differentiable at $x$, then $G'(x)$ denotes its derivative at $x$.
We say that $G$ has a closed graph $\gra{G}$ if for any sequence $\sets{(x^k, u^k)}$ such that $(x^k, u^k) \in \gra{G}$ and $(x^k, u^k) \to (x, u)$ as $k \to \infty$, then $(x, u) \in \gra{G}$.
For further details on these concepts and properties, we refer the reader to \cite{Bauschke2011}.
\begin{definition}\label{de:monotonicity}
Let $G : \R^p \to 2^{\R^p}$ be a given single-valued or multivalued mapping.
\begin{compactitem}
\item We say that $G$ is monotone if for all $(x, u), (y, v) \in \gra{G}$, we have $\iprods{u - v, x - y} \geq 0$.
\item We say that $G$ is $\rho$-co-monotone if there exists $\rho \in \R$ such that $\iprods{u - v, x - y} \geq \rho\norms{u - v}^2$ for all  $(x, u), (y, v) \in \gra{G}$.
Clearly, if $\rho = 0$, then $G$ reduces to a monotone mapping.
See \cite{bauschke2020generalized} for the definitions and further properties.
\begin{compactitem}[$\diamond$]
	\item If $\iprods{u - v, x - y} \geq \rho\norms{u - v}^2$ for $\rho < 0$, then we say that $G$ is $\vert\rho\vert$-co-hypomonotone.
	A $\vert\rho\vert$-co-hypomonotone mapping can be nonmonotone (see examples, e.g., in \cite{TranDinh2025a}).
	\item If $\iprods{u - v, x - y} \geq \rho\norms{u - v}^2$ for $\rho > 0$, then we say that $G$ is $\rho$-co-coercive. 
	Moreover, if $G$ is $\rho$-co-coercive, then it is also monotone and $L$-Lipschitz continuous with $L := \frac{1}{\rho}$.
\end{compactitem}
\end{compactitem}
\end{definition}
If $f$ is a convex and $L$-smooth (i.e., $\nabla{f}$ is $L$-Lipschitz continuous) function, then $\nabla{f}$ is $\frac{1}{L}$-co-coercive.
The $\rho$-co-coercivity of a mapping $G$ is also equivalent to the nonexpansiveness of $F = \Id - 2\rho G$ in the fixed-point problem \eqref{eq:fixed_point_prob}, see \cite{Bauschke2011}.
All the concepts defined above are global.
If $G$ satisfies any of the above properties for all $(x, u), (y, v)\in\gra{G}\cap\Wc$, where $\Wc$ is a neighborhood of $(\bar{x}, \bar{u})$ in $\R^p\times\R^p$, then we say that $G$ has such a local property around $(\bar{x}, \bar{u}) \in \gra{G}$.

\beforesubsec
\subsection{\textbf{The past-extragradient method.}}
\label{subsec:OG_method}
\aftersubsec
A classical method for solving \eqref{eq:GE} in the variational inequality setting $T=\Nc_{\Xc}$ (the normal cone of $\Xc$), where $\Xc$ is nonempty, closed, and convex and $G$ is monotone and Lipschitz continuous, is Korpelevich's extragradient method, introduced by Korpelevich in \cite{Korpelevic1976} in 1976 (also by Antipin in \cite{antipin1976} (1976)).
This method has inspired various variants, among which the past-extragradient method, a notable variant proposed by Popov in \cite{popov1980modification}, is particularly well-known.
When applied to solving \eqref{eq:GE}, this variant is equivalent to the optimistic gradient method studied in \cite{daskalakis2018training,mertikopoulos2019optimistic,mokhtari2020convergence}.

This method starts from $x^0 \in \dom{\Phi}$, chooses $y^{-1} := x^0$, and at each iteration $k\geq 0$, updates 
\begin{equation}\label{eq:popov_method}
\arraycolsep=0.2em
\left\{\begin{array}{lcl}
y^k &:= & J_{\eta T}( x^k - \eta Gy^{k-1}), \vspace{1ex}\\
x^{k+1} &:= & J_{\eta T}(x^k - \eta Gy^k),
\end{array}\right.
\tag{PEG}
\end{equation}
where $\eta > 0$ is a given stepsize.
This scheme has been well studied in the literature, including \cite{malitsky2020forward,popov1980modification,tran2024revisiting}, mostly in the deterministic setting.
It also has close connections to other methods such as the projected reflected gradient scheme \cite{malitsky2015projected}, the forward-reflected-backward splitting algorithm \cite{malitsky2020forward}, and the golden-ratio method \cite{malitsky2019golden}. 
As mentioned, when $T = 0$, it is also equivalent to the optimistic gradient method \cite{daskalakis2018training,mertikopoulos2019optimistic,mokhtari2020convergence}.
Inspired by \eqref{eq:popov_method}, we develop an accelerated algorithmic framework with variance reduction in Section~\ref{sec:SFOG} to solve \eqref{eq:GE} under Assumptions~\ref{as:A0}-\ref{as:A3}.

\beforesubsec
\subsection{\textbf{Examples of co-hypomonotone operators.}}
\label{subsec:co-hypomonotone_example}
\aftersubsec
\revise{
We provide here a number of mappings that are co-hypomonotone, but not necessarily monotone.
}

\revise{
\vspace{0.75ex}
\noindent\textbf{$\mathrm{(a)}$~Example 1.}
As shown in \cite{TranDinh2025a}, given a $p\times p$ symmetric and invertible matrix $\mbf{G}$ and a vector $\mbf{g} \in \R^p$, an affine mapping $Gx := \mbf{G}x + \mbf{g}$ is $\rho$-co-monotone with $\rho := \lambda_{\min}(\mbf{G}^{-1})$.
If $\rho := \lambda_{\min}(\mbf{G}^{-1}) < 0$, then $G$ is nonmonotone, but $\vert\rho\vert$-co-hypomonotone.
More generally, for an arbitrary matrix $\mbf{G} \in\R^{p\times p}$, the affine mapping $Gx := \mbf{G}x + \mbf{g}$ is $\rho$-co-hypomonotone for some $\rho\geq 0$ if and only if $\mbf{G} + \mbf{G}^\top + 2\rho \mbf{G}^{\top} \mbf{G} \succeq 0$.
As shown in \textbf{Example 4} below, the subdifferentials of certain nonconvex regularizers are not linear but piecewise linear, and they are also locally $\rho$-co-hypomonotone.

\revise{
\vspace{0.75ex}
\noindent\textbf{$\mathrm{(b)}$~Example 2.}
As shown in \cite[Proposition~4.17]{evens2025convergence}, consider a twice continuously differentiable function $\Psi:\R^{p_1}\times\R^{p_2}\to\R$ and its associated saddle mapping $\Phi(u, v) := [\nabla_u\Psi(u,v), -\nabla_v\Psi(u,v)]$.
If $\Psi$ satisfies the \textit{$\alpha$-interaction-dominance} condition with $\alpha\geq0$ and parameter $\rho>0$, then $\Phi$ is $\rho$-co-hypomonotone.
The interaction-dominance condition was introduced in \cite{grimmer2023landscape} for studying smooth nonconvex-nonconcave minimax problems.
Recently, \cite{evens2025convergence} generalized co-hypomonotonicity to \textit{semi-monotonicity}, which encompasses both hypomonotonicity and co-hypomonotonicity while covering broader parameter regimes.
}

\revise{
\vspace{0.75ex}
\noindent\textbf{$\mathrm{(c)}$~Example 3.}
Let $A : \R^p\rightrightarrows \R^p$ be a monotone operator and $\rho > 0$ be given.
Then, the mapping $Tx := -\rho^{-1}x - Ax$ is $\rho$-co-hypomonotone and, in general, nonmonotone (see Appendix~\ref{apdx:subsec:co-hypomonotone_example}).
}

\vspace{0.75ex}
\noindent\textbf{$\mathrm{(d)}$~Example 4.}
As shown in \cite{TranDinh2025a}, consider the [generalized] subdifferential $T x := \partial{\mathrm{SCAD}_a}(x)$ of the SCAD (smoothly clipped absolute deviation) regularizer from \cite{fan2001variable} with $a=3.7$, which is widely used in statistics.
The operator $T$ is not globally $\rho$-co-hypomonotone. 
However, on each smooth concave region of $\mathrm{SCAD}_a$, it is $\rho$-co-hypomonotone, where $\rho := a - 1$.
More precisely, $\iprods{u-v,x-y} = -(a-1)\norms{u-v}^2$ for all $(x,u), (y, v) \in \gra{T}$ belonging to the same concave region.
Similar local co-hypomonotonicity properties hold for several other well-known nonconvex regularizers, such as MCP (minimax concave penalty) and the log-sum penalty.
}

\beforesubsec
\subsection{\textbf{Examples of co-hypomonotone operators satisfying Assumption~\ref{as:A3}.}}
\label{subsec:nonmonotone_example}
\aftersubsec
The Lipschitz continuity-type condition \eqref{eq:G_Lipschitz} in Assumption~\ref{as:A1} is common, especially when $\alpha = 0$ or $\alpha = 1$.
Let us provide concrete examples that satisfy our ``co-hypomonotonicity-type'' condition \eqref{eq:G_cohypomonotonicity} in Assumption~\ref{as:A3} and condition  \eqref{eq:G_cohypomonotonicity2} (i.e., $T=0$).

\vspace{0.5ex}
\noindent$\mathrm{(a)}$ \textbf{Example 1}~(\textit{\textbf{Concrete linear mappings}}).
We first consider a concrete linear example.
Consider the following three square matrices $\mbf{G}_1$, $\mbf{G}_2$, and $\hat{\mbf{G}}$ and two vectors $\mbf{g}_1$ and $\mbf{g}_2$:
\begin{equation}\label{eq:exam1}
\mbf{G}_1 = \begin{bmatrix}-2 & 0\\ 0 & 0 \end{bmatrix}, \quad \mbf{G}_2 = \begin{bmatrix}0 & 0\\ 0 & 1 \end{bmatrix}, \quad  \hat{\mbf{G}} = \begin{bmatrix}0 & 0\\ 0 & 1/2 \end{bmatrix}, \quad \mbf{g}_1 = \begin{pmatrix} 0 \\ 1 \end{pmatrix},   \quad \textrm{and} \quad \mbf{g}_2 = \begin{pmatrix} -1 \\ 1 \end{pmatrix}.
\end{equation}
We define three mappings $G_1x := \mbf{G}_1x + \mbf{g}_1$, $G_2x := \mbf{G}_2x + \mbf{g}_2$, and $Tx := \hat{\mbf{G}}x$.
Then, we have $Gx := \frac{1}{2}(G_1x + G_2x)$ and $\Phi{x} := Gx + Tx$.
Since $\mbf{\Phi} := \frac{1}{2}(\mbf{G}_1 + \mbf{G}_2) + \hat{\mbf{G}} = \begin{bmatrix}-1 & 0\\ 0 & 1\end{bmatrix}$ is not positive semidefinite, $\Phi$ is nonmonotone.

Now, for any vector $u \in \R^2$, we can easily check that  for some $\rho_n  \geq \rho_c > 0$, we have 
\begin{equation}\label{eq:exam1_cond1}
\begin{array}{lcl}
u^{\top}\mbf{\Phi}u \geq -\rho_n \norms{\mbf{\Phi}u}^2 + \frac{\rho_c}{2}\big(\norms{\mbf{G}_1u}^2 + \norms{\mbf{G}_2u}^2\big).
\end{array}
\end{equation}
Since $\mbf{\Phi}$ is symmetric, the condition \eqref{eq:exam1_cond1} is equivalent to 
\begin{equation*}
\begin{array}{lcl}
\mbf{\Phi} + \rho_n\mbf{\Phi}^{\top}\mbf{\Phi} - \frac{\rho_c}{2}\big(\mbf{G}_1^{\top}\mbf{G}_1 + \mbf{G}_2^{\top}\mbf{G}_2\big) \succeq 0.
\end{array}
\end{equation*}
By a direct calculation,  there exist $\rho_n := 2\epsilon + 1$ and $\rho_c := \epsilon > 0$ for any $\epsilon > 0$ such that this inequality holds.
Therefore, with these given $\rho_n$ and $\rho_c$, for all $x, y \in \R^2$, we have
\begin{equation*}
\begin{array}{lcl}
\iprods{\Phi{x} - {\Phi}y, x - y} \geq -\rho_n\norms{\Phi{x} - \Phi{y}}^2 + \frac{\rho_c}{2}\big(\norms{G_1x - G_1y}^2 + \norms{G_2x - G_2y}^2\big).
\end{array}
\end{equation*}
Overall, we have shown that the mappings $G$ and $T$ satisfy Assumption~\ref{as:A3}.

\vspace{0.5ex}
\noindent$\mathrm{(b)}$~\textbf{Example 2}~(\textit{\textbf{A class of \eqref{eq:GE} studied in \cite{cai2023variance,tran2024accelerated,TranDinh2025a}}}).
In \cite{cai2023variance,tran2024accelerated,TranDinh2025a}, the authors studied a class of \eqref{eq:GE} in which $G$ is $\frac{1}{L}$-co-coercive in expectation, i.e., for all $x, y \in \dom{G}$, there exists $L > 0$ such that
\begin{equation*}
\begin{array}{l}
	\iprods{Gx - Gy, x - y} = \Expsn{\xi}{\iprods{\mbf{G}(x, \xi) - \mbf{G}(y, \xi), x - y}} \geq \frac{1}{L} \Expsn{\xi}{\norms{\mbf{G}(x, \xi) - \mbf{G}(y, \xi) }^2},
\end{array}
\end{equation*} 
and $T$ is maximally $\rho$-co-hypomonotone \cite{tran2024accelerated,TranDinh2025a} or monotone \cite{cai2023variance}. 

As shown in Appendix~\ref{apdx:subsec:example_verification}, one can check that $\Phi := G + T$ satisfies \eqref{eq:G_cohypomonotonicity} in Assumption \ref{as:A3} with $\rho_n := \frac{(1+\tau)\rho}{\tau} \geq 0$ and $0 \leq \rho_c \leq \frac{1}{L} - (1+\tau)\rho$ for any $\tau > 0$.
Here, if $0 \leq L\rho \leq \frac{1}{1+\tau} < 1$, then  $\rho_c \geq 0$. 

\vspace{0.5ex}
\noindent$\mathrm{(c)}$ \textbf{Example 3}~(\textit{\textbf{General linear mappings}}).
Let us consider the following linear mappings:
\begin{equation}\label{eq:exam2}
\begin{array}{lcl}
G_ix := \mbf{G}_ix + \mbf{g}_i \quad (i \in [n]) \quad \textrm{and} \quad Gx := \mbf{G}x + \mbf{g} = \frac{1}{n}\sum_{i=1}^nG_ix,
\end{array}
\end{equation}
where $\mbf{G}_i \in \R^{p\times p}$ and $\mbf{g}_i \in \R^p$ for $i \in [n]$ are given, $\mbf{G} := \frac{1}{n}\sum_{i=1}^n\mbf{G}_i$, and $\mbf{g} := \frac{1}{n}\sum_{i=1}^n\mbf{g}_i$.
Suppose that there exist $\rho_n > 0$ and $\rho_c > 0$ such that
\begin{equation}\label{eq:exam2_cond}
\begin{array}{lcl}
\frac{1}{2}(\mbf{G} + \mbf{G}^{\top}\big) + \rho_n\mbf{G}^{\top}\mbf{G} - \frac{\rho_c}{n}\sum_{i=1}^n\mbf{G}_i^{\top}\mbf{G}_i \succeq 0.
\end{array}
\end{equation}
Then, one can easily show that $G$ satisfies the condition~\eqref{eq:G_cohypomonotonicity2}, a special case of the condition \eqref{eq:G_cohypomonotonicity} in Assumption~\ref{as:A3} when $T = 0$.

If $\lambda_{\min}(\frac{1}{2}(\mbf{G} + \mbf{G}^{\top})) < 0$, then $G$ is nonmonotone.
Additionally, if we assume that $\lambda_{\min}(\mbf{G}^{\top}\mbf{G}) > 0$ (i.e., $\mbf{G}$ is nonsingular), then we can choose 
\begin{equation}\label{eq:exam2_cond2}
\begin{array}{lcl}
\rho_n \geq -\frac{\lambda_{\min}(\frac{1}{2}(\mbf{G} + \mbf{G}^{\top}))}{\lambda_{\min}(\mbf{G}^{\top}\mbf{G})} > 0 \quad \textrm{and} \quad 0 \leq \rho_c \leq \frac{n \left[\lambda_{\min}(\frac{1}{2}(\mbf{G} + \mbf{G}^{\top})) + \rho_n\lambda_{\min}(\mbf{G}^{\top}\mbf{G}) \right]}{\sum_{i=1}^n\lambda_{\max}(\mbf{G}_i^{\top}\mbf{G}_i)}.
\end{array}
\end{equation}
In this case, $\rho_n$ and $\rho_c$ chosen as in \eqref{eq:exam2_cond2} satisfy \eqref{eq:exam2_cond}.
Consequently, the condition \eqref{eq:G_cohypomonotonicity2} (a special case of Assumption~\ref{as:A3} when $T=0$) is satisfied.

\vspace{0.5ex}
\noindent$\mathrm{(d)}$~\textbf{Example 4}~(\textit{\textbf{Bounded differential operator diversity}}).
For given $n$ operators $G_i$, we define their difference diversity ratio:
\begin{equation}\label{eq:exam3_ratio}
R(x, y) := \frac{n\sum_{i=1}^n\norms{G_ix - G_iy}^2}{\norms{\sum_{i=1}^n(G_ix - G_iy)}^2}, \quad \forall x, y \in \dom{G}, \ Gx \neq Gy.
\end{equation}
This concept is similar to the differential gradient diversity concept in \cite{yin2018gradient} often used in federated learning and distributed optimization.
Here, we extend it to operators.
Since $\norms{\sum_{i=1}^n(G_ix - G_iy)}^2 \leq n\sum_{i=1}^n\norms{G_ix - G_iy}^2$ by Jensen's inequality, we have $R(x, y) \geq 1$ for all $x, y\in\dom{G}$.
Next, we assume that $R(x, y)$ is upper bounded for all $x, y \in \dom{G}$, i.e., there exists $M \geq 1$ such that $1 \leq R(x, y) \leq M$ for all $x, 
y \in\dom{G}$.

Now, suppose that $G$ is also $\rho$-co-hypomonotone with $\rho > 0$.
Then it satisfies the condition~\eqref{eq:G_cohypomonotonicity2} in Assumption~\ref{as:A3} with $\rho_n > \rho$ and  $\rho_c := \frac{\rho_n - \rho}{M} > 0$  (see Appendix~\ref{apdx:subsec:example_verification}).

\vspace{0.5ex}
\noindent$\mathrm{(e)}$~\textbf{Example 5}~(\textit{\textbf{Bounded derivatives}}).
Suppose that $\dom{G_i}$ is convex and $G_i$ is continuously differentiable for all $i \in [n]$.
Let $\bar{\sigma}_i := \sup_{x \in \dom{G_i}}\sigma_{\max}(G_i'(x))$ be the upper bound of the largest singular value of $G_i'(x)$ for $i \in [n]$ and $\underline{\sigma} > 0$ be  the smallest constant such that $\norms{Gx - Gy} \geq \underline{\sigma}\norms{x - y}$ for all $x, y\in\dom{G}$.
Denote $\bar{\sigma}^2 := \frac{1}{n}\sum_{i=1}^n\bar{\sigma}_i^2$.
Suppose that  $ 0 < \underline{\sigma} \leq \bar{\sigma} < +\infty$.
Then, if $G$ is additionally $\rho$-co-hypomonotone with $\rho > 0$, then it satisfies the condition~\eqref{eq:G_cohypomonotonicity2} in Assumption~\ref{as:A3} with  $\rho_n > \rho$ and $\rho_c := \frac{(\rho_n - \rho)\underline{\sigma}^2 }{\bar{\sigma}^2} > 0$ (see Appendix \ref{apdx:subsec:example_verification}).

\beforesec
\section{Variance Reduced Fast  Optimistic Gradient Method for the General Case.}\label{sec:SFOG}
\aftersec
In this section, we exploit three main ideas: Popov's past-EG scheme, Nesterov's acceleration, and variance reduction techniques to develop a novel algorithmic framework for solving a class of \eqref{eq:GE} under the \textbf{Lipschitz continuity of $G$} and the \textbf{$\rho_n$-co-hypomonotonicity of $\Phi$}.
That is: we solve  a class of \eqref{eq:GE} covered by Assumption~\ref{as:A1} with $\alpha = 1$ and Assumption~\ref{as:A3} with $\rho_n \geq 0$ and $\rho_c = 0$.

\beforesubsec
\subsection{\textbf{The proposed algorithm.}}
\label{subsec:iAEG4GE}
\aftersubsec
The following variant of Nesterov's accelerated method from \cite{Beck2009,Nesterov1983} for minimizing the sum $\phi :=  f + g$ of a convex and $L$-smooth function $f$ and a possibly nonsmooth and convex function $g$ can be expressed as follows:
\vspace{-0.5ex}
\begin{equation}\label{eq:NAG4cvx}
\arraycolsep=0.2em
\left\{\begin{array}{lcl}
\hat{x}^k &:= &  \frac{t_k-s}{t_k}x^k + \frac{s}{t_k}z^k, \vspace{1ex}\\
x^{k+1} &:= & \hat{x}^k - \eta(\nabla{f}(\hat{x}^k) + v^{k+1}), \vspace{1ex}\\
z^{k+1} &:= & z^k - \frac{\eta t_k}{s}(\nabla{f}(\hat{x}^k) + v^{k+1}),
\end{array}\right.
\tag{NAG}
\vspace{-0.5ex}
\end{equation}
where $v^{k+1} \in \partial{g}(x^{k+1})$, $0 < \eta \leq \frac{1}{L}$ is a given stepsize, and $t_k = k + s + 1$ for some $s > 1$.
\revise{Here, we have used the identity  $x^{k+1} = \prox_{\eta g}(x^{k+1} + \eta v^{k+1})$ in the second line.}
This algorithm is known to achieve an $\BigOs{1/k^2}$ convergence rate on the objective residual $\phi(x^k) - \phi(x^{\star})$, significantly improving over the $\BigOs{1/k}$ rate in the standard proximal gradient method.
Under an appropriate choice of $s$, \eqref{eq:NAG4cvx} can achieve up to an $\SmallOs{1/k^2}$ rate \cite{attouch2016rate}.

\vspace{1ex}
\noindent\textbf{Main idea.}
Our idea is to exploit this \eqref{eq:NAG4cvx} scheme and combine it with Popov's past-EG mechanism to develop an accelerated past-extragradient scheme for \eqref{eq:GE}.
However, since we can only access an unbiased stochastic oracle $\mbf{G}$ of  $G$, we will approximate the value $Gy^k$ by a stochastic variance-reduced estimator $\widetilde{G}y^k$ constructed from $\mbf{G}(y^k, \xi)$.
This combination leads to the following new \textit{\textbf{V}ariance-reduced \textbf{A}ccelerated \textbf{P}ast \textbf{E}xtra\textbf{G}radient} method, which we call \ref{eq:VFOG}.

\vspace{0.75ex}
\noindent\textbf{The proposed method.}
Given an initial point $x^0 \in \R^p$, we set $z^0 = y^{-1} := x^0$, \revise{$\widetilde{G}y^{-1} = \widetilde{G}y^{-1} := \widetilde{G}x^0$,} and at each iteration $k \geq 0$, we construct a stochastic estimator $\widetilde{G}y^k$ of $Gy^k$, and then update
\begin{equation}\label{eq:VFOG}
\arraycolsep=0.2em
\left\{\begin{array}{lcl}
\hat{x}^k &:= & \frac{t_k-s}{t_k}x^k + \frac{s}{t_k}z^k, \vspace{1ex}\\
y^k &:= & \hat{x}^k - (\eta - \beta_k)(\widetilde{G}y^{k-1} + v^k), \vspace{1ex}\\
x^{k+1} &:= & \hat{x}^k - \eta (\widetilde{G}y^k + v^{k+1})  +  \beta_k (\widetilde{G}y^{k-1} + v^k), \vspace{1ex}\\
z^{k+1} &:= & z^k - \frac{\gamma_k}{s}(\widetilde{G}y^{k-1} + v^k),
\end{array}\right.
\tag{VAPEG}
\end{equation}
where $v^k \in Tx^k$,  $t_k \geq s > 1$, $\eta > 0$, $\beta_k \geq 0$, $\gamma_k > 0$ are given parameters, determined later.
There are two key components of \eqref{eq:VFOG} that require further discussion.
\begin{compactitem}
\item First, we will clearly characterize the conditions on the  estimator $\widetilde{G}y^k$ of $Gy^k$ to guarantee the convergence of \eqref{eq:VFOG}.
We will also give concrete examples of $\widetilde{G}y^k$.
\item Second, we will specify the choice of parameters $t_k$, $s$, $\eta$, $\beta_k$, and $\gamma_k$ in our convergence analysis to achieve better convergence than its non-accelerated counterparts.
\end{compactitem}
The scheme \eqref{eq:VFOG} is sufficiently general to cover the following special cases.
\begin{compactitem}
\item[$\mathrm{(a)}$] If $\widetilde{G}y^k := Gy^k$ and $\widetilde{G}y^{k-1} := Gx^k$, then we obtain a deterministic fast extragradient method for solving \eqref{eq:GE}.
It is similar to the one in  \cite{tran2025accelerated,yuan2024symplectic}, but the last step is different.
\item[$\mathrm{(b)}$] 
If $\widetilde{G}y^k := Gy^k$ and $\widetilde{G}y^{k-1} := Gy^{k-1}$, then we get a deterministic accelerated past-extragradient or fast optimistic gradient (OG) method for solving \eqref{eq:GE}.
This deterministic variant is also new under Assumptions~\ref{as:A1}-\ref{as:A3}, and also different from \cite{tran2025accelerated}.
\end{compactitem}

\vspace{0.75ex}
\noindent\textbf{Comparison with \eqref{eq:popov_method} and \eqref{eq:NAG4cvx}.}
In fact, \eqref{eq:VFOG} mimics Popov's past-extragradient scheme, Nesterov's acceleration, and variance reduction techniques to create such a new algorithm.
Formally, if $T = 0$, $\widetilde{G}y^k = Gy^k$, $\beta_k = 0$, and $s = 0$, then this scheme reduces to \eqref{eq:popov_method} discussed in Section~\ref{sec:background}.
At the same time, it incorporates an intermediate step $\hat{x}^k = \frac{t_k-s}{t_k}x^k +  \frac{s}{t_k}z^k$ and uses the varying parameters $t_k$ and $\beta_k$, similarly to Nesterov's accelerated method \eqref{eq:NAG4cvx} in convex optimization.
However, the first fundamental difference between \eqref{eq:NAG4cvx} and \eqref{eq:VFOG} is the use of the stepsize $\frac{\gamma_k}{s}$ instead of $\frac{\eta t_k}{s}$ in the last line of \eqref{eq:VFOG}.
The second one is the search direction $\eta \widetilde{G}y^k -  \beta_k \widetilde{G}y^{k-1}$ in  \eqref{eq:VFOG} combining two consecutive approximations $\widetilde{G}y^k$ and $\widetilde{G}y^{k-1}$ of $G$, which acts like $\eta\nabla{f}(\hat{x}^k)$ in \eqref{eq:NAG4cvx} for smooth and convex optimization.

\vspace{0.75ex}
\noindent\textbf{The implementable version.}
Since $v^{k+1}$ is on the right-hand side of the third line of \eqref{eq:VFOG}, we can use the resolvent $J_{\eta T}$ of $\eta T$ to rewrite it equivalently to
\begin{equation}\label{eq:VFOG_impl}
\arraycolsep=0.2em
\left\{\begin{array}{lcl}
\hat{x}^k &:= & \frac{s}{t_k}z^k + \frac{t_k-s}{t_k}x^k, \vspace{1ex}\\
y^k &:= & \hat{x}^k - (\eta - \beta_k)(\widetilde{G}y^{k-1} + v^k), \vspace{1ex}\\
x^{k+1} &:= & J_{\eta T} \big( \hat{x}^k - \eta \widetilde{G}y^k  +  \beta_k (\widetilde{G}y^{k-1} + v^k) \big), \vspace{1ex}\\
z^{k+1} &:= & z^k - \frac{\gamma_k}{s}(\widetilde{G}y^{k-1} + v^k), \vspace{1ex}\\
v^{k+1} &:= & \frac{1}{\eta}\big(\hat{x}^k - x^{k+1}   +  \beta_k (\widetilde{G}y^{k-1} + v^k) \big) - \widetilde{G}y^k.
\end{array}\right.
\end{equation}
Clearly, each iteration of \eqref{eq:VFOG_impl} requires one construction of the estimator $\widetilde{G}y^k$ and one evaluation of $J_{\eta T}$.
\revise{
Note that, for our theoretical analysis, we require the exact evaluation of the resolvent $J_{\eta T}$.
This requirement is satisfied in many cases where $T$ has a ``simple form,'' such as the subdifferential of a proximally tractable function (e.g., $\norms{\cdot}_1$, coordinate-wise separable functions, or the indicator function of a simple convex set).
In practice, however, the resolvent $J_{\eta T}$ can also be evaluated approximately to a prescribed accuracy, thereby broadening the applicability of our method \eqref{eq:VFOG_impl}, see, e.g., \cite{Rockafellar1976b}.
}

\beforesubsec
\subsection{\textbf{Main result 1: Convergence properties of \ref{eq:VFOG}.}}\label{subsec:convergence1}
\aftersubsec
There are different strategies to approximate $Gy^k$ by $\widetilde{G}y^k$.
For instance, one can use increasing mini-batch stochastic approximation strategies as in \cite{boct2021minibatch,pethick2023solving} or  control variate techniques as in \cite{Defazio2014,SVRG,nguyen2017sarah}.
In what follows, we do not specify how to approximate $Gy^k$, but introduce a general recursive bound for $\Expsn{k}{\norms{\widetilde{G}y^k - Gy^k}^2}$.

\beforesubsubsec
\subsubsection{Error approximation condition between $\widetilde{G}y^k$ and $Gy^k$.}\label{subsubsec:error_bound_def}
\aftersubsubsec
We assume that $\widetilde{G}y^k$ is constructed by accessing the unbiased stochastic oracle $\mbf{G}(y^k, \xi)$ of $Gy^k$ adapted to the filtration $\sets{\Fc_k}$ such that for $k \geq \revise{0}$, there exist two constants $\kappa \in (0, 1]$ and $\Theta \geq 0$, a nonnegative sequence $\sets{\delta_k}$,  and a sequence of nonnegative random variables $\sets{\Delta_k}$ such that the following condition holds almost surely for the error term $e^k := \widetilde{G}y^k - Gy^k$: 
\myeq{eq:error_cond3}{
\arraycolsep=0.2em
\left\{\begin{array}{lcl}
\Expsn{k}{ \norms{e^k}^2 } & \leq & \Expsn{k}{\Delta_k}, \vspace{1ex}\\
\Expsn{k}{ \Delta_k } & \leq & (1-\kappa)\Delta_{k-1} + \Theta L^2 \norms{x^k - y^{k-1} }^2 + \frac{\delta_k}{t_k^2},
\end{array}\right.
}
where $\sets{(x^k, y^k)}$ and $\sets{t_k}$ are updated by our scheme \eqref{eq:VFOG}.

The condition \eqref{eq:error_cond3} is rather general.
Indeed, if $\kappa=1$ and $\Theta = 0$, then it reduces to $\Expsn{k}{ \norms{e^k}^2 } \leq \frac{\delta_k}{t_k^2}$.
If $\kappa = 1$ and $\delta_k = 0$, then it becomes $\Expsn{k}{ \norms{e^k}^2 } \leq \Theta L^2 \norms{x^k - y^{k-1}}^2$, which can be viewed as an adaptive error bound condition.
One way to choose $\delta_k$ is $\delta_k := \frac{\delta}{(k+r+1)^{1+\nu}}$ for any $\nu \geq 0$ and given $\delta \geq 0$ and $r > 2$.
In the sequel, we will  denote by $\mcal{S}_K := \sum_{k=0}^{K}\delta_k$ and $\mcal{S}_{\infty} := \sum_{k=0}^{\infty}\delta_k$.

\vspace{0.75ex}
\noindent\textbf{\textit{Adaptive mini-batch estimator.}}
One way to construct $\widetilde{G}y^k$ is to use a mini-batch estimator with adaptively increasing size.
Given a mini-batch $\Bc_k$ of size $b_k$, we construct 
\myeq{eq:mini_batch_Gy}{
\arraycolsep=0.2em
\begin{array}{lcl}
\widetilde{G}y^k & := & \frac{1}{b_k}\sum_{\xi_k \in \Bc_k}\mbf{G}(y^k, \xi_k).
\end{array}
}
It is well-known that $\Expsn{\Bc_k}{ \norms{\widetilde{G}y^k - Gy^k}^2 } \leq \frac{\sigma^2}{b_k}$, where $\sigma^2$ is the variance of the stochastic oracle $\mbf{G}$ in Assumption~\ref{as:A0}.
To guarantee the condition \eqref{eq:error_cond3}, we need to choose $b_k$ adaptively.
For example, we can choose $b_k$ using the following condition:
\myeq{eq:mini_batch_size_cond}{
\arraycolsep=0.2em
\begin{array}{lcl}
\frac{\sigma^2}{b_k} \leq \frac{(1-\kappa)\sigma^2}{b_{k-1}} + \Theta L^2\norms{x^k - y^{k-1}}^2 + \frac{\delta_k}{t_k^2}.
\end{array}
}
Clearly, by setting $\Delta_k := \frac{\sigma^2}{b_k}$, $\widetilde{G}y^k$ defined by \eqref{eq:mini_batch_Gy} satisfies \eqref{eq:error_cond3}.
Moreover,  the condition \eqref{eq:mini_batch_size_cond} can be guaranteed by choosing an appropriate increasing mini-batch size $b_k$.
This technique has been used in several previous works such as \cite{boct2021minibatch,pethick2023solving}.
Nevertheless,  \eqref{eq:error_cond3} can cover a wide range of estimators, including unbiased and biased instances, beyond the classical mini-batch instance \eqref{eq:mini_batch_Gy}.

\beforesubsubsec
\subsubsection{Convergence rates in expectation.}\label{subsubsec:convergence_result1}
\aftersubsubsec
Given \revise{$s > 2$}, $\kappa$ and $\Theta$ in \eqref{eq:error_cond3}, we first define
\begin{equation}\label{eq:fixed_constants}
\arraycolsep=0.2em
\left\{\begin{array}{lcl}
\Lambda_0 &:= & \frac{s^2}{s^2 - (1-\kappa)(s+1)^2} \revise{\left[\frac{25s-34}{2(s-2)} + \frac{9s^2 + 740s - 1192}{64(s-1)}\right]} \vspace{1ex}\\
\omega & := & \revise{\frac{9s^2 + 754s - 1220}{256(s-1)}} + \frac{2\Lambda_0\Theta(s+1)^2}{s^2}, \vspace{1ex}\\
\lambda &:= & \frac{1}{\sqrt{2(1+\omega)(s+1)}}, \vspace{1ex}\\
\mu &:= & \frac{\revise{s-2}}{\revise{8}(s-1)}\lambda.
\end{array}\right.
\end{equation}
It is easy to check that these constants are positive and $\lambda \in (0, 1)$.

Now, let us state the convergence of \eqref{eq:VFOG} in Theorem~\ref{th:VFOG1_convergence} below for the recursive error bound condition \eqref{eq:error_cond3}.
The proof of Theorem~\ref{th:VFOG1_convergence} is given in Subsection~\ref{subsec:th:VFOG1_convergence}. 

\begin{theorem}\label{th:VFOG1_convergence}
Let  $\Lambda_0$, $\mu$, and $\lambda$ be defined by \eqref{eq:fixed_constants}.
For \eqref{eq:GE}, suppose that $G$ satisfies Assumption~\ref{as:A0} and Assumption~\ref{as:A1} with $\alpha = 1$ $($i.e., $G$ is $L$-Lipschitz continuous$)$, and $\Phi$ satisfies Assumption \ref{as:A3} with $\rho_n \geq 0$ and $\rho_c = 0$ $($i.e., $\Phi$ is $\rho_n$-co-hypomonotone$)$ such that $L\rho_n ~\revise{<}~ \mu$.

Let $\sets{(x^k, y^k, z^k)}$ be generated by \eqref{eq:VFOG} using an estimator $\widetilde{G}y^k$ of $Gy^k$ satisfying the condition \eqref{eq:error_cond3} with $\frac{2s+1}{(s+1)^2} < \kappa \leq 1$ for a given \revise{$s > 2$} 
and the following parameters:
\begin{equation}\label{eq:SFOG_para_update2}
\hspace{-0.5ex}
\arraycolsep=0.2em
\begin{array}{ll}
\frac{\revise{8}(s-1)\rho_n}{\revise{s-2}} \leq \eta <  \frac{\lambda}{L}, \quad 
\gamma_k := \frac{\revise{(s-2)}\eta(k+s)}{\revise{32(s-1)}(k+s+1)},  \quad \textrm{and} \quad \beta_k :=  \big[\frac{\revise{3}(s-2)\eta}{\revise{16}(s-1)} + 2\rho_n]\frac{k+1}{k+s+1} - \frac{\gamma_k}{k+s+1}.
\end{array}
\hspace{-0ex}
\end{equation}
Then, the following statement holds:
\vspace{0.75ex}
\begin{compactitem}
	\item[$\mathrm{(a)}$] $($\textbf{Summability bounds}$)$
	There exist $v^k \in Tx^k$ for all $k \geq 0$ such that
	\begin{equation}\label{eq:SFOG_summability_bounds}
		\arraycolsep=0.2em
		\begin{array}{lcl}
			\sum_{k=0}^K \revise{\frac{(s-2)(4s-7)}{s-1}}   (k+1)    \Exp{\norms{Gx^k + v^k }^2} & \leq & \revise{64} \big( \mcal{R}_0^2 + \Lambda_0\mcal{S}_K \big), \vspace{1.75ex}\\
			\sum_{k=0}^K s (k+s+1) \Exp{\norms{Gy^{k} + v^{k+1} }^2} & \leq & 2  \mcal{R}_0^2 + 2\Lambda_0\mcal{S}_K, 
		\end{array}
	\end{equation}
	where $\mcal{S}_K := \sum_{k=0}^K \delta_k$ and $\mcal{R}_0^2$ is given by 
	\begin{equation}\label{eq:SFOG_BigO_R0}
	\arraycolsep=0.2em
	\begin{array}{lcl}
		\mcal{R}_0^2 & := &  \revise{\frac{2(13s-10)s^2 - (s-1)(s-2)}{64(s-1)}} \norms{Gx^0 + v^0 }^2 + \revise{\frac{16s^3(s-1)}{(s-2)\eta^2}} \norms{x^0 - x^{\star} }^2.
	\end{array}
	\end{equation}
\end{compactitem}
\begin{compactitem}	
	\item[$\mathrm{(b)}$]$($\textbf{The $\BigO{1/k^2}$ convergence rates}$)$ We have the following bounds:
	\begin{equation}\label{eq:SFOG_BigO_rates}
	\arraycolsep=0.2em
	\begin{array}{lcl}
	\Expn{ \norms{Gx^k + v^k }^2 } &\leq &  \dfrac{  C_0^2 (\mcal{R}_0^2 + \Lambda_0 \mcal{S}_k )}{ (k+s)^2 }  \vspace{1ex}\\
	\Expn{ \norms{Gy^{k-1} + v^k }^2 } & \leq &  \dfrac{2C_0^2 ( \mcal{R}_0^2 + \Lambda_0 \mcal{S}_k)}{ (k+s)^2 },
	\end{array}
	\end{equation}
	where $C_0^2 := \revise{\frac{128(s-1)}{19s+10}}$.
	\item[$\mathrm{(c)}$]$($\textbf{The $\SmallO{1/k^2}$-convergence rates}$)$
	If additionally $\mcal{S}_{\infty} := \sum_{k=0}^{\infty}\delta_k < +\infty$, then 
	\begin{equation}\label{eq:SFOG_small_o_rates}
	\arraycolsep=0.2em
	\begin{array}{lcl}
	\lim_{k\to\infty}  k^2  \Expn{\norms{Gx^k + v^k }^2} = 0 \ \quad \textrm{and} \ \quad \lim_{k\to\infty}  k^2  \Expn{\norms{Gy^{k-1} + v^k }^2} = 0.
	\end{array}
	\end{equation}
	
	\item[$\mathrm{(d)}$]$($\textbf{Iteration-complexity}$)$
	For a given $\epsilon > 0$, if $\mcal{S}_{\infty} < +\infty$, then the number of iterations $k$ to achieve $\Expn{ \norms{Gx^k + v^k}^2 } \leq \epsilon^2$ is at most $\BigOs{1/\epsilon}$.
	This  is also the number of evaluations of $J_{\eta T}$. 
\end{compactitem}
\end{theorem}

\beforesubsubsec
\subsubsection{Almost sure convergence.}\label{subsubsec:almost_sure_convergence1}
\aftersubsubsec
Next, we state the almost sure convergence rates and the almost sure convergence of iterates in the following theorem, whose proof is given in Subsection~\ref{apdx:th:VFOG1_almost_sure_convergence}.

\begin{theorem}\label{th:VFOG1_almost_sure_convergence}
Under the same assumptions and settings as in Theorem~\ref{th:VFOG1_convergence}, assuming that $\Sc_{\infty} := \sum_{k=0}^{\infty}\delta_k < +\infty$, we have the following statements.
\begin{compactitem}
\item[$\mathrm{(a)}$] The following summability results hold:
\begin{equation}\label{eq:SFOG_as_summability}
\arraycolsep=0.2em
\begin{array}{lcl}
\sum_{k=0}^{\infty}   (k+1)  \norms{Gx^k + v^k }^2 & < & +\infty, \quad \textrm{almost surely}, \vspace{1.5ex}\\
\sum_{k=0}^{\infty} (k+s)   \norms{Gy^{k-1} + v^k }^2 & < & +\infty, \quad \textrm{almost surely}.
\end{array}
\end{equation}
\item[$\mathrm{(b)}$] The following limits hold almost surely:
\begin{equation}\label{eq:SFOG_as_limits}
\arraycolsep=0.2em
\begin{array}{lcl}
\lim_{k\to\infty}   k^2 \norms{Gx^k + v^k }^2 = 0 \quad \textrm{and} \quad \lim_{k\to\infty}   k^2  \norms{Gy^{k-1} + v^k }^2 = 0.
\end{array}
\end{equation}
\item[$\mathrm{(c)}$] 
If, additionally, $\gra{\Phi}$ is closed, then $\sets{x^k}$ almost surely converges to a $\zer{\Phi}$-valued random variable $x^{\star} \in \zer{\Phi}$.
In addition,  the other two sequences of iterates $\sets{y^k}$ and $\sets{z^k}$ also almost surely converge to the same random variable $x^{\star} \in \zer{\Phi}$.
\end{compactitem}
\end{theorem}

Note that if Assumption~\ref{as:A3} holds with $\rho_c = 0$, then Assumption~\ref{as:A2} also automatically holds  with $\rho_{*} = \rho_n$.
Therefore, we do not require Assumption~\ref{as:A2} in Theorems~\ref{th:VFOG1_convergence} and \ref{th:VFOG1_almost_sure_convergence}.

If $T$ is $\rho_T$-co-hypomonotone for some $\rho_T \geq 0$, then,  as proven in \cite{TranDinh2025a}, the resolvent $J_{\lambda T}$ is well-defined, single-valued, and nonexpansive if $\lambda > \rho_T$.
Let us define the following forward-backward splitting (FBS) residual of \eqref{eq:GE}, widely used in monotone operator theory:
\begin{equation}\label{eq:FBS_residual}
\arraycolsep=0.2em
\begin{array}{lcl}
\mcal{G}_{\lambda}x := \frac{1}{\lambda}\big(x - J_{\lambda T}(x - \lambda Gx) \big).
\end{array}
\end{equation}
It is well-known that $x^{\star} \in \zer{\Phi}$ iff $\mcal{G}_{\lambda}x^{\star} = 0$.
Therefore, for a given tolerance $\epsilon > 0$, we can say that $x^k$ is an $\epsilon$-solution of \eqref{eq:GE} (in expectation) if $\Expn{ \norms{\mcal{G}_{\lambda}x^k}^2} \leq \epsilon^2$.
Since $\norms{\mcal{G}_{\lambda}x^k} \leq \norms{Gx^k + v^k}$ for any $v^k \in Tx^k$ (see \cite{tran2024revisiting}), if we obtain $\Expn{ \norms{Gx^k + v^k}^2 } \leq \epsilon^2$ as in Theorem~\ref{th:VFOG1_convergence}, then we conclude that $x^k$ is also an $\epsilon$-solution of \eqref{eq:GE}, i.e., $\Expn{ \norms{\mcal{G}_{\lambda}x^k}^2} \leq \epsilon^2$.

\begin{remark}\label{re:closedness_of_graph}
Note that if $\Phi$ is maximally monotone, then $\gra{\Phi}$ is closed, see \cite[Proposition 20.38]{Bauschke2011}.
The closed-graph assumption is weaker than maximal monotonicity and is satisfied by many classes of continuous or closed set-valued mappings.
In the single-valued case, one can replace this closedness assumption of $\gra{\Phi}$ by a demiclosedness of $\Phi$, see, e.g., \cite{davis2022variance}.
If $T = 0$, then $\Phi = G$, which is continuous, and thus $\gra{G}$ is closed.
\end{remark}

\beforesubsec
\subsection{\textbf{Oracle complexity of \ref{eq:VFOG} with adaptive mini-batching estimators.}}
\label{subsec:minibatching_estimator}
\aftersubsec
Now, we estimate the oracle complexity of \eqref{eq:VFOG} using the adaptive mini-batch estimator $\widetilde{G}y^k$ in \eqref{eq:mini_batch_Gy} in the following corollary, whose proof can be found in Appendix~\ref{apdx:co:complexity_of_VFOG}.

\begin{corollary}\label{co:complexity_of_VFOG}
Under the same conditions and settings as in Theorem~\ref{th:VFOG1_convergence}, if $\widetilde{G}y^k$ constructed by \eqref{eq:mini_batch_Gy} with increasing mini-batch size $b_k := \big\lceil \frac{\sigma^2(k+s+1)^{3+\nu}}{\delta} \big\rceil$ for some constants $\nu > 0$ and $\delta > 0$ is used in \eqref{eq:VFOG}, then the expected total number of oracle calls $\mbf{G}$ to achieve $x^K$ such that $\Expn{\norms{Gx^K + \revise{v^K}}^2} \leq \epsilon^2$ for a given tolerance $\epsilon > 0$ and $\revise{v^K} \in Tx^K$ is at most $\Expn{\Tc_K} = \BigOs{ \frac{\sigma^2}{\epsilon^{4 + \nu}} }$.
The expected total number of $J_{\eta T}$ evaluations is at most  $\BigOs{\epsilon^{-1}}$.

Alternatively, if we choose $b_k := \big\lceil \frac{\sigma^2(k+s+1)^3\log(k+s+1)}{\delta} \big\rceil$, then the above complexity is improved to $\Expn{\Tc_K} = \widetilde{\mathcal{O}}\big( \frac{\sigma^2}{\epsilon^{4}} \big)$.
In this case, the expected total number of $J_{\eta T}$ evaluations is at most $\widetilde{\mcal{O}}(\epsilon^{-1})$.
\end{corollary}

The work \cite{boct2021minibatch} studied VIPs, a special case of \eqref{eq:GE}, and established a $\BigOs{\epsilon^{-4}\ln(\epsilon^{-1})}$ oracle complexity using increasing mini-batch estimators as \eqref{eq:mini_batch_Gy}.
Our complexity in Corollary~\ref{co:complexity_of_VFOG} is either $\BigOs{\epsilon^{-4+\nu}}$ for any sufficiently small $\nu > 0$ or $\BigOs{\epsilon^{-4}\ln(\ln(\epsilon^{-1}))}$ depending on the choice of $b_k$, which is comparable with \cite{boct2021minibatch}.
Nevertheless, our assumptions, Assumptions~\ref{as:A2} and \ref{as:A3}, are different from the pseudo-monotonicity in \cite{boct2021minibatch}.
Moreover, as mentioned earlier, our method can work with any unbiased and biased estimators under the condition \eqref{eq:error_cond3}, beyond \eqref{eq:mini_batch_Gy}. 
Increasing mini-batches were also used in \cite{pethick2023solving}, but the authors did not explicitly derive any oracle complexity, and also used different assumptions and metrics.

Note that our mini-batch size $b_k$ in Corollary~\ref{co:complexity_of_VFOG} is loosely chosen from the condition \eqref{eq:error_cond3} with $\kappa = 1$ and $\Theta = 0$.
However, we can directly use the condition~\eqref{eq:mini_batch_size_cond} to estimate $b_k$, which gives us a smaller batch size than the worst case $b_k$ chosen in Corollary~\ref{co:complexity_of_VFOG}.

Corollary~\ref{co:complexity_of_VFOG} also shows that if we use standard stochastic approximation \eqref{eq:mini_batch_Gy} for $Gy^k$, then our accelerated method \eqref{eq:VFOG} does not achieve better complexity than $\widetilde{\mcal{O}}(\epsilon^{-4})$.
This is due to the accumulation of errors in accelerated methods compared to non-accelerated ones as studied in \cite{Devolder2010}.
To achieve better complexity bounds, other variance-reduction techniques should be exploited.
For example, in Section~\ref{sec:VFOG} below, we will use control variate techniques.

\beforesec
\section{Convergence and Complexity of \ref{eq:VFOG} using Control Variate Estimators.}\label{sec:VFOG}
\aftersec
In this section, we will employ control variate techniques to develop a variance-reduction variant of \eqref{eq:VFOG}.
Unlike Section~\ref{sec:SFOG}, we let $\alpha \in [0, 1)$ in Assumption~\ref{as:A1} and require $\rho_c > 0$ in Assumption~\ref{as:A3}.

\beforesubsec
\subsection{\textbf{A class of stochastic variance-reduced  estimators.}}
\label{subsec:Vr_Estimator}
\aftersubsec
We propose to use the following class of variance-reduced estimators $\widetilde{G}y^k$ of $Gy^k$, which was introduced in our recent work \cite{TranDinh2025a}.
This class of estimators is relatively broad and covers various well-known instances, including Loopless SVRG, SAGA, and Loopless SARAH.

\begin{definition}\label{de:VR_Estimators}
\textit{
Given $\underline{\kappa} \in (0, 1)$ and $\Theta > 0$, let $\sets{y^k}$ be the sequence of iterates generated by \eqref{eq:VFOG} and $\widetilde{G}y^k$ be a stochastic estimator of $Gy^k$ constructed from an i.i.d. sample $\Bc_k \subseteq [n]$, adapted to the filtration $\sets{\Fc_k}$.
We say that $\widetilde{G}y^k$ satisfies a $\textbf{VR}(\Delta_k; \kappa_k, \Theta_k, \delta_k )$  $($variance-reduction$)$  property if there exist two parameters $\kappa_k \in [\underline{\kappa}, 1]$ and $\Theta_k \in [\underline{\Theta}, \bar{\Theta}]$, a nonnegative sequence $\set{\delta_k}$, and a nonnegative random quantity  $\Delta_k$  such that $($almost surely$)$:
\myeq{eq:VR_property}{
\arraycolsep=0.2em
\left\{\begin{array}{lcl}
\Expsn{k}{ \norms{\widetilde{G}y^k - Gy^k }^2 } & \leq & \Expsn{k}{ \Delta_k }, \vspace{1ex}\\
\Expsn{k}{ \Delta_k }  & \leq & (1- \kappa_k)\Delta_{k-1}  + \Theta_k \Expsn{\xi}{ \norms{\mbf{G}(y^k,\xi) - \mbf{G}(y^{k-1},\xi) }^2 } + \frac{\delta_k}{t_k^2},
\end{array}\right.
}
where $\sets{y^k}$ and $\sets{t_k}$ are updated by our scheme \eqref{eq:VFOG}.
}
\end{definition}

Definition~\ref{de:VR_Estimators} covers both unbiased and biased estimators $\widetilde{G}y^k$ of $Gy^k$ since we do not impose the unbiased condition $\Expsn{\Bc_k}{\widetilde{G}y^k} = Gy^k$.
It is also different from existing works, including \cite{driggs2019accelerating,tran2024accelerated,TranDinh2024}, as we only require the condition~\eqref{eq:VR_property}, making it broader and allowing it to cover several existing estimators that may violate the definitions in \cite{driggs2019accelerating,tran2024accelerated,TranDinh2024}.
Our class is also different from other generalizations such as \cite{beznosikov2023stochastic,gorbunov2022stochastic,loizou2021stochastic} for stochastic optimization since they required the unbiasedness and  additional or different conditions.

Now we provide three examples: Loopless SVRG, SAGA, and Loopless SARAH, that satisfy Definition~\ref{de:VR_Estimators}.
Note that the Loopless SVRG and SAGA are unbiased, whereas the Loopless SARAH is biased.
In this paper, we only focus on these three estimators, but we believe that other variance-reduced estimators such as Hybrid SGD \cite{Tran-Dinh2019a},  SAG \cite{LeRoux2012,schmidt2017minimizing}, SARGE \cite{driggs2019bias}, SEGA \cite{hanzely2018sega}, and JacSketch \cite{gower2021stochastic} could also be used in our methods.

\beforesubsubsec
\subsubsection{\textbf{Example 1: L-SVRG estimator.}}
\label{subsubsec:svrg_estimator}
\aftersubsubsec
SVRG was introduced in \cite{SVRG}, and its Loopless version, called L-SVRG, was proposed in \cite{kovalev2019don}.
We show that this estimator satisfies Definition~\ref{de:VR_Estimators}.

For given $G$ in \eqref{eq:GE}, two iterates $y^k$ and $\bar{y}^k$, and i.i.d. sample $\Bc_k$ of size $b_k$, we define  $\mbf{G}_{\Bc_k}(\cdot) := \frac{1}{b_k}\sum_{\xi_k \in \Bc_k}\mbf{G}(\cdot, \xi_k)$ is a mini-batch estimator, and then construct $\widetilde{G}y^k$ as follows:
\myeq{eq:loopless_svrg}{
\widetilde{G}y^k := G\bar{y}^k  + \mbf{G}_{\Bc_k}y^k  - \mbf{G}_{\Bc_k}\bar{y}^k,
\tag{L-SVRG}
}
where, for all $k\geq 1$,  $\bar{y}^k$ is updated by 
\myeq{eq:xy_hat}{
\bar{y}^{k}  := \begin{cases}
y^{k-1} &\text{with probability}~\mbf{p}_{k} \vspace{1ex} \\
\bar{y}^{k-1} &\text{with probability}~1-\mbf{p}_{k}.
\end{cases}
}
Here, $\mbf{p}_k\in(0,1)$ is a given Bernoulli probability parameter, and $\bar{y}^0:=y^0$.

The following lemma shows that $\widetilde{G}y^k$ satisfies Definition~\ref{de:VR_Estimators}, whose proof is given in \cite{TranDinh2025a}.

\vspace{-1ex}
\begin{lemma}\label{le:loopless_svrg_bound}
The estimator $\widetilde{G}y^k$ constructed by \eqref{eq:loopless_svrg} such that $b_k \geq b_{k-1}$ satisfies Definition~\ref{de:VR_Estimators} with 
$\kappa_k := \frac{\mbf{p}_k}{2}$,  $\Theta_k := \frac{4}{b_k\mbf{p}_k}$, $\delta_k = 0$, and $\Delta_k := \frac{1}{b_k} \Expsn{\xi}{ \norms{\mbf{G}(y^k, \xi) - \mbf{G}(\bar{y}^k, \xi)}^2 }$.
\end{lemma}

\beforesubsubsec
\subsubsection{\textbf{Example 2: SAGA estimator.}}
\label{subsubsec:saga_estimator}
\aftersubsubsec
SAGA was introduced in \cite{Defazio2014} for solving finite-sum convex optimization.
We apply it to \eqref{eq:GE} when $G$ has a finite-sum structure as in \eqref{eq:finite_sum}.

For given $G$ in \eqref{eq:finite_sum}, an iterate sequence $\sets{y^k}_{k\geq 0}$, and a random mini-batch $\Bc_k \subseteq [n] := \sets{1, \cdots, n}$ of size $b_k$ sampled independently across iterations, for $k\geq 1$, we update $\widehat{G}_i^k$ for all $i \in [n]$ as follows:
\myeq{eq:SAGA_ref_points}{
\widehat{G}_i^k := \left\{\begin{array}{lll}
G_iy^{k-1}  &\text{if}~ i \in \Bc_k, \vspace{1ex} \\
\widehat{G}_i^{k-1} & \text{if}~i \notin \Bc_k \quad \textrm{(no update)}.
\end{array}\right.
}
Then, we  construct a SAGA estimator $\widetilde{G}y^k$ for $Gy^k$ as follows:
\myeq{eq:SAGA_estimator}{
\arraycolsep=0.2em
\begin{array}{lcl}
\widetilde{G}y^k & := & \frac{1}{n} \sum_{i =1}^n \widehat{G}_i^k  + G_{\Bc_k}y^k  - \widehat{G}_{\Bc_k}^k, 
\end{array}
\tag{SAGA}
}
where $G_{\Bc_k}y^k := \frac{1}{b_k}\sum_{i \in \Bc_k}G_iy^k$ and $\widehat{G}_{\Bc_k}^{k} := \frac{1}{b_k}\sum_{i \in\Bc_k}\widehat{G}_i^{k}$.
In this case, we need to store $n$ component $\widehat{G}^k_{i}$ computed so far for all $i \in [n]$ in a table $\mathbb{T}_k := [\widehat{G}_1^k, \cdots, \widehat{G}_n^k]$ initialized at $\widehat{G}_i^0 := G_iy^0$ for all $i \in [n]$.
SAGA requires significant memory to store $\mathbb{T}_k$ if $n$ and $p$ are both large.
We have the following result, whose proof can be found in \cite{TranDinh2025a}.

\vspace{-1ex}
\begin{lemma}\label{le:SAGA_estimator_full}
The estimator $\widetilde{G}y^k$ constructed by \eqref{eq:SAGA_estimator} such that $b_{k-1} - \frac{b_kb_{k-1}}{4n} \leq b_k \leq b_{k-1}$ satisfies Definition~\ref{de:VR_Estimators} with $\kappa_k := \frac{b_k }{2n}$,  $\Theta_k := \frac{5n}{b_k^2}$, $\delta_k = 0$, and $\Delta_k :=   \frac{1}{nb_k}{\sum_{i=1}^n } \norms{ G_iy^k -  \widehat{G}_i^k }^2$.
\end{lemma}

\beforesubsubsec
\subsubsection{\textbf{Example 3: L-SARAH estimator.}}
\label{subsubsec:sarah_estimator}
\aftersubsubsec
SARAH was introduced in \cite{nguyen2017sarah} for finite-sum convex optimization.
Its loopless variant, called L-SARAH, was studied in \cite{driggs2019accelerating,li2020page} for convex and nonconvex optimization, and recently in \cite{cai2023variance,cai2022stochastic,TranDinh2025a} for co-coercive root-finding problems.
This estimator is constructed as follows.

Given $G$, $\sets{y^k}$, and an i.i.d. sample $\Bc_k \subseteq [n]$ of size $b_k$, we construct  $\widetilde{G}y^k$ as
\myeq{eq:loopless_sarah_estimator}{
\widetilde{G}y^k := \begin{cases} 
\widetilde{G}y^{k-1}  + \mbf{G}_{\Bc_k}y^k - \mbf{G}_{\Bc_k}y^{k-1} &\text{with probability}~1 - \mbf{p}_k, \vspace{1ex}\\
Gy^k &\text{with probability}~\mbf{p}_k,
\end{cases} 
\tag{L-SARAH}
}
where $\mbf{G}_{\Bc_k}(\cdot) := \frac{1}{b_k}\sum_{\xi_k \in\Bc_k}\mbf{G}(\cdot, \xi_k)$, $\widetilde{G}y^0 := Gy^0$, and $\mbf{p}_k \in [\underline{\mbf{p}}, 1]$ is a given Bernoulli switching probability and $\underline{\mbf{p}} \in (0, 1)$ is a given lower bound.
The following lemma shows that $\widetilde{G}y^k$ satisfies Definition~\ref{de:VR_Estimators}, whose proof is presented in \cite{TranDinh2025a}.

\vspace{-1ex}
\begin{lemma}\label{le:loopless_sarah_bound}
The estimator \revise{$\widetilde{G}y^k$} constructed by \eqref{eq:loopless_sarah_estimator} satisfies Definition~\ref{de:VR_Estimators} with 
$\kappa_k :=  \mbf{p}_k$, $\Theta_k  := \frac{1}{b_k}$, $\delta_k = 0$, and $\Delta_k :=  \norms{\widetilde{G}y^k - Gy^k }^2$.
\end{lemma}

\beforesubsec
\subsection{\textbf{Main result 2: Convergence properties  of \ref{eq:VFOG} using $\widetilde{G}y^k$ in Definition~\ref{de:VR_Estimators}.}}
\label{subsec:convergence_of_VFOG}
\aftersubsec
For given \revise{$s > 2$}, $\underline{\Theta}$ in Definition~\ref{de:VR_Estimators}, and $L$ and $\alpha \in [0, 1)$ in Assumption~\ref{as:A1}, we define
\begin{equation}\label{eq:omega_quantity}
\arraycolsep=0.2em
\left\{\begin{array}{lcl}
\revise{\varphi_s} &\revise{:=}& \revise{\frac{9(85s-134)}{64(s-1)} + \frac{9(s-2)}{64},} \vspace{1ex}\\
\hat{\omega} &:= & \frac{\revise{s-2}}{\revise{2}(1-\alpha)(s-1)} + \revise{\frac{\varphi_s}{4} + \frac{s-2}{16(s-1)}}, \vspace{1ex}\\
\hat{\lambda} & := & \frac{1}{\sqrt{2(s+1)(1+ \hat{\omega})}}, \vspace{1ex}\\
\hat{\mu} & := & \frac{\revise{(s-2)}\hat{\lambda}}{\revise{8}(s-1)}, \vspace{1ex}\\ 
\Gamma & := & \frac{3s^2}{s+1} \revise{\big[\varphi_s + \frac{25s-34}{2(s-2)}\big]}, \vspace{1ex}\\ 
\hat{\Lambda}_0 & := &  \revise{ \frac{(s-2)\rho_c}{24(s^2-1)\rho_n\underline{\Theta}}}.
\end{array}\right.
\end{equation}
It is easy to check that these constants are positive and $\hat{\lambda} \in (0, 1)$. 

Now, we are ready to state the convergence and convergence rates of \eqref{eq:VFOG}, whose proof is deferred to Subsection~\ref{subsection:th:VFOG_convergence} for the sake of presentation.

\begin{theorem}\label{th:VFOG_convergence}
Let $\hat{\mu}$, $\hat{\lambda}$, and $\Gamma$  be defined by \eqref{eq:omega_quantity}.
Suppose that $G$ and $T$ in \eqref{eq:GE} satisfy Assumptions~\ref{as:A0}, \ref{as:A1}, \ref{as:A2}, and \ref{as:A3} with $\alpha \in [0, 1)$, $\rho_n \geq \rho_c > 0$, $0 \leq \rho_{*} \leq \rho_n$, and $L\rho_n < \hat{\mu}$.

Let $\sets{(x^k, y^k, z^k)}$ be generated by \eqref{eq:VFOG} using a variance-reduced estimator $\widetilde{G}y^k$ of $Gy^k$ satisfying Definition~\ref{de:VR_Estimators} and $\eta$, $\gamma_k$, and $\beta_k$ such that:
\begin{equation}\label{eq:SAEG_para_update2}
\hspace{-0.25ex}
\arraycolsep=0.2em
\begin{array}{lcl}
\frac{\revise{8}(s-1)\rho_n}{\revise{s-2}} \leq \eta < \frac{\hat{\lambda}}{L},  \quad   \gamma_k := \frac{\revise{(s-2)}\eta(k+s)}{\revise{32(s-1)}(k+s+1)}, \quad\textrm{and} \quad  \beta_k :=  \big[\frac{\revise{3}(s-2)\eta}{\revise{16}(s-1)} + 2\rho_n]\frac{ k+1 }{k+s+1} - \frac{\gamma_k}{k+s+1}.
\end{array}
\hspace{-0.25ex}
\end{equation}
Suppose further that both $\kappa_k$ and $\Theta_k$ in Definition~\ref{de:VR_Estimators} satisfy the following condition:
\begin{equation}\label{eq:SAEG_kappa_Theta_cond}
\arraycolsep=0.2em
\begin{array}{lcl}
\kappa_k \geq \frac{\eta \Gamma}{\rho_c}\Theta_k  + \frac{2}{k+s+1}. 
\end{array}
\end{equation}
Then, the following statements hold:
\vspace{0.75ex}
\begin{compactitem}
\item[$\mathrm{(a)}$]$($\textbf{The $\BigO{1/k^2}$ convergence rates}$)$ We have the following bounds:
\begin{equation}\label{eq:SAEG_BigO_rates}
\arraycolsep=0.2em\hspace{-3ex}
\begin{array}{lcl}
\Expn{ \norms{Gx^k + v^k }^2 } \leq \dfrac{\hat{C}_0^2 (\hat{\mcal{R}}_0^2 + \hat{\Lambda}_0 \Sc_k)}{ (k+s)^2} \quad \textrm{and} \quad \Expn{ \norms{Gy^{k-1} + v^k }^2 } \leq  \dfrac{2 \hat{C}_0^2 (\hat{\mcal{R}}_0^2 + \hat{\Lambda}_0 \Sc_k )}{(k+s)^2 },
\end{array}\hspace{-3ex}
\end{equation}
where $v^k \in Tx^k$, $\Sc_k := \sum_{l=0}^k \delta_l$, and $\hat{C}_0^2$ and $\hat{\mcal{R}}_0^2$ are respectively given by
\begin{equation}\label{eq:R0_quantity}
\arraycolsep=0.2em
\begin{array}{lcl}
\hat{C}_0^2 & := & \revise{\frac{128(s-1)}{19s+10}} 
\quad \text{and} \quad 
\hat{\mcal{R}}_0^2 := \revise{\frac{ 2(13s-10)s^2 - (s-1)(s-2)  }{64(s-1)}} \norms{Gx^0 + v^0 }^2 + \revise{\frac{16s^3(s-1)}{(s-2)\eta^2 }} \norms{x^0 - x^{\star} }^2.
\end{array}
\end{equation}
\item[$\mathrm{(b)}$]$($\textbf{Summability bounds}$)$
The following bounds also hold:
\begin{equation}\label{eq:SAEG_summability_bounds}
\arraycolsep=0.2em
\begin{array}{lcl}
\sum_{k=0}^K s (k+s+1) \Exp{\norms{Gy^k + v^{k+1} }^2} & \leq & 2 (\hat{\mcal{R}}_0^2 + \hat{\Lambda}_0 \Sc_K ), \vspace{1.75ex}\\
\sum_{k=0}^K \revise{\frac{(s-2)(4s-7)}{s-1}} (k+1)   \Exp{\norms{Gx^k + v^k }^2} & \leq & \revise{64} (\hat{\mcal{R}}_0^2 + \hat{\Lambda}_0 \Sc_K ).
\end{array}
\end{equation}
\item[$\mathrm{(c)}$]$($\textbf{The $\SmallO{1/k^2}$-convergence rates}$)$
If, in addition, $\Sc_{\infty} := \sum_{k=0}^\infty \delta_k < +\infty$, then we also have the following limits:
\begin{equation}\label{eq:SAEG_small_o_rates}
\arraycolsep=0.2em
\begin{array}{lcl}
\lim_{k\to\infty} k^2  \Expn{\norms{Gx^k + v^k }^2} = 0 \quad \textrm{and} \quad \lim_{k\to\infty}  k^2 \Expn{\norms{Gy^{k-1} + v^k }^2} = 0.
\end{array}
\end{equation}
\item[$\mathrm{(d)}$]$($\textbf{Iteration-complexity}$)$
For a given $\epsilon > 0$, if $\Sc_{\infty} < +\infty$, then the number of iterations $k$ to achieve $\Expn{ \norms{Gx^k + v^k}^2 } \leq \epsilon^2$ is at most $\BigOs{1/\epsilon}$.
This  is also the number of evaluations of $J_{\eta T}$. 

\item[$\mathrm{(e)}$]$($\textbf{Almost sure convergence}$)$
The conclusions of Theorem~\ref{th:VFOG1_almost_sure_convergence} still hold in this case.
\end{compactitem}
\end{theorem}

The constant $s$ in Theorem~\ref{th:VFOG_convergence} was left free, provided that $s > 2$. 
It provides flexibility in choosing its value for adjusting other parameters and theoretical bounds.
For example, if we assume that $\alpha = 0$ in \eqref{eq:G_Lipschitz} and choose $s = \revise{3}$, then $\hat{\omega} = \revise{2.44336}$.
Consequently, we have  $\hat{\lambda} := \frac{1}{\sqrt{\revise{27.54687}}} \approx \revise{0.19053}$ and $\hat{\mu} :=  \frac{\revise{(s-2)}\hat{\lambda}}{\revise{8}(s-1)} \approx \revise{0.01191}$.
In this case, the conditions $L\rho_n \leq \hat{\mu}$ and $\revise{16}\rho_n \leq \eta \leq \frac{\hat{\lambda}}{L}$ in Theorem~\ref{th:VFOG_convergence} hold.
\revise{The supremum of the admissible stepsize is attained asymptotically as $s\downarrow 2$, in which case $\eta \lesssim  \frac{0.27123}{L}$.}
Note that these parameter bounds are not tight, and their range could potentially be improved by further refining our analysis.

Theorem~\ref{th:VFOG_convergence}  only establishes convergence rates for \eqref{eq:VFOG}, but does not provide its oracle complexity.
It also requires the condition \eqref{eq:SAEG_kappa_Theta_cond} to hold for $\widetilde{G}y^k$ in Definition~\ref{de:VR_Estimators}.
The iteration-complexity of \eqref{eq:VFOG} for achieving an $\epsilon$-solution $x^k$ in expectation, i.e., $\Expn{\norms{Gx^k + v^k}^2 } \leq \epsilon^2$ is $\BigOs{1/\epsilon}$, improving by a factor of $\frac{1}{\epsilon}$ over  its non-accelerated counterparts, including \cite{TranDinh2024}.

\beforesubsec
\subsection{\textbf{Oracle complexity of \eqref{eq:VFOG} using L-SVRG, SAGA, and L-SARAH}.}\label{subsec:VFOG_oracle_complexity}
\aftersubsec
In this subsection, we specify the oracle complexity of \eqref{eq:VFOG} using three concrete estimators: L-SVRG, SAGA, and L-SARAH  in Subsection \ref{subsec:Vr_Estimator}.
The key point is to choose an appropriate distribution $\mbf{p}_k$ and a mini-batch size $b_k$ to fulfill the condition \eqref{eq:SAEG_kappa_Theta_cond}.
Note that from Lemmas \ref{le:loopless_svrg_bound}, \ref{le:SAGA_estimator_full}, and \ref{le:loopless_sarah_bound}, since $\delta_k = 0$, we always have $\Sc_{\infty} := \sum_{k=0}^\infty \delta_k = 0$ for these three estimators.
Now, let us start with the Loopless SVRG estimator in \eqref{eq:loopless_svrg}.

\begin{corollary}\label{co:SVRG_complexity}
Suppose that Assumptions~\ref{as:A1}-\ref{as:A3} hold for \eqref{eq:GE} under the finite-sum setting \eqref{eq:finite_sum}.
Let $\sets{(x^k, y^k, z^k)}$ be generated by \eqref{eq:VFOG} using the SVRG estimator $\widetilde{G}y^k$ from \eqref{eq:loopless_svrg} for $Gy^k$ and the parameters as in Theorem~\ref{th:VFOG_convergence}.
For given $\nu \in (0, 1/2)$ and $c_2 > 0$, if  $n \geq \max\big\{ \big[\frac{c_2(s+2)}{s-2}\big]^{1/\nu}, \ \big[\frac{8\Gamma \eta}{c_2^2 \rho_c }\big]^{1/(1-2\nu)}\big\}$ for $\Gamma$ in \eqref{eq:omega_quantity}, and  $b_k$ and $\mbf{p}_k$ in \eqref{eq:xy_hat} are chosen such that
\begin{equation}\label{eq:SVRG_choice_of_params}
\arraycolsep=0.2em
\begin{array}{lcl}
b_k = b :=  \big\lceil \frac{8\Gamma \eta n^{2\nu}}{ c_2^2 \rho_c } \big\rceil  \quad\textrm{and} \quad \mbf{p}_k := \begin{cases}
\frac{c_2}{n^{\nu}} + \frac{4}{k+s + 1} & \textrm{if}~k \leq \lfloor \frac{4n^{\nu}}{c_2} - s\rfloor, \vspace{1ex}\\
\frac{2c_2}{n^{\nu}} & \textrm{otherwise}.
\end{cases}
\end{array}
\end{equation}
Then, we have $1 \leq b_k \leq n$ and $\mbf{p}_k \in (0, 1]$.
For any $\epsilon > 0$, the expected total number of oracle calls  $\mbf{G}$ for $k\geq 0$ in \eqref{eq:VFOG} is at most 
\begin{equation}\label{eq:SVRG_complexity_bound}
\arraycolsep=0.2em
\begin{array}{lcl}
\Expn{\Tc_K} & = & \BigO{ n\ln(n^{\nu}) + n^{3\nu} + \frac{\mcal{R}_0(n^{1-\nu} + n^{2\nu})}{\epsilon}}
\end{array}
\end{equation}
to achieve $\Expn{\norms{Gx^K + v^K}^2} \leq \epsilon^2$.
If we choose $\nu = \frac{1}{3}$, then $\Expn{\revise{\Tc_K}} = \BigOs{ n \ln(n^{1/3}) + \frac{n^{2/3}}{\epsilon}}$.
\end{corollary}

\proof{\textbf{Proof.}}
Since $\widetilde{G}y^k$ is updated by  \eqref{eq:loopless_svrg}, 
by Lemma~\ref{le:loopless_svrg_bound}, we have $\kappa_k = \frac{\mbf{p}_k}{2}$ and $\Theta_k = \frac{4}{b_k\mbf{p}_k}$.
For $\Gamma$ given in \eqref{eq:omega_quantity}, since $\kappa_k = \frac{\mbf{p}_k}{2}$ and $\Theta_k = \frac{4}{b_k\mbf{p}_k}$, the condition \eqref{eq:SAEG_kappa_Theta_cond} holds if
\begin{equation*}
\arraycolsep=0.2em
\begin{array}{lcl}
\mbf{p}_k & \geq & \frac{8\Gamma \eta}{\rho_c b_k\mbf{p}_k} + \frac{4}{k+s+1}.
\end{array}
\end{equation*}
Let $\nu \in [0, 1]$, $c_1 > 0$, and $c_2 > 0$ be given.
If we choose $\mbf{p}_k \geq \frac{c_2}{n^{\nu}}$ and $b_k \geq c_1n^{2\nu}$, then  the last condition holds if 
\begin{equation}\label{eq:SAEG_co31_proof1}
\arraycolsep=0.2em
\begin{array}{lcl}
\mbf{p}_k & \geq &  \frac{8\Gamma\eta }{\rho_c\mred{c_1c_2n^{\nu}}} + \frac{4}{k+s+1} \geq  \frac{8\Gamma\eta }{\rho_c\mred{c_1c_2n^{\nu}}} =: \underline{\mbf{p}} > 0.
\end{array}
\end{equation}
Therefore, if we choose $c_1 = \frac{8\Gamma\eta}{\rho_c\mred{c_2^2}}$, then $ \frac{8\Gamma\eta }{\rho_c \mred{c_1c_2n^{\nu}}} = \frac{c_2}{n^{\nu}}$.
Moreover, if  $k \geq \frac{4n^{\nu}}{c_2} - s$, then $\frac{4}{k+s+1} \leq \frac{c_2}{n^{\nu}}$.
Overall, if we choose $\mbf{p}_k$ as in \eqref{eq:SVRG_choice_of_params}, then it satisfies \eqref{eq:SAEG_co31_proof1}.

On the other hand, to guarantee $\mbf{p}_k \leq 1$, we require $\frac{2c_2}{n^{\nu}} \leq 1$ and $\frac{c_2}{n^{\nu}} + \frac{4}{s+2} \leq 1$, which are guaranteed if $n \geq \big[\frac{c_2(s+2)}{s-2}\big]^{1/\nu}$.
In addition, we need $\frac{8\Gamma \eta}{\rho_c\mred{c_2^2}}n^{2\nu} = c_1n^{2\nu} \leq b \leq n$, which holds if $n \geq \big[\frac{8\Gamma}{\mred{c_2^2}}\big]^{1/(1-2\nu)}$.
Combining these two conditions, we get $n \geq \max\big\{ \big[\frac{c_2(s+2)}{s-2}\big]^{1/\nu}, \ \big[\frac{8\Gamma \eta}{\rho_c\mred{c_2^2}}\big]^{1/(1-2\nu)}\big\}$. 
In this case, we can choose $b := \big\lceil \frac{8\Gamma \eta n^{2\nu}}{\rho_c \mred{c_2^2}} \big\rceil$ as in \eqref{eq:SVRG_choice_of_params}.
We can also easily verify that $b \geq 1$.

Now, from \eqref{eq:SAEG_BigO_rates} of Theorem~\ref{th:VFOG_convergence} and $S_K = 0$, to guarantee $\Expn{\norms{Gx^K + v^K}^2} \leq \epsilon^2$, we impose $\frac{\mred{\hat{C}_0^2\hat{\mcal{R}}_0^2}}{(K+s)^2} \leq \epsilon^2$ for \revise{$\hat{C}_0^2  := \frac{128(s-1)}{(19s+10)}$}.
This leads to the choice $K := \big\lfloor \frac{ \mred{\hat{C}_0\hat{\mcal{R}}_0} }{\epsilon} \big\rfloor$.

The expected total number of evaluations of $\mbf{G}(y^k, \xi_k)$ is $\Expn{ \Tc_K } = n + \Expn{ \hat{\Tc}_K }$, where $n$ is the total evaluations of the first epoch ($k=0$), and $\hat{\Tc}_K$ is the evaluations of the next $K$ iterations.
Using \eqref{eq:SVRG_choice_of_params}, we can evaluate the expected total number  $\Expn{ \hat{\Tc}_K }$  of oracle calls $\mbf{G}$ as follows:
\myeqn{
\arraycolsep=0.2em
\begin{array}{lcl}
\Expn{\hat{\Tc}_K } & := &  \sum_{k=1}^K[n\mbf{p}_k + 2(1-\mbf{p}_k)b] \leq \sum_{k=1}^K(n\mbf{p}_k + 2b) \vspace{1ex}\\
& \leq & \sum_{k=1}^{\frac{4n^{\nu}}{c_2} - s}\big(c_2n^{1-\nu} + \frac{4n}{k+s+1} + 2c_1n^{2\nu} \big) + \big(K - \frac{ 4n^{\nu}}{c_2} + s)(2c_2n^{1-\nu} + 2c_1n^{2\nu}) \vspace{1ex}\\
& \leq & 4n + 4n\ln\big( \frac{4n^{\nu}}{c_2} + 2\big) + c_1n^{3\nu} + \frac{(2c_2n^{1-\nu} +  2c_1n^{2\nu})\mred{\hat{C}_0\hat{\mcal{R}}_0}}{\epsilon} \vspace{1ex}\\
& = & \BigO{ n\ln(n^{\nu}) + n^{3\nu} + \frac{\mred{\hat{\mcal{R}}_0}(n^{1-\nu} + n^{2\nu})}{\epsilon}}.
\end{array}
}
This expression, together with $\Expn{\Tc_K} = n + \Expn{ \hat{\Tc}_K}$, leads to \eqref{eq:SVRG_complexity_bound}.
Clearly, if $\nu = \frac{1}{3}$, then we get $\Expn{\Tc_K} = \BigOs{n\ln(n^{1/3}) + \frac{n^{2/3}}{\epsilon}}$.
\Eproof
\endproof

Our next result is the complexity of \eqref{eq:VFOG} using the SAGA estimator \eqref{eq:SAGA_estimator}.

\begin{corollary}\label{co:SAGA_complexity}
Suppose that Assumptions~\ref{as:A1}-\ref{as:A3} hold for \eqref{eq:GE} under the finite-sum setting \eqref{eq:finite_sum}.
Let $\sets{(x^k, y^k, z^k)}$ be generated by \eqref{eq:VFOG} using the SAGA estimator $\widetilde{G}y^k$ from \eqref{eq:SAGA_estimator} for $Gy^k$ and the parameters as in Theorem~\ref{th:VFOG_convergence}.
Suppose that $n \geq \frac{10\Gamma\eta}{\rho_c}\big(\frac{s+2}{s-2}\big)^3$ and we choose $b_k$  in \eqref{eq:SAGA_estimator} as 
\begin{equation}\label{eq:SAGA_choice_of_params}
\arraycolsep=0.2em
\begin{array}{lcl}
b_k := \begin{cases} 
\big\lceil \sqrt[3]{\frac{10\Gamma\eta}{\rho_c}} n^{2/3} + \frac{4n }{k+s + 1} \big\rceil & \textrm{if}~k \leq \big\lfloor \frac{4\sqrt[3]{\rho_c}}{\sqrt[3]{10\Gamma\eta}}n^{1/3}  - s \big\rfloor, \vspace{1ex}\\
\big\lceil 2\sqrt[3]{\frac{10\Gamma \eta}{\rho_c}} n^{2/3} \big\rceil & \textrm{otherwise}.
\end{cases}
\end{array}
\end{equation}
Then, we have $1 \leq b_k \leq n$, and  for any $\epsilon > 0$, the expected total number of oracle calls $\mbf{G}$ in \eqref{eq:VFOG} is at most 
\begin{equation}\label{eq:SAGA_complexity_bound}
\arraycolsep=0.2em
\begin{array}{lcl}
\Expn{\Tc_K} & = & \BigO{ n\ln(n^{1/3}) + \frac{\mcal{R}_0 n^{2/3} }{ \epsilon} }
\end{array}
\end{equation}
to achieve $\Expn{\norms{Gx^K + v^K }^2} \leq \epsilon^2$.
\end{corollary}

\proof{\textbf{Proof.}}
As indicated by Lemma~\ref{le:SAGA_estimator_full}, $\widetilde{G}y^k$ generated by \eqref{eq:SAGA_estimator} satisfies Definition~\ref{de:VR_Estimators} with $\kappa_k := \frac{b_k}{2n}$ and $\Theta_k :=  \frac{5n}{b_k^2}$.
For given $\nu \in [0, 1/2]$ and $c_1 > 0$, let us choose $b_k$ as
\begin{equation}\label{eq:SAEG_saga_proof1}
\arraycolsep=0.2em
\begin{array}{lcl}
b_k & = & \begin{cases}
c_1 n^{2\nu} + \frac{4n}{k+s+1} & \textrm{if $k \leq \lfloor \frac{4n^{1-2\nu}}{c_1} - s -1 \rfloor$}, \\
2c_1 n^{2\nu} & \textrm{otherwise}.
\end{cases}
\end{array}
\end{equation}
Since $\kappa_k := \frac{b_k}{2n}$ and $\Theta_k := \frac{5n}{b_k^2}$, the condition  \eqref{eq:SAEG_kappa_Theta_cond} holds if $b_k \geq \frac{10n^2\Gamma \eta}{\rho_c b_k^2} + \frac{4n}{k+s+1}$.
Because $b_k$ is updated by \eqref{eq:SAEG_saga_proof1}, the last condition is satisfied if $b_k \geq \frac{ 10\Gamma \eta n^{2(1-2\nu)}}{\rho_c c_1^2} + \frac{4n}{k+s+1}$, which automatically holds if we impose $\frac{ 10n^{2(1-2\nu)}\Gamma \eta}{\rho_c c_1^2} = c_1n^{2\nu}$.
The last equality holds if $\nu = \frac{1}{3}$ and $c_1 := \sqrt[3]{\frac{10\Gamma\eta}{\rho_c}}$.
Substituting these values into \eqref{eq:SAEG_saga_proof1}, we obtain $b_k$ as in \eqref{eq:SAGA_choice_of_params}.
From \eqref{eq:SAEG_kappa_Theta_cond}, we can observe that $\Gamma\eta > \rho_c$, leading to $b \geq 1$.

Since $b_k \leq n$, we need to impose $\sqrt[3]{\frac{10\Gamma\eta}{\rho_c}} n^{2/3} \leq \frac{n(s-2)}{s+2}$, leading to $n \geq \frac{ 10\Gamma\eta}{\rho_c}\big(\frac{s+2}{s-2}\big)^3$.
In addition, for $b_k$ in \eqref{eq:SAGA_choice_of_params}, we can easily check that it satisfies $b_{k-1} - \frac{b_kb_{k-1}}{4n} \leq b_k \leq b_{k-1}$ in Lemma~\ref{le:SAGA_estimator_full}.

Now, for a given tolerance $\epsilon > 0$, from \eqref{eq:SAEG_BigO_rates} of Theorem~\ref{th:VFOG_convergence}, to guarantee that $\Expn{\norms{Gx^K + v^K }^2} \leq \epsilon^2$, we need to impose $\frac{\mred{\hat{C}_0^2\hat{\mcal{R}}_0^2}}{(K+s)^2} \leq \epsilon^2$ for \revise{$\hat{C}_0^2  := \frac{128(s-1)}{(19s+10)}$}.
This requirement leads to the choice $K := \big\lfloor \frac{ \mred{\hat{C}_0 \hat{\mcal{R}}_0} }{\epsilon} \big\rfloor$.

Let us estimate the expected total number  $\Expn{\Tc_K}$ of oracle calls $\mbf{G}$ for \eqref{eq:VFOG}, which is $\Expn{\Tc_K} = n + \Expn{\hat{\Tc}_K}$, where
\begin{equation*} 
\arraycolsep=0.2em
\begin{array}{lcl}
\Expn{ \hat{\Tc}_K } &= & 2\sum_{k=1}^Kb_k = \sum_{k=1}^{ \frac{4n^{1/3}}{c_1} - s}\big(c_1n^{2/3} + \frac{4n}{k+s+1}\big) + 4c_1n^{2/3}\big(K - \frac{4n^{1/3}}{c_1} + s \big) \vspace{1ex}\\
& \leq &  4n + 4n\ln\big( \frac{4n^{1/3}}{c_1} + 1 \big) + \frac{4c_1 \mred{\hat{C}_0 \hat{\mcal{R}}_0} n^{2/3}}{\epsilon} \vspace{1ex}\\
& = & \BigO{ n \ln(n^{1/3}) + \frac{\mred{\hat{\mcal{R}}_0}n^{2/3}}{\epsilon}}.
\end{array}
\end{equation*}
This implies that $\Expn{\Tc_K} = n + \Expn{ \hat{\Tc}_K} \leq \BigOs{ n \ln(n^{1/3}) + \frac{\mred{\hat{\mcal{R}}_0} n^{2/3}}{\epsilon}}$ as stated in \eqref{eq:SAGA_complexity_bound}.
\Eproof
\endproof

Finally, we estimate the complexity of \eqref{eq:VFOG} using the  \ref{eq:loopless_sarah_estimator} estimator.

\begin{corollary}\label{co:SARAH_complexity}
Suppose that Assumptions~\ref{as:A1}-\ref{as:A3} hold for \eqref{eq:GE} under the finite-sum setting \eqref{eq:finite_sum}.
Let $\sets{(x^k, y^k, z^k)}$ be generated by \eqref{eq:VFOG} using the SARAH estimator $\widetilde{G}y^k$ from \eqref{eq:loopless_sarah_estimator} for $Gy^k$ and the parameters as in Theorem~\ref{th:VFOG_convergence}.
For given $\nu \in (0, 1]$ and $c_1 \in (0, n^{1-\nu}]$, if $n \geq \big[\frac{2\Gamma\eta}{c_1\rho_c}\big]^{1/\nu}$ for $\Gamma$ given in \eqref{eq:omega_quantity},  and  $b$ and $\mbf{p}_k$ in \eqref{eq:loopless_sarah_estimator} are chosen as 
\begin{equation}\label{eq:SARAH_choice_of_params}
\arraycolsep=0.2em
\begin{array}{lcl}
b_k = b :=  \lceil c_1 n^{\nu} \rceil   \quad\textrm{and} \quad \mbf{p}_k := \begin{cases}
\frac{\Gamma \eta}{c_1 \rho_c n^{\nu}} + \frac{2}{k+s + 1} & \textrm{if}~k \leq \lfloor \frac{2c_1\rho_c n^{\nu}}{\Gamma \eta} - s\rfloor, \vspace{1ex}\\
\frac{2\Gamma \eta}{c_1\rho_c n^{\nu}} & \textrm{otherwise}.
\end{cases}
\end{array}
\end{equation}
Then, we have $1 \leq b_k \leq n$ and $\mbf{p}_k \in (0, 1]$.
For any given $\epsilon > 0$, the expected total number of oracle calls $\mbf{G}$ is at most 
\begin{equation}\label{eq:SARAH_complexity_bound}
\arraycolsep=0.2em
\begin{array}{lcl}
\Expn{\Tc_K} & = & \BigO{ n \ln(n^{\nu}) + n^{2\nu} + \frac{\mcal{R}_0(n^{1-\nu} + n^{\nu})}{\epsilon}}
\end{array}
\end{equation}
to achieve $\Expn{\norms{Gx^K + v^K }^2} \leq \epsilon^2$.
If we choose $\nu = \frac{1}{2}$, then $\Expn{\Tc_K} = \BigOs{ n\ln(n^{1/2}) + \frac{\sqrt{n}}{\epsilon}}$.
\end{corollary}

\proof{\textbf{Proof.}}
By Lemma~\ref{le:loopless_sarah_bound}, $\widetilde{G}y^k$ generated by \eqref{eq:loopless_sarah_estimator} satisfies Definition~\ref{de:VR_Estimators} with $\kappa_k = \mbf{p}_k$ and $\Theta_k = \frac{1}{b_k}$.
Using these values, the condition \eqref{eq:SAEG_kappa_Theta_cond} holds if
\begin{equation*} 
\arraycolsep=0.2em
\begin{array}{lcl}
\mbf{p}_k \geq \frac{\Gamma \eta }{\rho_c b_k} + \frac{2}{k+s+1}.
\end{array}
\end{equation*}
Let us choose $b := \lceil c_1n^{\nu} \rceil$ as in \eqref{eq:SARAH_choice_of_params} for some $\nu \in (0, 1)$ and $c_1 > 0$.
Then, the last condition holds if we choose $\mbf{p}_k$ as in \eqref{eq:SARAH_choice_of_params}.
Moreover, to ensure $0 < \mbf{p}_k \leq 1$, we need to impose $b_k \geq c_1n^{\nu} \geq \frac{ 2\Gamma \eta}{\rho_c} \geq 1$, leading \mred{to} $n \geq \big[\frac{2\Gamma\eta}{c_1\rho_c}\big]^{1/\nu}$ as stated.
In addition, we also have $b_k \leq n$.

Now, from \eqref{eq:SAEG_BigO_rates} of Theorem~\ref{th:VFOG_convergence}, to guarantee that $\Expn{\norms{Gx^K + v^K }^2} \leq \epsilon^2$, we need to impose $\frac{\mred{\hat{C}_0^2\hat{\mcal{R}}_0^2}}{(K + s)^2} \leq \epsilon^2$ for \revise{$\hat{C}_0^2  := \frac{128(s-1)}{(19s+10)}$}.
This requirement leads to the choice of $K := \big\lfloor \frac{ \mred{\hat{C}_0\hat{\mcal{R}}_0} }{\epsilon} \big\rfloor$.
In this case, the expected total number of oracle calls $\mbf{G}$ is $\Expn{ \Tc_K } = n + \Expn{ \hat{\Tc}_K }$, where $n$ is the number of evaluations in the first epoch $k=0$, and $\hat{\Tc}_K$ is the number of evaluations over the next $K$ iterations.
Then, we can evaluate $\Expn{ \hat{\Tc}_K }$ as 
\begin{equation*}
\arraycolsep=0.2em
\begin{array}{lcl}
\Expn{ \hat{\Tc}_K } & := &  \sum_{k=1}^K[n\mbf{p}_k + 2(1-\mbf{p}_k)b] \leq \sum_{k=1}^K(n\mbf{p}_k + 2b) \vspace{1ex}\\
& \leq & \sum_{k=1}^{\frac{2c_1n^{\nu}}{\Gamma} - s}\big( \frac{\Gamma n^{1-\nu}}{c_1} + \frac{2n}{k+s+1} + 2c_1n^{\nu} \big) + \big( K - \frac{2c_1n^{\nu}}{\Gamma} \big)\big( \frac{2\Gamma n^{1-\nu}}{c_1}  + 2c_1n^{\nu} \big) \vspace{1ex}\\
& \leq & 2n + 2n \ln\big( \frac{2c_1n^{\nu}}{\Gamma} + 2 \big) + \frac{4c_1^2n^{2\nu}}{\Gamma} + \frac{\mred{\hat{C}_0\hat{\mcal{R}}_0}}{\epsilon}\big( \frac{2 n^{1-\nu}}{\Gamma} + 2c_1n^{\nu} \big) \vspace{1ex}\\
& = & \BigO{ n \ln( n^{\nu} ) + n^{2\nu} + \frac{\mred{\hat{\mcal{R}}_0}(n^{1-\nu} + n^{\nu})}{\epsilon}}.
\end{array}
\end{equation*}
This expression together with $\Expn{\Tc_K} = n + \Expn{ \hat{\Tc}_K }$ leads to \eqref{eq:SARAH_complexity_bound}.
Finally, if we choose $\nu = \frac{1}{2}$, then $\Expn{\Tc_K} = \BigOs{n\ln(n^{1/2}) + \frac{n^{1/2}}{\epsilon}}$.
\Eproof
\endproof

\begin{remark}
\label{re:choice_of_estimator_params}
Note that both $\mbf{p}_k$ and $b_k$ selected in Corollaries \ref{co:SVRG_complexity}, \ref{co:SAGA_complexity}, and \ref{co:SARAH_complexity} aim at achieving the best oracle complexity.
Other choices of these parameters still work as long as they satisfy the condition $\kappa_k \geq \frac{\Gamma\eta}{\rho_c}\Theta_k + \frac{2}{k+s+1}$ in Theorem~\ref{th:VFOG_convergence}.
For instance, the probability $\mbf{p}_k$ in Corollary \ref{co:SARAH_complexity} can be fixed at $\mbf{p}_k := \frac{\Gamma \eta}{c_1\rho_c \sqrt{n}} + \frac{2}{s+1}$ for all $k \geq 1$ and $s > 2$.
This choice is more convenient for the implementation of \eqref{eq:VFOG} in practice. 
\end{remark}

\begin{remark}
\label{re:comparison}
Compared with \cite{TranDinh2024}, Assumption~\ref{as:A3} is more restrictive than the weak-Minty condition in  \cite{TranDinh2024}.
However, \eqref{eq:VFOG} achieves a better convergence rate by a factor of $\frac{1}{k}$ and better oracle complexity by a factor $\frac{1}{\epsilon}$.
Moreover, \eqref{eq:VFOG} also works with biased estimators.
For SARAH, it improves by a factor $n^{1/6}$ over SVRG and SAGA.
Compared to \cite{alacaoglu2021stochastic,alacaoglu2021forward,cai2022stochastic}, \eqref{eq:VFOG} can work with a class of nonmonotone operators, and its convergence guarantee is in $\Expn{\norms{Gx^K + v^K}^2}$,
whereas the gap function guarantees in  \cite{alacaoglu2021stochastic,alacaoglu2021forward} apply only to the monotone case.
Therefore, a direct comparison of their oracle complexities is not meaningful because the convergence metrics differ.
Nevertheless, since \revise{the} methods in \cite{alacaoglu2021stochastic,alacaoglu2021forward} are non-accelerated, we believe that their convergence rate is at most $\BigOs{1/k}$, leading to an oracle complexity worse by a factor of $\frac{1}{\epsilon}$ than ours.
Other works such as \cite{davis2022variance,tran2024accelerated,TranDinh2025a} require stronger assumptions such as [star] co-coercivity and strong quasi-monotonicity, excluding our problems of interest.
\end{remark}

\begin{remark}\label{re:expected_case}
We only consider the finite-sum case \eqref{eq:finite_sum} in Corollaries \ref{co:SVRG_complexity}, \ref{co:SAGA_complexity}, and \ref{co:SARAH_complexity}.
If we consider the general case $G(\cdot) = \Expsn{\xi}{ \mbf{G}(\cdot, \xi)  }$, then we can also prove that the expected total number of oracle calls $\Expn{\Tc_K}$ is at most $\widetilde{\mathcal{O}} (\epsilon^{-3})$ to achieve an $\epsilon$-solution.
However, we skip this analysis.
\end{remark}
\beforesec
\section{The Detailed Proofs of Theorems~\ref{th:VFOG1_convergence}, \ref{th:VFOG1_almost_sure_convergence}, and \ref{th:VFOG_convergence}.}\label{sec:convergence_analysis}
\aftersec
For readability, we present the technical proofs of Theorems \ref{th:VFOG1_convergence}, \ref{th:VFOG1_almost_sure_convergence}, and \ref{th:VFOG_convergence} in this section.

\beforesubsec
\subsection{\textbf{The full proof of Theorem \ref{th:VFOG1_convergence}.}}\label{subsec:th:VFOG1_convergence}
\aftersubsec
The proof of Theorem~\ref{th:VFOG1_convergence} is divided into several steps.

\beforesubsubsec
\subsubsection{\textbf{Key technical lemmas for convergence analysis.}}
\label{subsec:VFOG_key_estimate}
\aftersubsubsec
For given $x^k$, $y^{k-1}$, $y^k$, and $v^k \in Tx^k$, we introduce the following quantities of our analysis:
\begin{equation}\label{eq:SFOG_quantities}
\arraycolsep=0.2em
\begin{array}{lclclcl}
\hat{w}^k & := & Gy^{k-1} + v^k \quad & \textrm{and} & \quad w^k & := &  Gx^k + v^k, \vspace{1ex}\\
e^k & := &  \widetilde{G}y^k - Gy^k \quad & \textrm{and} & \quad \hat{e}^k & := &  \widetilde{G}y^{k-1} - Gx^k. 
\end{array}
\end{equation}
Our first step is to establish the following key lemma, whose proof is in Appendix~\ref{apdx:le:VFOG_key_estimate1}.

\begin{lemma}\label{le:VFOG_key_estimate1}
Assume that Assumptions~\ref{as:A1} and \ref{as:A3} hold for \eqref{eq:GE}.
Let $\sets{(x^k, y^k, z^k)}$ be updated by \eqref{eq:VFOG} using the following parameters:
\begin{equation}\label{eq:SFOG_para_cond1}
\arraycolsep=0.2em
\begin{array}{lcl}
t_{k+1} := t_k + 1, \quad \gamma_k := \frac{\gamma(t_k - 1)}{t_k}, \quad \text{and} \quad \beta_k := \frac{(c_1\eta + c_2\eta + 2\rho_n)(t_k-s) - \gamma_k}{t_k},
\end{array}
\end{equation}
for some constants $\gamma ~\revise{>}~ 0$,  $c_1 > 0$, and $c_2 > 0$, determined later, \revise{and $t_0 \geq s + 1$ for a given $s > 2$}.
Let $w^k$, $\hat{w}^k$, $e^k$, and $\hat{e}^k$ be defined by \eqref{eq:SFOG_quantities}, $v^k \in Tx^k$, and $\Pc_k$ be defined by 
\begin{equation}\label{eq:SFOG_Ek_func}
\arraycolsep=0.2em
\begin{array}{lcl}
\Pc_k := \frac{a_k}{2}\norms{w^k}^2 + st_{k-1}\iprods{w^k, x^k - z^k} + \frac{s^2(s-1)}{2\gamma_{\revise{k-1}}}\norms{z^k - x^{\star}}^2.
\end{array}
\end{equation}
Then, for any constant $\omega \geq 0$ and \revise{$\theta > 0$}, we have
\begin{equation}\label{eq:SFOG_key_est1}
\hspace{-4ex}
\arraycolsep=0.2em
\begin{array}{lcl}
\Pc_k - \Pc_{k+1}  & \geq &  \frac{\eta t_k^2}{2}\big[  (1-\alpha) \big( \tnorms{ w^{k+1} - \hat{w}^{k+1}  }^2 - \norms{w^{k+1} - \hat{w}^{k+1} }^2 \big) + \omega L^2\norms{x^{k+1} - y^k}^2 \big]  \vspace{1ex}\\
&& + {~} \rho_ct_k(t_k-s) \tnorms{ Gx^{k+1} - Gx^k }^2 + \frac{\eta t_k(t_k-s - M^2\eta^2t_k)}{2}\norms{\hat{w}^{k+1} - w^k}^2  \vspace{1ex}\\
&& + {~} \frac{\eta st_k}{2} \norms{ \hat{w}^{k+1} }^2 + \frac{a_k - \hat{a}_k}{2} \norms{ w^k }^2 +   s(s-1)\iprods{w^k, x^k - x^{\star}} \vspace{1ex}\\
&&  - {~} b_k \norms{e^k}^2 - \hat{b}_k \norms{\hat{e}^k}^2,
\end{array}
\hspace{-4ex}
\end{equation}
where $M^2 := 2(1+\omega)L^2$, \revise{$\psi := (1-c_1 - c_2)\eta - 2\rho_n$}, and 
\begin{equation}\label{eq:SFOG_coeffs}
\hspace{-0ex}
\arraycolsep=0.1em
\left\{\begin{array}{lcl}
a_{k+1} &:= &  \psi t_k^2 + 2s\rho_n t_k - \gamma(t_k - 1), \vspace{1ex}\\
\hat{a}_k &:= & \psi(t_k - s)^2  + s(\eta + 2\rho_n)(t_k-s)  + \frac{2\gamma\mred{(t_k-1)^2}}{t_k}, \vspace{1ex}\\ 
b_k &:= & \big(\frac{1}{2c_1} +  \frac{\revise{(1+\theta)}M^2\eta^2}{\revise{2}} \big)\eta t_k^2, \vspace{1ex}\\
\hat{b}_k &:= & \big( \frac{\beta_k^2t_k}{2 c_2\eta} + \frac{\gamma(t_k - 1)}{2 t_k} +  \frac{\revise{(1+\theta^{-1})}M^2\eta^3t_k}{\revise{2}} \big)t_k + \frac{\gamma(s-1)(t_k-1)}{2}\big( \frac{2}{t_k} + t_k - 2  \big).
\end{array}\right.
\hspace{-0ex}
\end{equation}
\end{lemma}

Next, we can further simplify Lemma~\ref{le:VFOG_key_estimate1} by choosing concrete values for the constants $c_1$ and $c_2$ to achieve a simplified bound for $\Pc_k$ as stated in the following lemma.

\begin{lemma}\label{le:VFOG_key_estimate2}
Under the same conditions and parameters as in Lemma~\ref{le:VFOG_key_estimate1}, suppose that 
\begin{equation}\label{eq:SFOG_para_cond2}
\hspace{-1ex}
\begin{array}{lcl}
s \ \revise{>} \ 2, \qquad  t_k \geq s + 1, \qquad \frac{M^2\eta^2t_k}{t_k-s} \leq 1, \qquad 0 < \gamma \leq \revise{\frac{9}{16}}\eta, \quad \textrm{and} \quad \revise{(s-2)}\eta \geq \revise{8(s-1)}\rho_n.
\end{array}
\hspace{-3ex}
\end{equation}
Then, $\beta_k$ in \eqref{eq:SFOG_para_cond1} satisfies $\beta_k  \leq \big[\frac{\revise{3}(s-2)\eta}{\revise{16}(s-1)} + 2\rho_n]\frac{(t_k-s)}{t_k}$ and $\Pc_k$ defined by \eqref{eq:SFOG_Ek_func} satisfies
\begin{equation}\label{eq:SFOG_key_est2}
\hspace{-0.5ex}
\arraycolsep=0.2em
\begin{array}{lcl}
\Pc_k - \Pc_{k+1}  & \geq &  \frac{  \eta t_k^2}{2} \big[  (1-\alpha) \big( \tnorms{w^{k+1} - \hat{w}^{k+1} }^2  - \norms{w^{k+1} - \hat{w}^{k+1} }^2 \big) + \omega L^2\norms{x^{k+1} - y^k }^2 \big] \vspace{1ex}\\
&& + {~}  \rho_c t_k(t_k-s)\tnorms{ Gx^{k+1} - Gx^k}^2 +  \frac{\eta t_k(t_k-s)}{2}\big(1  - \frac{M^2\eta^2 t_k}{t_k - s} \big) \norms{ \hat{w}^{k+1} - w^k}^2 \vspace{1ex}\\
&& + {~}   \frac{\eta st_k}{2} \norms{ \hat{w}^{k+1} }^2   + \frac{\phi_k}{2} \norms{ w^k }^2 +  s(s-1)\iprods{ w^k, x^k - x^{\star}}   \vspace{1ex}\\
&& - {~}  \frac{\varphi \eta t_{k-1}^2}{\revise{8}}\norms{ w^k - \hat{w}^k }^2 - \frac{ \revise{(25s-34)}\eta t_k^2}{\revise{2}(s-2)}\norms{e^k}^2 - \varphi \eta t_{k-1}^2 \norms{e^{k-1}}^2,
\end{array}
\hspace{-2ex}
\end{equation}
where
\begin{equation}\label{eq:SFOG_phik_and_varphi} 
\arraycolsep=0.2em
\left\{\begin{array}{lcl}
\phi_k & := & \big[ \revise{\frac{5(s-2)}{8}}\eta - 4(s-1)\rho_n - 3\gamma \big](t_k - 1) + \frac{(s-1)\revise{(3s+10)}\eta}{\revise{16}} + 2 (2s^2 - 3s + 1) \rho_n + \gamma,  \vspace{1ex}\\
\varphi & := & \revise{\frac{9(85s - 134)}{64(s-1)}} + \frac{\revise{9}(s-1)\gamma}{\revise{2}\eta}.
\end{array}\right.
\end{equation}
\end{lemma}

\begin{remark}\label{re:choice_of_parameters}
Lemma~\ref{le:VFOG_key_estimate2} only provides one choice for the constants $c_1$ and $c_2$.
This choice affects the ranges of parameters $\eta$, $\beta_k$, and $\gamma_k$ in \eqref{eq:VFOG} as well as the range of $L\rho_n$.
In this paper, we do not try to optimize the choice of these constants.
We only provide one instance to simplify our convergence analysis.
However, other choices are also possible.
\end{remark}

Finally, for $\Pc_k$ defined by \eqref{eq:SFOG_Ek_func}, we can lower bound it as in the following lemma.

\begin{lemma}\label{le:Ek_lowerbound}
Under the same settings as in Lemma~\ref{le:VFOG_key_estimate2} and Assumption~\ref{as:A2}, $\Pc_k$ defined by  \eqref{eq:SFOG_Ek_func} is lower bounded by
\begin{equation}\label{eq:Ek_lowerbound}
\arraycolsep=0.2em
\begin{array}{lcl}
\Pc_k & \geq & \frac{A_k}{2}\norms{w^k }^2 + \revise{\frac{s^2}{2}\big(\frac{s-1}{\gamma} - \frac{4}{\eta}\big)}\norms{z^k - x^{\star} }^2,
\end{array}
\end{equation}
where $A_k := \big[\revise{\frac{3(3s-2)}{16(s-1)}}\eta - 2\rho_n\big]t_{k-1}^2 + 2s(\rho_n - \rho_{*}) t_{k-1} - \gamma(t_{k-1} - 1)$.
\end{lemma}

\beforesubsubsec
\subsubsection{\textbf{The Lyapunov function and additional supporting lemmas.}}\label{apdx:subsubsec:case1:supporting_lemmas}
\aftersubsubsec
We consider the following Lyapunov function:
\begin{equation}\label{eq:SFOG_th1_Lyfunc} 
\hspace{-0.5ex}
\arraycolsep=0.2em
\begin{array}{lcl}
\Lc_k & := &  \Pc_k + \frac{\omega \eta L^2 t_{k-1}^2}{2} \norms{ x^k - y^{k-1} }^2 + \frac{\Lambda_0\eta (1-\kappa)(s+1)^2t_{k-1}^2}{s^2} \Delta_{k-1} +  \revise{\big(\varphi + \frac{s-2}{32(s-1)}\big) \eta} t_{k-1}^2 \norms{e^{k-1}}^2,
\end{array}
\hspace{-1ex}
\end{equation}
where $\Lambda_0$ is given in \eqref{eq:fixed_constants}, $\kappa$ is in \eqref{eq:error_cond3}, \revise{$\varphi$ is in \eqref{eq:SFOG_phik_and_varphi},} and $\Pc_k$ is the function defined by \eqref{eq:SFOG_Ek_func}.

Next, we need the following supporting lemmas to simplify the proof of Theorem~\ref{th:VFOG1_convergence}.
The proof of these lemmas can be found in Appendix~\ref{apdx:sec:convergence_analysis}.

\begin{lemma}\label{le:SFOG_supporting_lemma0}
Under the same settings as in Theorem~\ref{th:VFOG1_convergence}, $\Lc_k$ given in \eqref{eq:SFOG_th1_Lyfunc} satisfies
\begin{equation}\label{eq:SFOG_lmB0_bound}  
\arraycolsep=0.2em
\begin{array}{lcl}
\Lc_k - \Expsn{k}{ \Lc_{k+1} }  & \geq &  \frac{\epsilon_1\eta t_k(t_k-s)}{2} \Expsn{k}{ \norms{ \hat{w}^{k+1} - w^k}^2 } +  \frac{\eta st_k}{2} \Expsn{k}{ \norms{ \hat{w}^{k+1} }^2 }  \vspace{1ex}\\
&& + {~} \revise{\frac{(s-2)(4s-7)\eta (k+1)}{64(s-1)}} \norms{ w^k }^2 + \frac{ \revise{(s-2)} \eta L^2 t_{k-1}^2}{\revise{32(s-1)}}   \norms{x^k -  y^{k-1} }^2  \vspace{1ex}\\
&& + {~}  \frac{\mred{(2s+1)\Lambda_0\eta(1-\kappa)t_{k-1}(t_{k-1}-s)} }{s^2}\Delta_{k-1}  +  \frac{\revise{(s-2)}\eta t_{k-1}^2}{\revise{32(s-1)}}  \norms{ e^{k-1} }^2 - \Lambda_0\eta \delta_k,
\end{array}
\end{equation}
where $\epsilon_1 := 1 - (s+1)M^2\eta^2 \geq 0$ and $\Lambda_0$ is given in \eqref{eq:fixed_constants}.
\end{lemma}

\begin{lemma}\label{le:SFOG_supporting_lemma1}
Let $\sets{(x^k, y^k, z^k)}$ be generated by \eqref{eq:VFOG} and the conditions of Lemmas~\ref{le:VFOG_key_estimate2} and \ref{le:Ek_lowerbound} hold.
Suppose further that
\begin{equation}\label{eq:SFOG_proof1_summability_bounds}
\arraycolsep=0.2em
\begin{array}{lcl}
\sum_{k=0}^{\infty} t_k \Exp{\norms{\hat{w}^k }^2} & < & +\infty, \vspace{1.75ex}\\
\sum_{k=0}^{\infty} t_k^2 \Expn{ \tnorms{ \hat{w}^{k+1} - \hat{w}^k }^2 } & < &  +\infty, \vspace{1ex}\\
\sum_{k=0}^{\infty} t_k^2 \Expn{\norms{ \hat{e}^k }^2 }  & < &  +\infty, \vspace{1ex}\\
\sum_{k=0}^{\infty} t_k^2 \Expn{\norms{ e^k }^2 }  & < &  +\infty.
\end{array}
\end{equation}
Then, we have 
\begin{equation}\label{eq:SFOG_proof1_limits}
\lim_{k\to\infty} \Expn{ \norms{x^k - z^k}^2 } = 0, \ \ \textrm{and} \ \ \sum_{k=0}^{\infty}\tfrac{1}{t_k} \Expn{ \norms{z^k - x^k}^2 } < +\infty.
\end{equation}
In addition, we also have
\begin{equation}\label{eq:SFOG_proof1_summability_bounds2}
\arraycolsep=0.2em
\begin{array}{lcl}
\sum_{k=0}^{\infty} t_k^2 \Expn{ \norms{x^{k+1} - y^k }^2 } & < &  +\infty, \vspace{1ex}\\
\sum_{k=0}^{\infty} t_k \Expn{ \norms{y^k - x^k}^2 } & < & +\infty, \vspace{1ex}\\
\sum_{k=0}^{\infty} t_k \Expn{ \norms{x^{k+1} - x^k}^2} & < & +\infty.
\end{array}
\end{equation}
\end{lemma}

\begin{lemma}\label{le:SFOG_supporting_lemma2}
Under the same settings as in Lemma~\ref{le:SFOG_supporting_lemma1}, for $\Pc_k$ given in \eqref{eq:SFOG_Ek_func}, if $\Expn{\Pc_k} \leq \bar{\Pc}_0$ for all $k\geq 0$ and a given $\bar{\Pc}_0 \geq 0$, then $\lim_{k\to\infty}\Expn{ \norms{x^k - x^{\star}}^2 }$ exists for any $x^{\star} \in \zer{\Phi}$.

If, additionally, $\lim_{k\to\infty}\revise{\Expn{\Pc_k}}$ exists, $\sets{ t_{k-1}^2\Expn{\norms{w^k}^2} } $ is bounded, and $\sum_{k=0}^{\infty} t_{k-1} \Expn{\norms{w^k}^2 } < +\infty$, then 
\begin{equation}\label{eq:SFOG_lmB2_limits}
\lim_{k\to\infty}t_k^2\Expn{ \norms{w^k}^2 } = 0  \quad \textrm{and} \quad \lim_{k\to\infty}t_k^2\Expn{ \norms{ \hat{w}^k}^2 } = 0.
\end{equation}
\end{lemma}

\beforesubsubsec
\subsubsection{\textbf{The proof of Theorem~\ref{th:VFOG1_convergence}.}}\label{apdx:th:convergence1}
\aftersubsubsec
For clarity, we divide the proof into the following steps.

\vspace{0.75ex}
\noindent\textbf{\underline{Step 1} (\textit{The summability results in \eqref{eq:SFOG_summability_bounds}}).}
First, \mred{applying} Lemma~\ref{le:SFOG_supporting_lemma0}, it follows from \eqref{eq:SFOG_lmB0_bound} after taking the full expectation \mred{that}
\begin{equation}\label{eq:SFOG_th31_proof1}
\arraycolsep=0.2em
\begin{array}{lcl}
\Expn{\Lc_k} -  \Expn{ \Lc_{k+1}}  & \geq & \frac{\eta s(k+s+1)}{2} \Exp{\norms{ \hat{w}^{k+1} }^2}  + \frac{\revise{(s-2)(4s-7)\eta}(k+1)}{\revise{64(s-1)}} \Expn{ \norms{ w^k }^2} \vspace{1ex}\\
&& + {~}  \frac{(2s+1) \Lambda_0\eta \mred{(1-\kappa)}t_{k-1}(t_{k-1}-s)}{s^2} \Expn{ \Delta_{k-1}} \vspace{1ex}\\
&& + {~} \frac{\epsilon_1\eta t_k(t_k-s)}{2} \Expn{ \norms{ \hat{w}^{k+1} - w^k}^2} + \frac{\revise{(s-2)} \eta L^2 t_{k-1}^2}{\revise{32(s-1)}}\Expn{ \norms{ x^k - y^{k-1} }^2 } \vspace{1ex}\\
&& + {~}   \frac{\revise{(s-2)}\eta t_{k-1}^2}{\revise{32(s-1)}} \Expn{ \norms{e^{k-1} }^2 }   - \Lambda_0\eta\delta_k.
\end{array}
\end{equation}
Moreover, from Lemma~\ref{le:Ek_lowerbound}, \eqref{eq:SFOG_th1_Lyfunc}, and the nonnegativity of the coefficients, we have $\Lc_k \geq 0$ for all $k \geq 0$.
Summing up \eqref{eq:SFOG_th31_proof1} from $k=0$ to $k := K$ and using $\Lc_k \geq 0$ and $\Sc_K := \sum_{k=0}^K\delta_k
 $, we obtain
\begin{equation}\label{eq:SFOG_th31_proof2}
\arraycolsep=0.2em
\begin{array}{lclcl}
\sum_{k=0}^K \eta s(k+s+1) \Exp{\norms{ \hat{w}^{k+1} }^2} & \leq & \mred{2 \Lc_0 + 2\Lambda_0\eta \Sc_{K}}, \vspace{1ex}\\ 
\sum_{k=0}^K \frac{\revise{(s-2)(4s-7)\eta}}{\revise{s-1}} (k+1) \Exp{\norms{ w^k }^2} & \leq & \revise{64} \Lc_0 + \revise{64}\Lambda_0\eta \Sc_{K}, \vspace{1ex}\\ 
\end{array}
\end{equation}
Since $z^0 = x^0 = y^{-1}$ and the conventions $e^{-1} = \mred{\Delta_{-1}} = 0$, we have from \eqref{eq:SFOG_th1_Lyfunc} that
\begin{equation*} 
\arraycolsep=0.2em
\begin{array}{lcl}
\Lc_0 & = & \Pc_0 =  \frac{a_0}{2}\norms{ w^0 }^2  + \frac{s^2(s-1)}{2\gamma_{\revise{-1}}}\norms{x^0 - x^{\star}}^2 \vspace{1ex}\\
& = &  \revise{\frac{[2(13s-10)s^2 - (s-1)(s-2)] \eta}{64(s-1)}} \norms{ w^0 }^2 + \revise{\frac{16s^3(s-1)}{(s-2)\eta}}\norms{x^0 - x^{\star} }^2.
\end{array}
\end{equation*}
Here, we have used $a_0 = \left[\revise{\frac{13s-10}{16(s-1)}}\eta - 2\rho_n \right]s^2 + 2s^2\rho_n - \gamma(s-1) =  \revise{\frac{[2(13s-10)s^2 - (s-1)(s-2)] \eta}{32(s-1)}} > 0$ for all $s > 2$ and $\gamma_{\revise{-1}} = \frac{\revise{(s-2)}\eta(t_{\revise{-1}}-1)}{\revise{32(s-1)}t_{\revise{-1}}} = \revise{\frac{{(s-2)\eta}}{32s}}$.
Substituting this $\Lc_0$ into \eqref{eq:SFOG_th31_proof2}, we obtain \eqref{eq:SFOG_summability_bounds}.

\vspace{0.75ex}
\noindent\textbf{\underline{Step 2} (\textit{The proof of the bounds in \eqref{eq:SFOG_BigO_rates}}).}
Now, from \eqref{eq:SFOG_th31_proof1}, by induction, we also have $\Expn{\Lc_k } \leq \Lc_0 + \Lambda_0\eta \Sc_k$.
However, from \eqref{eq:Ek_lowerbound} and \eqref{eq:SFOG_th1_Lyfunc}, we also have
\begin{equation}\label{eq:SFOG_th31_proof1_100} 
\arraycolsep=0.2em
\begin{array}{lcl}
\Lc_k \geq \Pc_k & \geq & \frac{A_k}{2}\norms{ w^k }^2 + \revise{\frac{s^2}{2}\big(\frac{s-1}{\gamma} - \frac{4}{\eta}\big)}\norms{z^k - x^{\star} }^2 \geq \frac{A_k}{2}\norms{ w^k }^2,
\end{array}
\end{equation}
where \revise{we have used $\gamma := \frac{(s-2)\eta}{32(s-1)}$ to evaluate $\frac{s^2}{2}\big(\frac{s-1}{\gamma} - \frac{4}{\eta}\big) = \frac{2s^2(8s^2 - 17s + 10)}{(s-2)\eta} > 0$ for all $s>2$ in the last inequality.
Moreover, using $\rho_{*} \leq \rho_n$, $2\rho_n \leq \frac{s-2}{4(s-1)}\eta$ from \eqref{eq:SFOG_para_update2}, and $\gamma t_{k-1} \leq \frac{\gamma}{2}t_{k-1}^2$ due to $t_{k-1} \geq s > 2$, we also have}
\revise{
\begin{equation*}
\begin{array}{lcl}
A_k & := & \big[ \frac{3(3s-2)}{16(s-1)} \eta - 2\rho_n\big]t_{k-1}^2 + 2s(\rho_n - \rho_{*}) t_{k-1} -  \gamma t_{k-1} + \gamma \vspace{1ex}\\
& \geq & \big[\frac{3(3s-2)}{16(s-1)}\eta - 2\rho_n - \frac{\gamma}{2}\big] t_{k-1}^2 
\geq \big[\frac{3(3s-2)}{16(s-1)} - \frac{s-2}{4(s-1)} - \frac{(s-2)}{64(s-1)}\big]\eta(k+s)^2 \vspace{1ex}\\
& = & \frac{(19s+10)\eta}{64(s-1)}(k+s)^2.
\end{array}
\end{equation*}}
Therefore, we obtain 
\begin{equation*} 
\arraycolsep=0.2em
\begin{array}{lcl}
\Expn{ \norms{ w^k }^2 } \leq \frac{ \revise{128(s-1)}( \Lc_0 + \Lambda_0\eta \Sc_k )}{\revise{(19s+10)\eta} (k+s)^2 } 
= \frac{ \revise{128(s-1)} (\mcal{R}_0^2 +  \Lambda_0 \Sc_k )}{\revise{(19s+10)} (k+s)^2 },
\end{array}
\end{equation*}
which is the first bound of \eqref{eq:SFOG_BigO_rates}, where \revise{$\Rc_0$ is defined in \eqref{eq:SFOG_BigO_R0}.}

\revise{Next, since
\begin{equation}\label{eq:SFOG_th31_proof2.1}
\begin{array}{lcl}
	A_k & = & \frac{3(3s-2)}{16(s-1)}\eta t_{k-1}^2 - 2\rho_n t_{k-1}(t_{k-1} - s) - 2s\rho_{*} t_{k-1} -  \gamma (t_{k-1} - 1) 
	\leq \frac{3(3s-2)}{16(s-1)}\eta t_{k-1}^2,
\end{array}
\end{equation}
from \eqref{eq:fixed_constants}, we have}
\begin{equation*} 
\arraycolsep=0.2em
\begin{array}{lcl}
\omega& = & \revise{\frac{\varphi}{4} + \frac{2\Lambda_0\Theta(s+1)^2}{s^2} + \frac{s-2}{16(s-1)} \geq \frac{\varphi}{4} \geq \frac{9(85s - 134)}{256(s-1)} \geq \frac{3(3s-2)}{16(s-1)} \geq \frac{A_k}{\eta t_{k-1}^2}}.
\end{array}
\end{equation*}
By \mred{the $L$-Lipschitz continuity of $G$ and} Young's inequality, we have
\mred{\begin{equation*}
\arraycolsep=0.2em
\begin{array}{lcl}
\frac{A_k}{2}\norms{w^k}^2 + \frac{\omega\eta L^2 t_{k-1}^2}{2}\norms{x^k - y^{k-1}}^2 &\geq& \frac{A_k}{2}\norms{w^k}^2 + \frac{\omega\eta t_{k-1}^2}{2}\norms{Gx^k - Gy^{k-1}}^2 \vspace{1ex}\\
&=& \frac{A_k}{2}\norms{w^k}^2 + \frac{\omega\eta t_{k-1}^2}{2}\norms{w^k - \hat{w}^k  }^2 \vspace{1ex} \\
&\geq& \frac{A_k}{4}\norms{ \hat{w}^k }^2 \vspace{1ex}\\
&\geq& \revise{\frac{(19s+10)\eta}{256(s-1)}} (k+s)^2\norms{ \hat{w}^k }^2.
\end{array}
\end{equation*}}
Combining this inequality, \eqref{eq:SFOG_th1_Lyfunc}, and \eqref{eq:SFOG_th31_proof1_100}, we can show that
\begin{equation*}
\arraycolsep=0.2em
\begin{array}{lcl}
\Lc_k &\geq & \Pc_k +  \frac{\omega\eta\revise{L^2} t_{k-1}^2}{2}\norms{\revise{x^k - y^{k-1}} }^2 \geq \revise{\frac{(19s+10)\eta}{256(s-1)}} (k+s)^2\norms{ \hat{w}^k }^2.
\end{array}
\end{equation*}
Hence, utilizing this inequality, with  similar arguments as  the proof of the first line of \eqref{eq:SFOG_BigO_rates}, we obtain the second line of \eqref{eq:SFOG_BigO_rates} from this expression. 

\vspace{0.75ex}
\noindent\textbf{\underline{Step 3} (\textit{Other intermediate summability results}).}
Next, since \mred{$\epsilon_1 \mred{:= 1 - (s+1)M^2\eta^2} > 0$ (due to the choice $\eta < \mred{\frac{\lambda}{L}}$ as in \eqref{eq:SFOG_para_update2}) and $\Sc_{\infty} := \sum_{k=0}^\infty \delta_k < + \infty$, we can rearrange \eqref{eq:SFOG_th31_proof1} and then use Lemma \ref{le:A1_descent} to prove that $\lim_{k \to \infty} \Expn{\Lc_k}$ exists and }
\begin{equation}\label{eq:SFOG_summability_proof1}
\arraycolsep=0.2em
\begin{array}{lcl}
	\mred{\sum_{k=0}^{\infty} t_k \Expn{ \norms{w^k}^2 }}  & < &  +\infty, \vspace{1ex}\\
	\mred{\sum_{k=0}^{\infty} t_k \Expn{ \norms{\hat{w}^k}^2 }}  & < &  +\infty, \vspace{1ex}\\
	\sum_{k=0}^{\infty} t_{k}^2 \Expn{ \norms{ x^{k+1} -  y^k }^2 } & < &  +\infty, \vspace{1ex}\\
	\sum_{k=0}^{\infty} t_k(t_k-s)  \Expn{ \norms{  \hat{w}^{k+1} - w^k}^2 }  & < & +\infty, \vspace{1ex}\\
	\sum_{k=0}^{\infty} t_k^2 \Expn{ \Delta_k }  & < &  +\infty, \vspace{1ex}\\
	\sum_{k=0}^{\infty} t_k^2 \Expn{\norms{e^k }^2 }  & < &  +\infty.
\end{array}
\end{equation}
By Assumption~\ref{as:A1} and Young's inequality, we easily get
\begin{equation*} 
\arraycolsep=0.2em
\begin{array}{lcl}
\norms{w^{k+1} - \hat{w}^{k+1}}^2 &\leq & L^2\norms{x^{k+1} - y^k}^2, \vspace{1ex}\\
\norms{w^{k+1} - w^k }^2 & \leq & 2\norms{w^{k+1} - \hat{w}^{k+1} }^2 + 2\norms{ \hat{w}^{k+1} - w^k}^2, \vspace{1ex}\\
\norms{ \hat{w}^{k+1} - \hat{w}^k }^2 & \leq & 2\norms{w^k - \hat{w}^{k+1} }^2 + 2\norms{w^k - \hat{w}^k }^2.
\end{array}
\end{equation*}
Moreover, we also have $\Expn{\norms{\hat{e}^k}^2} \leq 2 \Expn{\norms{ e^{k-1} }^2 }  + 2\Expn{ \norms{w^k - \hat{w}^k  }^2 }$.
Using these relations and $t_k(t_k-s) \sim t_k^2$, we can derive from \eqref{eq:SFOG_summability_proof1} that
\begin{equation}\label{eq:SFOG_summability_proof2}
\arraycolsep=0.2em
\begin{array}{lcl}
\sum_{k=0}^{\infty} t_k^2 \Expn{ \norms{ w^{k+1} - \hat{w}^{k+1} }^2 } & < &  +\infty, \vspace{1ex}\\
\sum_{k=0}^{\infty} t_k^2 \Expn{ \norms{ w^{k+1} - w^k }^2 } & < &  +\infty, \vspace{1ex}\\
\sum_{k=0}^{\infty} t_k^2 \Expn{ \norms{ \hat{w}^{k+1} - \hat{w}^k }^2 } & < &  +\infty, \vspace{1ex}\\
\sum_{k=0}^{\infty} t_k^2 \Expn{\norms{ \hat{e}^k }^2 }  & < &  +\infty.
\end{array}
\end{equation}
\vspace{0.75ex}
\noindent\textbf{\underline{Step 4} (\textit{Applying Lemmas~\ref{le:SFOG_supporting_lemma1} and \ref{le:SFOG_supporting_lemma2}}).}
First, from \mred{\eqref{eq:SFOG_summability_proof1}} and \eqref{eq:SFOG_summability_proof2}, we see that the conditions \eqref{eq:SFOG_proof1_summability_bounds} of Lemma~\ref{le:SFOG_supporting_lemma1} are fulfilled. 
Next, from \eqref{eq:SFOG_th31_proof1}, by induction, we get $\Expn{\Lc_k} \leq \Lc_0 + \eta\Lambda_0\Sc_{\infty}$ for all $k\geq 0$.
Moreover,  by \eqref{eq:SFOG_th31_proof1_100}, we have $\Pc_k \leq \Lc_k$.
Therefore, we conclude that $\Expn{\Pc_k} \leq \Lc_0 + \eta\Lambda_0\Sc_{\infty} \equiv \bar{\Pc}_0$.

\mred{Next}, by \eqref{eq:SFOG_th1_Lyfunc}, we have
\begin{equation*} 
\arraycolsep=0.2em
\begin{array}{lcl}
\Lc_k & := &  \Pc_k + \frac{\omega \eta L^2 t_{k-1}^2}{2} \norms{ x^k - y^{k-1} }^2 + \frac{\Lambda_0\eta (1-\kappa)(s+1)^2t_{k-1}^2}{s^2} \Delta_{k-1} +  \revise{\big[ \varphi + \frac{s-2}{32(s-1)} \big]} \eta t_{k-1}^2 \norms{e^{k-1}}^2.
\end{array}
\end{equation*}
\mred{Since t}he last three terms are summable due to \eqref{eq:SFOG_summability_proof1}, their limit is zero.
\mred{Using these facts and the existence of $\lim_{k \to \infty} \Expn{\Lc_k}$ proved above,} we conclude that $\lim_{k\to\infty}\Expn{\Pc_k}$ also exists.

From the first line of \eqref{eq:SFOG_BigO_rates}, we also say that $\sets{ t_{k-1}^2\Expn{\norms{w^k}^2} } $ is bounded, and from the \revise{first} line of \eqref{eq:SFOG_summability_bounds}, we have $\sum_{k=0}^{\infty} t_{k-1} \Expn{\norms{w^k}^2 } < +\infty$.
Overall, we have verified all the conditions of Lemma~ \ref{le:SFOG_supporting_lemma2}.
Hence, we obtain \eqref{eq:SFOG_small_o_rates} from \eqref{eq:SFOG_lmB2_limits}.

\vspace{0.75ex}
\noindent\textbf{\underline{Step 5} (\textit{Iteration-complexity}).}
For a given tolerance $\epsilon > 0$, from \eqref{eq:SFOG_BigO_rates} and $\Sc_k \leq \Sc_{\infty} < +\infty$, to guarantee $\Expn{ \norms{Gx^k + v^k}^2 } \leq \epsilon^2$, we need to impose $\mred{\frac{C_0^2 (\mcal{R}_0^2 + \Lambda_0 S_{\infty})}{(k+s)^2 }} \leq \epsilon^2$.
This relation implies that $k = \BigOs{1/\epsilon}$, which proves our statement (d).
\Eproof

\beforesubsec
\subsection{\textbf{The proof of Theorem~\ref{th:VFOG1_almost_sure_convergence}: Almost sure convergence.}}\label{apdx:th:VFOG1_almost_sure_convergence}
\aftersubsec
Similar to the proof of Theorem~\ref{th:VFOG1_convergence}, we also divide this proof into several steps as follows.

\vspace{0.75ex}
\noindent\textbf{\underline{Step 1}~(\textit{Almost sure summability bounds)}.}
For $\Lc_k$ given in \eqref{eq:SFOG_th1_Lyfunc}, from \eqref{eq:SFOG_lmB0_bound}, we have
\begin{equation}\label{eq:SFOG_th2_proof1a}  
\arraycolsep=0.2em
\begin{array}{lcl}
\Lc_k - \Expsn{k}{ \Lc_{k+1} }  & \geq &  \frac{\epsilon_1\eta t_k(t_k-s)}{2} \Expsn{k}{ \norms{ \hat{w}^{k+1} - w^k}^2 } +  \frac{\eta st_k}{2} \Expsn{k}{ \norms{ \hat{w}^{k+1} }^2 }  \vspace{1ex}\\
&& + {~} \frac{\revise{(s-2)(4s-7)}\eta(k+1)}{\revise{64(s-1)}} \norms{ w^k }^2 + \frac{ \revise{(s-2)} \eta L^2 t_{k-1}^2 }{\revise{32(s-1)}}   \norms{x^k - y^{k-1} }^2   \vspace{1ex}\\
&& + {~} \frac{(2s+1) \Lambda_0\eta\mred{(1-\kappa)} t_{k-1}(t_{k-1}-s)}{s^2}\Delta_{k-1} +  \frac{\revise{(s-2)}\eta t_{k-1}^2}{\revise{32(s-1)}} \norms{ e^{k-1} }^2 - \Lambda_0\eta\delta_k.
\end{array}
\end{equation}
Let us define 
\begin{equation}\label{eq:Uk_def_proof1} 
\arraycolsep=0.2em
\begin{array}{lcl}
U_k & := & \Lc_k +    \frac{\epsilon_1\eta t_{k-1}(t_{k-1} - s)}{2} \norms{ \hat{w}^k - w^{k-1} }^2 +  \frac{\eta st_{k-1} }{2}   \norms{ \hat{w}^k }^2.
\end{array}
\end{equation}
Then, we can rearrange \eqref{eq:SFOG_th2_proof1a} and \mred{use} $U_k$ from \eqref{eq:Uk_def_proof1} to get
\begin{equation}\label{eq:SFOG_th2_proof1}  
\arraycolsep=0.2em
\begin{array}{lcl}
\Expsn{k}{ U_{k+1} }  & \leq &  U_k - \frac{\epsilon_1\eta t_{k-1}(t_{k-1} - s)}{2} \norms{ \hat{w}^k - w^{k-1} }^2 - \frac{\eta st_{k-1} }{2}   \norms{ \hat{w}^k }^2  \vspace{1ex}\\
&& - {~} \frac{\revise{(s-2)(4s-7)}\eta(k+1)}{\revise{64(s-1)}} \norms{ w^k }^2 - \frac{(2s+1) \Lambda_0\eta \mred{(1-\kappa)} t_{k-1}(t_{k-1}-s)}{s^2}\Delta_{k-1} \vspace{1ex}\\
&& - {~} \frac{ \revise{(s-2)} \eta L^2 t_{k-1}^2 }{\revise{32(s-1)}}   \norms{x^k - y^{k-1} }^2  -  \frac{\revise{(s-2)}\eta t_{k-1}^2}{\revise{32(s-1)}}  \norms{ e^{k-1} }^2 + \Lambda_0\eta\delta_k.
\end{array}
\end{equation}
Note that $\Lc_k \geq 0$ almost surely for all $k \geq 0$ and $\Sc_{\infty} := \sum_{k=0}^{\infty}\delta_k < +\infty$ by our assumption,  applying Lemma~\ref{le:RS_lemma} to the last inequality, we conclude that $\lim_{k\to\infty} U_k$ exists almost surely.
Moreover, we also have the following almost sure summability bounds:
\begin{equation}\label{eq:SFOG_th32_proof2}
\arraycolsep=0.2em
\begin{array}{lcl}
\sum_{k=0}^{\infty}(k+s)\norms{\hat{w}^k}^2  & < &  +\infty, \vspace{1ex}\\
\sum_{k=0}^{\infty}(k+1)\norms{w^k}^2  & < &  +\infty, \vspace{1ex}\\
\sum_{k=0}^{\infty}(k+s+1)^2 \norms{e^k}^2  & < &  +\infty, \vspace{1ex}\\
\sum_{k=0}^{\infty}(k+1)(k+s+1) \Delta_k  & < &  +\infty, \vspace{1ex}\\
\sum_{k=0}^{\infty}(k+1)(k+s+1)\norms{\hat{w}^{k+1} - w^k}^2  & < &  +\infty, \vspace{1ex}\\
\sum_{k=0}^{\infty}(k+s)^2\norms{x^{k+1} - y^k }^2  & < &  +\infty.
\end{array}
\end{equation}
The first two lines of \eqref{eq:SFOG_th32_proof2} are exactly \eqref{eq:SFOG_as_summability}.
Similar to the proof of \eqref{eq:SFOG_summability_proof2}, but using \eqref{eq:SFOG_th32_proof2}, we can also show that
\begin{equation}\label{eq:SFOG_th32_proof2b}
\arraycolsep=0.2em
\begin{array}{lcl}
\sum_{k=0}^{\infty} t_k^2 \norms{\hat{w}^{k+1} - \hat{w}^k}^2  & < &  +\infty \quad \textrm{almost surely}, \vspace{1ex}\\
\sum_{k=0}^{\infty} t_k^2 \norms{w^{k+1} - \hat{w}^{k+1}}^2  & < &  +\infty \quad \textrm{almost surely}.
\end{array}
\end{equation}
From \eqref{eq:SFOG_th32_proof2} we also have $\lim_{k\to\infty}  t_{k-1}(t_{k-1} - s) \norms{ \hat{w}^k - w^{k-1} }^2 = 0$ and $\lim_{k \to\infty}  t_{k-1} \norms{ \hat{w}^k }^2 = 0$ almost surely.
Combining these facts and the existence of $\lim_{k\to\infty}U_k$ almost surely, \eqref{eq:Uk_def_proof1} implies that $\lim_{k\to\infty} \Lc_k$ also exists almost surely.

\vspace{0.75ex}
\noindent\textbf{\underline{Step 2} (\textit{The bounds and almost sure convergence of $\revise{u}^k := s(z^k - x^k)$}).}
Let $u^k :=  s(z^k - x^k)$ \revise{and $\hat{u}^k := s(z^k - x^k) - s\eta\hat{w}^k$}.
Then, similar to the proof of \revise{Lemma~\ref{le:SFOG_supporting_lemma1}}, we can show that
\myeqn{
	\arraycolsep=0.15em
	\begin{array}{lcl}
		\revise{\hat{u}^{k+1}}
		& = &  \big(1 - \frac{s}{t_k}\big)\revise{\hat{u}^k}+ \frac{s}{t_k} \cdot \frac{t_k}{s} \big[  [s(\eta-\beta_k)-\gamma_k - \frac{s^2\eta}{t_k}]\hat{w}^k + s\eta e^k - (s\beta_k + \gamma_k)e^{k-1}\big].
	\end{array}
}
Since $\frac{s}{t_k} \in (0,1)$, by convexity of $\norms{\cdot}^2$, Young's inequality, $\beta_k \leq \eta$, and $\gamma_k \leq \eta$, we obtain the following inequality after taking the conditional expectation $\Exps{k}{\cdot}$:
\begin{equation*} 
\arraycolsep=0.15em
\begin{array}{lcl}
	\revise{\Exps{k}{\norms{\hat{u}^{k+1}}^2}}
	& \leq &  \big(1 - \frac{s}{t_k}\big)\revise{\norms{\hat{u}^k}^2} + \frac{3t_k}{s}[s(\eta-\beta_k)-\gamma_k - \frac{s^2\eta}{t_k}]^2\norms{\hat{w}^k}^2 \vspace{1ex} \\
	&& + {~} 3s\eta^2 t_k \Exps{k}{\norms{e^k}^2} + \frac{3t_k}{s}(s\beta_k + \gamma_k)^2 \norms{e^{k-1}}^2 \vspace{1ex}\\
	&\leq& \big(1 - \frac{s}{t_k}\big)\revise{\norms{\hat{u}^k}^2} + \frac{12}{s}\big[(2s^2+1) + \frac{s^4}{t_k^2}\big]\eta^2 t_k\norms{\hat{w}^k}^2 \vspace{1ex}\\
	&& + {~} 3s\eta^2 t_k \Exps{k}{\norms{e^k}^2} + \frac{3(s+1)^2\eta^2 t_k}{s} \norms{e^{k-1}}^2.
\end{array}
\end{equation*}
By \eqref{eq:error_cond3} and $1-\kappa \leq 1$, we have 
\begin{equation*} 
\arraycolsep=0.2em
\begin{array}{lcl}
\Expsn{k}{\norms{e^k}^2 } \leq \Expsn{k}{\Delta_k } &\leq &  \Delta_{k-1} + \Theta L^2 \norms{x^k - y^{k-1}}^2 + \frac{\delta_k}{t_k^2}.
\end{array}
\end{equation*}
Combining the last two inequalities, and noting that $\kappa \in (0, 1]$, we finally get
\begin{equation}\label{eq:SAEG_th32_proof5a} 
\arraycolsep=0.2em
\begin{array}{lcl}
	\revise{\Exps{k}{\norms{\hat{u}^{k+1}}^2}}
	& \leq & \norms{\revise{\hat{u}^k}}^2 - \frac{s}{t_k}\norms{\revise{\hat{u}^k}}^2 +\frac{12}{s}\big[(2s^2+1)\eta^2 t_k + \frac{s^4\eta}{t_k}\big]\norms{\hat{w}^k}^2  + \frac{3(s+1)^2\eta^2 t_k}{s} \norms{e^{k-1}}^2 \vspace{1ex}\\
	&& + {~} 3s\eta^2 t_k \Delta_{k-1} + 3s\eta^2 L^2 \Theta  t_k \norms{x^k - y^{k-1}}^2 + \frac{3s\eta^2 \delta_k}{t_k}.
\end{array}
\end{equation}
Note that the third to the sixth terms of this expression are almost surely summable due to \eqref{eq:SFOG_th32_proof2} and \eqref{eq:SFOG_th32_proof2b}.
Moreover, we also have $\sum_{k=0}^{\infty}\frac{\delta_k}{t_k} \leq \sum_{k=0}^{\infty}\delta_k < +\infty$.
Applying again Lemma~\ref{le:RS_lemma} to the last inequality and using these facts,  we conclude that almost surely
\begin{equation}\label{eq:SAEG_th32_proof5} 
\arraycolsep=0.0em
\begin{array}{lcl}
\lim_{k\to\infty} \norms{z^k - x^k - \eta\hat{w}^k}^2  \ \ \textrm{exists \quad and \quad} \ \sum_{k=0}^{\infty}\frac{\norms{z^k - x^k - \eta\hat{w}^k}^2 }{t_k} <+\infty.
\end{array}
\end{equation}
In addition, we have $\sum_{k=0}^{\infty}\frac{1}{t_k} = +\infty$.
Utilizing this fact and \eqref{eq:SAEG_th32_proof5},  by a contradiction argument, we must have $\lim_{k\to\infty} \norms{z^k - x^k - \eta\hat{w}^k}^2 = 0$ almost surely.
Moreover, from \eqref{eq:SFOG_th32_proof2}, it follows that
\begin{equation*} 
\arraycolsep=0.0em
\begin{array}{lcl}
	\lim_{k \to +\infty} (k+s)\norms{\hat{w}^k}^2 = 0 \quad \text{almost surely}
\end{array}
\end{equation*}
By Young's inequality, we have
\begin{equation*} 
	\norms{s(z^k - x^k)}^2 \leq 2s^2 \norms{z^k - x^k - \eta\hat{w}^k}^2 + 2s^2\eta^2\norms{\hat{w}^k}^2.
\end{equation*}
Combining the last three facts and \eqref{eq:SAEG_th32_proof5}, and using $t_k = k + s + 1$, we conclude that
\begin{equation}\label{eq:SAEG_th32_proof5.3} 
\arraycolsep=0.0em
\begin{array}{lcl}
	\lim_{k\to\infty} \norms{z^k - x^k}^2  = 0 \textrm{\quad and \quad} \ \sum_{k=0}^{\infty}\frac{\norms{z^k - x^k}^2}{t_k} <+\infty.
\end{array}
\end{equation}
\vspace{0.75ex}
\noindent\textbf{\underline{Step 3} (\textit{The proof of \eqref{eq:SFOG_as_limits}}).}
First, from \eqref{eq:SFOG_quantities} we have $\hat{e}^k = \widetilde{G}y^{k-1} - Gx^k = e^{k-1} + Gy^{k-1} - Gx^k$.
By Young's inequality and \eqref{eq:G_Lipschitz} in Assumption~\ref{as:A1}, the last expression leads to
$\norms{\hat{e}^k}^2 \leq 2\norms{e^{k-1}}^2 + 2\norms{Gx^k - Gy^{k-1}}^2 \leq 2\norms{e^{k-1}}^2 + 2L^2\norms{x^k - y^{k-1}}^2$.
Exploiting this inequality and $\iprods{w^k, x^k - x^{\star}} \geq -\rho_n\norms{w^k}^2$ for any $x^{\star} \in \zer{\Phi}$ from \eqref{eq:G_cohypomonotonicity} of Assumption~\ref{as:A3}, similar to the proof of \eqref{eq:SFOG_summability_proof5a}, we can show that
\begin{equation}\label{eq:SFOG_th32_proof6} 
\arraycolsep=0.2em
\begin{array}{lcl}
\Expsn{k}{ \norms{z^{k+1} - x^{\star}}^2 }
& \leq &  \big(1 + \frac{\eta}{st_k^2}\big) \norms{z^k - x^{\star}}^2 + \big(\frac{\eta t_k}{s} + \frac{\revise{2}\eta \rho_n }{s} +  \frac{2\eta^2}{s^2} \big)\norms{ w^k }^2  \vspace{1ex}\\
&& + {~}   \frac{\eta }{st_k}\norms{z^k - x^k }^2 + \big( \frac{\eta t_k^2}{s} + \frac{2\eta^2}{s^2} \big)[ 2\norms{ e^{k-1} }^2 + 2L^2\norms{x^k - y^{k-1}}^2].
\end{array}
\end{equation}
Since $\sum_{k=0}^{\infty}\frac{\eta}{st_k^2} < +\infty$ and the last three terms of \eqref{eq:SFOG_th32_proof6} are  almost surely summable due to $t_k = k + s \revise{~+~1}$, \eqref{eq:SFOG_th32_proof2}, and \eqref{eq:SAEG_th32_proof5.3}, 
applying Lemma~\ref{le:RS_lemma} to \eqref{eq:SFOG_th32_proof6}, we conclude that  $\lim_{k\to\infty}\norms{z^k - x^{\star}}^2$ exists almost surely.

Next, from \eqref{eq:SFOG_th1_Lyfunc}, we have 
\begin{equation*} 
\arraycolsep=0.2em
\begin{array}{lcl}
\Lc_k & := &  \Pc_k + \frac{\omega \eta L^2 t_{k-1}^2}{2} \norms{ x^k - y^{k-1} }^2 + \frac{\Lambda_0\eta (1-\kappa)(s+1)^2t_{k-1}^2}{s^2} \Delta_{k-1} +  \revise{\big(\varphi + \frac{s-2}{32(s-1)}\big)} \eta t_{k-1}^2 \norms{e^{k-1}}^2.
\end{array}
\end{equation*}
However, the last three terms are almost surely summable due to \eqref{eq:SFOG_th32_proof2}, their limit is zero.
Hence, we conclude that $\lim_{k\to\infty} \Pc_k$ also exists almost surely.

Since $\lim_{k\to\infty} \Pc_k$ exists almost surely, and \revise{$\Pc_k \geq \frac{A_k}{2}\norms{w^k}^2 \geq \frac{C_1 t_{k-1}^2}{2}\norms{w^k}^2$ due to the proof of Theorem~\ref{th:VFOG1_convergence}, where $C_1 := \frac{(19s+10)\eta}{64(s-1)}$}, there exists $\bar{M} > 0$ such that $t_{k-1}^2 \norms{w^k}^2  \leq \bar{M}^2$ almost surely.
Using this inequality and $\lim_{k\to\infty} \norms{x^k-z^k} = 0$ almost surely, we can show that almost surely
\begin{equation*} 
\arraycolsep=0.2em
\begin{array}{lcl}
t_{k-1} \vert  \iprods{w^k, x^k - z^k}  \vert &\leq & t_{k-1}  \norms{w^k}  \norms{x^k - z^k}  \leq \bar{M} \norms{x^k - z^k} \to 0 \quad \text{as}~k \to \infty.
\end{array}
\end{equation*}
This expression implies that $\lim_{k\to\infty} t_{k-1} \iprods{w^k, x^k - z^k }   = 0$ almost surely.

Now, recall $\Pc_k$ from \eqref{eq:SFOG_Ek_func} as
\begin{equation*} 
\arraycolsep=0.2em
\begin{array}{lcl}
\Pc_k = \frac{a_k}{2}\norms{w^k}^2 + st_{k-1}\iprods{w^k, x^k - z^k} + \frac{s^2(s-1)}{2\gamma_{k-1}}\norms{z^k - x^{\star}}^2.
\end{array}
\end{equation*}
Since $\lim_{k\to\infty} \Pc_k $ and $\lim_{k\to\infty} \norms{z^k - x^{\star}}^2$ exist almost surely, $\lim_{k\to\infty} t_{k-1} \iprods{w^k, x^k - z^k }   = 0$ almost surely, $\lim_{k\to\infty}\gamma_{k-1} = \gamma > 0$, and $\lim_{k\to\infty}\frac{t_{k-1}^2}{a_k} = \frac{1}{\psi}$, we conclude that $\lim_{k\to\infty}t_{k-1}^2  \norms{w^k}^2 $ exists almost surely.
However, since $\sum_{k=0}^{\infty}t_{k-1} \norms{w^k}^2  < +\infty$ due to \eqref{eq:SFOG_th32_proof2}, applying Lemma~\ref{le:A2_sum}, we can easily argue that $\lim_{k\to\infty}t_{k-1}^2  \norms{w^k}^2  = 0$ almost surely.
This limit implies the first limit of \eqref{eq:SFOG_as_limits}.
Note that 
\begin{equation*} 
\arraycolsep=0.2em
\begin{array}{lcl}
t_{k-1}^2\norms{\hat{w}^k }^2 & \leq & 2t_{k-1}^2\norms{w^k}^2 + 2t_{k-1}^2\norms{Gx^k - Gy^{k-1}}^2 \leq 2t_{k-1}^2\norms{w^k}^2 + 2L^2t_{k-1}^2\norms{x^k - y^{k-1}}^2.
\end{array}
\end{equation*}
Since $\lim_{k\to\infty}t_{k-1}^2 \norms{x^k - y^{k-1}}^2  = 0$ and $\lim_{k\to\infty} t_{k-1}^2   \norms{w^k}^2  = 0$ almost surely, the last inequality implies that $\lim_{k\to\infty} t_{k-1}^2  \norms{ \hat{w}^k }^2  = 0$, which leads to the second limit of \eqref{eq:SFOG_as_limits}.

\vspace{0.75ex}
\noindent\textbf{\underline{Step 4} (\textit{The almost sure convergence of iterates}).}
Note that $\vert \norms{x^k - x^{\star}} - \norms{z^k - x^{\star}}\vert \leq \norms{x^k - z^k} $.
Since $\lim_{k\to\infty}\norms{x^k - z^k}  = 0$ almost surely and $\lim_{k\to\infty} \norms{z^k - x^{\star}}$ exists almost surely, we conclude that $\lim_{k\to\infty}\norms{x^k - x^{\star}}$ exists almost surely.

Let $\Dc = \sets{x_1^{\star}, x_2^{\star}, \cdots}$ be a countable dense set in $\zer{\Phi}$.
For each $x_j^{\star} \in \Dc$, we have proven $\mathbb{P}\big( \lim_{k\to\infty}\norms{x^k - x^{\star}_j} \ \text{exists}\big) = 1$.
If we define $\Omega_j := \sets{\omega : \lim_{k\to\infty}\norms{x^k(\omega) - x^{\star}_j} \ \text{exists}}$, then $\mathbb{P}(\Omega_j) = 1$.
Since $\Dc$ is countable, $\Omega := \bigcap_{j=1}^{\infty}\Omega_j$ also satisfies $\mathbb{P}(\Omega) = 1$.

Now, fix $\omega \in \Omega$.
For every $\hat{x}^{\star} \in \Dc$, $\lim_{k \to\infty}\norms{x^k(\omega) - \hat{x}^{\star}}$ exists.
For every $x^{\star} \in \zer{\Phi}$ and $\hat{x}^{\star} \in \Dc$, we have $\vert \norms{x^k - x^{\star}} - \norms{x^k - \hat{x}^{\star}}\vert \leq \norms{x^{\star} - \hat{x}^{\star}}$.
Therefore, for any $m, n \geq 0$, we have 
\begin{equation*}
\vert \norms{x^m - x^{\star}} - \norms{x^n - x^{\star}} \vert \leq \vert \norms{x^m - \hat{x}^{\star}} - \norms{x^n - \hat{x}^{\star}} \vert + 2\norms{x^{\star} - \hat{x}^{\star}}.
\end{equation*} 
For any $\epsilon > 0$, since $\Dc$ is dense in $\zer{\Phi}$, we can choose $\hat{x}^{\star} \in \Dc$ such that $\norms{x^{\star} - \hat{x}^{\star}} \leq \frac{\epsilon}{4}$.
On the other hand, since $\set{\norms{x^k - \hat{x}^{\star}}}$ converges, it is Cauchy.
Consequently, for sufficiently large $m$ and $n$, we have $\vert \norms{x^m - \hat{x}^{\star}} - \norms{x^n - \hat{x}^{\star}} \vert \leq \frac{\epsilon}{2}$.
Using these facts into the last inequality, we get $\vert \norms{x^m - x^{\star}} - \norms{x^n - x^{\star}} \vert \leq \epsilon$.
Thus, $\lim_{k\to\infty}\norms{x^k - x^{\star}}$ exists for every $x^{\star} \in \zer{\Phi}$.

From the second line of \eqref{eq:SFOG_th32_proof2}, we have $\lim_{k\to\infty}\norms{Gx^k + v^k} = 0$ almost surely, for $v^k \in Tx^k$.
This means that $\lim_{k\to\infty}\norms{w^k} = 0$ almost surely for $(x^k, w^k) \in \gra{\Phi}$.
Moreover, $\gra{\Phi}$ is closed by our assumption.
Applying Lemma~\ref{le:Opial_lemma} from the appendix, we conclude that  $\sets{x^k}$ converges almost surely to a $\zer{\Phi}$-valued random variable $x^{\star} \in \zer{\Phi}$.
Finally, since $\lim_{k\to\infty}\norms{x^k - z^k} = 0$ \revise{almost surely} and, \revise{by the last line of \eqref{eq:SFOG_th32_proof2}, $\lim_{k\to\infty}\norms{x^{k+1} - y^k} = 0$ almost surely,} we can also conclude that both $\sets{z^k}$ and $\sets{y^k}$ also converge almost surely to the same $x^{\star}$.
\Eproof

\beforesubsec
\subsection{\textbf{The full proof of Theorem~\ref{th:VFOG_convergence}.}}\label{subsection:th:VFOG_convergence}
\aftersubsec
Before presenting the full proof of Theorem~\ref{th:VFOG_convergence}, we need the following supporting lemma, whose proof is given in Appendix~\ref{apdx:le:VrAEG4GE_descent_property}.

\begin{lemma}\label{le:VrAEG4GE_descent_property}
Let $\hat{\omega}$, $\hat{\lambda}$, and $\hat{\mu}$ be defined by \eqref{eq:omega_quantity}.
For \eqref{eq:GE}, suppose that $G$ satisfies Assumptions~\ref{as:A1} and \ref{as:A3} with $\alpha \in [0, 1)$ and $\rho_c > 0$ such that $L\rho_n \leq \hat{\mu}$.
Let $\sets{(x^k, y^k, z^k)}$ be generated by \eqref{eq:VFOG} using a variance-reduced estimator $\widetilde{G}y^k$ for $Gy^k$ satisfying Definition~\ref{de:VR_Estimators}.
\revise{For a given} $s > 2$, \revise{let $\varphi_s := \frac{9(85s-134)}{64(s-1)} + \frac{9(s-2)}{64}$} and suppose that $\eta$, $\beta_k$, and $\gamma_k$ are respectively updated by
\begin{equation}\label{eq:SAEG_para_update_v1}
\arraycolsep=0.2em\hspace{-3ex}
\begin{array}{ll}
& \frac{\revise{8}(s-1)\rho_n}{\revise{s-2}} \leq \eta \leq  \frac{\hat{\lambda}}{L}, \quad \gamma_k := \frac{\revise{(s-2)}\eta(k+s)}{\revise{32(s-1)}(k+s+1)}, \quad\textrm{and} \quad  \beta_k :=  \big[\frac{\revise{3}(s-2)\eta}{\revise{16}(s-1)} + 2\rho_n]\frac{ k+1 }{k+s+1} - \frac{\gamma_k}{k+s+1}.
\end{array}\hspace{-3ex}
\end{equation}
Suppose further that   there exist $\bar{c} \in (0, 1)$ and a sufficiently small $\epsilon_2 \geq 0$ such that for every $k \geq 0$, the parameters $\kappa_k$ and $\Theta_k$ in \eqref{eq:VR_property} of Definition~\ref{de:VR_Estimators} satisfy the following condition:
\begin{equation}\label{eq:SAEG_error_cond3}
\arraycolsep=0.2em
\begin{array}{lcl}
\revise{\big[\varphi_s + \frac{25s-34}{2(s-2)}\big]}\frac{ 3\eta \Theta_k t_k}{\rho_c(t_k - s)} \leq \bar{c} \leq \frac{1}{1 + \epsilon_2} \Big[ 1 - \frac{(1- \kappa_k)t_k^2}{t_{k-1}^2} \Big]. 
\end{array}
\end{equation}
For $\bar{c}$ in \eqref{eq:SAEG_error_cond3} and $a_k$ in Lemma~\ref{le:VFOG_key_estimate2},  define the following \textbf{Lyapunov function}:
\begin{equation}\label{eq:SAEG_Lyapunov_func1}
\arraycolsep=0.2em
\begin{array}{lcl}
\hat{\Lc}_k &:= & \frac{a_k}{2}\norms{w^k}^2 + st_{k-1}\iprods{w^k, x^k - z^k} + \frac{s^2(s-1)}{2\gamma_{\revise{k-1}}}\norms{z^k - x^{\star}}^2 + \rho_c  t_{k-1}^2 \tnorms{ w^k -  \hat{w}^k }^2  \vspace{1ex}\\
&& + {~} \frac{\revise{\varphi_s}\eta L^2\mred{t_{k-1}^2}}{\revise{8}} \norms{x^k - y^{k-1} }^2 +  \revise{\big[\varphi_s + \frac{25s-34}{2(s-2)}\big]}\frac{(1-\bar{c})\eta t_{k-1}^2}{\bar{c}} \Delta_{k-1} +  \revise{\varphi_s}\eta t_{k-1}^{\revise{2}}\norms{e^{k-1}}^2.
\end{array}
\end{equation}
Then, the following inequality holds:
\begin{equation}\label{eq:SAEG_key_est5}
\hspace{-0ex}
\arraycolsep=0.1em
\begin{array}{lcl}
\hat{\Lc}_k -  \Expsn{k}{ \hat{\Lc}_{k+1}}  & \geq & \frac{\eta s(k+s+1)}{2} \Expsn{k}{ \norms{\hat{w}^{k+1} }^2 }  + \frac{\phi_k}{2} \norms{ w^k }^2 +  s(s-1)\iprods{w^k, x^k - x^{\star}} \vspace{1ex}\\
&& + {~} \frac{\eta t_k(t_k-s)}{2}\big(1  - \frac{ M^2\eta^2 t_k }{t_k - s} \big) \Expsn{k}{ \norms{ \hat{w}^{k+1} - w^k}^2 } \vspace{1ex}\\
&& + {~}  \frac{2\alpha\rho_c t_k^2}{1-\alpha} \Expsn{k}{ \norms{ w^{k+1} - \hat{w}^{k+1} }^2 } + \frac{\revise{(s-2)}\eta L^2 \revise{t_k^2}}{\revise{32(s-1)}} \Expsn{k}{\norms{ \revise{x^{k+1} - y^k}  }^2}  \vspace{1ex}\\ 
&& + {~} \revise{\big[\varphi_s + \frac{25s-34}{2(s-2)}\big]} \epsilon_2 \eta t_{k-1}^2\Delta_{k-1} - \revise{\big[\varphi_s + \frac{25s-34}{2(s-2)}\big] \frac{\eta \delta_k}{\bar{c}}},  
\end{array}
\hspace{-0ex}
\end{equation}
where $\phi_k := \big[ \revise{\frac{(20s-23)(s-2)}{32(s-1)}}\eta - 4(s-1)\rho_n \big](k + s) + (s-1)\big[ \frac{(\revise{3s+10})\eta}{\revise{16}} + 2(2s-1)\rho_n\big] + \frac{\revise{(s-2)}\eta}{\revise{32(s-1)}}$.
\end{lemma}

\vspace{1ex}
\noindent{\textbf{The proof of Theorem~\ref{th:VFOG_convergence}}.}\label{subsec:VR_method}
For the sake of clarity, we divide this proof into several steps.

\vspace{0.75ex}
\noindent\textbf{\underline{Step 1} (\textit{Verifying condition~\eqref{eq:SAEG_kappa_Theta_cond}}).}
First, the update rules of parameters in \eqref{eq:SAEG_para_update2} come from \eqref{eq:SAEG_para_update_v1} of Lemma~\ref{le:VrAEG4GE_descent_property}.
Next, we verify the condition \eqref{eq:SAEG_error_cond3} of Lemma~\ref{le:VrAEG4GE_descent_property}.
Since we can choose $\epsilon_2$ sufficiently small, \eqref{eq:SAEG_error_cond3} holds if 
\begin{equation*} 
\arraycolsep=0.2em
\begin{array}{lcl}
\revise{\big[\varphi_s + \frac{25s-34}{2(s-2)}\big]}\frac{3\eta \Theta_k t_k}{\rho_c(t_k - s)} < 1 - \frac{(1- \kappa_k)t_k^2}{t_{k-1}^2}.
\end{array}
\end{equation*}
If we denote $\hat{C} := 3\revise{\big[\varphi_s + \frac{25s-34}{2(s-2)}\big]}$, then since $t_k = k+s+1$, this condition becomes
\begin{equation*} 
\arraycolsep=0.2em
\begin{array}{lcl}
\frac{\hat{C}\eta \Theta_k (k+s+1)}{\rho_c(k+1)} <   1 - \frac{(1- \kappa_k)(k+s+1)^2}{(k+s)^2} = \frac{\kappa_k(k+s+1)^2}{(k+s)^2} - \frac{2(k+s) + 1}{(k+s)^2}.
\end{array}
\end{equation*}
Because $k\geq 0$ and $2(k+s) + 1 < 2(k+s+1)$, this condition holds if 
\begin{equation*} 
\arraycolsep=0.2em
\begin{array}{lcl}
\kappa_k \geq \frac{ s^2 \hat{C} \eta \Theta_k }{(s+1) \rho_c } + \frac{2}{k+s+1} =  \frac{\Gamma \eta \Theta_k}{\rho_c} + \frac{2}{k+s+1},
\end{array}
\end{equation*}
where $\Gamma  := \revise{\frac{3s^2}{s+1}\big[\varphi_s + \frac{25s-34}{2(s-2)}\big]}$ is given in \eqref{eq:omega_quantity}.
This condition is exactly \eqref{eq:SAEG_kappa_Theta_cond}.

Note that if there exist $\underline{\kappa}$ and $\Theta > 0$ such that $0 < \underline{\kappa} \leq \kappa_k < 1$ and $0 < \Theta_k \leq \Theta$, then under the condition \eqref{eq:SAEG_kappa_Theta_cond}, we can easily show that there exists $\epsilon_2 > 0$ such that 
\begin{equation*} 
\arraycolsep=0.2em
\begin{array}{lcl}
\frac{\rho_c(t_k - s)}{3 \revise{\big[\varphi_s + \frac{25s-34}{2(s-2)}\big]} \eta \Theta_k t_k}\big[1 - \frac{(1- \kappa_k)t_k^2}{t_{k-1}^2}\big] - 1 \geq \epsilon_2 >  0.
\end{array}
\end{equation*}
\vspace{0.75ex}
\noindent\textbf{\underline{Step 2} (\textit{The summability results in \eqref{eq:SAEG_summability_bounds}}).}
By \eqref{eq:weak_Minty} in Assumption~\ref{as:A2}, we have $\iprods{ w^k, x^k - x^{\star}} \geq - \rho_{*}\norms{w^k}^2$.
Moreover, \revise{similar to the proof of Lemma~\ref{le:SFOG_supporting_lemma0}}, we also have
\begin{equation*}
\arraycolsep=0.2em
\begin{array}{lcl}
\phi_k - 2 \rho_{*}s(s-1) & = &  \big[ \revise{\frac{(20s-23)(s-2)}{32(s-1)}}\eta - 4(s-1)\rho_n \big](k + s ) + \frac{(s-1)\revise{(3s+10)}\eta}{16} \vspace{1ex}\\
&& + {~} 2(2s-1)(s-1)\rho_n + \revise{\frac{(s-2)\eta}{32(s-1)}}  - 2\rho_{*} s(s-1) \vspace{1ex}\\

& \geq & \revise{\frac{(s-2)(4s-7)\eta(k + 1)}{32(s-1)}} > 0.
\end{array}
\end{equation*}
Substituting these facts into  \eqref{eq:SAEG_key_est5} in Lemma~\ref{le:VrAEG4GE_descent_property}, we get
\begin{equation}\label{eq:SAEG_th31_proof1a}
\hspace{-2ex}
\arraycolsep=0.2em 
\begin{array}{lcl}
\hat{\Lc}_k -  \Expsn{k}{ \hat{\Lc}_{k+1}}  & \geq & \frac{\eta s(k+s+1)}{2} \Expsn{k}{\norms{ \hat{w}^{k+1} }^2}  + \revise{\frac{(s-2)(4s-7)\eta(k + 1)}{64(s-1)}} \norms{ w^k }^2 \vspace{1ex}\\
&& + {~} \frac{\eta t_k(t_k-s)}{2}\big(1  - \frac{ M^2\eta^2 t_k }{t_k - s} \big) \Expsn{k}{ \norms{ \hat{w}^{k+1} - w^k}^2} + \frac{\revise{(s-2)}\eta L^2 \revise{t_k^2}}{\revise{32(s-1)}} \Expsn{k}{ \norms{ \revise{x^{k+1} - y^k}  }^2} \vspace{1ex}\\
&& +  {~}  \frac{2\alpha\rho_c t_k^2}{1-\alpha} \Expsn{k}{ \norms{ w^{k+1} - \hat{w}^{k+1} }^2 } + \revise{\big[\varphi_s + \frac{25s-34}{2(s-2)}\big]} \epsilon_2 \eta t_{k-1}^2\Delta_{k-1} - \revise{\big[\varphi_s + \frac{25s-34}{2(s-2)}\big] \frac{\eta \delta_k}{\bar{c}}}. 
\end{array}
\hspace{-4ex}
\end{equation}
Now, taking the full expectation $\Expn{\cdot}$ on both sides of \eqref{eq:SAEG_th31_proof1a}, we obtain
\begin{equation}\label{eq:SAEG_th31_proof1}
\hspace{-3ex}
\arraycolsep=0.1em
\begin{array}{lcl}
\Expn{\hat{\Lc}_k } -  \Expn{ \hat{\Lc}_{k+1}}  & \geq & \frac{\eta s(k+s+1)}{2} \Expn{\norms{ \hat{w}^{k+1} }^2}  + \revise{\frac{(s-2)(4s-7)\eta (k + 1)}{64(s-1)}} \Expn{ \norms{ w^k }^2 } \vspace{1ex}\\
&& + {~} \frac{\revise{(s-2)}\eta L^2 \revise{t_k^2}}{\revise{32(s-1)}}\Expn{ \norms{ \revise{x^{k+1} - y^k} }^2 } + \revise{\big[\varphi_s + \frac{25s-34}{2(s-2)}\big]} \epsilon_2 \eta t_{k-1}^2 \Expn{ \Delta_{k-1} } \vspace{1ex}\\
&& + {~} \frac{\eta t_k(t_k-s)}{2}\big(1  - \frac{ M^2\eta^2 t_k }{t_k - s} \big) \Expn{ \norms{ \hat{w}^{k+1} - w^k}^2} \vspace{1ex}\\
&& + {~} \frac{2\alpha\rho_c t_k^2}{1-\alpha} \Expn{ \norms{ w^{k+1} - \hat{w}^{k+1} }^2 } - \revise{\big[\varphi_s + \frac{25s-34}{2(s-2)}\big]}\frac{\eta \delta_k}{\bar{c}}.  
\end{array}
\hspace{-6ex}
\end{equation}
Summing \eqref{eq:SAEG_th31_proof1} from $k=0$ to $K$ and noting that $\hat{\Lc}_k \geq 0$ due to \eqref{eq:SAEG_Lyapunov_func1} and Lemma~\ref{le:Ek_lowerbound}, we have
\begin{equation}\label{eq:SAEG_th31_proof2}
\arraycolsep=0.2em
\begin{array}{lcl}
\sum_{k=0}^K \eta s(k+s+1) \Exp{\norms{ \hat{w}^{k+1} }^2} & \leq & 2 (\hat{\Lc}_0 \mred{+ \hat{\Lambda}_0\eta \Sc_K}), \vspace{1ex}\\
\sum_{k=0}^K \revise{\frac{(s-2)(4s-7)\eta}{s-1}} (k+1) \Exp{\norms{ w^k }^2} & \leq & \revise{64} (\hat{\Lc}_0 \mred{+ \hat{\Lambda}_0\eta \Sc_K}),
\end{array}
\end{equation}
\mred{where $\hat{\Lambda}_0$ is defined in \eqref{eq:omega_quantity} and $\Sc_K := \sum_{k=0}^K \delta_k$.}
Here, we have used the fact that $\bar{c} \geq \frac{24(s^2-1)\rho_n\underline{\Theta}}{(s-2)\rho_c} \big[ \varphi_s + \frac{25s - 34}{2(s-2)} \big]$ derived from  \eqref{eq:SAEG_error_cond3}, leading to $\hat{\Lambda}_0$ as in  \eqref{eq:omega_quantity}.

Finally, since $z^0 = x^0 = y^{-1}$ and the conventions $e^{-1} = 0$, and $\Delta_{-1} = 0$, we have from \eqref{eq:SAEG_Lyapunov_func1} that
\begin{equation*} 
\arraycolsep=0.2em
\begin{array}{lcl}
\hat{\Lc}_0 & = & \frac{a_0}{2}\norms{ w^0 }^2  + \frac{s^2(s-1)}{2\gamma_{\revise{-1}}}\norms{x^0 - x^{\star}}^2 
= \frac{\revise{[2(13s-10)s^2 - (s-1)(s-2)]}\eta}{\revise{64(s-1)}} \norms{ w^0 }^2 + \revise{\frac{16s^3(s-1)}{(s-2)\eta}}\norms{x^0 - x^{\star} }^2.
\end{array}
\end{equation*}
Here, we have used $a_ 0 = [\frac{\revise{(13s-10)}\eta}{\revise{16}(s-1)} - 2\rho_n]s^2 + 2s^2\rho_n  -  \frac{\revise{(s - 2)} \eta}{\revise{32}} = \frac{\revise{[2(13s-10)s^2 - (s-1)(s-2)]}\eta}{\revise{32}(s-1)} > 0$ for all $s > 2$ and $\gamma_{\revise{-1}} = \frac{\revise{(s-2)}\eta(t_{\revise{-1}}-1)}{\revise{32(s-1)}t_{\revise{-1}}} = \revise{\frac{{(s-2)\eta}}{32s}}$.
Substituting $\hat{\Lc}_0$ into \eqref{eq:SAEG_th31_proof2}, we obtain \eqref{eq:SAEG_summability_bounds}.

\vspace{0.75ex}
\noindent\textbf{\underline{Step 3} (\textit{The proof of the bounds in \eqref{eq:SAEG_BigO_rates}}).}
Now, from \eqref{eq:SAEG_th31_proof1}, by induction, we also have $\mred{\Expn{\hat{\Lc}_k } \leq \hat{\Lc}_0 + \hat{\Lambda}_0 \eta \Sc_k}$.  
However, from \mred{\eqref{eq:SAEG_Lyapunov_func1} and} \eqref{eq:Ek_lowerbound}, we also have
\begin{equation}\label{eq:SAEG_th31_proof1_100} 
\arraycolsep=0.2em
\begin{array}{lcl}
\mred{\hat{\Lc}}_k \geq \Pc_k & \geq & \frac{A_k}{2}\norms{ w^k }^2 + \revise{\frac{s^2}{2}\big(\frac{s-1}{\gamma} - \frac{4}{\eta}\big)}\norms{z^k - x^{\star} }^2 \geq \frac{A_k}{2}\norms{ w^k }^2,
\end{array}
\end{equation}
where, similarly to the proof of Theorem~\ref{th:VFOG1_convergence}, 
\begin{equation*}
\begin{array}{lcl}
	A_k & := & \big[\revise{\frac{3(3s-2)\eta}{16(s-1)}} - 2\rho_n\big]t_{k-1}^2 + 2s(\rho_n - \rho_{*}) t_{k-1} -  \gamma t_{k-1} + \gamma 
		 \geq \revise{\frac{(19s+10)\eta}{64(s-1)}(k+s)^2,}
\end{array}
\end{equation*}
Therefore, we obtain 
\begin{equation*} 
\arraycolsep=0.2em
\begin{array}{lcl}
\Expn{ \norms{ w^k }^2 } \leq \frac{\revise{128}(s-1)(\hat{\Lc}_0 + \mred{\hat{\Lambda}_0 \eta \Sc_k})}{\revise{(19s+10)\eta} (k+s)^2 } = \frac{\revise{128}(s-1) (\hat{\mcal{R}}_0^2 + \mred{\hat{\Lambda}_0 \Sc_k})}{\revise{(19s+10)} (k+s)^2 },
\end{array}
\end{equation*}
which is the first bound of \eqref{eq:SAEG_BigO_rates}, where 
\begin{equation*}
\arraycolsep=0.2em
\begin{array}{lcl}
\hat{\mcal{R}}_0^2 & := & \frac{\revise{2(13s-10)s^2 - (s-1)(s-2)}}{\revise{64(s-1)}} \norms{ w^0 }^2 + \revise{\frac{16s^3(s-1)}{(s-2)\eta^2}}\norms{x^0 - x^{\star} }^2.
\end{array}
\end{equation*}
From \eqref{eq:SAEG_Lyapunov_func1}, using
\begin{equation*}
\arraycolsep=0.2em
\begin{array}{lcl}
\revise{\frac{\varphi_s\eta}{4} + 2\rho_c \geq \frac{\varphi_s \eta}{4} \geq \frac{9(85s - 134)\eta}{256(s-1)} \geq \frac{3(3s-2)\eta}{16(s-1)} \stackrel{\eqref{eq:SFOG_th31_proof2.1}}{\geq} \frac{A_k}{t_{k-1}^2}},
\end{array}
\end{equation*}
\mred{the $L$-Lipschitz continuity of $G$, the relation $\tnorms{\cdot} \geq \norms{\cdot}$, and Young's inequality,} we have
\begin{equation*}
\arraycolsep=0.2em
\begin{array}{lcl}
\Tc_{[1]} &:= &\frac{A_k}{2}\norms{w^k}^2 + \frac{\revise{\varphi_s} \eta L^2 t_{k-1}^2}{8} \norms{x^k - y^{k-1}}^2 + \rho_c t_{k-1}^2 \tnorms{w^k - \hat{w}^k}^2 \vspace{1ex}\\
&\geq& \frac{A_k}{2}\norms{w^k}^2 +  \big( \frac{\revise{\varphi_s}\eta}{8} + \rho_c \big) t_{k-1}^2\norms{w^k - \hat{w}^k}^2  \vspace{1ex} \\
& \geq  & \frac{A_k}{4}\norms{ \hat{w}^k }^2 \vspace{1ex}\\
& \geq &  \revise{\frac{(19s+10)\eta}{256(s-1)}} (k+s)^2\norms{ \hat{w}^k }^2.
\end{array}
\end{equation*}
Combining this inequality, \eqref{eq:Ek_lowerbound}, and \eqref{eq:SAEG_Lyapunov_func1}, we can show that
\begin{equation*}
\arraycolsep=0.2em
\begin{array}{lcl}
\hat{\Lc}_k &\geq & \Pc_k + \frac{\revise{\varphi_s}\eta L^2 t_{k-1}^2}{8} \norms{x^k - y^{k-1}}^2 + \rho_c t_{k-1}^2 \tnorms{w^k - \hat{w}^k}^2 \vspace{1ex}\\
&\geq& \revise{\frac{(19s+10)\eta}{256(s-1)}} (k+s)^2\norms{ \hat{w}^k }^2.
\end{array}
\end{equation*}
Hence, utilizing this inequality, with  similar arguments as  the proof of the first bound of \eqref{eq:SAEG_BigO_rates}, we obtain the second line of \eqref{eq:SAEG_BigO_rates} from this expression. 

\vspace{0.75ex}
\noindent\textbf{\underline{Step 4} (\textit{Other intermediate summability results}).}
\mred{Next, rearranging \eqref{eq:SAEG_th31_proof1}, we have
\begin{equation*}
	\arraycolsep=0.1em
	\begin{array}{lcl}
		\Expn{ \hat{\Lc}_{k+1}}  & \leq & \Expn{\hat{\Lc}_k } - \frac{\eta s(k+s+1)}{2} \Expn{\norms{ \hat{w}^{k+1} }^2}  - \revise{\frac{(s-2)(4s-7)\eta(k + 1)}{64(s-1)}}  \Expn{ \norms{ w^k }^2 } \vspace{1ex}\\
		&& - {~} \frac{\revise{(s-2)}\eta L^2 \revise{t_k^2}}{\revise{32(s-1)}}\Expn{ \norms{ \revise{x^{k+1} - y^k} }^2 } - \revise{\big[\varphi_s + \frac{25s - 34}{2(s-2)}\big]}\epsilon_2 \eta t_{k-1}^2 \Expn{ \Delta_{k-1} }  - \frac{2\alpha\rho_c t_k^2}{1-\alpha} \Expn{ \norms{ w^{k+1} - \hat{w}^{k+1} }^2 } \vspace{1ex}\\
		&& - {~} \frac{\eta t_k(t_k-s)}{2}\big(1  - \frac{ M^2\eta^2 t_k }{t_k - s} \big) \Expn{ \norms{ \hat{w}^{k+1} - w^k}^2} + \revise{\big[\varphi_s + \frac{25s - 34}{2(s-2)}\big]}\frac{\eta \delta_k}{\bar{c}}. 
	\end{array}
\end{equation*}
Since $\epsilon_2 > 0$ and $\eta < \mred{\frac{\hat{\lambda}}{L}}$ as in \eqref{eq:SAEG_para_update2}, the terms from the second to seventh on the right-hand side of this inequality are all nonnegative.
Under the assumption $\Sc_{\infty} := \sum_{k=0}^\infty \delta_k < +\infty$, we can apply Lemma \ref{le:A1_descent} to the last inequality and obtain that 
\begin{equation}\label{eq:SAEG_summability_proof1}
\arraycolsep=0.2em
\begin{array}{lcl}
\sum_{k=0}^{\infty} t_k^2 \Expn{ \norms{ x^{k+1} - y^k }^2 } & < &  +\infty, \vspace{1ex}\\
\sum_{k=0}^{\infty} t_k(t_k-s)  \Expn{ \norms{  \hat{w}^{k+1} - w^k}^2 }  & < & +\infty, \vspace{1ex}\\
\sum_{k=0}^{\infty} t_k^2 \Expn{\Delta_k}  & < &  +\infty.
\end{array}
\end{equation}
\revise{From \eqref{eq:G_Lipschitz} of Assumption~\ref{as:A1}, for any $\alpha \in [0,1)$, we have $\tnorms{ w^{k+1} - \hat{w}^{k+1} }^2 \leq \frac{L^2}{1-\alpha} \norms{x^{k+1} - y^k}^2$.}
Using \revise{this bound,} the relations $\tnorms{\cdot} \geq \norms{\cdot}$, and $\Expn{\norms{e^k}^2} \leq \Expn{\Delta_k}$ from the first line of \eqref{eq:VR_property} after taking total expectations, we also obtain from \eqref{eq:SAEG_summability_proof1} that
\begin{equation}\label{eq:SAEG_summability_proof1b}
	\arraycolsep=0.2em
	\begin{array}{lcl}
		\sum_{k=0}^{\infty} t_k^2 \Expn{ \tnorms{ w^{k+1} - \hat{w}^{k+1} }^2 } & < &  +\infty, \vspace{1ex}\\
		\sum_{k=0}^{\infty} t_k^2 \Expn{ \norms{ w^{k+1} - \hat{w}^{k+1} }^2 } & < &  +\infty, \vspace{1ex}\\
		\sum_{k=0}^{\infty} t_k^2 \Expn{ \norms{e^k}^2 }  & < &  +\infty.
	\end{array}
\end{equation}}
By Young's inequality, we easily get
\begin{equation*} 
\arraycolsep=0.2em
\begin{array}{lcl}
\norms{w^{k+1} - w^k }^2 & \leq & 2\norms{w^{k+1} - \hat{w}^{k+1} }^2 + 2\norms{ \hat{w}^{k+1} - w^k}^2, \vspace{1ex}\\
\norms{ \hat{w}^{k+1} - \hat{w}^k }^2 & \leq & 2\norms{w^k - \hat{w}^{k+1} }^2 + 2\norms{w^k - \hat{w}^k }^2.
\end{array}
\end{equation*}
Moreover, we also have $\Expn{\norms{\hat{e}^k}^2} \leq 2 \Expn{\norms{ e^{k-1} }^2 }  + 2\Expn{ \norms{w^k - \hat{w}^k  }^2 }$.
Using these relations, we can obviously derive from \eqref{eq:SAEG_summability_proof1} \mred{and \eqref{eq:SAEG_summability_proof1b}} that
\begin{equation}\label{eq:SAEG_summability_proof2}
\arraycolsep=0.2em
\begin{array}{lcl}
\sum_{k=0}^{\infty} t_k^2 \Expn{ \norms{ w^{k+1} - w^k }^2 } & < &  +\infty, \vspace{1ex}\\
\sum_{k=0}^{\infty} t_k^2 \Expn{ \norms{ \hat{w}^{k+1} - \hat{w}^k }^2 } & < &  +\infty, \vspace{1ex}\\
\sum_{k=0}^{\infty} t_k^2 \Expn{\norms{ \hat{e}^k }^2 }  & < &  +\infty.
\end{array}
\end{equation}

\vspace{0.75ex}
\noindent\textbf{\underline{Step 5} (\textit{Applying Lemmas~\ref{le:SFOG_supporting_lemma1} and \ref{le:SFOG_supporting_lemma2}}).}
The proof of \eqref{eq:SAEG_small_o_rates} follows the same arguments as the proof of \textbf{Step 4} in Theorem~\ref{th:VFOG1_convergence} and we omit the details.
The iteration-complexity stated in Statement (d) is also proven as in Theorem~\ref{th:VFOG1_convergence}, and we omit it here.

\vspace{0.75ex}
\noindent\textbf{\underline{Step 6} (\textit{Almost sure convergence}).}
The almost sure convergence results in Statement~(e) follow from an argument similar to the proof of Theorem \ref{th:VFOG1_almost_sure_convergence}, and hence we omit the details.
\Eproof

\beforesec
\section{Numerical Experiments.}\label{sec:numerical_experiments}
\aftersec
In this section, we present three numerical examples to validate our algorithm, \eqref{eq:VFOG}, and compare it with closely related methods.
Our Python implementation was run on an Apple M4 MacBook Pro (24 GB RAM, 1 TB SSD).

\beforesubsec
\subsection{\textbf{Bilinear matrix game: Policeman vs. Burglar problem.}}\label{subsec:bilinear_matrix_game}
\aftersubsec
Consider a city with $m^2$ houses distributed on a grid $m \times m$, where the $j$-th house has wealth $w_j$ for all $j = 1, \cdots, m^2$. 
Each night, a Burglar selects a house $j$ to rob, while a Policeman chooses a location to station himself near house $k$. 
Once the burglary begins, the Policeman is immediately informed of the burglary location, and the probability that he apprehends the Burglar is given by $\mbf{p}_c := \exp\{-\theta\, \mathrm{dist}(j, k)\}$, where $\mathrm{dist}(j, k)$ is the distance between houses $j$ and $k$. 
However, neither the Burglar nor the Policeman has access to the true wealth $w_j$; instead, they need to estimate it using a set of $n$ observations ${\hat{w}^{(i)}_{j}}$ for $i \in [n]$.
The Burglar aims to maximize his expected reward, which is $\mcal{R}_{jk} := \frac{1}{n} \sum_{i=1}^n {\hat{w}^{(i)}_{j}} \left(1 - \exp\{-\theta\, \mathrm{dist}(j, k)\}\right)$,
while the Policeman’s objective is to minimize $\mcal{R}_{jk}$.

\vspace{1ex}
\noindent\textbf{\textit{$\mathrm{(a)}$~Mathematical model.}}
The Policeman vs. Burglar problem described above can be viewed as a zero-sum game.
It was studied in \cite{nemirovski2013mini} and we slightly modify it to obtain a stochastic version.
Following  \cite{nemirovski2013mini}, we define a payoff matrix $\mbf{L}$ as follows:
\begin{equation*}
\begin{array}{lcl}
	\mbf{L} := \frac{1}{n}\sum_{i=1}^n \mbf{L}^{(i)}, \quad \mbf{L}^{(i)}_{j,k} = {\hat{w}^{(i)}_{j}} \left(1 - \exp\{-\theta\, \mathrm{dist}(j, k)\}\right), \quad 1 \le j, k \le m^2.
\end{array}
\end{equation*}
Then, we can formulate the above problem as the following classical bilinear matrix game: 
\begin{equation}\label{eq:bilinear_matrix_game}
\min_{u \in \Delta_{p_1}}\max_{v \in \Delta_{p_2}}\Big\{ \Lc(u, v) := \iprods{\mbf{L}u, v} \Big\},
\end{equation}
where $u$ and $v$ represent mixed strategies of the Policeman and Burglar,  and $\Delta_{p_1}$ and $\Delta_{p_2}$ are the standard simplexes in $\R^{p_1}$ and $\R^{p_2}$, respectively with $p_1=p_2=m^2$.

\mred{Let us} define $x := [u;v] \in \R^p$ with $p := p_1+p_2$, $G_i x := \mbf{G}_ix =  [\mbf{L}^{(i)\top}v; -\mbf{L}^{(i)}u]$ for $i \in [n]$, and $Tx := [\partial \delta_{\Delta_{p_1}}(u); \partial \delta_{\Delta_{p_2}}(v)]$, where $\delta_{\Xc}(\cdot)$ is the indicator of $\Xc$.
Then, we can write the optimality condition of \eqref{eq:bilinear_matrix_game} as $0 \in \frac{1}{n}\sum_{i=1}^n G_ix + Tx$, a special case of \eqref{eq:GE}.
Since $\mbf{G} := \frac{1}{n}\sum_{i=1}^n \mbf{G}_i$ is a skew-symmetric matrix, $G$ is monotone and $L$-Lipschitz continuous.

\vspace{1ex}
\noindent\textit{\textbf{$\mathrm{(b)}$~Input data.}} 
Following \cite{nemirovski2013mini}, we first choose $\mathrm{dist}(j, k) := \vert j - k\vert$ and $\theta := 0.8$. 
Next, we generate the nominal wealth $w_j$ from a standard normal distribution, and then take $\vert w_j\vert$ to ensure its nonnegativity.
Then, the noisy wealth $w^{(i)}_j$ is generated as $w^{(i)}_j := |w_j + \epsilon^{(i)}_j|$, where  $\epsilon^{(i)}_j$ is a normal random variable with zero mean and variance $\sigma^2 = 0.05$.
Finally, we produce two sets of experiments, each containing $10$ problem instances, as follows.
\begin{compactitem}
\item \textbf{Experiment 1.} $m = 10$ and $n = 1000$, corresponding to $p_1=100$ and $p = 2p_1 = 200$.
\item \textbf{Experiment 2.} $m = 15$ and $n = 2000$, corresponding to $p_1=225$ and $p = 2p_1 = 450$.
\end{compactitem}

\vspace{1ex}
\noindent\textit{\textbf{$\mathrm{(c)}$~Algorithms and parameters.}} 
We choose three competitors as follows.
\begin{compactitem}
\item The deterministic optimistic gradient method, e.g., in \cite{daskalakis2018training}, abbreviated by \texttt{OG}.
We consider this method as a baseline for our comparison.
\item The variance-reduced forward-reflected-backward splitting algorithm in \cite{alacaoglu2021forward}, called \texttt{VrFRBS}.
\item The variance-reduced extragradient (EG) algorithm in \cite{alacaoglu2021stochastic}, abbreviated by \texttt{VrEG}.
\end{compactitem}
Both \texttt{VrFRBS} and \texttt{VrEG} are variants of the EG method and use a Loopless SVRG estimator.

We implement four variants of our method \eqref{eq:VFOG}: \texttt{VAPEG-Sgd},  \texttt{VAPEG-Svrg}, \texttt{VAPEG-Saga}, and \texttt{VAPEG-Sarah}, and compare them with \texttt{OG}, \texttt{VrFRBS}, and \texttt{VrEG}.
Here, \texttt{VAPEG-Sgd}  uses the mini-batch estimator \eqref{eq:mini_batch_Gy}, while the last three variants use L-SVRG, SAGA, and L-SARAH estimators, respectively. 
We run all algorithms for 200 epochs, where each epoch corresponds to $n$ calls of $\mbf{G}$.

The parameters of these algorithms are chosen as follows.
For all four variants of our methods, we choose  $\eta = \frac{1}{8L}$, (i.e., $s = 3$ and $\omega = 7$).
For \texttt{OG}, we choose $\eta = \frac{1}{L}$.
 For \texttt{VrFRBS}, we choose $\eta = \frac{0.95(1 - \sqrt{1 - \mbf{p}})}{2L}$, and for \texttt{VrEG}, we choose $\eta = \frac{0.95\sqrt{1 - \alpha}}{L}$ with $\alpha := 1 - \mbf{p}$.

For the probability of reference points and mini-batch size,  we choose:
\begin{compactitem}
\item $\mbf{p} = 0.5n^{-1/3}$  for \texttt{VAPEG-Svrg}, \texttt{VrFRBS}, and \texttt{VrEG}.
For example, if $n=1000$, then $\mbf{p} = 0.05$, and if $n=2000$, then $\mbf{p} = 0.04$.

\item $\mbf{p} = 0.5n^{-1/2}$  \mred{for \texttt{VAPEG-Sarah}.}
If $n= 1000$, then $\mbf{p} = 0.016$, while if $n=2000$, then $\mbf{p} = 0.011$.

\item $b_k = \lfloor 0.5n^{2/3} \rfloor$   for \texttt{VAPEG-Svrg}, \texttt{VAPEG-Saga}, \texttt{VrFRBS}, and \texttt{VrEG}.
In this case, if $n=1000$, then $b_k = 50$, and if $n=2000$, then $b_k = 79$.

\item $b_k = \lfloor 0.5n^{1/2} \rfloor$ for \texttt{VAPEG-Sarah} (e.g., $b_k = \mred{15}$ for $n=1000$ and $b_k = 22$ for $n=2000$).

\item For \texttt{VAPEG-Sgd}, we choose increasing mini-batch size $b_k = \max\sets{5, \min\sets{0.05(l+1)^3, n}}$, where $l$ is the epoch counter.
This choice makes $b_k$ increasing, but not too fast since $b_k = \BigOs{(k+s)^3}$.
\end{compactitem}
The initial point is chosen as $x^0 := [\frac{1}{p_1}\cdot\texttt{ones($p_1$)}; \frac{1}{p_2}\cdot\texttt{ones($p_2$)}]$ for all methods.

\begin{figure}[hpt!]
	\vspace{-0ex}
	\centering
	\includegraphics[width=\textwidth]{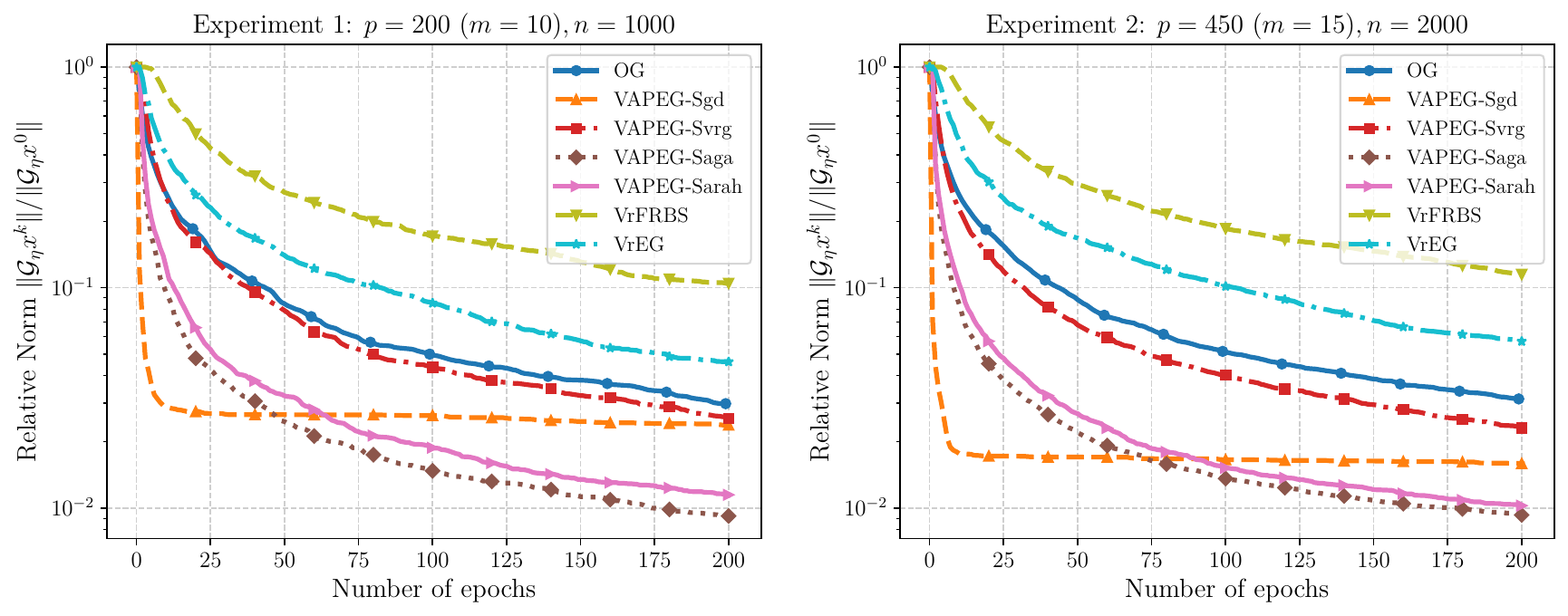}
	\vspace{-3ex}
	\caption{
	The performance of the seven algorithms for solving the bilinear matrix game \eqref{eq:bilinear_matrix_game} in two experiments, averaging over $10$ problem instances each.
	Here, $\mcal{G}_{\eta}$ is the FBS residual defined by \eqref{eq:FBS_residual}.
	}
	\label{fig:exam1_matrix_games_results1}
	\vspace{-2ex}
\end{figure}

\vspace{1ex}
\noindent\textit{\textbf{$\mathrm{(d)}$~Numerical results.}}
The  results of these experiments are presented in Figure \ref{fig:exam1_matrix_games_results1}.

From Figure \ref{fig:exam1_matrix_games_results1}, we can see that with the theoretical parameters, two \ref{eq:VFOG} variants: \texttt{VAPEG-Saga} and \texttt{VAPEG-Sarah} significantly outperform the competitors in both experiments.
However, \texttt{VAPEG-Saga} is slightly better than \texttt{VAPEG-Sarah} since \texttt{VAPEG-Saga} uses only one batch $\Bc_k$ as opposed to two in \texttt{VAPEG-Sarah} per iteration.
As a compensation, \texttt{VAPEG-Saga} requires storing a table of $n$ operator values.
Our \texttt{VAPEG-Svrg} also works well, but it is worse than \texttt{VAPEG-Saga} and \texttt{VAPEG-Sarah}.
With its theoretical parameters, \texttt{VrFRBS} exhibits the poorest performance among the tested methods.
\texttt{VrEG} is still worse than the baseline scheme \texttt{OG}.
Our \texttt{VAPEG-Sgd} performs well and better than other methods in early epochs, but then it plateaus at a certain accuracy level, perhaps due to the effect of variance. 

\begin{figure}[hpt!]
	\vspace{-2ex}
	\centering
	\includegraphics[width=\textwidth]{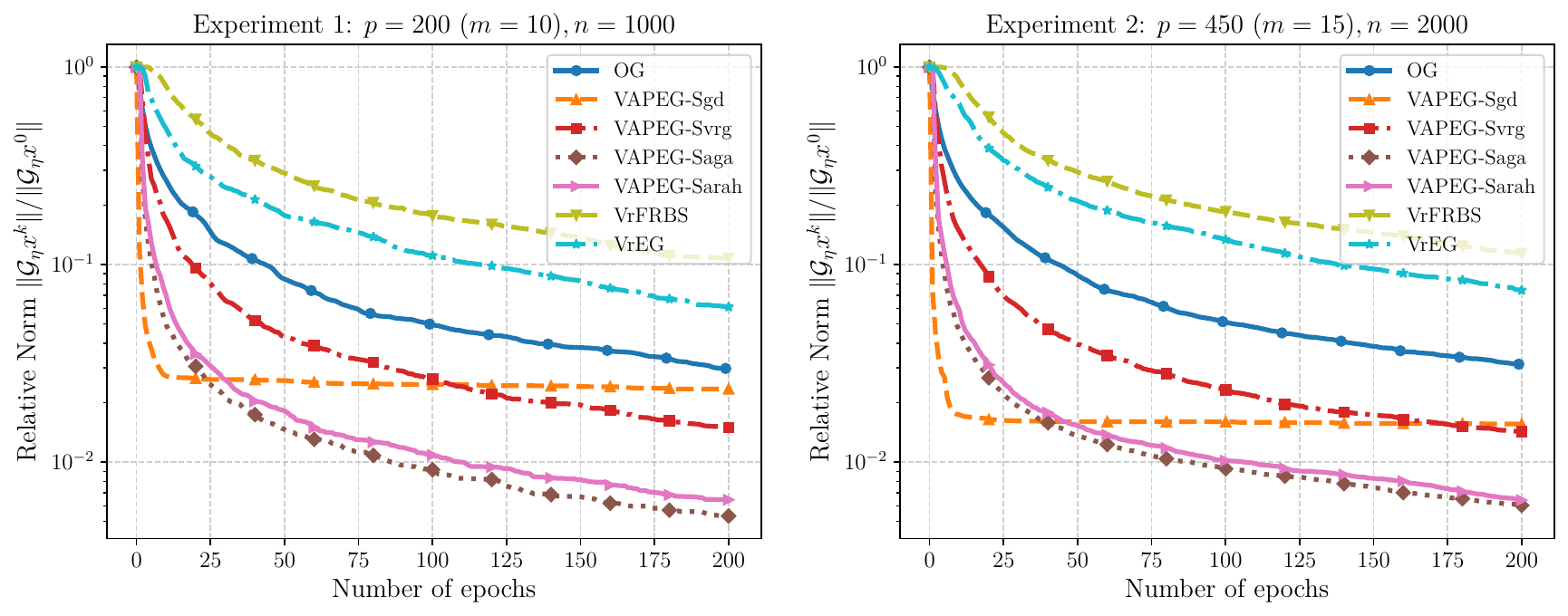}
	\vspace{-3ex}
	\caption{The performance of seven algorithms with a smaller probability $\mbf{p}_k$ and \revise{mini-batch} size $b_k$ to solve the bilinear matrix game \eqref{eq:bilinear_matrix_game} in two experiments, averaging over $10$ problem instances each.}
	\label{fig:exam1_matrix_games_results2}
	\vspace{-2ex}
\end{figure}

Next, we test the effect of the probability $\mbf{p}_k$ and the mini-batch size $b_k$ of $\Bc_k$ using the same experiments as in Figure~\ref{fig:exam1_matrix_games_results1}.
For the methods involving a reference-refresh probability, we halve \(p_k\), and for the applicable control-variate methods we also halve \(b_k\); the VAPEG-Sgd batch schedule is left unchanged.
Then, the results are presented in Figure~\ref{fig:exam1_matrix_games_results2}. 

The performance of control-variate methods improves because of the increased number of inner iterations per epoch.
This improvement is seen across all variants and consistent with the first test.

\beforesubsec
\subsection{\textbf{Markov decision processes.}}\label{subsec:MDP_exam}
\aftersubsec
Markov decision processes (MDPs) are commonly used to model sequential decision-making problems \cite{puterman1994markov}, with applications in numerous areas such as reinforcement learning \cite{mnih2015human} and healthcare \cite{alagoz2010markov,grandclem2019solving}.  

\noindent\textit{\textbf{$\mathrm{(a)}$~Mathematical model.}} 
A finite MDP consists of a state space $\Sc := \sets{1, 2, \cdots, n}$ with $n$ states and a finite action space $\mcal{A} := \sets{1, 2, \cdots, m}$ with $m$ actions. 
Each state-action pair $(s, a)$ is associated with a reward $r_{sa}$ and a transition probability distribution $P_{sa} \in \Delta_n$ over the next states. 
Define $r_\infty := \max_{s,a} r_{sa}$. 
Without loss of generality, we assume $r_{sa} \geq 0$ for all $(s,a) \in \Sc \times \mcal{A}$.
Given a discount factor $\bar{\gamma} \in (0,1)$ and an initial state distribution $\mbf{p}_0 \in \Delta_n$, a common objective in MDPs is to maximize the expected infinite-horizon discounted reward. 

Following the same procedure as in \cite{jin2020efficient,puterman1994markov}, this MDP problem can be reformulated into the following \revise{convex}-concave saddle-point problem:
\myeq{eq:MDP_minimax}{
	\min_{v \geq 0, \norms{v}_2 \leq \frac{\sqrt{n}r_\infty}{1- \bar{\gamma}}} \max\limits_{\mu \in \Delta_{nm}} \Big\{ \Lc(v, \mu) := (1 - \bar{\gamma}) \mbf{p}_0^\top v + \sum\limits_{s=1}^n \sum\limits_{a=1}^m \mu_{sa} \left( r_{sa} + \bar{\gamma} P_{sa}^\top v - v_s \right) \Big\},
}
where $v \in \R^n$ is the value-function variable, whose optimal value represents the expected discounted return from each state $s \in \mcal{S}$, and  $\mu_{sa}$ ($s \in \mcal{S}$ and $a \in \mcal{A}$) are the dual variables.

Next, define $x := [v; \mu] \in \R^{n(m+1)}$ as the concatenation of $v$ and $\mu = [\mu_1; \cdots; \mu_n] \in \R^{nm}$ (the concatenation of $\mu_s$'s with $\mu_s = (\mu_{s1}, \cdots, \mu_{sm})^\top$, $s = 1, \cdots, n$).
For each $s \in \Sc$, let $B_s := \left[ \bar{\gamma} P_{s1} - e_s, \cdots, \bar{\gamma} P_{sm} - e_s \right] \in \R^{n \times m}$, $r_s := (r_{s1}, \cdots, r_{sm})^\top$, $\Bc^+_n(0, R) := \sets{v \in \R^n_+: \norms{v}_2 \leq R}$, and
\myeqn{
\begin{array}{lcl}
	G_s x := \begin{bmatrix}
		n(1- \bar{\gamma})\mbf{p}_{0,s} e_s + nB_s \mu_s\\
		\mbf{0}_{(s-1) m}\\
		- nr_s - nB_s^\top v \\
		\mbf{0}_{(n-s) m}
	\end{bmatrix}, \qquad
	T x := \begin{bmatrix}
		\partial \delta_{\Bc^+_n(0, \frac{\sqrt{n}r_{\infty}}{1-\bar{\gamma}})}(v) \\
		\partial \delta_{\Delta_{nm}}(\mu)
	\end{bmatrix}.
\end{array}
}
Then, the optimality condition of \eqref{eq:MDP_minimax} can be written as \eqref{eq:GE}: $0 \in Gx + Tx$, where $Gx := \frac{1}{n}\sum_{s=1}^n G_s x$ is a finite-sum representation over the state space $\mcal{S}$.

\vspace{1ex}
\noindent\textit{\textbf{$\mathrm{(b)}$~Experimental setup and input data.}}
We evaluate the performance of \eqref{eq:VFOG} on solving problem \eqref{eq:MDP_minimax} using randomly generated MDPs from the \emph{garnet} testbed \cite{archibald1995generation,bhatnagar2012garnet}, which is widely used for benchmarking algorithms in sequential decision-making.
Garnet MDPs are defined by a branching factor $n_b/n$, representing the proportion of reachable next states from each $(s,a)$ pair.
Our setup is characterized by a tuple $(|\Sc|, |\mcal{A}|, n_b, \bar{\gamma})$, and then we consider two experiments: \textit{Experiment 1} with $(|\Sc|, |\mcal{A}|, n_b, \bar{\gamma}) = (2000, 5, 1000, 0.9)$, leading to $p=12000$ and \textit{Experiment 2} with $(|\Sc|, |\mcal{A}|, n_b, \bar{\gamma}) = (4000, 10, 2000, 0.9)$, leading to $p=44000$.
For each setting, we generate 10 random MDP instances where rewards are sampled uniformly at random from $[0,1]$. 

\vspace{1ex}
\noindent\textit{\textbf{$\mathrm{(c)}$~Algorithms and parameters.}}
In this experiment, we use three variants of our \eqref{eq:VFOG}: \texttt{VAPEG-Svrg}, \texttt{VAPEG-Saga}, and \texttt{VAPEG-Sarah} since they use control variate technique.
We again compare them with the baseline method \texttt{OG}, and two other competitors: \texttt{VrFRBS} and \texttt{VrEG} using SVRG estimators.
Unlike the previous example, we decide to select the parameters for all algorithms using a grid search procedure on the interval $\frac{1}{L_B}[10^{-5}, 10]$ to search for the best stepsize $\eta$ in all algorithms, where $L_B = \norms{B}$.
Note that $L = L_B$ is the Lipschitz constant of the aggregate operator $G$.
Applying this procedure, we obtain the following values of $\eta$ for each algorithm.
\begin{compactitem}
	\item For our variants \texttt{VAPEG-Svrg}, \texttt{VAPEG-Saga}, and \texttt{VAPEG-Sarah}, we find $\eta = \frac{1}{10^3L}$.
	\item For the baseline method: \texttt{OG}, we get $\eta = \frac{1}{10^2L}$.
	\item For \texttt{VrFRBS}, we obtain $\eta = \frac{0.95(1 - \sqrt{1 - \mbf{p}})}{10^2 \times 2L}$.
	\item For \texttt{VrEG}, we find $\eta = \frac{0.95\sqrt{1 - \alpha}}{10^3L}$ with $\alpha := 1 - \mbf{p}$.
\end{compactitem}
As in the first example, we still choose $\mbf{p} = 0.5n^{-1/3}$ and $b = \lfloor 0.5n^{2/3} \rfloor$ for \texttt{VAPEG-Svrg}, \texttt{VrFRBS}, and \texttt{VrEG}, $\mbf{p} = 0.5n^{-1/2}$ and $b = \lfloor 0.5n^{1/2} \rfloor$ for \texttt{VAPEG-Sarah}.
The initial point is chosen as $x^0 := [\frac{1-\bar{\gamma}}{r_{\infty}}\texttt{ones(n)}; \frac{1}{n\mred{m}}\texttt{ones(n\mred{m})}]$ for all methods.

\begin{figure}[hpt!]
	\vspace{-0ex}
	\centering
	\includegraphics[width=\textwidth]{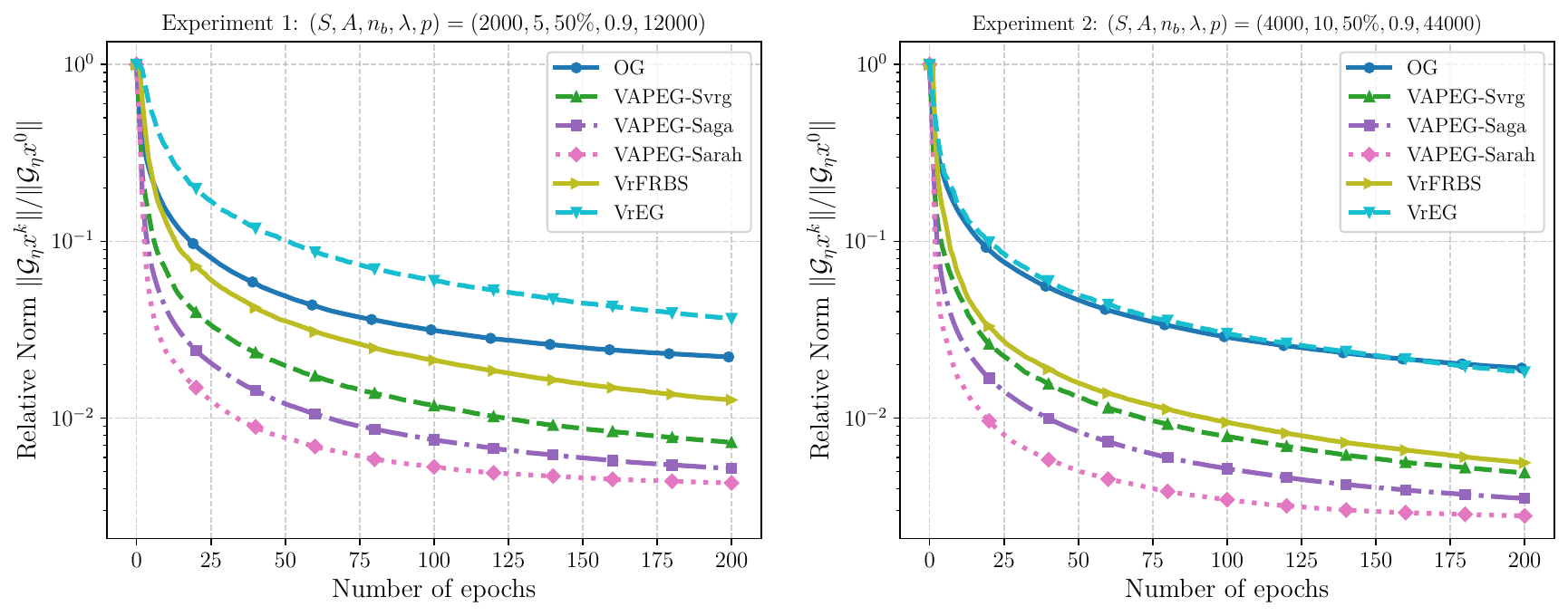}
	\vspace{-3ex}
	\caption{The performance of \revise{six} algorithms to solve the infinite horizon discounted reward MDPs reformulation \eqref{eq:MDP_minimax} in two experiments, averaging over 10 problem instances each.}
	\label{fig:exam3_MDP_results1}
	\vspace{-2ex}	
\end{figure}

\vspace{1ex}
\noindent\textit{\textbf{$\mathrm{(d)}$~Numerical results.}}
We run all the algorithms using the above setting on two experiments, and the results are shown in Figure~\ref{fig:exam3_MDP_results1}.

As shown in Figure~\ref{fig:exam3_MDP_results1}, the three VAPEG variants: \texttt{VAPEG-Svrg}, \texttt{VAPEG-Saga}, and \texttt{VAPEG-Sarah}, consistently outperform both the baseline method (\texttt{OG}) and competing algorithms across the entire training horizon.
Among them, \texttt{VAPEG-Svrg} achieves the worst residual norm, while \texttt{VAPEG-Sarah} demonstrates the best overall performance.
These results suggest that acceleration indeed provides a significant performance boost.

\begin{figure}[hpt!]
	\centering
	\includegraphics[width=\textwidth]{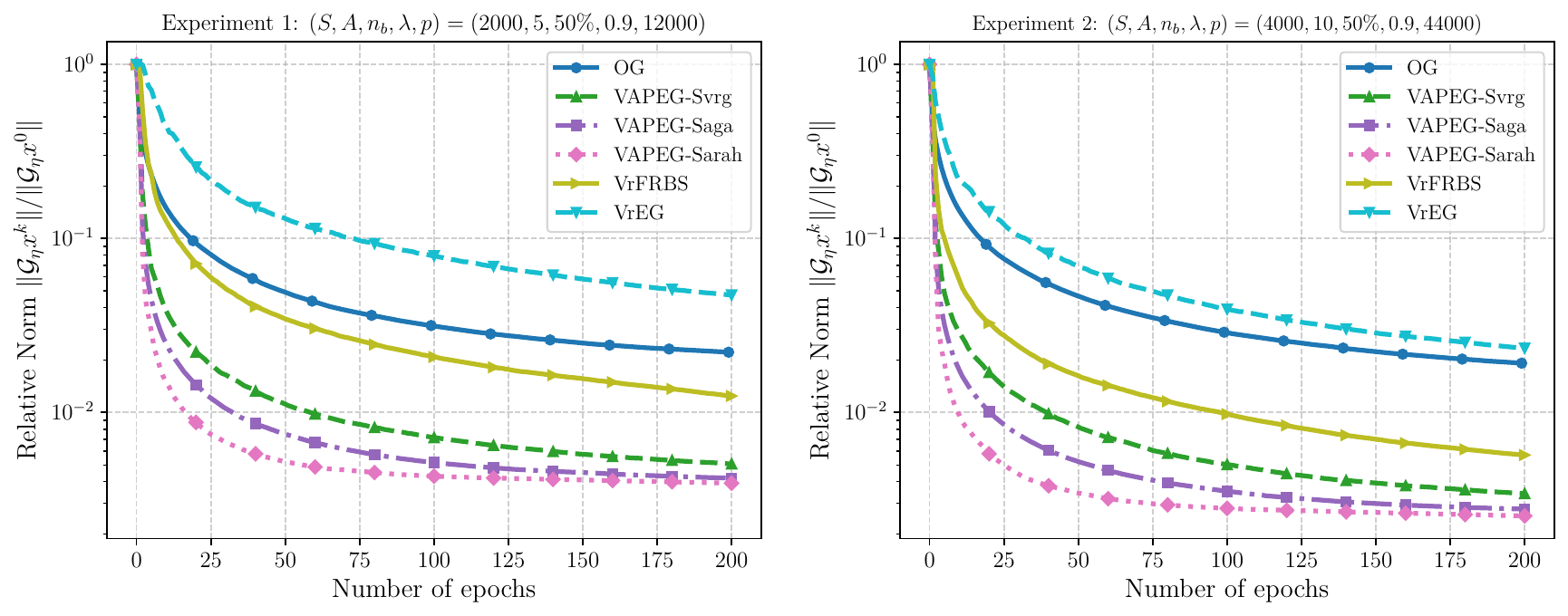}
	\vspace{-3ex}
	\caption{The performance of \revise{six} algorithms to solve the infinite horizon discounted reward MDPs reformulation \eqref{eq:MDP_minimax} in two experiments using smaller $\mbf{p}$ and batch size $b$, averaging over 10 problem instances each.}
	\label{fig:exam3_MDP_results2}
	\vspace{-1ex}	
\end{figure}

Finally, we examine the effect of the probability $\mbf{p}_k$ and the mini-batch size $b_k$ on the performance of our methods by halving their values to obtain more inner iterations in each epoch.
The corresponding results are presented in Figure \ref{fig:exam3_MDP_results2}.

We observe that the three variants of \eqref{eq:VFOG} outperform their non-accelerated counterparts.
A smaller mini-batch size $b_k$ enables our algorithms to run more iterations, while a smaller probability $\mbf{p}_k$ reduces the number of full evaluations of $G$.
Consequently, both factors contribute to the speedup of the algorithms observed in Figure \ref{fig:exam3_MDP_results2}.

\revise{
\beforesubsec
\subsection{\textbf{Sparse soft-robust bond portfolio construction.}}\label{subsec:robust_portfolio_exam}
\aftersubsec
In this example, we compare our algorithms with three competing methods on a sparse soft-robust bond portfolio construction problem. 
The model is adapted from the robust bond portfolio formulation in \cite[Section 6.2]{bian2026nonsmooth}.

\vspace{1ex}
\noindent\textbf{$\mathrm{(a)}$~Mathematical model.}
Consider a portfolio of $m$ bonds with quantities $z = (z_1, \cdots, z_m) \in \Zc \subseteq \R^m_{+}$ over time periods $t = 1, \cdots, T$.
A market state is modeled by $q = (q^u, q^s) \in \R^{T+m}$, where $q^u = (q^u_1, \cdots, q^u_T) \in \R^T$ is the yield curve and $q^s = (q^s_1, \cdots, q^s_m) \in \R^m_{+}$ is the bond spread.
For a market state $q = (q^u, q^s)$, let $\pi \in \R^m_{+}$ denote the bond prices, where the price $\pi_j$ of bond $j$ is modeled by $\pi_j(q) := \sum_{t=1}^T \alpha_{j,t}\exp\big(-t(q^u_t + q^s_j)\big)$,
with $\alpha_{j,t} \geq 0$ denoting the cash flow paid by bond $j$ at time $t$.

Next, fix a reference market state $\bar{q}$ to evaluate local bond-price sensitivities, and let $y := (u^\top, s^\top)^\top \in \R^{T + m}$ denote changes in yields and spreads relative to this state.
Thus, the stressed state is $\bar{q} + y$.
Since $\pi$ is highly nonlinear, we approximate it around $\bar{q}$ by $\pi(\bar{q} + y) \approx \pi(\bar{q}) + J_{\pi}(\bar{q})y$, where $J_\pi$ is the Jacobian of the bond-price map at $\bar q$.
Denoting $B := -J_{\pi}(\bar{q})$, we obtain $\pi(\bar{q} + y)^\top z \approx \pi(\bar{q})^\top z - z^\top By$.
Let $\phi: \R^m \to \R$ be a smooth convex nominal portfolio objective, which may include tracking error against a benchmark, risk, or transaction cost terms.
Absorbing the reference-state portfolio value, we define $f(z) := \phi(z)-\pi(\bar q)^\top z$, which remains smooth and convex.

Now, instead of using a bounded uncertainty set together with a nonsmooth uncertainty penalty, we penalize the overall magnitude of the market shock and its movements along prescribed market-shock patterns through
\begin{equation}\label{eq:portfolio_h_function}
\begin{array}{lcl}
	h(y) &:=& \frac{\gamma}{2}\norms{y}^2 + \nu \sum_{\ell=1}^r \log \cosh(c_{\ell}^\top y),
\end{array}
\end{equation}
where $\gamma>0$, $\nu>0$, and $C\in\mathbb R^{r\times(T+m)}$ has rows $c_\ell^\top$, $\ell \in [r]$.
Each vector $c_\ell$ represents a prescribed pattern of joint movements in yields and spreads, while $c_\ell^\top y$ measures the signed projection of the market shock onto that pattern.
Since $\log\cosh(\cdot)$ is even and increasing in the absolute value of its argument, the second term penalizes large movements along these patterns.
The parameter $\gamma$ controls the overall magnitude penalty, while $\nu$ controls the strength of the pattern-specific penalty.

Finally, to promote a sparse bond portfolio, we use the SCAD penalty $\Rc_{\lambda, \theta}(z) := \sum_{j=1}^m R_{\lambda, \theta} (z_j)$ proposed in \cite{fan2001variable} as a continuous surrogate for cardinality.
The resulting \textit{sparse soft-robust bond portfolio model} is the following nonconvex-concave minimax problem:
\begin{equation}\label{eq:portfolio_minimax}
	\min_{z \in \Zc}\max_{y \in \R^{T+m}} \set{\Lc(z,y) := f(z) + z^\top By - h(y) + \tau \Rc_{\lambda, \theta}(z)}
\end{equation}
Here, $\lambda>0$ determines the SCAD threshold, \(\theta>2\) is the SCAD shape parameter, and $\tau>0$ controls the strength of the sparsity regularization.
The SCAD penalty can be decomposed as $R_{\lambda, \theta}(t) = \lambda |t| - \tilde{R}_{\lambda, \theta}(t)$, where
\begin{equation}\label{eq:portfolio_SCAD_complementary}
	\tilde{R}_{\lambda, \theta}(t) = \begin{cases} 
		0, &|t|\leq\lambda,\\[1mm]
		\dfrac{(|t|-\lambda)^2}{2(\theta-1)}, &\lambda<|t|\leq \theta\lambda,\\[2mm] 
		\lambda|t|-\dfrac{(\theta+1)\lambda^2}{2}, 	&|t|>\theta\lambda.
	\end{cases}
\end{equation}
Let $\tilde{\Rc}_{\lambda, \theta}(z) := \sum_{j=1}^m \tilde{R}_{\lambda, \theta} (z_j)$.
Then, the saddle-point problem \eqref{eq:portfolio_minimax} can be written equivalently as
\begin{equation}\label{eq:portfolio_minimax2}
	\min_{z \in \R^m}\max_{y \in \R^{T+m}} \set{\Lc(z,y) := f(z) + z^\top By - h(y) - \tau \tilde{\Rc}_{\lambda, \theta}(z) + \tau \lambda \norms{z}_1 + \delta_{\Zc}(z)},
\end{equation}
where $\delta_{\Zc}$ is the indicator function of the closed convex set $\Zc$.
Suppose $f(z) = \frac{1}{n}\sum_{i=1}^n f_i(z)$,
where every $f_i$ is continuously differentiable.
Now, let us define
\begin{equation}\label{eq:portfolio_inclusion_elements}
	x = \begin{bmatrix}  z \\ y \end{bmatrix}, \qquad 
	G_ix := \begin{bmatrix}  \nabla f_i(z) + By - \tau\nabla \tilde{\Rc}_{\lambda, \theta}(z) \\ \nabla h(y) - B^\top z \end{bmatrix}, \quad \text{and} \quad
	Tx := \begin{bmatrix} \tau \lambda \partial \norms{z}_1 + \partial \delta_{\Zc}(z) \\ 0 \end{bmatrix}.
\end{equation}
Then, the optimality condition of \eqref{eq:portfolio_minimax2} can be written as $0 \in Gx + Tx$, where $Gx := \frac{1}{n}\sum_{i=1}^n G_ix$.

Proposition~\ref{pro:properties_of_G_and_T} summarizes the properties of $G$ and $T$, with its proof deferred to Appendix~\ref{apdx:prop:h_function_smooth}.

\begin{proposition}\label{pro:properties_of_G_and_T}
Suppose that $f$ in \eqref{eq:portfolio_minimax2} is convex and $L_f$-average smooth, i.e.:
\begin{equation*}
\begin{array}{lcl}
	\frac{1}{n}\sum_{i=1}^n \norms{\nabla f_i(z) - \nabla f_i(z')}^2 \leq L_f^2 \norms{z - z'}^2, \qquad \forall z, z' \in \R^m.
\end{array}
\end{equation*}
Suppose further that $B \in \R^{m \times (T+m)}$ in \eqref{eq:portfolio_minimax2} has full row rank such that its smallest singular value $\sigma := \sigma_{\min}(B^\top) > L_h\big(\frac{\tau}{\gamma(\theta - 1)}\big)^{1/2}$, where $L_h := \gamma + \nu \norms{C}^2$.
Then the following statements hold:
\begin{compactenum}
\item[$\mathrm{(i)}$]
	The mapping $G$ in \eqref{eq:portfolio_inclusion_elements} is $L$-Lipschitz continuous with $L := \frac{1}{2} \big[ \tilde{L}_f + L_h + \big[(\tilde{L}_f - L_h)^2 + 4\norms{B}^2 \big]^{1/2} \big]$, where $\tilde{L}_f := L_f + \frac{\tau}{\theta - 1}$.
	
\item[$\mathrm{(ii)}$]
	The mapping $\Phi := G+T$ is globally $\rho$-co-hypomonotone, with $\rho := \frac{\tau}{(\theta - 1)\big(\sigma - L_h\sqrt{\tau/((\theta - 1)\gamma)} \big)^2}$.
	
\item[$\mathrm{(iii)}$]
	The mapping $\Phi := G+T$ satisfies Assumption~\ref{as:A3} with $\rho_n := (1 + \kappa)\rho$ and $\rho_c := \frac{\kappa\rho}{1+\hat{L}^2}$, where $\hat{L}^2 := \frac{L^2\big\{\gamma^2 L_h^2(\theta - 1)^2\left(\gamma^2 + \norms{B}^2\right) + \big[\gamma\left(\sigma^2(\theta - 1) - \tau L_h\right) + L_h(\theta - 1)\norms{B}\sqrt{\gamma^2 + \norms{B}^2}\big]^2\big\}}{\gamma^4 [\sigma^2(\theta - 1) - \tau L_h ]^2}$ and $\kappa > 0$ arbitrarily.
\end{compactenum}
\end{proposition}

\vspace{1ex}
\noindent\textbf{$\mathrm{(b)}$~Data generation.}
This example aims to illustrate the performance of our methods and their competitors on a nonmonotone problem.
The data is generated synthetically as follows.

\smallskip
\noindent\textbf{\textit{$\mathrm{(i)}$~Bond cash flows and reference market state.}}
For each bond $j \in [m]$, we independently generate its maturity and coupon rate according to
$M_j \sim \operatorname{Unif}\sets{2,\cdots,T}$ and $\kappa_j \sim \operatorname{Unif}[0.018,0.042]$, respectively.
The face value is normalized to one, and the cash flows are
$\alpha_{j,t} := \kappa_j \boldsymbol{1}_{\sets{t \leq M_j}} + \boldsymbol{1}_{\sets{t = M_j}}$.
The reference yield curve is deterministic, $\bar q_t^u := 0.018 + \frac{0.012}{T-1}(t-1)$, $t \in [T]$, while the bond-specific reference spreads are independently generated from $\bar q_j^s \sim \operatorname{Unif}[0.004,0.007]$, $j \in [m]$.
Once the reference market state $\bar{q}$ is generated, the matrix $B=-J_\pi(\bar q) \in \R^{m \times (T+m)}$ is constructed analytically. Specifically, the yield-curve and spread blocks of $B$ are given by
\begin{equation*}
\begin{array}{lcll}
	B_{j,t}	&=&	t\alpha_{j,t} \exp\bigl(-t(\bar q_t^u+\bar q_j^s)\bigr),
	\qquad &\qquad t \in [T], \vspace{1ex}\\
	B_{j,T+\ell} &=& \mathbf{1}_{\sets{j=\ell}} 	\sum_{t=1}^{T} t\alpha_{j,t} \exp\bigl(-t(\bar q_t^u+\bar q_j^s) \bigr),
	\qquad &\qquad \ell \in [m].
\end{array}
\end{equation*}
Since the spread block is diagonal with strictly positive diagonal entries, every generated matrix \(B\) has full row rank, and thus $\sigma := \sigma_{\min}(B^\top) > 0$.

\smallskip
\noindent\textbf{\textit{$\mathrm{(ii)}$~Market-stress directions.}}
For each setting, define $\boldsymbol{t} := (1,2,\cdots,T)^\top$, $\boldsymbol{1}_{T} := (1,\cdots,1)^\top\in\R^T$, $\bar{t} := \frac{T+1}{2}$, and
$\widetilde{\boldsymbol{t}} := \boldsymbol{t} - \bar t\,\boldsymbol{1}_{T}$.
The level, slope, and curvature directions and the common credit-spread direction are defined as
\begin{equation*}
\begin{array}{lcl}
	c^{\mathrm{lev}} := \frac{\boldsymbol{1}_{T}}{\sqrt{T}}, 
	\qquad
	c^{\mathrm{slp}} := \frac{\widetilde{\boldsymbol{t}}}{\norms{\widetilde{\boldsymbol{t}}}},
	\qquad
	c^{\mathrm{curv}} := \frac{\widetilde{\boldsymbol{t}}^{\odot 2} - \frac{ \norms{ \widetilde{\boldsymbol{t}} }^2}{T}\boldsymbol{1}_{T}}{
	\big\Vert \widetilde{\boldsymbol{t}}^{\odot 2} - \frac{ \norms{ \widetilde{\boldsymbol{t}} }^2 }{T} \boldsymbol{1}_{T} \big\Vert},
	\qquad \text{and} \qquad
	c^{\mathrm{spr}} := \frac{\boldsymbol{1}_{m}}{\sqrt{m}},
\end{array}
\end{equation*}
respectively, where $^{\odot 2}$ denotes the elementwise square.
We then set
\begin{equation*}
C := \begin{bmatrix}
	(c^{\mathrm{lev}})^\top & 0\\
	(c^{\mathrm{slp}})^\top & 0\\
	(c^{\mathrm{curv}})^\top & 0\\
	0 & (c^{\mathrm{spr}})^\top
\end{bmatrix}
\in \mathbb{R}^{4\times(T+m)}.
\end{equation*}
The four rows of $C$ are orthonormal, and hence $\norms{C}=1$.
We use $\nu=1$ and $\gamma := \nu\norms{C}^2 = 1$.
Thus, $L_h = \gamma+\nu\norms{C}^2 = 2$.

\smallskip
\noindent\textbf{\textit{$\mathrm{(iii)}$~Finite-sum scenarios.}}
The shifted portfolio objective is chosen as $f(z) := \frac{1}{n}\sum_{i=1}^n f_i(z)$, where $f_i(z) := \frac{1}{2}(g_i^\top z-r_i)^2
+ \frac{\mu}{2}\norms{z}^2 + b^\top z$.
Here, $g_i\in\R^m$ represents the modeled bond responses in scenario $i$, so that $g_i^\top z$ is the aggregate portfolio response, while $r_i \in \R$ is the corresponding benchmark response.
Thus, the first term is the average squared tracking error, while the remaining terms provide quadratic regularization and a linear portfolio adjustment.

Suppose that each market scenario is characterized by a $6$-dimensional factor model.
We first define the standardized bond features
\begin{equation*}
\begin{array}{lcl}
	\widehat M_j := \frac{M_j-\bar M}{\mathrm{sd}_M}, 	\qquad
	\widehat\kappa_j  := \frac{\kappa_j-\bar\kappa}{\mathrm{sd}_\kappa},  \qquad
	\widehat s_j := 	\frac{\bar q_j^s-\bar s}{\mathrm{sd}_s},
\end{array}
\end{equation*}
where $(\bar{M}, \mathrm{sd}_M)$, $(\bar{\kappa}, \mathrm{sd}_{\kappa})$, and $(\bar{s}, \mathrm{sd}_s)$ are the empirical mean and standard deviation of maturity, coupon, and spread, respectively.
We form
\(\widehat F\in\mathbb{R}^{m\times6}\) row-wise as
\begin{equation*}
\begin{array}{lcl}
	\widehat F_{j,:} := \big[ 1,\  \widehat M_j,\  \widehat M_j^2 - \frac{1}{m} \sum_{\ell=1}^{m}\widehat M_\ell^2,\  \widehat\kappa_j,\  \widehat s_j,\  \widehat M_j\widehat s_j \big].
\end{array}
\end{equation*}
The six columns represent a common movement, linear and quadratic maturity effects, a coupon effect, a credit-spread effect, and a maturity--credit interaction.
Let $D_F := \operatorname{diag}\bigl(\Vert \widehat F_{:,1} \Vert, \cdots, \Vert \widehat F_{:,6}\Vert\bigr)$.
If any empirical standard deviation or column norm of $\widehat F$ is zero, we discard the draw and resample the corresponding bond characteristics, ensuring that $D_F$ is nonsingular.
We then define $F := a_F\widehat F D_F^{-1}$ so that all factors have equal magnitude $a_F$.
For every $i \in [n]$, we independently generate $\omega_i\sim \Nc(0,\Id_6)$ and set $g_i:=F\omega_i$.

To construct a low-curvature direction and verify nonmonotonicity, we impose two bonds with identical observable characteristics.
This imposition is applied before constructing the cash flows, $B$, $F$, and the scenario vectors $g_i$.
For each problem instance, we select
$j_1 := \left\lfloor\frac{m}{4}\right\rfloor$, $j_2 := j_1+1$, and set $M_{j_2}:=M_{j_1}$, $\kappa_{j_2} := \kappa_{j_1}$, and $\bar{q}_{j_2}^s := \bar{q}_{j_1}^s$.
It follows that $F_{j_2,:} = F_{j_1,:}$ and hence $g_{i,j_1}=g_{i,j_2}$ for every $i$.
Therefore, defining $d_{\mathrm{sw}} := \frac{e_{j_1}-e_{j_2}}{\sqrt{2}}$, where $e_j$ is the $j^{\mathrm{th}}$-unit vector, we have $g_i^\top d_{\mathrm{sw}}=0$ for every $i$.

\smallskip
\noindent\textbf{\textit{$\mathrm{(iv)}$~Selection of the theoretical parameters.}}
We use the same SCAD and box-constraint parameters for all instances: $\mathcal Z=[0,1]^m$, $\lambda=0.25$, and $\theta=2.5$.
The parameters $\tau$ and $\mu$ are selected as follows.
For each instance, we initialize $\tau := (\theta-1) \min \sets{0.05, 0.1\gamma \big(\frac{\sigma}{L_h}\big)^2}$ and set $\mu := \frac{\tau}{4(\theta-1)}$.
This choice satisfies $\sigma^2 > L_h^2\frac{\tau}{(\theta-1)\gamma}$.
Next, we compute the average smoothness parameter of $f$ as $L_f := \big[ \lambda_{\max}\big( \frac{1}{n}\sum_{i=1}^n(g_ig_i^{\top} + \mu \Id)^2 \big) \big]^{1/2}$.
If we denote $\tilde{L}_f := L_f+\frac{\tau}{\theta-1}$, then the Lipschitz constant of $G$ as $L = \frac{\tilde{L}_f + L_h + \sqrt{(\tilde{L}_f - L_h)^2 + 4\norms{B}^2}}{2}$.
For a desired tolerance $\epsilon \in (0,1)$, if $L\rho_n > \epsilon$, we replace $\tau$ by $\frac{\tau}{2}$, update $\mu = \frac{\tau}{4(\theta-1)}$, and recompute $L_f$, $L$, and $\rho_n$.
This procedure terminates once $L\rho_n \leq \epsilon$.
Consequently, every generated instance satisfies the sufficient global co-hypomonotonicity condition and the stronger numerical requirement $L\rho_n \leq \epsilon < 1$.
In this case, Assumption~\ref{as:A3} also holds.

\smallskip
\noindent\textbf{\textit{$\mathrm{(v)}$~Sparse stationary point.}}
For each problem instance, we sample an active set $\Sc \subseteq \sets{1,\cdots,m}$ uniformly such that $j_1,j_2\in\Sc$ (i.e., we uniformly sample $\tilde{\Sc} \subseteq \sets{1, \cdots, m} \setminus \sets{j_1, j_2}$ of cardinality $\vert \tilde{\Sc} \vert = \vert \Sc \vert - 2$, then define $\Sc := \tilde{\Sc} \cup \sets{j_1, j_2}$).
We choose $z_j^\star :=\frac{\theta+1}{2}\lambda \mbf{1}_{\sets{j \in \Sc}}$ for all $j \in [m]$.
Hence, all active coordinates satisfy $\lambda < z_j^\star < \theta\lambda < 1$.
We compute $y^\star$ as the solution of
$\nabla h(y^\star)=B^\top z^\star$, whose existence and uniqueness are guaranteed by the strong convexity of $h$. 
We then set $r_i:=g_i^\top z^\star$, $i \in [n]$, so that $z^\star$ reproduces the generated benchmark response in every scenario.
Define $x^{\star} := [z^{\star}, y^{\star}]$, $\xi^\star := \sum_{j=1}^m \tau \lambda \mbf{1}_{\sets{j \in \Sc}}e_j$, and choose $b := -\mu z^\star -By^\star + \tau \nabla\widetilde R_{\lambda,\theta}(z^\star) - \xi^\star$.
By construction, $\xi^\star \in \tau \lambda \partial\norms{z^\star}_1 + \Nc_{\Zc}(z^\star)$
and $0 = Gx^\star +	[\xi^\star, 0]$.

\smallskip
\noindent\textbf{\textit{$\mathrm{(vi)}$~Nonmonotonicity verification.}}
Since $g_i^\top d_{\mathrm{sw}} = 0$ for every $i$, the Hessian of $f$ satisfies $H_f d_{\mathrm{sw}} = \mu d_{\mathrm{sw}}$.
Choose $\varepsilon \in \big( 0, \frac{(\theta-1)\lambda}{\sqrt{2}}\big)$ arbitrarily, and define $z^\varepsilon := z^\star+\varepsilon d_{\mathrm{sw}}$ and
$x^\varepsilon := [z^\varepsilon, y^\star]$.
Then, the two perturbed coordinates $z_{j_1}^{\varepsilon} = z_{j_1}^{\star} + \frac{\varepsilon}{\sqrt{2}}$ and $z_{j_2}^{\varepsilon} = z_{j_2}^{\star} - \frac{\varepsilon}{\sqrt{2}}$ remain strictly in the middle SCAD region $(\lambda, \theta \lambda)$.
Define $w^{\star} := Gx^{\star} + [\xi^{\star}, 0]$ and $w^{\varepsilon} := Gx^{\varepsilon} + [\xi^{\varepsilon}, 0]$.
It is straightforward to verify that $\xi^{\varepsilon} = \xi^{\star} = \sum_{j=1}^m \tau \lambda \mathbf{1}_{\sets{j \in \Sc}} e_j$ satisfy $\xi^{\varepsilon} \in \tau \lambda \partial\norms{z^{\varepsilon}}_1 + \Nc_{\Zc}(z^{\varepsilon})$ and $\xi^{\star} \in \tau \lambda \partial\norms{z^{\star}}_1 + \Nc_{\Zc}(z^{\star})$.
Then, similar to \eqref{eq:portfolio_cohypomonotone_proof1}, we have
\begin{equation*} 
\begin{array}{lcl}
	\iprods{w^\varepsilon-w^\star, x^\varepsilon-x^\star} &=& \iprods{\nabla f(z^\varepsilon) - \nabla f(z^\star), z^\varepsilon - z^\star} - \tau \iprods{\nabla \tilde{\Rc}_{\lambda, \theta}(z^\varepsilon) - \nabla \tilde{\Rc}_{\lambda, \theta}(z^\star), z^\varepsilon - z^\star} \vspace{1ex}\\
	&& + {~} \iprods{\xi^\varepsilon - \xi^\star, z^\varepsilon - z^\star} + \iprods{\nabla h(y^\star)-\nabla h(y^\star), y^\star - y^\star}.
\end{array}
\end{equation*}
The last two terms on the right-hand side are zero.
Moreover, using the gradient of $\tilde{\Rc}_{\lambda, \theta}$, $z^{\varepsilon} - z^{\star} = \varepsilon d_{\mathrm{sw}}$, and $g_i^\top d_{\mathrm{sw}} = 0$ for all $i$, we obtain
\begin{equation*}
\begin{array}{lcl}
	\nabla \tilde{\Rc}_{\lambda, \theta}(z^\varepsilon) - \nabla \tilde{\Rc}_{\lambda, \theta}(z^\star) & = & \frac{1}{\theta - 1}(z^\varepsilon - z^\star),  \quad \text{and} \vspace{1ex}\\
	\nabla f(z^\varepsilon) - \nabla f(z^\star) & = & \big(\frac{1}{n}\sum_{i=1}^n g_i g_i^\top + \mu\Id\big)(z^\varepsilon - z^\star) = \frac{1}{n}\sum_{i=1}^n \varepsilon g_i (g_i^\top d_{\mathrm{sw}}) + \mu(z^\varepsilon - z^\star) = \mu(z^\varepsilon - z^\star).
\end{array}
\end{equation*}
Therefore, we get
\begin{equation*} 
\begin{array}{lcl}
	\iprods{w^\varepsilon-w^\star, x^\varepsilon-x^\star}  	&=& \big(\mu-\frac{\tau}{\theta-1}\big)	\norms{z^\varepsilon-z^\star}^2	< 0
\end{array}
\end{equation*}
due to the choice $\mu = \frac{\tau}{4(\theta - 1)} < \frac{\tau}{\theta-1}$.
Since $\varepsilon$ can be chosen arbitrarily small, the last inequality shows that $\Phi$ fails to be monotone in every neighborhood of $x^\star$.

\vspace{1ex}
\noindent\textbf{$\mathrm{(c)}$~Experiment setup and algorithms.}
Our setup is characterized by a tuple $(n, m, T)$, with two experiments: $(n, m, T) = (1000, 200, 20)$, leading to $p = 420$, and $(n, m, T) = (2000, 460, 30)$, leading to $p = 950$.
For each setting, we generate 10 random problem instances as described above.
We test four variants of our \eqref{eq:VFOG} and compare them with the same three benchmarks as before.
The parameters are chosen as in Subsection~\ref{subsec:bilinear_matrix_game}.
For all algorithms, the initial point is chosen as $x^0 = [z^0, \mbf{0}]$, where $z^0$ is sampled i.i.d. from $\Uc(1.1\lambda, 0.9\theta\lambda)$, so that every portfolio component lies strictly inside the concave region of the SCAD penalty.

\vspace{1ex}
\noindent\textbf{$\mathrm{(d)}$~Numerical results.}
The results of these experiments are reported in Figure~\ref{fig:exam2_bond_portfolio}.
\begin{figure}[hpt!]
	\vspace{-2ex}
	\centering
	\includegraphics[width=\textwidth]{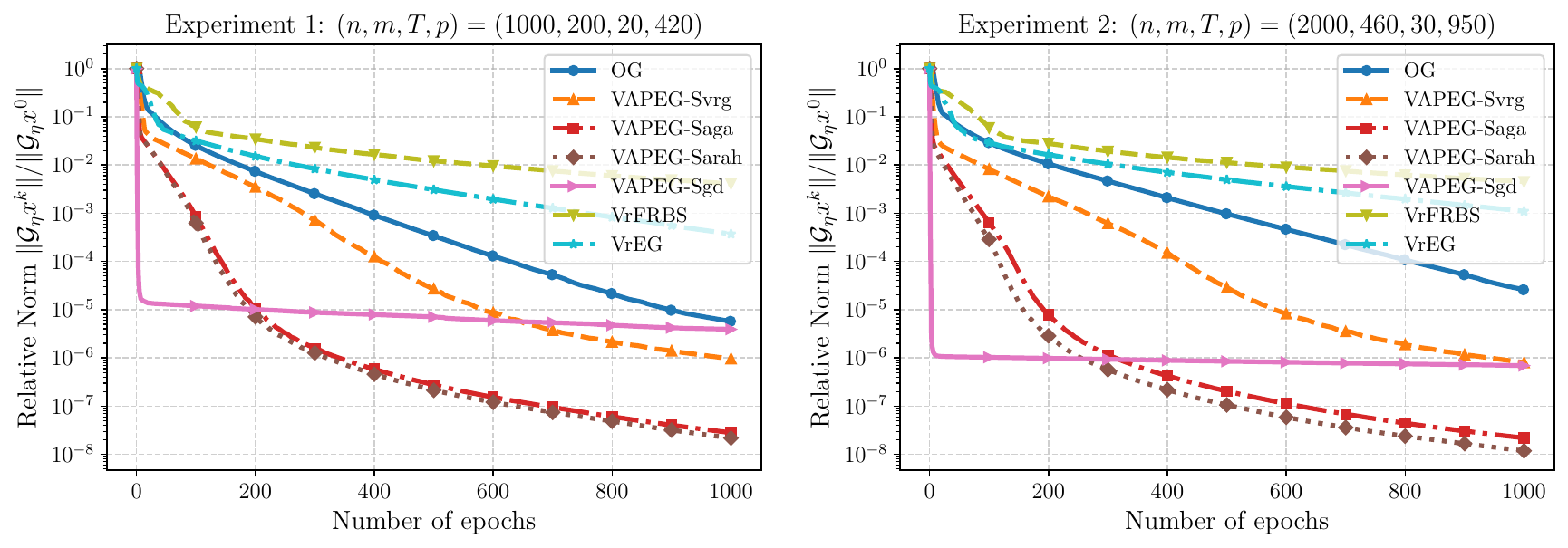}
	\vspace{-3ex}
	\caption{
		The performance of the seven algorithms for solving the soft-robust bond-portfolio construction problem \eqref{eq:portfolio_minimax2} in two experiments, averaged over $10$ problem instances for each.
		Here, $\mcal{G}_{\eta}$ is the FBS residual defined by \eqref{eq:FBS_residual}.
	}
	\label{fig:exam2_bond_portfolio}
	\vspace{-1ex}
\end{figure}

Unlike the bilinear matrix game problem considered in Subsection~\ref{subsec:bilinear_matrix_game}, whose associated operator is monotone, the soft-robust bond-portfolio problem is genuinely nonmonotone due to the concave part of the SCAD penalty.
Moreover, by construction, the problem remains nonmonotone around the stationary point.
Nevertheless, the performance of the algorithms remains largely consistent across the two experiments.
In both experiments, \texttt{VAPEG-Saga} and \texttt{VAPEG-Sarah} converge substantially faster and attain significantly higher accuracy than the other methods, reaching relative residuals of approximately $10^{-8}$, with \texttt{VAPEG-Sarah} performing slightly better than \texttt{VAPEG-Saga}.
\texttt{VAPEG-Svrg} also performs relatively well, reaching a final residual of roughly $10^{-6}$, whereas
\texttt{VAPEG-Sgd} also reaches the magnitude of $10^{-6}$. 
Meanwhile, \texttt{OG} converges more slowly, and the two benchmark methods, \texttt{VrFRBS} and \texttt{VrEG}, exhibit the slowest residual decay.
These results suggest that the advantages of the proposed variance-reduced schemes persist beyond the monotone setting and extend to structured nonmonotone finite-sum inclusions.
}

\beforesec
\section*{Acknowledgements.}
This work was partly supported by the National Science Foundation (NSF): NSF-RTG grant No. DMS-2134107 and the Office of Naval Research (ONR), grant No. N00014-23-1-2588.

\begin{APPENDICES}
\beforesec
\section{Preliminary Results and Example Verification.}\label{apdx:tech_results}
This appendix provides useful lemmas.
It also verifies the examples in Subsection~\ref{subsec:nonmonotone_example}.

\beforesubsec
\subsection{\textbf{Useful lemmas.}}
\label{apdx:subsec:useful_lemmas}
\aftersubsec
We recall some useful results which will be used in our analysis.

\begin{lemma}[\cite{Bauschke2011}, Lemma 5.31]\label{le:A1_descent}
Let $\sets{u_k}$, $\sets{v_k}$, $\sets{\gamma_k}$, and $\sets{\varepsilon_k}$ be nonnegative sequences such that $A_{\infty} := \sum_{k=0}^{\infty}\gamma_k < +\infty$ and $B_{\infty} := \sum_{k=0}^{\infty} \varepsilon_k < +\infty$.
In addition, we assume that
\begin{equation}\label{eq:lm_A1_cond}
u_{k+1} \leq (1 + \gamma_k)u_k - v_k + \varepsilon_k, \quad \textrm{for all $k \geq 0$}.
\end{equation}
Then, we conclude that $\lim_{k\to\infty}u_k$ exists and $\sum_{k=0}^{\infty}v_k < +\infty$.
\end{lemma}

\begin{lemma}[\cite{TranDinh2025a}]\label{le:A2_sum}
Given a nonnegative sequence $\sets{u_k}$ and $\omega \geq 0$ such that $\lim_{k\to\infty} k^{\omega + 1}u^k$ exists and $\sum_{k=0}^{\infty}k^{\omega}u_k < +\infty$.
Then, we conclude that $\lim_{k\to\infty}k^{\omega + 1}u^k = 0$.
\end{lemma}

We also need the well-known Robbins-Siegmund  supermartingale theorem \cite{robbins1971convergence}.

\begin{lemma}[\cite{robbins1971convergence}]\label{le:RS_lemma}
Let $\sets{U_k}$, $\sets{\alpha_k}$, $\sets{V_k}$ and $\sets{R_k}$ be  sequences of nonnegative integrable random variables on some arbitrary probability space and adapted to the filtration $\set{\Fc_k}$ with $\sum_{k=0}^{\infty}\alpha_k < +\infty$ and $\sum_{k=0}^{\infty}R_k < +\infty$ almost surely, and 
\begin{equation}\label{eq:RS_martingale_cond}
\Expn{ U_{k+1} \mid \Fc_k} \leq (1 + \alpha_k)U_k - V_k + R_k, 
\end{equation}
almost surely for all $k \geq 0$.
Then, almost surely, $\sets{U_k}$ converges to a nonnegative random variable $U$ and $\sum_{k=0}^{\infty}V_k <+\infty$ almost surely.
\end{lemma}

The following lemma is proven similarly to \cite[Theorem 3.2.]{combettes2015stochastic} and \cite[Proposition 4.1.]{davis2022variance}.
We recall it here without repeating the proof to cover the composite mapping $\Phi$ in \eqref{eq:GE}.

\begin{lemma}\label{le:Opial_lemma}
Suppose that $\gra{\Phi}$ in \eqref{eq:GE} is closed and $\zer{\Phi} \neq\emptyset$.
Let $\sets{x^k}$ be a sequence of random vectors such that for all $x^{\star} \in \zer{\Phi}$, $\sets{\norms{x^k - x^{\star}}^2}$ almost surely converges to a $[0, \infty)$-valued random variable. 
In addition, assume that $\sets{\norms{w^k}}$ for $(x^k, w^k) \in \gra{\Phi}$ also almost surely converges to zero.
Then, $\sets{x^k}$ almost surely converges to a $\zer{\Phi}$-valued random variable $x^{\star} \in \zer{\Phi}$.
\end{lemma}

\revise{
\beforesubsec
\subsection{\textbf{Verifying Example 3 in Subsection~\ref{subsec:co-hypomonotone_example}}}\label{apdx:subsec:co-hypomonotone_example}
\aftersubsec
We verify the co-hypomonotonicity of the operator in \textbf{Example 3} from Subsection~\ref{subsec:co-hypomonotone_example}.
For any $x,y\in\dom{A}$, let $w_x\in Ax$ and $w_y\in Ay$.
Then, $u := -\rho^{-1}x - w_x \in Tx$ and $v := -\rho^{-1}y - w_y \in Ty$.
Therefore, we have
\begin{equation*}
\arraycolsep=0.2em
\begin{array}{lcl}
	\iprods{u - v, x - y} & = & -\rho^{-1}\norms{x - y}^2 - \iprods{w_x - w_y, x - y}, \vspace{1ex}\\
	\rho\norms{u-v}^2 & = & \rho\norms{ \rho^{-1}(x - y) + (w_x - w_y) }^2 = \rho^{-1}\norms{x-y}^2 + 2\iprods{w_x - w_y, x - y} + \rho\norms{w_x - w_y}^2.
\end{array}
\end{equation*}
Adding these two identities gives
$\iprods{u - v, x - y} + \rho\norms{u-v}^2 = \iprods{w_x - w_y, x - y} + \rho\norms{w_x - w_y}^2 \geq 0$,
where the inequality follows from the monotonicity of $A$ and $\rho>0$.
Consequently, $\iprods{u - v, x - y} \geq -\rho\norms{u - v}^2$ for all $(x,u),(y,v)\in\gra{T}$, showing that $T$ is $\rho$-co-hypomonotone.
\Eproof
}
\beforesubsec
\subsection{\textbf{Verifying examples in Subsection~\ref{subsec:nonmonotone_example}.}}\label{apdx:subsec:example_verification}
\aftersubsec
We present more details of Subsection~\ref{subsec:nonmonotone_example}.

\vspace{0.75ex}
\noindent\textit{\textbf{Example 2}}.
We provide a more detailed explanation of \textbf{Example 2} in Subsection~\ref{subsec:nonmonotone_example}.
First, using Jensen's inequality, it is easy to check that $G$ satisfies \eqref{eq:G_Lipschitz} in Assumption \ref{as:A1} with $\alpha = 1$.
Second, for all $x, y \in \dom{\Phi}$, let $w_x := Gx + v_x$ and $w_y := Gy + v_y$, where $v_x \in Tx$ and $v_y \in Ty$, respectively. 
Then, from Young's inequality, we have $(1+\tau)\norms{Gx - Gy}^2 - \norms{v_x - v_y}^2 \geq - \frac{(1+\tau)}{\tau}\norms{w_x - w_y}^2$ for any $\tau > 0$.
Utilizing the assumptions in \cite{tran2024accelerated,TranDinh2025a}, Jensen's inequality, and the last relation, we can easily prove that
\begin{equation*}
\begin{array}{lcl}
	\iprods{w_x - w_y, x - y} &=& \iprods{Gx - Gy, x - y} + \iprods{v_x - v_y, x - y} \vspace{1ex}\\
	&\geq& \frac{1}{L}\Expsn{\xi}{\norms{\mbf{G}(x, \xi) - \mbf{G}(y, \xi) }^2} - \rho \norms{v_x - v_y}^2 \vspace{1ex}\\
	&\geq& \big[ \frac{1}{L} - (1+\tau)\rho \big] \Expsn{\xi}{\norms{\mbf{G}(x, \xi) - \mbf{G}(y, \xi)}^2} + \rho \big[(1+\tau) \norms{Gx - Gy}^2 - \norms{v_x - v_y}^2 \big] \vspace{1ex}\\
	&\geq&  \big[\frac{1}{L} - (1+\tau)\rho \big] \Expsn{\xi}{\norms{\mbf{G}(x,\xi) - \mbf{G}(y, \xi)}^2} - \frac{(1+\tau)\rho}{\tau}\norms{w_x - w_y}^2.
\end{array}
\end{equation*}
This inequality shows that $\Phi := G + T$ studied in \cite{tran2024accelerated,TranDinh2025a} satisfies \eqref{eq:G_cohypomonotonicity} in Assumption \ref{as:A3} with $\rho_n := \frac{(1+\tau)\rho}{\tau} \geq 0$, $0 \leq \rho_c \leq \frac{1}{L} - (1+\tau)\rho$.
Clearly, to guarantee $\rho_c \geq 0$, we require $0 \leq L\rho \leq \frac{1}{1+\tau} < 1$.

\vspace{1ex}
\noindent\textit{\textbf{Example 4}}.
Since $R(x,y) \leq M$, we have $n \sum_{i=1}^n \norms{G_i x - G_i y}^2 \leq M \norms{\sum_{i=1}^n (G_i x - G_i y)}^2$.
It implies that $\Expsn{\xi}{\norms{\mbf{G}(x,\xi) - \mbf{G}(y,\xi) }^2} = \frac{1}{n} \sum_{i=1}^n \norms{G_i x - G_i y}^2 \leq \frac{M}{n^2} \norms{\sum_{i=1}^n (G_i x - G_i y)}^2 = M\norms{Gx - Gy}^2$.\\
Now, since $G$ is $\rho$-co-hypomonotone, we have
\begin{equation*}
\begin{array}{lcl}
	\iprods{Gx - Gy, x - y} &\geq& -\rho \norms{Gx - Gy}^2 = -\rho_n \norms{Gx - Gy}^2 + (\rho_n - \rho)\norms{Gx - Gy}^2 \vspace{1ex}\\
	&\geq& -\rho_n \norms{Gx - Gy}^2 + \frac{\rho_n - \rho}{M}\Expsn{\xi}{\norms{\mbf{G}(x, \xi) - \mbf{G}(y, \xi) }^2}.
\end{array}
\end{equation*}
Note that in order for the last inequality to hold (and to match the condition~\eqref{eq:G_cohypomonotonicity2} in Assumption~\ref{as:A3}), we require $\rho_n > \rho$.
In this case, we get $\rho_c := \frac{\rho_n - \rho}{M}$.

\vspace{1ex}
\noindent\textit{\textbf{Example 5}}.
By the definitions of $\bar{\sigma}_i$ and $\underline{\sigma}$, we have $\norms{G_ix - G_iy}\leq \bar{\sigma}_i\norms{x - y}$ and $\norms{Gx - Gy} \geq \underline{\sigma}\norms{x - y}$, respectively.
Combining these two relations, we get $\norms{G_ix - G_iy} \leq \frac{\bar{\sigma}_i}{\underline{\sigma}} \norms{Gx - Gy}$.
This implies $\Expsn{\xi}{\norms{\mbf{G}(x,\xi) - \mbf{G}(y, \xi) }^2} := \frac{1}{n} \sum_{i=1}^n \norms{G_i x - G_i y}^2 \leq \frac{1}{\underline{\sigma}^2} \left(\frac{1}{n}\sum_{i=1}^n \bar{\sigma}_i^2\right) \norms{Gx - Gy}^2 = \frac{\bar{\sigma}^2}{\underline{\sigma}^2}\norms{Gx - Gy}^2$, where we redefine $\bar{\sigma}^2 := \frac{1}{n}\sum_{i=1}^n \bar{\sigma}_i^2$.\\
Now, since $G$ is $\rho$-co-hypomonotone, we have
\myeqn{
	\begin{array}{lcl}
	\iprods{Gx - Gy, x - y} &\geq& -\rho \norms{Gx - Gy}^2 = -\rho_n \norms{Gx - Gy}^2 + (\rho_n - \rho)\norms{Gx - Gy}^2 \vspace{0.5ex}\\
	&\geq& -\rho_n \norms{Gx - Gy}^2 + \frac{\underline{\sigma}^2}{\bar{\sigma}^2}(\rho_n - \rho)\Expsn{\xi}{\norms{\mbf{G}(x, \xi) - \mbf{G}(y, \xi) }^2}.
	\end{array}
}	
Note that in order for the last inequality to hold (and to match the condition~\eqref{eq:G_cohypomonotonicity2} in Assumption~\ref{as:A3}), we require $\rho_n > \rho$.
In this case, we also get $\rho_c := \frac{\underline{\sigma}^2}{\bar{\sigma}^2}(\rho_n - \rho)$.

\beforesec
\section{The Proof of Lemmas~\ref{le:VFOG_key_estimate1}, \ref{le:VFOG_key_estimate2}, and \ref{le:Ek_lowerbound}.}\label{apdx:sec:technical_lemmas_proofs}
\aftersec
This appendix presents the full proof of the three lemmas: Lemmas~\ref{le:VFOG_key_estimate1}, \ref{le:VFOG_key_estimate2}, and \ref{le:Ek_lowerbound}.

\beforesubsec
\subsection{\textbf{The proof of Lemma~\ref{le:VFOG_key_estimate1}.}}\label{apdx:le:VFOG_key_estimate1}
\aftersubsec
First, using $e^k$ and $\hat{e}^k$ from \eqref{eq:SFOG_quantities}, we can write $\widetilde{G}y^k = Gy^k + e^k$ and $\widetilde{G}y^{k-1} = Gx^k + \hat{e}^k$.
Denote $w^k := Gx^k + v^k$ and $\hat{w}^k := Gy^{k-1} + v^k$ for a given $v^k \in Tx^k$.
Then, we can rewrite the first and third lines of \eqref{eq:VFOG} equivalently to
\begin{equation*} 
\arraycolsep=0.2em
\begin{array}{lcl}
t_kx^{k+1} & = & sz^k + (t_k-s)x^k - \eta t_k(\hat{w}^{k+1}  + e^k) + \beta_kt_k(w^k + \hat{e}^k).
\end{array}
\end{equation*}
Rearranging this expression in two different ways, and then using $sz^k = sz^{k+1} + \gamma_k(\widetilde{G}y^{k-1} + v^k) = sz^{k+1} + \gamma_k (w^k + \hat{e}^k)$ from the last line of \eqref{eq:VFOG}, we have
\begin{equation*} 
\arraycolsep=0.2em
\begin{array}{lcl}
t_k(t_k-s)(x^{k+1} - x^k) &= & s(t_k-s)(z^k - x^k) - \eta t_k(t_k-s)(\hat{w}^{k+1} + e^k) + \beta_kt_k(t_k-s)(w^k + \hat{e}^k), \vspace{1ex}\\
t_k(t_k-s)(x^{k+1} - x^k) &= & st_k (z^k - x^{k+1}) - \eta t_k^2(\hat{w}^{k+1} + e^k) + \beta_kt_k^2(w^k + \hat{e}^k), \vspace{1ex}\\
& = & st_k(z^{k+1} - x^{k+1}) -  \eta t_k^2(\hat{w}^{k+1} + e^k) + t_k(\beta_kt_k + \gamma_k )(w^k + \hat{e}^k).
\end{array}
\end{equation*}
Second, by the condition \eqref{eq:G_cohypomonotonicity} of $\Phi$ in Assumption~\ref{as:A3} and using the notion $\tnorms{\cdot}$, we have
\begin{equation*} 
\arraycolsep=0.2em
\begin{array}{lcl}
t_k(t_k-s)\iprods{w^{k+1}, x^{k+1} - x^k} - t_k(t_k-s)\iprods{w^k, x^{k+1} - x^k} &\geq & -\rho_n t_k(t_k-s)\norms{w^{k+1} - w^k}^2 \vspace{1ex}\\
&& + {~} \rho_c t_k(t_k-s)\tnorms{Gx^{k+1} - Gx^k}^2.
\end{array}
\end{equation*}
Combining the last three expressions, we get
\begin{equation*} 
\arraycolsep=0.2em
\begin{array}{lcl}
\Tc_{[1]} &:=& s(t_k-s)\iprods{w^k, x^k - z^k} - st_k\iprods{w^{k+1}, x^{k+1} - z^{k+1}} \vspace{1ex}\\
& \geq & \eta t_k^2\iprods{w^{k+1}, \hat{w}^{k+1} } - t_k( \beta_kt_k + \gamma_k ) \iprods{w^{k+1}, w^k}  -  \eta t_k(t_k-s)\iprods{w^k, \hat{w}^{k+1} } \vspace{1ex}\\
&& + {~}  \beta_kt_k(t_k-s)\norms{w^k}^2 -  \rho_n t_k(t_k-s)\norms{w^{k+1} - w^k}^2  \vspace{1ex}\\
&& +  {~} \rho_c t_k(t_k-s)\tnorms{Gx^{k+1} - Gx^k}^2  + \eta t_k\iprods{t_k w^{k+1} - (t_k-s)w^k, e^k}\vspace{1ex}\\
&& - {~} \beta_kt_k \iprods{t_k w^{k+1} - (t_k-s)w^k, \hat{e}^k} - \gamma_k t_k\iprods{w^{k+1}, \hat{e}^k}.
\end{array}
\end{equation*}
Now, by Young's inequality, for any $c_1 > 0$ and $c_2 > 0$, we can show that
\begin{equation*} 
\arraycolsep=0.2em
\begin{array}{lcl}
\Tc_{[2]} &:=&  \iprods{t_k w^{k+1} - (t_k-s)w^k, e^k}\vspace{1ex}\\
& \geq & -\frac{c_1}{2 t_k}\norms{t_k w^{k+1} - (t_k-s)w^k }^2 - \frac{t_k}{2c_1}\norms{e^k}^2 \vspace{1ex}\\
& = & -\frac{c_1t_k}{2}\norms{w^{k+1}}^2 - \frac{c_1(t_k-s)^2}{2t_k}\norms{w^k}^2 + c_1(t_k-s) \iprods{w^{k+1}, w^k} - \frac{ t_k}{2c_1}\norms{e^k}^2, \vspace{1ex}\\
\Tc_{[3]} &:=& -  \iprods{t_k w^{k+1} - (t_k-s)w^k, \hat{e}^k}\vspace{1ex}\\
& \geq & -\frac{c_2\eta}{2\beta_k t_k}\norms{t_k w^{k+1} - (t_k-s)w^k }^2 - \frac{\beta_k t_k}{2c_2\eta}\norms{\hat{e}^k}^2 \vspace{1ex}\\
& = & -\frac{c_2\eta t_k}{2\beta_k}\norms{w^{k+1}}^2 - \frac{c_2\eta (t_k-s)^2}{2\beta_k t_k}\norms{w^k}^2 + \frac{c_2\eta(t_k-s)}{\beta_k} \iprods{w^{k+1}, w^k} - \frac{\beta_k t_k}{2c_2\eta}\norms{\hat{e}^k}^2, \vspace{1ex}\\
\Tc_{[4]}& := &  - \iprods{w^{k+1}, \hat{e}^k} \geq -\frac{1}{2}\norms{w^{k+1}}^2 - \frac{1}{2}\norms{\hat{e}^k}^2.
\end{array}
\end{equation*}
Substituting $\Tc_{[2]}$, $\Tc_{[3]}$, and $\Tc_{[4]}$ into $\Tc_{[1]}$, we can further derive that
\begin{equation*} 
\arraycolsep=0.2em
\begin{array}{lcl}
\Tc_{[1]} &:=& s(t_k-s)\iprods{w^k, x^k - z^k} - st_k\iprods{w^{k+1}, x^{k+1} - z^{k+1}} \vspace{1ex}\\
& \geq & \eta t_k^2\iprods{w^{k+1}, \hat{w}^{k+1} } - \eta t_k(t_k-s)\iprods{ w^k, \hat{w}^{k+1} } + \rho_c t_k(t_k-s)\tnorms{Gx^{k+1} - Gx^k}^2 \vspace{1ex}\\
&& {~} - \frac{t_k(t_k-s)}{2}\big[2(\rho_n - \beta_k)  + \frac{c_1\eta(t_k-s)}{t_k} + \frac{c_2\eta(t_k-s)}{t_k}  \big] \norms{ w^k }^2 \vspace{1ex}\\
&& {~} - \frac{t_k}{2}\big[ 2\rho_n(t_k-s) + c_1\eta t_k + c_2\eta t_k + \gamma_k \big] \norms{ w^{k+1} }^2 \vspace{1ex}\\
&& + {~} t_k\big[ 2\rho_n(t_k-s) + c_1\eta(t_k - s) + c_2\eta(t_k-s) - \beta_kt_k - \gamma_k \big]\iprods{ w^{k+1}, w^k } \vspace{1ex}\\
&& - {~} \frac{\eta t_k^2}{2c_1}\norms{e^k}^2  - \frac{t_k}{2}\big(\frac{\beta_k^2t_k}{c_2\eta} + \gamma_k \big)\norms{\hat{e}^k}^2.
\end{array}
\end{equation*}
Next, assume that $(c_1\eta + c_2\eta + 2\rho_n)(t_k-s)  - t_k\beta_k - \gamma_k  = 0$.
This condition is guaranteed if we update $\beta_k$ as in  \eqref{eq:SFOG_para_cond1}, i.e.:
\begin{equation*} 
\arraycolsep=0.2em
\begin{array}{lcl}
\beta_k = \frac{(c_1\eta + c_2\eta + 2\rho_n)(t_k-s) - \gamma_k}{t_k}. 
\end{array}
\end{equation*}
Partly substituting this $\beta_k$ into $\Tc_{[1]}$ and using the following identity
\begin{equation*} 
\arraycolsep=0.2em
\begin{array}{lcl}
t_k(t_k-s)\iprods{ w^k, \hat{w}^{k+1} } &= & \frac{t_k(t_k-s)}{2}\norms{ \hat{w}^{k+1} }^2 + \frac{t_k(t_k-s)}{2}\norms{ w^k }^2 - \frac{t_k(t_k-s)}{2}\norms{ \hat{w}^{k+1}  - w^k }^2
\end{array}
\end{equation*}
we can simplify the last expression $\Tc_{[1]}$ as
\begin{equation*} 
\arraycolsep=0.2em
\begin{array}{lcl}
\Tc_{[1]} &:=& s(t_k-s)\iprods{ w^k, x^k - z^k} - st_k\iprods{ w^{k+1}, x^{k+1} - z^{k+1}} \vspace{1ex}\\
& \geq & \eta t_k^2\iprods{ w^{k+1}, \hat{w}^{k+1} } - \frac{\eta t_k(t_k-s)}{2}\norms{ \hat{w}^{k+1} }^2 + \frac{\eta t_k(t_k-s)}{2}\norms{ \hat{w}^{k+1}  - w^k}^2 \vspace{1ex}\\
&& {~} - \frac{t_k(t_k-s)}{2}\big[2\rho_n + \eta + \frac{(c_1+c_2)\eta(t_k-s)}{t_k} - \frac{2(c_1\eta + c_2\eta + 2\rho_n)(t_k-s)}{t_k} + \frac{2\gamma_k}{t_k}  \big] \norms{ w^k }^2 \vspace{1ex}\\
&& {~} - \frac{t_k}{2}\big[ (c_1+c_2)\eta t_k + \gamma_k + 2\rho_n(t_k-s) \big] \norms{ w^{k+1} }^2 \vspace{1ex}\\
&& - {~} \frac{\eta t_k^2}{2c_1}\norms{e^k}^2 - \frac{t_k}{2}\big(\frac{\beta_k^2t_k}{c_2\eta} + \gamma_k \big)\norms{\hat{e}^k}^2 
+ {~} \rho_c t_k(t_k-s)\tnorms{ Gx^{k+1} - Gx^k}^2.
\end{array}
\end{equation*}
For a given $\omega \geq 0$, we define $M^2 := 2(1+\omega)L^2$.
Then, from the second and third lines of \eqref{eq:VFOG}, we have 
\begin{equation*} 
\arraycolsep=0.2em
\begin{array}{lcl}
t_k(x^{k+1} - y^k) = -\eta t_k (\hat{w}^{k+1}  + e^k) + \eta t_k(w^k + \hat{e}^k) = -\eta t_k( \hat{w}^{k+1} - w^k) - \eta t_k(e^k -\hat{e}^k).
\end{array}
\end{equation*}
Applying the condition \eqref{eq:G_Lipschitz} from Assumption~\ref{as:A1}, we can show that
\begin{equation*} 
\arraycolsep=0.2em
\begin{array}{lcl}
\norms{w^{k+1} - \hat{w}^{k+1}}^2 + E_k = \norms{Gx^{k+1} -  Gy^k}^2 + E_k \leq (1+\omega)L^2\norms{x^{k+1} - y^k}^2,
\end{array}
\end{equation*}
where 
\begin{equation*} 
\arraycolsep=0.2em
\begin{array}{lcl}
E_k & := &  (\alpha-1)\norms{Gx^{k+1} -  Gy^k}^2 + (1-\alpha)\tnorms{Gx^{k+1} -  Gy^k}^2 + \omega L^2\norms{x^{k+1} - y^k}^2 \vspace{1ex}\\
& = &  (1-\alpha)\big[ \tnorms{w^{k+1} -  \hat{w}^{k+1} }^2 - \norms{w^{k+1} -  \hat{w}^{k+1} }^2\big] + \omega L^2\norms{x^{k+1} - y^k}^2.
\end{array}
\end{equation*}
Combining the last two expressions, we get
\begin{equation*} 
\arraycolsep=0.2em
\begin{array}{lcl}
t_k^2\norms{w^{k+1} -  \hat{w}^{k+1} }^2 +   t_k^2E_k &\leq & (1+\omega)L^2t_k^2\norms{x^{k+1} - y^k}^2 \vspace{1ex}\\
& = & (1+\omega)L^2\eta^2t_k^2\norms{  \hat{w}^{k+1} - w^k + e^k - \hat{e}^k }^2 \vspace{1ex}\\

&\leq& \revise{M^2 \eta^2 t_k^2 \norms{\hat{w}^{k+1} - w^k}^2 + M^2\eta^2 t_k^2 \norms{e^k - \hat{e}^k}^2} \vspace{1ex}\\

&\leq & M^2\eta^2 t_k^2\norms{  \hat{w}^{k+1} - w^k }^2 + \revise{(1+\theta)}M^2\eta^2t_k^2 \norms{e^k}^2  + \revise{(1+\theta^{-1})}M^2\eta^2t_k^2 \norms{\hat{e}^k}^2,
\end{array}
\end{equation*}
\revise{for some $\theta > 0$.}
Multiplying both sides of this inequality by $\frac{\eta}{2}$ and partly expanding the resulting expression and rearranging the result, we can deduce that
\myeqn{
\arraycolsep=0.2em
\begin{array}{lcl}
0 &\geq & \frac{\eta t_k^2}{2}\norms{w^{k+1}}^2 + \frac{\eta t_k^2}{2}\norms{  \hat{w}^{k+1} }^2 - \eta t_k^2 \iprods{ w^{k+1},  \hat{w}^{k+1} } + \frac{\eta t_k^2}{2}E_k \vspace{1ex}\\
&& - {~} \frac{M^2\eta^3t_k^2}{2}\norms{ \hat{w}^{k+1} - w^k}^2 - \frac{\revise{(1+\theta)}M^2\eta^3t_k^2}{\revise{2}} \norms{e^k}^2 - \frac{\revise{(1+\theta^{-1})}M^2\eta^3t_k^2}{\revise{2}} \norms{\hat{e}^k}^2.
\end{array}
}
Adding this inequality to $\Tc_{[1]}$ yields
\myeqn{
\arraycolsep=0.2em
\begin{array}{lcl}
\Tc_{[1]} &:=& s(t_k-s)\iprods{w^k, x^k - z^k} - st_k\iprods{w^{k+1}, x^{k+1} - z^{k+1}} \vspace{1ex}\\
& \geq &  \frac{\eta t_k^2}{2}E_k  + \frac{\eta t_k[ t_k - s - M^2\eta^2t_k]}{2}\norms{  \hat{w}^{k+1}  - w^k }^2 +  \rho_ct_k(t_k-s)\tnorms{ Gx^{k+1} - Gx^k }^2 \vspace{1ex}\\
&& - {~}  \frac{t_k(t_k-s)}{2}\big[ 2\rho_n + \eta + \frac{(c_1+c_2)\eta(t_k-s)}{t_k}  - \frac{2(c_1\eta + c_2\eta + 2\rho_n)(t_k-s)}{t_k} + \frac{2\gamma_k}{t_k}  \big] \norms{ w^k }^2 \vspace{1ex}\\
&& + {~}  \frac{t_k}{2}\Big(  \eta t_k - \big[ (c_1+c_2)\eta t_k  + \gamma_k + 2\rho_n(t_k-s) \big] \Big) \norms{ w^{k+1} }^2 + \frac{\eta st_k}{2}\norms{  \hat{w}^{k+1} }^2 \vspace{1ex}\\
&& - {~} \big( \frac{\eta }{2c_1} + \frac{\revise{(1+\theta)}M^2\eta^3}{\revise{2}} \big)t_k^2 \norms{e^k}^2 
-  \big( \frac{\beta_k^2t_k}{2 c_2\eta} + \frac{\gamma_k}{2} + \frac{\revise{(1+\theta^{-1})}M^2\eta^3t_k}{\revise{2}} \big)t_k \norms{\hat{e}^k}^2.
\end{array}
}
Next, utilizing $s(z^{k+1} - z^k) = -\gamma_k (w^k + \hat{e}^k)$ from the last line of \eqref{eq:VFOG} and Young's inequality, for any $c_k > 0$, we can derive that
\myeqn{
\arraycolsep=0.2em
\begin{array}{lcl}
\bar{\Tc}_{[4]} &:=& \frac{s^2(s-1)}{2\gamma_k} \norms{ z^k - x^{\star}}^2 - \frac{s^2(s-1)}{2\gamma_k}\norms{z^{k+1} - x^{\star}}^2 \vspace{1ex}\\
& = &  - \frac{s^2(s-1)}{\gamma_k}\iprods{z^{k+1} - z^k, z^k - x^{\star}} -  \frac{s^2(s-1)}{2\gamma_k}\norms{z^{k+1} - z^k}^2 \vspace{1ex}\\
& = & s(s-1) \iprods{w^k + \hat{e}^k, z^k - x^{\star}} - \frac{(s-1)\gamma_k}{2}\norms{w^k + \hat{e}^k}^2 \vspace{1ex}\\
& \geq & s(s-1)\iprods{w^k, z^k - x^{\star}} - \frac{s^2(s-1)c_k}{2}\norms{z^k - x^{\star}}^2 - \mred{(s-1)\gamma_k}\norms{w^k}^2 -  \frac{s-1}{2}\big(  2\gamma_k + \frac{1}{c_k} \big) \norms{\hat{e}^k}^2.
\end{array}
}
Since $\gamma_k := \frac{\gamma(t_k - 1)}{t_k}$ and $t_k = t_{k-1} + 1$ from \eqref{eq:SFOG_para_cond1} for some $\gamma \geq 0$, we can easily prove that
\myeqn{
\arraycolsep=0.2em
\begin{array}{lcl}
\frac{1}{2\gamma_{k-1}} - \frac{1}{2\gamma_k} & = & \frac{1}{2\gamma}\big(\frac{t_{k-1}}{t_{k-1}- 1} - \frac{t_k}{t_k- 1}\big)  = \frac{(t_k- 1)(t_k-1) - t_k(t_k-2 )}{2\gamma(t_k- 1)(t_k-2)} = \frac{1}{2\gamma(t_k-1)(t_k-2)}.
\end{array}
}
Therefore, if we choose $c_k := \frac{1}{\gamma_{k-1}} - \frac{1}{\gamma_k} = \frac{1}{\gamma(t_k- 1)(t_k-2)}$, then we have $\frac{1}{2\gamma_k} + \frac{c_k}{2} = \frac{1}{2\gamma_{k-1}}$.
Using this relation, the expression $\bar{\Tc}_{[4]}$ implies that
\begin{equation*} 
\arraycolsep=0.2em
\begin{array}{lcl}
\Tc_{[4]} &:=& \frac{s^2(s-1)}{2\gamma_{k-1}} \norms{ z^k - x^{\star}}^2 - \frac{s^2(s-1)}{2\gamma_k}\norms{z^{k+1} - x^{\star}}^2 \vspace{1ex}\\
& \geq & \mred{s(s-1)\iprods{w^k, z^k - x^{\star}} - \frac{\gamma(s-1)(t_k - 1)}{t_k}\norms{w^k}^2 - \frac{(s-1)\gamma}{2}\big[ \frac{2(t_k-1)}{t_k} + (t_k-1)(t_k-2)  \big] \norms{\hat{e}^k}^2.}
\end{array}
\end{equation*}
Since $t_k - t_{k-1} = 1$ due to \eqref{eq:SFOG_para_cond1}, we  have $s(t_{k-1} - t_k + s) = s(s-1)$.
Adding $\Tc_{[4]}$ to $\Tc_{[1]}$ and using $s(t_{k-1} - t_k + s) = s(s-1)$, we have
\begin{equation*} 
\arraycolsep=0.2em
\begin{array}{lcl}
\Tc_{[5]} &:=& st_{k-1}\iprods{w^k, x^k - z^k} + \frac{s^2(s-1)}{2\gamma_{k-1}} \norms{ z^k - x^{\star}}^2 \vspace{1ex}\\
&& - {~} st_k\iprods{ w^{k+1}, x^{k+1} - z^{k+1}} - \frac{s^2(s-1)}{2\gamma_k}\norms{z^{k+1} - x^{\star}}^2  \vspace{1ex}\\
& \geq &  \frac{\eta t_k^2}{2}E_k  + \frac{\eta t_k(t_k-s - M^2\eta^2t_k)}{2}\norms{  \hat{w}^{k+1}  - w^k }^2 
+ \rho_ct_k(t_k-s)\tnorms{Gx^{k+1} - Gx^k}^2 \vspace{1ex}\\
&& + {~}  s(s-1)\iprods{w^k, x^k - x^{\star}} +  \frac{\eta st_k}{2} \norms{  \hat{w}^{k+1} }^2  \vspace{1ex}\\
&& + {~} \frac{t_k}{2}\big[ (1 - c_1 - c_2)\eta t_k - \gamma_k -  2\rho_n(t_k-s)   \big] \norms{ w^{k+1} }^2 \vspace{1ex}\\ 
&& - {~}  \Big\{ \frac{t_k(t_k-s)}{2}\big[2\rho_n + \eta   - \frac{[(c_1 + c_2)\eta + 4\rho_n](t_k-s)}{t_k} + \frac{2\gamma(t_k-1)}{t_k^2}   \big]  + \frac{\gamma (s-1) (t_k-1)}{t_k} \Big\} \norms{ w^k }^2 \vspace{1ex}\\
&& - {~}  \Big[ \big( \frac{\beta_k^2 t_k}{2 c_2\eta} + \frac{\gamma (t_k-1)}{2 t_k} + \frac{\revise{(1+\theta^{-1})}M^2\eta^3t_k}{\revise{2}}  \big)t_k + \gamma (s-1) \big( \frac{\mred{t_k-1}}{t_k} + \frac{(t_k-1)(t_k-2)}{2} \big) \Big] \norms{\hat{e}^k}^2 \vspace{1ex} \\
&& - {~} \big( \frac{\eta }{2c_1} + \frac{\revise{(1+\theta)}M^2\eta^3}{2} \big)t_k^2 \norms{e^k}^2.
\end{array}
\end{equation*}
Now, we recall from \eqref{eq:SFOG_coeffs} that
\myeqn{
\arraycolsep=0.2em
\left\{\begin{array}{lcl}
a_{k+1} &:= &   (1 - c_1 - c_2)\eta t_k^2  -  2\rho_n t_k(t_k-s) - \gamma(t_k-1), \vspace{1ex}\\
\hat{a}_k &:= & \big[(1 - c_1 - c_2)\eta - 2\rho_n\big](t_k - s)^2  + s(\eta + 2\rho_n)(t_k-s)  + \frac{2\gamma(t_k-1)^2 }{t_k}, \vspace{1ex}\\ 
b_k &:= & \big(\frac{1}{2c_1} + \frac{\revise{(1+\theta)}M^2\eta^2}{\revise{2}} \big)\eta t_k^2, \vspace{1ex}\\
\hat{b}_k &:= & \big( \frac{\beta_k^2t_k}{2 c_2\eta} + \frac{\gamma(t_k-1)}{2t_k} + \frac{\revise{(1+\theta^{-1})}M^2\eta^3t_k  }{\revise{2}}\big)t_k + \frac{(s-1)\gamma}{2}\big[ \frac{2(t_k-1)}{ t_k } + (t_k-1)(t_k - 2)  \big].
\end{array}\right.
}
Then, we can rewrite the last expression $\Tc_{[5]}$ as 
\myeqn{
\arraycolsep=0.2em
\begin{array}{lcl}
\Tc_{[6]} &:=& \frac{a_k}{2}\norms{w^k}^2 +  st_{k-1}\iprods{w^k, x^k - z^k} + \frac{s^2(s-1)}{2\gamma_{k-1}} \norms{ z^k - x^{\star}}^2 \vspace{1ex}\\
&& - {~} \frac{a_{k+1}}{2}\norms{w^{k+1} }^2 -  st_k\iprods{w^{k+1}, x^{k+1} - z^{k+1}} - \frac{s^2(s-1)}{2\gamma_k}\norms{z^{k+1} - x^{\star}}^2  \vspace{1ex}\\
& \geq &  \frac{\eta t_k^2}{2}E_k + \rho_ct_k(t_k-s)\tnorms{Gx^{k+1} - Gx^k}^2 + \frac{\eta t_k(t_k-s - M^2\eta^2t_k)}{2}\norms{ \hat{w}^{k+1} - w^k}^2 \vspace{1ex}\\
&& + {~} \frac{\eta st_k}{2} \norms{  \hat{w}^{k+1} }^2 +  s(s-1)\iprods{w^k, x^k - x^{\star}}  + \frac{a_k - \hat{a}_k}{2} \norms{ w^k }^2 - b_k \norms{e^k}^2 - \hat{b}_k \norms{\hat{e}^k}^2.
\end{array}
}
Substituting $\Pc_k$ from \eqref{eq:SFOG_Ek_func} into this new expression $\Tc_{[6]}$, we eventually get \eqref{eq:SFOG_key_est1}.
\Eproof

\beforesubsec
\subsection{\textbf{The proof of Lemma~\ref{le:VFOG_key_estimate2}.}}\label{apdx:le:VFOG_key_estimate2}
\aftersubsec
First, since $s > 2$ due to \eqref{eq:SFOG_para_cond2} and $t_k = t_{k-1} + 1$ due to \eqref{eq:SFOG_para_cond1},  we have $t_k(t_k-s) \leq (t_k-1)^2$.
Thus, for $\psi := (1 - c_1 - c_2)\eta - 2\rho_n$, from \eqref{eq:SFOG_coeffs}, we can show that
\myeqn{
\arraycolsep=0.1em
\left\{\begin{array}{lcl}
a_k &:= & \psi t_{k-1}^2 + 2s\rho_n t_{k-1} - \gamma(t_{k-1} - 1) =  \psi (t_k - 1)^2 + (2s\rho_n - \gamma)(t_k - 1) + \gamma, \vspace{1ex}\\
\hat{a}_k &:= & \psi(t_k - s)^2  + s(\eta + 2\rho_n)(t_k-s)  + \frac{2\gamma(t_k-1)^2}{t_k} \vspace{1ex}\\ 
& \leq &  \psi(t_k - s)^2  + [s(\eta + 2\rho_n) + 2\gamma](t_k - 1) - s(s-1)(\eta + 2\rho_n).
\end{array}\right.
}
Using these relations, we can lower bound
\myeqn{
\arraycolsep=0.2em
\begin{array}{lcl}
a_k - \hat{a}_k & \geq &  \psi(s-1)(2t_k - s - 1) - (s\eta + 3\gamma)(t_k - 1) + \gamma + s(s-1)(\eta + 2\rho_n) \vspace{1ex}\\
& = & [2\psi(s-1) - s\eta - 3\gamma](t_k - 1) - \psi(s-1)^2 + \gamma + s(s-1)(\eta + 2\rho_n) \vspace{1ex}\\
& = & \big[ [s-2 - 2(c_1+c_2)(s-1)]\eta - 4\rho_n(s-1) - 3\gamma \big](t_k - 1) + \gamma + s(s-1)(\eta + 2\rho_n) \vspace{1ex}\\
&& - {~} (s-1)^2[(1-c_1-c_2)\eta - 2\rho_n].
\end{array}
}
Next, let us choose \revise{$c_1 := \frac{s-2}{16(s-1)} > 0$ and $c_2 := \frac{s-2}{8(s-1)} > 0$}.
Then, we obtain from the above expressions that
\myeqn{
\arraycolsep=0.2em
\begin{array}{lcl}
a_k &= & \big[\revise{\frac{(13s-10)\eta}{16(s-1)}} - 2\rho_n\big]t_{k-1}^2 + 2s\rho_nt_{k-1} - \gamma (t_{k-1} - 1), \vspace{1ex}\\
a_k - \hat{a}_k & \geq &  \big[ \revise{\frac{5(s-2)}{8}}\eta - 4(s-1)\rho_n - 3\gamma \big](t_k - 1) + \revise{\frac{(s-1)(3s+10)\eta}{16}} + 2(2s^2-3s+1)\rho_n + \gamma  =: \phi_k.
\end{array}
}
Using \revise{$c_1 + c_2 = \frac{3(s-2)}{16(s-1)}$} and \revise{$(s-2)\eta \geq 8(s-1)\rho_n$} from \eqref{eq:SFOG_para_cond2}, we can derive from \eqref{eq:SFOG_para_cond1} that
\myeqn{
\arraycolsep=0.2em
\begin{array}{lcl}
\beta_k &= & \frac{[(c_1+c_2)\eta + 2\rho_n](t_k-s) - \gamma_k}{t_k} \leq \frac{[(c_1+c_2)\eta + 2\rho_n](t_k-s)}{t_k} = \big[\frac{\revise{3}(s-2)}{\revise{16}(s-1)} + \frac{2\rho_n}{\eta}]\frac{\eta(t_k-s)}{t_k} \leq \frac{\revise{7(s-2)}\eta(t_k-s)}{\revise{16(s-1)}t_k}.
\end{array}
}
\revise{Next, let us choose $\theta := 8$.
Then,}
since $M^2\eta^2 \leq \frac{t_k-s}{t_k} \revise{\ \leq 1}$ and $s > 2$ from \eqref{eq:SFOG_para_cond2}, and $\beta_k \leq \frac{\revise{7(s-2)}\eta(t_k-s)}{\revise{16(s-1)}t_k}$, we can upper bound $b_k$ and $\hat{b}_k$ from \eqref{eq:SFOG_coeffs} respectively as follows:
\myeqn{ 
\arraycolsep=0.2em
\begin{array}{lcl}
b_k & = & \big(\frac{1}{2c_1} + \frac{\revise{(1+\theta)}M^2\eta^2}{\revise{2}} \big)\eta t_k^2 \leq \big[ \frac{8(s-1)}{s-2} + \frac{9}{2} \big]\eta t_k^2 = \frac{\revise{(25s - 34)}\eta t_k^2}{\revise{2}(s-2)}, \vspace{1ex}\\


\hat{b}_k &= & \big( \frac{\beta_k^2t_k}{2 c_2\eta} + \frac{\gamma(t_k-1)}{2t_k} + \frac{\revise{(1+\theta^{-1})}M^2\eta^3t_k}{\revise{2}} \big)t_k + \frac{(s-1)\gamma}{2}\big[ \frac{2(t_k-1)}{t_k} + (t_k-1)(t_k-2) \big] \vspace{1ex}\\

& \leq & \frac{\revise{49(s-2)}\eta (t_k-s)^2}{\revise{64(s-1)}} + \frac{\revise{9}\eta t_k(t_k-s)}{\revise{16}} + \frac{\gamma (t_k - 1)}{2} + \frac{\gamma (s-1) (t_k-1)}{t_k} + \frac{(s-1)\gamma (t_k-1)(t_k-2)}{2} \vspace{1ex}\\

& \leq &  \frac{\revise{49(s-2)}\eta}{\revise{64(s-1)}}[t_{k-1}^2 - 2(s-1)t_{k-1} + (s-1)^2 ] + \frac{\revise{9}\eta}{\revise{16}} [ t_{k-1}^2 - (s-2)t_{k-1} - (s-1)] \vspace{1ex}\\
&& + {~} \frac{(s-1)\gamma}{2}(t_{k-1}^2 - t_{k-1}) + \frac{\gamma t_{k-1}}{2} + (s-1)\gamma \vspace{1ex}\\

& = & \big[ \frac{\revise{49(s-2)}\eta}{\revise{64(s-1)}} + \frac{\revise{9}\eta}{\revise{16}} + \frac{(s-1)\gamma}{2} \big] t_{k-1}^2 - \big[ \frac{\revise{49}(s-2)\eta}{\revise{32}} + \frac{\revise{9}(s-2)\eta}{\revise{16}} + \frac{(s-1)\gamma}{2} - \frac{\gamma}{2} \big] t_{k-1}\vspace{1ex}\\
&& + {~} \eta(s-1)\big[ \frac{\revise{49}(s-2)}{\revise{64}} - \revise{\frac{9}{16}} \big]  + (s-1)\gamma \vspace{1ex}\\

& = & \big[ \frac{\revise{(85s - 134)}\eta}{\revise{64(s-1)}} + \frac{(s-1)\gamma}{2}\big] t_{k-1}^2 - \big[ \frac{\revise{67(s-2)}\eta}{\revise{32}} + \frac{(s-2)\gamma}{2}\big] t_{k-1} + \frac{(s-1)[\revise{(49s - 134)}\eta + \revise{64}\gamma]}{\revise{64}}\vspace{1ex}\\

& \leq &  \big[ \frac{\revise{(85s - 134)}\eta}{\revise{64(s-1)}} + \frac{(s-1)\gamma}{2}\big] t_{k-1}^2, 
\end{array}
}
provided that $ \big[ \frac{\revise{67(s-2)}\eta}{\revise{32}} + \frac{(s-2)\gamma}{2}\big] t_{k-1} \geq \frac{(s-1)[\revise{(49s - 134)}\eta + \revise{64}\gamma]}{\revise{64}}$.

Using these bounds into \eqref{eq:SFOG_key_est1}, we get
\begin{equation*} 
\hspace{-0.5ex}
\arraycolsep=0.2em
\begin{array}{lcl}
\Pc_k - \Pc_{k+1}  & \geq &  \frac{\eta t_k^2}{2}\big[ (1-\alpha) \big( \tnorms{ w^{k+1} - \hat{w}^{k+1} }^2 - \norms{ w^{k+1} - \hat{w}^{k+1} }^2 \big) + \omega L^2\norms{x^{k+1} - y^k }^2 \big] \vspace{1ex}\\
&& + {~} \rho_c t_k(t_k-s)\tnorms{ Gx^{k+1} - Gx^k }^2   + \frac{\eta t_k(t_k-s)}{2}\big(1  - \frac{M^2\eta^2 t_k}{t_k - s} \big) \norms{ \hat{w}^{k+1}  - w^k }^2 \vspace{1ex}\\
&& + {~}   \frac{\eta st_k}{2} \norms{ \hat{w}^{k+1} }^2 +   s(s-1)\iprods{ w^k, x^k - x^{\star}}  + \frac{\phi_k}{2} \norms{ w^k }^2 \vspace{1ex}\\
&& - {~}  \frac{\revise{(25s - 34)}\eta t_k^2}{\revise{2}(s-2)}\norms{e^k}^2 -   \big[ \revise{\frac{85s - 134}{64(s-1)}} + \frac{(s-1)\gamma}{2\eta}\big] \eta t_{k-1}^2 \norms{\hat{e}^k}^2.
\end{array}
\hspace{-2ex}
\end{equation*}
By Young's inequality, from \eqref{eq:SFOG_quantities} we get $\norms{\hat{e}^k}^2 \leq \revise{\big(1+\frac{1}{8}\big)}\norms{Gx^k - Gy^{k-1}}^2 + \revise{\big(1+8\big)}\norms{e^{k-1}}^2 = \revise{\frac{9}{8}}\norms{w^k - \hat{w}^k }^2 + \revise{9}\norms{e^{k-1}}^2$.
Utilizing this relation into the last estimate, we obtain \eqref{eq:SFOG_key_est2}.
Finally, since $t_{k-1} \geq s$, $s \geq 2$, and \revise{$0 < \gamma \leq \frac{9}{16}\eta$} from \eqref{eq:SFOG_para_cond2}, we can prove that the condition $ \big[ \frac{\revise{67(s-2)}\eta}{\revise{32}} + \frac{(s-2)\gamma}{2}\big] t_{k-1} \geq \frac{(s-1)[\revise{(49s - 134)}\eta + \revise{64}\gamma]}{\revise{64}}$ automatically holds.
\Eproof

\beforesubsec
\subsection{\textbf{The proof of Lemma~\ref{le:Ek_lowerbound}.}}\label{apdx:le:Ek_lowerbound}
\aftersubsec
Since $w^k := Gx^k + v^k$ for $v^k \in Tx^k$, from \eqref{eq:SFOG_Ek_func}, for any $c_k > 0$, we can rewrite $\Pc_k$ as 
\begin{equation*} 
\arraycolsep=0.2em
\begin{array}{lcl}
\Pc_k & = & \frac{a_k}{2}\norms{w^k}^2 - st_{k-1}\iprods{w^k, z^k - x^{\star}} + \frac{s^2(s-1)t_{\revise{k-1}}}{2\gamma (t_{\revise{k-1}}-1)}\norms{z^k - x^{\star}}^2 + st_{k-1}\iprods{w^k, x^k - x^{\star}} \vspace{1ex}\\
& = & \frac{a_k-c_k^2t_{k-1}^2}{2}\norms{w^k}^2 + \frac{1}{2c_k^2}\norms{c_k^2t_{k-1}w^k - s(z^k - x^{\star})}^2 + \frac{s^2}{2}\big[  \frac{(s-1)t_{\revise{k-1}}}{\gamma (t_{\revise{k-1}}-1)} - \frac{1}{c_k^2}\big] \norms{z^k - x^{\star}}^2 \vspace{1ex}\\
&& + {~} st_{k-1}\iprods{w^k, x^k - x^{\star}}.
\end{array}
\end{equation*}
Let us choose $c_k^2 := \revise{\frac{\eta}{4}}$.
Then, we have $a_k - c_k^2t_{k-1}^2 = \big[\revise{\frac{3(3s-2)}{16(s-1)}}\eta - 2\rho_n\big]t_{k-1}^2 + 2s\rho_nt_{k-1} -  \gamma(t_{k-1}-1)$.
Using this \revise{expression, $\frac{t_{\revise{k-1}}}{t_{\revise{k-1}} - 1} \geq 1$,} and noting that $\iprods{w^k, x^k - x^{\star}} \geq - \rho_{*}\norms{w^k}^2$ from \eqref{eq:weak_Minty} of Assumption~\ref{as:A2},  the last expression leads to
\begin{equation*} 
\arraycolsep=0.2em
\begin{array}{lcl}
\Pc_k & \geq & \frac{A_k}{2}\norms{w^k}^2 + \revise{\frac{s^2}{2}\big(\frac{s-1}{\gamma} - \frac{4}{\eta}\big)}\norms{z^k - x^{\star}}^2,
\end{array}
\end{equation*}
which proves \eqref{eq:Ek_lowerbound}, where $A_k := a_k - c_k^2t_{k-1}^2 - 2s\rho_{*}t_{k-1} = \big[\revise{\frac{3(3s-2)}{16(s-1)}}\eta - 2\rho_n\big]t_{k-1}^2 + 2s(\rho_n - \rho_{*})t_{k-1} - \gamma(t_{k-1}-1)$.
\Eproof

\beforesec
\section{The Missing Proofs of Lemmas in Section~\ref{sec:convergence_analysis}.}\label{apdx:sec:convergence_analysis}
\aftersec
This appendix presents the missing proofs of the technical lemmas in Section~\ref{sec:convergence_analysis}.

\beforesubsec
\subsection{\textbf{The proof of Lemma~\ref{le:SFOG_supporting_lemma0}.}}\label{apdx:subsec:le:SFOG_supporting_lemma0}
\aftersubsec
First, let us choose \revise{$\gamma := \frac{s-2}{32(s-1)}\eta$} in \eqref{eq:SFOG_para_cond1} of Lemma \ref{le:VFOG_key_estimate1}, \revise{which satisfies the condition $0 < \gamma \leq \frac{9}{16}\eta$ in \eqref{eq:SFOG_para_cond2}}.
This leads to the update of $\gamma_k$ as in \eqref{eq:SFOG_para_update2}.
In addition, \mred{$\varphi$ in \eqref{eq:SFOG_phik_and_varphi} of Lemma \ref{le:VFOG_key_estimate2} becomes} \revise{$\varphi := \frac{9(85s-134)}{64(s-1)} + \frac{9(s-2)}{64}$}.

Next, under the conditions of Theorem~\ref{th:VFOG1_convergence} and \mred{the choice 
of $\gamma$ as above}, all the conditions of Lemmas~\ref{le:VFOG_key_estimate2} and \ref{le:Ek_lowerbound} are fulfilled. 

Now, since $\rho_c = 0$ and $\alpha = 1$, taking the conditional expectation $\Expsn{k}{\cdot}$ on both sides of \eqref{eq:SFOG_key_est2} in Lemma~\ref{le:VFOG_key_estimate2}, we obtain from this inequality that
\begin{equation*} 
\arraycolsep=0.2em
\begin{array}{lcl}
\Pc_k - \Expsn{k}{ \Pc_{k+1} }  & \geq &  \frac{  \omega \eta L^2 t_k^2}{2} \Expsn{k}{ \norms{x^{k+1} - y^k }^2 } +  \frac{\eta t_k(t_k-s)}{2}\big(1  - \frac{M^2\eta^2 t_k}{t_k - s} \big) \Expsn{k}{ \norms{ \hat{w}^{k+1} - w^k}^2 } \vspace{1ex}\\
&& + {~}   \frac{\eta st_k}{2} \Expsn{k}{ \norms{ \hat{w}^{k+1} }^2 }  + \frac{\phi_k}{2} \norms{ w^k }^2 +  s(s-1)\iprods{ w^k, x^k - x^{\star}}   \vspace{1ex}\\
&& - {~} \revise{\frac{\varphi \eta  t_{k-1}^2}{8}}\norms{ w^k - \hat{w}^k }^2 - \frac{\revise{(25s-34)}\eta t_k^2}{\revise{2}(s-2)} \Expsn{k}{ \norms{e^k}^2 } - \revise{\varphi \eta  t_{k-1}^2} \norms{e^{k-1}}^2.
\end{array}
\end{equation*}
\mred{Rearranging the last two terms and using $\Expsn{k}{\norms{e^k}^2 } \leq \Expsn{k}{\Delta_k}$ \mred{from the first line of \eqref{eq:error_cond3}}}, we have 
\begin{equation}\label{eq:SFOG_supporting_lemma0_proof1}
\arraycolsep=0.2em\hspace{-3ex}
\begin{array}{lcl}
\Pc_k - \Expsn{k}{ \Pc_{k+1} }  & \geq &  \frac{  \omega \eta L^2 t_k^2}{2} \Expsn{k}{ \norms{x^{k+1} - y^k }^2 } +  \frac{\eta t_k(t_k-s)}{2}\big(1  - \frac{M^2\eta^2 t_k}{t_k - s} \big) \Expsn{k}{ \norms{ \hat{w}^{k+1} - w^k}^2 } \vspace{1ex}\\
&& + {~}   \frac{\eta st_k}{2} \Expsn{k}{ \norms{ \hat{w}^{k+1} }^2 }  + \frac{\phi_k}{2} \norms{ w^k }^2 +  s(s-1)\iprods{ w^k, x^k - x^{\star}}   \vspace{1ex}\\
&& - {~} \revise{\frac{\varphi \eta  t_{k-1}^2}{8}} \norms{ w^k - \hat{w}^k }^2 - \revise{\left[\frac{25s-34}{2(s-2)} + \varphi + \frac{s-2}{32(s-1)}\right]\eta t_k^2} \Expsn{k}{\Delta_k} \vspace{1ex}\\
&& + {~} \revise{\left[\varphi + \frac{s-2}{32(s-1)}\right] \eta t_k^2} \Expsn{k}{ \norms{e^k}^2 } - \revise{\left[\varphi + \frac{s-2}{32(s-1)}\right] \eta t_{k-1}^2} \norms{e^{k-1}}^2 +  \revise{\frac{(s-2)\eta t_{k-1}^2}{32(s-1)}} \norms{e^{k-1}}^2.
\end{array}\hspace{-3ex}
\end{equation}
Now, multiplying the second line of \eqref{eq:error_cond3} by $\Lambda_0\eta t_k^2$ for $\Lambda_0$ given in \eqref{eq:fixed_constants}, then rearranging the results, we get
\begin{equation*} 
	\arraycolsep=0.2em
	\begin{array}{llcl}
		& 0 & \geq &  \Lambda_0 \eta t_k^2 \Expsn{k}{ \Delta_k } - \Lambda_0\eta (1-\kappa)t_k^2 \Delta_{k-1} - \Lambda_0\Theta L^2 \eta t_k^2 \norms{x^k - y^{k-1} }^2 - \Lambda_0\eta \delta_k.
	\end{array}
\end{equation*}
Since $t_k = k + s + 1$, we can show that $t_{k-1}^2 \leq t_k^2 \leq \frac{(s+1)^2t_{k-1}^2}{s^2}$ and $\frac{(s+1)^2t_{k-1}^2}{s^2} - t_k^2 \geq \frac{(2s+1)t_{k-1}(t_{k-1}-s)}{s^2}$.
Using these relations, we obtain from the last inequality that
\begin{equation*} 
\arraycolsep=0.2em
\begin{array}{lcl}
0 & \geq & \frac{\Lambda_0(1-\kappa) (s+1)^2 }{s^2} \eta t_k^2 \Expsn{k}{ \Delta_k } -  \frac{\Lambda_0(1-\kappa)(s+1)^2}{s^2}\eta t_{k-1}^2 \Delta_{k-1}
-  \frac{\Lambda_0\Theta L^2(s+1)^2}{s^2} \eta t_{k-1}^2 \norms{x^k - y^{k-1} }^2 \vspace{1ex}\\
&& {~} - \Lambda_0\eta \delta_k + \frac{\mred{\Lambda_0 [s^2 - (1-\kappa)(s+1)^2 ]}}{s^2}\eta t_k^2 \Expsn{k}{ \Delta_k } + \frac{\mred{(2s+1)\Lambda_0(1-\kappa)} }{s^2}\eta t_{k-1}(t_{k-1}-s)\Delta_{k-1}.
\end{array}
\end{equation*}
From \eqref{eq:G_Lipschitz} of Assumption~\ref{as:A1} with $\alpha = 1$, we also have 
\begin{equation*} 
\arraycolsep=0.2em
\begin{array}{lcl}
\norms{w^k - \hat{w}^k}^2 = \norms{Gx^k - Gy^{k-1}}^2 \leq L^2\norms{x^k - y^{k-1}}^2.
\end{array}
\end{equation*}
Combining the last two expressions, we can derive from \eqref{eq:SFOG_supporting_lemma0_proof1} that
\begin{equation}\label{eq:SFOG_lmB0_proof1} 
\hspace{-1ex}
\arraycolsep=0.2em
\begin{array}{lcl}
\Pc_k - \Expsn{k}{ \Pc_{k+1} }  & \geq &  \frac{  \omega \eta L^2 t_k^2}{2} \Expsn{k}{ \norms{x^{k+1} - y^k }^2 } +  \frac{\eta t_k(t_k-s)}{2}\big(1  - \frac{M^2\eta^2 t_k}{t_k - s} \big) \Expsn{k}{ \norms{ \hat{w}^{k+1} - w^k}^2 } \vspace{1ex}\\
&& + {~}   \frac{\eta st_k}{2} \Expsn{k}{ \norms{ \hat{w}^{k+1} }^2 }  + \frac{\phi_k}{2} \norms{ w^k }^2 +  s(s-1)\iprods{ w^k, x^k - x^{\star}}   \vspace{1ex}\\
&& + {~}  \Big[ \frac{\Lambda_0 [s^2 - (1-\kappa)(s+1)^2 ]}{s^2} -  \revise{\frac{25s-34}{2(s-2)} - \varphi - \frac{s-2}{32(s-1)}}\Big] \eta t_k^2 \Expsn{k}{\Delta_k}  + \revise{\frac{(s-2)\eta t_{k-1}^2}{32(s-1)}} \norms{e^{k-1}}^2  \vspace{1ex}\\
&& - {~} \big( \revise{\frac{\varphi}{8}} + \frac{\Lambda_0\Theta(s+1)^2 }{s^2} \big) \eta L^2 t_{k-1}^2 \norms{ x^k - y^{k-1} }^2 \vspace{1ex}\\
&& + {~}  \frac{\Lambda_0 \eta (1-\kappa)  (s+1)^2t_k^2}{s^2} \Expsn{k}{ \Delta_k } -  \frac{\Lambda_0\eta (1-\kappa)(s+1)^2t_{k-1}^2}{s^2} \Delta_{k-1} \vspace{1ex}\\
&& + {~} \revise{\left[\varphi + \frac{s-2}{32(s-1)}\right] \eta t_k^2} \Expsn{k}{ \norms{e^k}^2 } - \revise{\left[\varphi + \frac{s-2}{32(s-1)}\right] \eta t_{k-1}^2} \norms{e^{k-1}}^2 \vspace{1ex}\\
&& + {~} \frac{\mred{(2s+1)\Lambda_0(1 - \kappa)} }{s^2}\eta t_{k-1}(t_{k-1}-s)\Delta_{k-1}  - \Lambda_0\eta \delta_k.
\end{array}
\hspace{-3ex}
\end{equation}
Next, since $\frac{2s+1}{(s+1)^2} < \kappa \leq 1$ from Theorem~\ref{th:VFOG1_convergence}, we have $\Lambda_0 := \revise{\frac{s^2}{s^2 - (1-\kappa)(s+1)^2}\left[\frac{25s-34}{2(s-2)} + \varphi + \frac{s-2}{32(s-1)}\right]} > 0$.
In addition, let us denote $\Lambda_1 :=  \revise{\frac{\varphi}{4}} + \frac{2\Lambda_0 \Theta (s+1)^2}{s^2}$ and choose $\omega := \Lambda_1 + \revise{\frac{s-2}{16(s-1)}}$ as in \eqref{eq:fixed_constants}.
Then, using this $\Lambda_0$ and $\Lambda_1$, the relation $\omega = \Lambda_1 + \revise{\frac{s-2}{16(s-1)}}$, and  $\iprods{w^k, x^k - x^{\star}} \geq -\rho_{*}\norms{w^k}^2$ from Assumption~\ref{as:A2}, we can show from \eqref{eq:SFOG_lmB0_proof1}  that
\myeq{eq:SFOG_lmB0_proof3}{ 
\hspace{-1ex}
\arraycolsep=0.2em
\begin{array}{lcl}
\Pc_k - \Expsn{k}{ \Pc_{k+1} }  & \geq &  \frac{  \omega \eta L^2 t_k^2}{2} \Expsn{k}{ \norms{x^{k+1} - y^k }^2 } +  \frac{\eta t_k(t_k-s)}{2}\big(1  - \frac{M^2\eta^2 t_k}{t_k - s} \big) \Expsn{k}{ \norms{ \hat{w}^{k+1} - w^k}^2 } \vspace{1ex}\\
&& + {~}   \frac{\eta st_k}{2} \Expsn{k}{ \norms{ \hat{w}^{k+1} }^2 }  + \big[ \frac{\phi_k}{2} - s(s-1)\rho_{*}\big] \norms{ w^k }^2   \vspace{1ex}\\
&& - {~}  \frac{\omega \eta L^2 t_{k-1}^2}{2} \norms{ x^k - y^{k-1} }^2 + \frac{\revise{(s-2)}\eta L^2 t_{k-1}^2}{\revise{32(s-1)}} \norms{ x^k - y^{k-1} }^2 \vspace{1ex}\\
&& + {~}  \frac{\Lambda_0(1-\kappa) \eta (s+1)^2t_k^2}{s^2} \Expsn{k}{ \Delta_k } -  \frac{\Lambda_0\eta (1-\kappa)(s+1)^2t_{k-1}^2}{s^2} \Delta_{k-1}  \vspace{1ex}\\
&& + {~} \revise{\left[\varphi + \frac{s-2}{32(s-1)}\right] \eta t_k^2} \Expsn{k}{ \norms{e^k}^2 } - \revise{\left[\varphi + \frac{s-2}{32(s-1)}\right] \eta t_{k-1}^2} \norms{e^{k-1}}^2 \vspace{1ex}\\
&& + {~}  \frac{\mred{(2s+1)\Lambda_0(1 - \kappa)}}{s^2}\eta t_{k-1}(t_{k-1}-s)\Delta_{k-1}  +   \revise{\frac{(s-2)\eta t_{k-1}^2}{32(s-1)}} \norms{e^{k-1}}^2 - \Lambda_0\eta \delta_k.
\end{array}
\hspace{-3ex}
}
Now, since \revise{$8(s-1)\rho_n \leq (s-2)\eta$} and $0 \leq \rho_{*} \leq \rho_n$, for $k \geq 0$, \revise{$s > 2$, and $\gamma := \frac{s-2}{32(s-1)}\eta$},  we have
\myeqn{
\arraycolsep=0.2em
\begin{array}{lcl}
\phi_k - 2 \rho_{*}s(s-1) & = &  \big[ \revise{\frac{5(s-2)}{8}}\eta - 4(s-1)\rho_n - 3\gamma \big](t_k - 1) +  \frac{(s-1)\revise{(3s+10)}\eta}{\revise{16}} \vspace{1ex}\\
&& + {~} 2 (2s-1)(s-1) \rho_n + \gamma - 2\rho_{*} s(s-1) \vspace{1ex}\\

& \geq & \revise{\frac{(s-2)(4s-7)\eta}{32(s-1)}}(k + 1)  +  \revise{\frac{(s-2)(4s-7)\eta}{32}} + \frac{(s-1)\revise{(3s+10)}\eta}{\revise{16}} \vspace{1ex}\\
&& + {~} 2 (2s-1)(s-1) \rho_n + \revise{\frac{s-2}{32(s-1)}\eta}  - 2\rho_n s(s-1)  \vspace{1ex}\\

&=&  \revise{\frac{(s-2)(4s-7)\eta}{32(s-1)}(k + 1) + \frac{10s^3 - 11s^2 - 4s + 4}{32(s-1)}\eta + 2(s-1)^2\rho_n} \vspace{1ex}\\

& \geq & \revise{\frac{(s-2)(4s-7)\eta}{32(s-1)}(k + 1)} > 0.
\end{array}
}
Furthermore, since $t_k \geq s + 1$, to guarantee $\frac{M^2\eta^2 t_k}{t_k - s} \leq 1$, we impose $(s+1)M^2\eta^2 \leq 1$.
However, since $\omega$ is given in \eqref{eq:fixed_constants}, the condition $(s+1)M^2\eta^2 \leq 1$ is guaranteed if we choose $\eta \leq \frac{1}{L\sqrt{2(1+\omega)(s+1)}}$.
In addition, since $\eta \geq \revise{\frac{8(s-1)\rho_n}{s-2}}$, we finally get 
\myeqn{
\arraycolsep=0.2em
\begin{array}{lcl}
\revise{\frac{8(s-1)\rho_n}{s-2}} \leq \eta \leq \frac{1}{L\sqrt{2(1+\omega)(s+1)}} = \frac{\lambda}{L}, \quad \textrm{where} \quad \lambda := \frac{1}{\sqrt{2(1+\omega)(s+1)}}.
\end{array}
}
This inequality holds if we choose $\eta$ as in \eqref{eq:SFOG_para_update2}.

Finally, if we recall $\Lc_k$ from \eqref{eq:SFOG_th1_Lyfunc}, then using the last relations and $\epsilon_1 := 1 - (s+1)M^2\eta^2 \geq 0$, \eqref{eq:SFOG_lmB0_proof3} leads to \eqref{eq:SFOG_lmB0_bound}.
\Eproof

\beforesubsec
\subsection{\textbf{The proof of Lemma~\ref{le:SFOG_supporting_lemma1}.}}\label{apdx:subsec:le:SFOG_supporting_lemma1}
\aftersubsec
First, from the first and second lines of \eqref{eq:VFOG}, we can easily derive that
\myeq{eq:SFOG_lmB1_proof1}{
	\arraycolsep=0.2em
	\begin{array}{lcl}
		\revise{u^k} & := &  s(z^k - x^k) = t_k(y^k - x^k) +   t_k(\eta - \beta_k)(\hat{w}^k  + e^{k-1}).
	\end{array}
}
Using this expression, the third line of \eqref{eq:VFOG}, and $x^{k+1} - y^k = - \eta(\hat{w}^{k+1} - \hat{w}^k  + e^k - e^{k-1})$, we can derive
\myeqn{
	\arraycolsep=0.15em
	\begin{array}{lcl}
		\revise{u^{k+1}} - \frac{(t_k-s)}{t_k}\revise{u^k} &= & \revise{u^{k+1}} - \revise{u^k} + \frac{s}{t_k}\revise{u^k} \vspace{1ex}\\
		&= & s(z^{k+1} - z^k) - s(x^{k+1} - x^k) + s(y^k - x^k) + s (\eta-\beta_k) ( \hat{w}^k + e^{k-1}) \vspace{1ex}\\
		& = &  s(z^{k+1} - z^k) - s(x^{k+1} - y^k) + s(\eta-\beta_k)( \hat{w}^k  + e^{k-1}) \vspace{1ex}\\
		& = & -\gamma_k( \hat{w}^k + e^{k-1})  + s \eta(\hat{w}^{k+1} - \hat{w}^k + e^k - e^{k-1}) +  s (\eta-\beta_k)( \hat{w}^k + e^{k-1}) \vspace{1ex}\\
		& = &  s\eta( \hat{w}^{k+1} - \hat{w}^k ) + [s(\eta-\beta_k)-\gamma_k]\hat{w}^k + s\eta e^k - (s\beta_k + \gamma_k)e^{k-1}.
	\end{array}
}
Denote $d^k := \frac{t_k}{s}\big[ s\eta( \hat{w}^{k+1} - \hat{w}^k ) + [s(\eta-\beta_k)-\gamma_k] \hat{w}^k + s\eta e^k - (s\beta_k + \gamma_k)e^{k-1} \big]$.
Then, the last estimate is rewritten as $\revise{u^{k+1}} = \big(1 - \frac{s}{t_k}\big)\revise{u^k} + \frac{s}{t_k}d^k$.
By convexity of $\norms{\cdot}^2$, $\frac{s}{t_k} \in (0, 1]$, and Young's inequality, we have 
\myeqn{
	\arraycolsep=0.2em
	\begin{array}{lcl}
		\norms{\revise{u^{k+1}}}^2 & \leq & \frac{t_k-s}{t_k}\norms{\revise{u^k}}^2 +  \frac{s}{t_k} \norms{d^k}^2 \vspace{1ex}\\
		& \leq &  \norms{\revise{u^k}}^2 - \frac{s}{t_k}\norms{\revise{u^k}}^2 + 4s\eta^2 t_k \norms{ \hat{w}^{k+1} - \hat{w}^k }^2 + \frac{4[s(\eta-\beta_k)-\gamma_k]^2t_k}{s}\norms{ \hat{w}^k }^2 \vspace{1ex}\\
		&& +  {~} 4 s\eta^2t_k\norms{e^k}^2 + \frac{4(s\beta_k + \gamma_k)^2t_k}{s} \norms{ e^{k-1}}^2.
	\end{array}
}
Note that from \revise{Lemmas \ref{le:VFOG_key_estimate2} and \ref{le:SFOG_supporting_lemma0}}, we have $\beta_k \leq \frac{\revise{7(s-2)}\eta(t_k-s)}{\revise{16(s-1)}t_k} \leq \eta$ and $\gamma_k = \revise{\frac{(s-2)\eta(t_k-1)}{ 32(s-1) t_k}} \leq \eta$.
Using these inequalities into the last estimate and taking the full expectation on both sides of the resulting inequality, we get
\myeq{eq:SFOG_lmB2_proof2}{
	\arraycolsep=0.2em
	\begin{array}{lcl}
		\Expn{ \norms{\revise{u^{k+1}}}^2} & \leq & \revise{ \Expn{ \norms{u^k}^2 }} - \frac{s}{t_k} \Expn{ \norms{\revise{u^k}}^2 } + 4s\eta^2 t_k \Expn{ \norms{ \hat{w}^{k+1} - \hat{w}^k  }^2 } \vspace{1ex}\\
		&& + {~} \frac{4 (2s+1)^2 \eta^2 t_k}{s} \Expn{ \norms{ \hat{w}^k }^2 }  + 4 s\eta^2t_k \Expn{ \norms{e^k}^2 } + \frac{4(s+1)^2\eta^2t_k}{s} \Expn{ \norms{ e^{k-1}}^2 }.
	\end{array}
}
Applying Lemma~\ref{le:A1_descent} to the last inequality by using the summability bounds in \eqref{eq:SFOG_proof1_summability_bounds}, we get
\myeqn{
	\lim_{k\to\infty} \Expn{ \norms{\revise{u^k}}^2 }  = \lim_{k\to\infty} s^2 \Expn{ \norms{x^k - z^k}^2 } = 0, \ \ \textrm{and} \ \ \sum_{k=0}^{\infty}\tfrac{s}{t_k} \Expn{ \norms{z^k - x^k}^2 } < +\infty.
}
These relations imply \eqref{eq:SFOG_proof1_limits}.

Next, since $t_k(x^{k+1} - y^k) = -\eta t_k\big(\hat{w}^{k+1} - \hat{w}^k + e^k - e^{k-1}\big)$ by \eqref{eq:VFOG},  by Young's inequality, we can show that
\myeqn{
	\arraycolsep=0.2em
	\begin{array}{lcl}
		t_k^2\norms{x^{k+1} - y^k}^2 &\leq & 3\eta^2t_k^2\norms{ \hat{w}^{k+1} - \hat{w}^k }^2 + 3\eta^2t_k^2\norms{e^k}^2 + 3\eta^2t_k^2\norms{e^{k-1}}^2.
	\end{array}
}
Combining this inequality and \eqref{eq:SFOG_proof1_summability_bounds}, we can prove that
\begin{equation*} 
	\arraycolsep=0.2em
	\begin{array}{lcl}
		\sum_{k=0}^{\infty} t_k^2 \Expn{ \norms{x^{k+1} - y^k }^2 } & < &  +\infty.
	\end{array}
\end{equation*}
This is exactly the first summability bound of \eqref{eq:SFOG_proof1_summability_bounds2}.

Finally, using again \eqref{eq:SFOG_lmB1_proof1} and then applying Young's inequality, we can show that 
\begin{equation*} 
	\arraycolsep=0.2em
	\begin{array}{lcl}
		t_k\norms{y^k - x^k}^2 & \leq & \frac{3s^2}{t_k}\norms{z^k - x^k}^2 +  6\eta^2 t_k \norms{ \hat{w}^k }^2 +  6\eta^2 t_k \norms{e^{k-1}}^2, \vspace{1ex}\\
		t_k\norms{x^{k+1} - x^k}^2 & \leq & 2t_k\norms{x^{k+1} - y^k}^2 +  2t_k\norms{y^k - x^k}^2.
	\end{array}
\end{equation*}
Utilizing the first line of \eqref{eq:SFOG_proof1_summability_bounds2}, \eqref{eq:SFOG_proof1_summability_bounds},  and the second summability bound of \eqref{eq:SFOG_proof1_limits},  we can also derive from these inequalities that
\begin{equation*} 
	\arraycolsep=0.0em
	\begin{array}{lcl}
		\sum_{k=0}^{\infty} t_k \Expn{ \norms{y^k - x^k}^2 } < +\infty \quad \text{and} \quad \sum_{k=0}^{\infty} t_k \Expn{ \norms{x^{k+1} - x^k}^2} < +\infty.
	\end{array}
\end{equation*}
These prove the last two lines of \eqref{eq:SFOG_proof1_summability_bounds2}, respectively.
\Eproof

\beforesubsec
\subsection{\textbf{The proof of Lemma~\ref{le:SFOG_supporting_lemma2}.}}\label{apdx:subsec:le:SFOG_supporting_lemma2}
\aftersubsec
Since $\Expn{ \Pc_k } \leq \bar{\Pc}_0$, from \eqref{eq:Ek_lowerbound} of Lemma~\ref{le:Ek_lowerbound}, we have
\begin{equation*} 
\arraycolsep=0.2em
\begin{array}{lcl}
\frac{A_k}{2}\Expn{ \norms{ w^k }^2 }  + \revise{\frac{s^2[(s-1)\eta - 4\gamma]}{2\gamma\eta}} \Expn{ \norms{z^k - x^{\star} }^2 } \leq \Expn{\Pc_k} \leq \bar{\Pc}_0,
\end{array}
\end{equation*}
where \revise{we have used $\frac{s^2[(s-1)\eta - 4\gamma]}{2\gamma\eta} = \frac{2s^2(8s^2 - 17s + 10)}{(s-2)\eta} > 0$ for all $s>2$}.
This inequality implies that $\Expn{ \norms{z^k - x^{\star}}^2 } \leq M_0 := \frac{\revise{2\gamma\eta}\bar{\Pc}_0}{\revise{s^2[(s-1)\eta - 4\gamma]}}$.

\noindent Now, from Young's inequality, we have
\begin{equation*} 
\arraycolsep=0.2em
\begin{array}{lcl}
	-\frac{2\gamma_k}{s}\iprods{\hat{e}^k, z^k - x^{\star}} &\leq& \frac{\gamma_k t_k^2}{s} \norms{\hat{e}^k}^2 + \frac{\gamma_k}{s t_k^2} \norms{z^k - x^{\star}}^2, \vspace{1ex}\\
	-\frac{2\gamma_k}{s}\iprods{w^k, z^k - x^k} &\leq& \frac{\gamma_k t_k}{s} \norms{w^k}^2 + \frac{\gamma_k}{s t_k} \norms{z^k - x^k}^2 \leq \frac{\eta t_k}{s} \norms{w^k}^2 + \frac{\eta}{s t_k} \norms{z^k - x^k}^2.
\end{array}
\end{equation*}
Using Young's inequality in {\tiny\textcircled{1}}, the condition \eqref{eq:weak_Minty} from Assumption \ref{as:A2} in {\tiny\textcircled{2}}, and the last two relations, we have
\begin{equation}\label{eq:SFOG_summability_proof5a} 
\arraycolsep=0.2em
\begin{array}{lcl}
\norms{z^{k+1} - x^{\star}}^2 & = & \norms{z^k - x^{\star}}^2 + 2\iprods{z^{k+1} - z^k, z^k - x^{\star}} + \norms{z^{k+1} - z^k}^2 \vspace{1ex}\\
& = & \norms{z^k - x^{\star}}^2 - \frac{2\gamma_k}{s}\iprods{ w^k + \hat{e}^k, z^k - x^{\star}} + \frac{\gamma_k^2}{s^2} \norms{w^k + \hat{e}^k }^2 \vspace{1ex}\\
& \stackrel{\tiny\textcircled{1}}{\leq} & \norms{z^k - x^{\star}}^2 - \frac{2\gamma_k}{s}\iprods{ w^k , z^k - x^{\star}} - \frac{2\gamma_k}{s}\iprods{ \hat{e}^k, z^k - x^{\star}} \vspace{1ex}\\
&& + {~} \frac{2\gamma_k^2}{s^2} \norms{ w^k }^2 +  \frac{2\gamma_k^2}{s^2} \norms{ \hat{e}^k }^2 \vspace{1ex}\\
&  \leq & \norms{z^k - x^{\star}}^2 - \frac{2\gamma_k}{s}\iprods{w^k, x^k - x^{\star}} - \frac{2\gamma_k}{s}\iprods{w^k, z^k - x^k }  \vspace{1ex}\\
&& + {~}  \frac{\gamma_kt_k^2}{s}\norms{ \hat{e}^k }^2 + \frac{\gamma_k}{s t_k^2} \norms{z^k - x^{\star}}^2 + \frac{2\gamma_k^2}{s^2} \norms{ w^k }^2 +  \frac{2\gamma_k^2}{s^2} \norms{ \hat{e}^k }^2 \vspace{1ex}\\
& \stackrel{\tiny\textcircled{2}}{\leq} & \big(1 + \frac{\eta}{st_k^2}\big) \norms{z^k - x^{\star}}^2 + \big(\frac{\eta t_k}{s} + \frac{\revise{2}\eta \rho_{*} }{s} +  \frac{2\eta^2}{s^2} \big)\norms{ w^k }^2  \vspace{1ex}\\
&& + {~}  \frac{\eta }{st_k}\norms{z^k - x^k }^2 + \big( \frac{\eta t_k^2}{s} + \frac{2\eta^2}{s^2} \big)\norms{ \hat{e}^k }^2.
\end{array}
\end{equation}
Taking the full expectation on both sides of this inequality and noting that the last three terms are summable, applying Lemma~\ref{le:A1_descent}, we conclude that $\lim_{k\to\infty}\Expn{\norms{z^k - x^{\star} }^2}$ exists.
However, since  $\lim_k \Expn{ \norms{x^k - z^k}^2 } = 0$ by \eqref{eq:SFOG_proof1_limits}, using Minkowski's inequality $\big\vert \sqrt{\Expn{ \norms{x^k - x^{\star}}^2} } - \sqrt{\Expn{\norms{z^k - x^{\star}}^2}} \big\vert \leq \sqrt{\Expn{\norms{x^k - z^k}^2} }$, we conclude that $\lim_{k\to\infty}\Expn{ \norms{x^k - x^{\star}}^2 } $ also exists.

Since $\sets{ t_{k-1}^2\Expn{ \norms{w^k}^2 } }$ is bounded, there exists $\bar{M} > 0$ such that $t_{k-1}^2\Expn{ \norms{w^k}^2 } \leq \bar{M}^2$.
Using this inequality and $\lim_{k\to\infty} \Expn{\norms{x^k - z^k}^2}  = 0$ from \eqref{eq:SFOG_proof1_limits}, we can show that when $k \to \infty$:
\begin{equation*} 
\arraycolsep=0.2em
\begin{array}{lcl}
t_{k-1}^2 \vert \Expn{  \iprods{w^k, x^k - z^k} } \vert^2 &\leq & t_{k-1}^2\Expn{ \norms{w^k}^2 } \Expn{ \norms{x^k - z^k}^2 } \leq \bar{M}^2\Expn{ \norms{x^k - z^k}^2} \to 0.
\end{array}
\end{equation*}
This expression implies that $\lim_{k\to\infty} t_{k-1}\Expn{  \iprods{ w^k, x^k - z^k }  } = 0$.
From \eqref{eq:SFOG_Ek_func}, we have
\begin{equation*} 
\arraycolsep=0.2em
\begin{array}{lcl}
\Pc_k = \frac{a_k}{2}\norms{w^k}^2 + st_{k-1}\iprods{ w^k, x^k - z^k} + \frac{s^2(s-1)}{2\gamma_{\revise{k-1}}}\norms{z^k - x^{\star}}^2.
\end{array}
\end{equation*}
Since $\lim_{k\to\infty}\Expn{\Pc_k}$ and $\lim_{k\to\infty} \Expn{ \norms{z^k - x^{\star}}^2}$ exist, $\gamma_{k-1} \to \gamma > 0$ as $k \to \infty$, and $\lim_{k\to\infty} t_{k-1}\Expn{  \iprods{ w^k, x^k - z^k }  } = 0$, we conclude that $\lim_{k\to\infty}a_k \Expn{ \norms{ w^k }^2 }$ exists.
Since $\lim_{k\to\infty}\frac{t_{k-1}^2}{a_k} = \frac{1}{\psi}$, we conclude that $\lim_{k\to\infty}t_{k-1}^2 \Expn{ \norms{ w^k}^2 }$ exists.
However, since $\sum_{k=0}^{\infty}t_{k-1}\Expn{ \norms{ w^k}^2 } < +\infty$ by our assumption, we can easily argue that $\lim_{k\to\infty}t_{k-1}^2 \Expn{ \norms{ w^k }^2 } = 0$, which proves the first limit of \eqref{eq:SFOG_lmB2_limits}.

Note that $\norms{w^k - \hat{w}^k}^2 = \norms{Gx^k - Gy^{k-1}}^2 \leq L^2\norms{x^k - y^{k-1}}^2$ by Assumption~\ref{as:A1}.
Utilizing this relation and the first line of \eqref{eq:SFOG_proof1_summability_bounds2}, we can conclude that $\lim_{k\to\infty} t_{k-1}^2\Expn{ \norms{w^k - \hat{w}^k}^2  } = 0$.

Finally, by Young's inequality, we have
\begin{equation*} 
\arraycolsep=0.2em
\begin{array}{lcl}
t_{k-1}^2\norms{ \hat{w}^k }^2 & \leq & 2t_{k-1}^2\norms{w^k}^2 + 2t_{k-1}^2\norms{w^k - \hat{w}^k }^2.
\end{array}
\end{equation*}
Since $\lim_{k\to\infty}t_{k-1}^2\Expn{ \norms{ w^k - \hat{w}^k }^2 } = 0$ and $\lim_{k\to\infty} t_{k-1}^2 \Expn{ \norms{ w^k }^2 } = 0$, the last inequality implies that $\lim_{k\to\infty} t_{k-1}^2 \Expn{ \norms{ \hat{w}^k }^2 } = 0$, which proves the second limit of \eqref{eq:SFOG_lmB2_limits}.
\Eproof

\beforesubsec
\subsection{\textbf{The proof of Lemma~\ref{le:VrAEG4GE_descent_property}.}}\label{apdx:le:VrAEG4GE_descent_property}
\aftersubsec
First, if we choose $\gamma := \revise{\frac{(s-2)\eta}{32(s-1)}}$, then $\varphi$ in \eqref{eq:SFOG_phik_and_varphi} becomes \revise{$\varphi = \varphi_s := \frac{9(85s - 134)}{64(s-1)} + \frac{9(s-2)}{64}$}.
Second, applying \eqref{eq:G_Lipschitz} from Assumption~\ref{as:A1} with $x = x^{k+1}$ and $y = y^k$, we can show that
\begin{equation*}
\arraycolsep=0.2em
\begin{array}{lcl}
\tau \omega L^2\norms{x^{k+1} - y^k}^2 & \geq &  \tau\omega (1-\alpha)\tnorms{w^{k+1} - \hat{w}^{k+1}}^2 + \tau\omega\alpha\norms{w^{k+1} - \hat{w}^{k+1}}^2, \quad \tau \in (0, 1].
\end{array}
\end{equation*}
Using $\tnorms{\cdot} \geq \norms{\cdot}$ and  this inequality, we can show that
\begin{equation*}
\arraycolsep=0.2em
\begin{array}{lcl}
\tilde{\Tc}_{[1]} &:= &   (1-\alpha) \big( \tnorms{w^{k+1} - \hat{w}^{k+1} }^2  - \norms{w^{k+1} - \hat{w}^{k+1} }^2 \big) + \omega L^2 \norms{x^{k+1} - y^k }^2  \vspace{1ex}\\
& \geq & \tau\omega(1-\alpha) \tnorms{w^{k+1} - \hat{w}^{k+1} }^2 + \tau\omega\alpha \norms{w^{k+1} - \hat{w}^{k+1}}^2  + (1-\tau)\omega L^2\norms{x^{k+1} - y^k }^2.
\end{array}
\end{equation*}
Next, substituting this relation into \eqref{eq:SFOG_key_est2} of Lemma~\ref{le:VFOG_key_estimate2}, we can derive that
\begin{equation}\label{eq:SAEG_key_est2_v2}
\hspace{-0.5ex}
\arraycolsep=0.2em
\begin{array}{lcl}
\Pc_k - \Pc_{k+1}  & \geq &  \frac{\tau\omega\alpha\eta t_k^2}{2} \norms{w^{k+1} - \hat{w}^{k+1} }^2 + \frac{\tau\omega(1-\alpha) \eta t_k^2}{2} \tnorms{ w^{k+1} - \hat{w}^{k+1} }^2  \vspace{1ex}\\
&& + {~} \frac{(1-\tau)\omega \eta L^2 t_k^2 }{2}\norms{x^{k+1} - y^k }^2  +   \rho_c t_k(t_k-s)\tnorms{ Gx^{k+1} - Gx^k }^2 \vspace{1ex}\\
&& + {~}  \frac{\eta t_k(t_k-s)}{2}\big(1  - \frac{ M^2\eta^2 t_k }{t_k - s} \big) \norms{ \hat{w}^{k+1} - w^k }^2 +  \frac{\eta st_k}{2} \norms{ \hat{w}^{k+1} }^2   + \frac{\phi_k}{2} \norms{ w^k }^2  \vspace{1ex}\\
&& + {~}   s(s-1)\iprods{ w^k, x^k - x^{\star}} - \frac{\revise{\varphi_s}\eta t_{k-1}^2}{\revise{8}} \norms{ w^k - \hat{w}^k }^2 \vspace{1ex}\\
&& - {~} \frac{\revise{(25s-34)}\eta t_k^2}{\revise{2}(s-2)}\norms{e^k}^2 - \revise{\varphi_s}  \eta t_{k-1}^2 \norms{e^{k-1}}^2.
\end{array}
\hspace{-2ex}
\end{equation}
Since $w^{k+1} - \hat{w}^{k+1} = Gx^{k+1} - Gy^k$ and $w^k - \hat{w}^k = Gx^k - Gy^{k-1}$, by Young's inequality, the relation $\tnorms{\cdot} \geq \norms{\cdot}$, and  \eqref{eq:G_Lipschitz} from Assumption~\ref{as:A1}, we can derive that
\begin{equation}\label{eq:three_estimates_proof1} 
\arraycolsep=0.2em
\begin{array}{lcl}
\tnorms{w^k - \hat{w}^k}^2  \geq \norms{ w^k - \hat{w}^k }^2, \vspace{1ex}\\
\tnorms{ w^{k+1} -  \hat{w}^{k+1} }^2 + \tnorms{Gx^{k+1} - Gx^k}^2 + \tnorms{ w^k - \hat{w}^k  }^2 \geq \frac{1}{3}\tnorms{ Gy^k  - Gy^{k-1} }^2, \vspace{1ex}\\
L^2\norms{x^k - y^{k-1} }^2 \geq \norms{w^k - \hat{w}^k}^2.
\end{array}
\end{equation}
Then, since $t_k^2 \geq t_k(t_k-s)$ and $t_{k-1}^2 \geq t_k(t_k - s)$ for $s \geq 2$, using the last two inequalities and the conditions 
\begin{equation}\label{eq:para_cond10_proof} 
\arraycolsep=0.2em
\begin{array}{lcl}
\tau\omega(1-\alpha)\eta \geq 4\rho_c \quad \textrm{and} \quad (1-\tau)\omega\eta \geq \revise{2\big[\frac{\varphi_s}{8} + \frac{s-2}{32(s-1)}\big]\eta},
\end{array}
\end{equation}
we can prove that 
\begin{equation*} 
\arraycolsep=0.2em
\begin{array}{lcl}
\bar{\Tc}_{[1]} &:= & \frac{\tau\omega(1 - \alpha ) \eta t_k^2}{2}\tnorms{w^{k+1} - \hat{w}^{k+1} }^2  +  \rho_c t_k(t_k-s) \tnorms{ Gx^{k+1} - Gx^k}^2  \vspace{1ex}\\
&& + {~} \frac{(1-\tau)\omega \eta L^2 t_k^2 }{2}\norms{x^{k+1} - y^k }^2 - \frac{\revise{\varphi_s} \eta t_{k-1}^2}{\revise{8}} \norms{ w^k - \hat{w}^k }^2 \vspace{1ex}\\

& \overset{\tiny\eqref{eq:para_cond10_proof}}{\geq} & \frac{2\rho_c t_k(t_k-s)}{2}\tnorms{ w^{k+1} -  \hat{w}^{k+1} }^2 +  \frac{2\rho_c t_k(t_k-s)}{2}\tnorms{ Gx^{k+1} - Gx^k }^2  +  \frac{2\rho_c t_k(t_k-s)}{2}\tnorms{ w^k - \hat{w}^k }^2 \vspace{1ex}\\
&& + {~}  \frac{2\rho_c t_k^2 }{2}\tnorms{w^{k+1} - \hat{w}^{k+1} }^2 - \frac{2 \rho_c  t_{k-1}^2}{2}\tnorms{ w^k - \hat{w}^k }^2 \vspace{1ex}\\
&& + {~} \revise{\big[\frac{\varphi_s}{8} + \frac{s-2}{32(s-1)}\big]} \eta L^2 t_k^2 \norms{x^{k+1} - y^k }^2 - \frac{\revise{\varphi_s} \eta t_{k-1}^2}{\revise{8}} \norms{ w^k - \hat{w}^k }^2 \vspace{1ex}\\

& \overset{\tiny\eqref{eq:three_estimates_proof1}}{\geq} & 
  \frac{\rho_c t_k(t_k-s)}{3} \tnorms{ Gy^k  -  Gy^{k-1} }^2  +  \frac{2\rho_c t_k^2 }{2}\tnorms{w^{k+1} - \hat{w}^{k+1} }^2 - \frac{2 \rho_c  t_{k-1}^2}{2}\tnorms{ w^k - \hat{w}^k }^2   \vspace{1ex}\\
&&  + {~} \frac{\revise{\varphi_s} \eta L^2 t_k^2}{\revise{8}} \norms{x^{k+1} - y^k }^2 - \frac{\revise{\varphi_s} \eta L^2 t_{k-1}^2}{\revise{8}} \norms{ x^k - y^{k-1} }^2 + \frac{ \revise{(s-2)}\eta L^2 t_k^2}{\revise{32(s-1)}}\norms{x^{k+1} - y^k }^2.
\end{array}
\end{equation*}
Substituting this expression into \eqref{eq:SAEG_key_est2_v2} and taking the conditional expectation $\Expsn{k}{\cdot}$ on both sides of the result, and splitting $\Expsn{k}{ \norms{e^k}^2 }$ \mred{and $\norms{e^{k-1}}^2$} into two parts, we can show that
\begin{equation}\label{eq:SAEG_key_est3}
\arraycolsep=0.2em
\begin{array}{lcl}
\Pc_k - \Expsn{k}{ \Pc_{k+1} } & \geq &
\frac{\eta t_k(t_k-s)}{2}\big(1  - \frac{M^2\eta^2 t_k }{t_k - s} \big) \Expsn{k}{ \norms{ \hat{w}^{k+1} - w^k}^2 } +  \frac{\tau\omega\alpha\eta t_k^2}{2} \Expsn{k}{ \norms{ w^{k+1} - \hat{w}^{k+1} }^2 } \vspace{1ex}\\
&& + {~} \frac{\revise{\varphi_s} \eta L^2 t_k^2}{\revise{8}} \Expsn{k}{\norms{x^{k+1} - y^k }^2} - \frac{\revise{\varphi_s} \eta L^2 t_{k-1}^2}{\revise{8}} \norms{ x^k - y^{k-1} }^2 \vspace{1ex}\\
&& + {~} \rho_c t_k^2 \Expsn{k}{\tnorms{w^{k+1} - \hat{w}^{k+1} }^2} - \rho_c  t_{k-1}^2\tnorms{ w^k - \hat{w}^k }^2 \vspace{1ex}\\
&&  + {~} \revise{\varphi_s} \eta t_k^2 \Expsn{k}{ \norms{e^k}^2 } - \revise{\varphi_s} \eta t_{k-1}^2 \norms{e^{k-1}}^2  + \frac{\revise{(s-2)}\eta L^2 \revise{t_k^2}}{\revise{32(s-1)}} \Expsn{k}{\norms{ \revise{x^{k+1} - y^k} }^2} \vspace{1ex}\\
&& + {~}  \frac{\rho_c t_k(t_k-s)}{3} \tnorms{Gy^k - Gy^{k-1} }^2  + \frac{\eta st_k}{2} \Expsn{k}{ \norms{ \hat{w}^{k+1} }^2 }  + \frac{\phi_k}{2} \norms{ w^k }^2 \vspace{1ex}\\
&& + {~} s(s-1)\iprods{ w^k, x^k - x^{\star}} - \revise{\big[\varphi_s + \frac{25s-34}{2(s-2)}\big]} \eta t_k^2 \Expsn{k}{ \norms{e^k}^2 }. 
\end{array}
\hspace{-5ex}
\end{equation}
Now, for a given $\bar{c} \in (0, 1)$, we can derive from \eqref{eq:VR_property} of Definition~\ref{de:VR_Estimators} that
\begin{equation}\label{eq:SAEG_th31_proof3}
\arraycolsep=0.2em
\begin{array}{lcl}
t_k^2\Expsn{k}{\Delta_k} &\leq &  \frac{(1- \kappa_k )t_k^2}{\bar{c}t_{k-1}^2} t_{k-1}^2 \Delta_{k-1}  - \frac{1-\bar{c}}{\bar{c}}t_k^2\Expsn{k}{ \Delta_k }  +  \frac{\Theta_k t_k^2}{\bar{c}} \tnorms{Gy^k - Gy^{k-1} }^2 \mred{~ + ~ \frac{\delta_k}{\bar{c}}}.
\end{array}
\end{equation}
Suppose that there exists $\epsilon_2 \geq 0$ such that  $\frac{(1 - \kappa_k)t_k^2}{\bar{c}t_{k-1}^2} \leq \frac{1-\bar{c}}{\bar{c}} - \epsilon_2$ for $\bar{c} \in (0, 1)$, which holds if $\bar{c} \leq \frac{1}{1+ \epsilon_2}\big[ 1 - \frac{(1 - \kappa_k)t_k^2}{t_{k-1}^2}\big]$.
Moreover, \mred{to handle the term $\tnorms{Gy^k - Gy^{k-1}}^2$,} we also need to impose $\frac{\rho_c t_k(t_k-s)}{3} - \revise{\big[\varphi_s + \frac{25s-34}{2(s-2)}\big]} \frac{\eta \Theta_k t_k^2}{\bar{c}} \geq 0$, leading to $\bar{c} \geq \revise{\big[\varphi_s + \frac{25s-34}{2(s-2)}\big]} \frac{3\eta \Theta_k t_k}{\rho_c(t_k - s)}$.
Combining both conditions on $\bar{c}$, we get
\begin{equation*} 
\arraycolsep=0.2em
\begin{array}{lcl}
\revise{\big[\varphi_s + \frac{25s-34}{2(s-2)}\big]} \frac{3\eta \Theta_k t_k}{\rho_c(t_k - s)} \leq \bar{c} \leq \frac{1}{1 + \epsilon_2}\big[ 1 - \frac{(1- \kappa_k)t_k^2}{t_{k-1}^2}\big].
\end{array}
\end{equation*}
This is exactly the condition~\eqref{eq:SAEG_error_cond3}.
In this case, using $\Expsn{k}{\norms{e^k}^2} \leq \Expsn{k}{ \Delta_k }$ from \eqref{eq:VR_property} of Definition~\ref{de:VR_Estimators},  \eqref{eq:SAEG_th31_proof3} becomes
\begin{equation*} 
\arraycolsep=0.2em
\begin{array}{lcl}
\eta t_k^2 \Expsn{k}{ \norms{e^k}^2 } &\leq & \frac{(1-\bar{c}) \eta t_{k-1}^2}{\bar{c}} \Delta_{k-1} - \frac{(1-\bar{c})\eta t_k^2}{\bar{c}}\Expsn{k}{\Delta_k } + \frac{\eta\delta_k}{\bar{c}} + \frac{\eta \Theta_k t_k^2}{\bar{c}} \tnorms{Gy^k - Gy^{k-1} }^2  -  \eta \epsilon_2 t_{k-1}^2\Delta_{k-1}.
\end{array}
\end{equation*}
\revise{Multiplying both sides of this inequality by $\revise{\big[\varphi_s + \frac{25s-34}{2(s-2)}\big]}$, then substituting the result into \eqref{eq:SAEG_key_est3}}
and using the condition $\frac{\rho_c t_k(t_k-s)}{3} - \revise{\big[\varphi_s + \frac{25s-34}{2(s-2)}\big]} \frac{\eta \Theta_k t_k^2}{\bar{c}} \geq 0$, we obtain
\begin{equation*}
\hspace{-0.5ex}
\arraycolsep=0.2em
\begin{array}{lcl}
\Pc_k - \Expsn{k}{ \Pc_{k+1} }  & \geq &  
\frac{\revise{\varphi_s}  \eta L^2 t_k^2}{\revise{8}} \Expsn{k}{\norms{x^{k+1} - y^k }^2} - \frac{\revise{\varphi_s}  \eta L^2 t_{k-1}^2}{\revise{8}} \norms{ x^k - y^{k-1} }^2  \vspace{1ex}\\
&& + {~} \rho_c t_k^2 \Expsn{k}{\tnorms{w^{k+1} - \hat{w}^{k+1} }^2} - \rho_c t_{k-1}^2\tnorms{ w^k - \hat{w}^k }^2 \vspace{1ex}\\
&&  + {~} \revise{\varphi_s} \eta t_k^2 \Expsn{k}{ \norms{e^k}^2 } - \revise{\varphi_s} \eta t_{k-1}^2 \norms{e^{k-1}}^2 \vspace{1ex}\\
&& + {~} \revise{\big[\varphi_s + \frac{25s-34}{2(s-2)}\big]} \frac{(1-\bar{c})\eta t_k^2}{\bar{c}}\Expsn{k}{\Delta_k } - \revise{\big[\varphi_s + \frac{25s-34}{2(s-2)}\big]}\frac{(1-\bar{c})\eta t_{k-1}^2}{\bar{c}} \Delta_{k-1} \vspace{1ex} \\
&& + {~} \frac{\eta st_k}{2} \Expsn{k}{ \norms{ \hat{w}^{k+1} }^2  } + \frac{\phi_k}{2} \norms{ w^k }^2 +  s(s-1)\iprods{ w^k, x^k - x^{\star}} \vspace{1ex}\\
&& + {~} \frac{\eta t_k(t_k-s)}{2}\big(1  - \frac{M^2\eta^2 t_k}{t_k - s} \big) \Expsn{k}{ \norms{ \hat{w}^{k+1} - w^k }^2 } + \frac{\tau\omega\alpha \eta t_k^2}{2} \Expsn{k}{ \norms{ w^{k+1} - \hat{w}^{k+1} }^2 } \vspace{1ex}\\
&& + {~} \frac{\revise{(s-2)}\eta L^2 \revise{t_k^2}}{\revise{32(s-1)}} \Expsn{k}{\norms{ \revise{x^{k+1} - y^k} }^2}   + \revise{\big[\varphi_s + \frac{25s-34}{2(s-2)}\big]} \epsilon_2 \eta t_{k-1}^2\Delta_{k-1} - \revise{\big[\varphi_s + \frac{25s-34}{2(s-2)}\big]} \frac{\eta\delta_k}{\bar{c}}.
\end{array}
\hspace{-1ex}
\end{equation*}
Recall the Lyapunov function $\hat{\Lc}_k$ from \eqref{eq:SAEG_Lyapunov_func1} as follows:
\begin{equation*} 
\arraycolsep=0.2em
\begin{array}{lcl}
	\hat{\Lc}_k &:= & \frac{a_k}{2}\norms{ w^k }^2 + st_{k-1}\iprods{ w^k, x^k - z^k} + \frac{s^2(s-1)}{2\gamma_{\revise{k-1}}}\norms{z^k - x^{\star}}^2 +  \rho_c  t_{k-1}^2  \tnorms{ w^k - \hat{w}^k }^2 \vspace{1ex}\\
	&& + {~} \frac{\revise{\varphi_s}  \eta L^2 t_{k-1}^2}{\revise{8}} \norms{ x^k - y^{k-1} }^2  +  \revise{\varphi_s} \eta t_{k-1}^2 \norms{e^{k-1}}^2 + \revise{\big[\varphi_s + \frac{25s-34}{2(s-2)}\big]}\frac{(1-\bar{c})\eta t_{k-1}^2}{\bar{c}}  \Delta_{k-1}.
\end{array}
\end{equation*}
Using this function, the last expression becomes
\begin{equation*} 
\arraycolsep=0.2em
\begin{array}{lcl}
\hat{\Lc}_k -  \Expsn{k}{ \hat{\Lc}_{k+1} }  & \geq & \frac{\eta st_k}{2} \Expsn{k}{ \norms{ \hat{w}^{k+1} }^2 }  + \frac{\phi_k}{2} \norms{ w^k }^2 +  s(s-1)\iprods{ w^k, x^k - x^{\star}} \vspace{1ex}\\
&& + {~} \frac{\eta t_k(t_k-s)}{2}\big(1  - \frac{M^2\eta^2 t_k }{t_k - s} \big) \Expsn{k}{ \norms{ \hat{w}^{k+1} - w^k }^2 } + \frac{\tau\omega\alpha \eta t_k^2}{2} \Expsn{k}{ \norms{ w^{k+1} - \hat{w}^{k+1} }^2 } \vspace{1ex}\\ 
&& + {~} \frac{\revise{(s-2)}\eta L^2 \revise{t_k^2}}{\revise{32(s-1)}} \Expsn{k}{\norms{ \revise{x^{k+1} - y^k} }^2}   + \revise{\big[\varphi_s + \frac{25s-34}{2(s-2)}\big]} \epsilon_2 \eta t_{k-1}^2\Delta_{k-1} - \revise{\big[\varphi_s + \frac{25s-34}{2(s-2)}\big]} \frac{\eta\delta_k}{\bar{c}}.
\end{array}
\end{equation*}
This exactly proves \eqref{eq:SAEG_key_est5}\mred{, with $\tau$ chosen as follows}.

If we choose $\tau := \frac{4\rho_c}{(1-\alpha)\omega\eta} \in (0, 1)$, then both conditions in \eqref{eq:para_cond10_proof} hold, provided that
\begin{equation}\label{eq:SAEG_proof10} 
\arraycolsep=0.2em
\begin{array}{lcl}
\omega \geq \frac{4\rho_c}{(1-\alpha)\eta} + \revise{2\big[\frac{\varphi_s}{8} + \frac{s-2}{32(s-1)}\big]}.
\end{array}
\end{equation}
Since $\rho_c \leq \rho_n$ from Assumption~\ref{as:A3} and $\revise{(s-2)}\eta \geq \revise{8}(s-1)\rho_n$ from \eqref{eq:SFOG_para_cond2}, the condition \eqref{eq:SAEG_proof10} holds if we choose $\omega = \hat{\omega} := \revise{\frac{s-2}{2(1-\alpha)(s-1)} + \frac{\varphi_s}{4} + \frac{s-2}{16(s-1)}}$ as in \eqref{eq:omega_quantity}.

To guarantee the nonnegativity of the \mred{coefficients on the} right-hand side of \eqref{eq:SAEG_key_est5}, we need to impose $\frac{M^2\eta^2 t_k}{t_k - s} \leq 1$.
Since $t_k \geq s + 1$, the last condition holds if $(s+1)M^2\eta^2 \leq 1$.
However, since \mred{$\hat{\omega}$} is given in \eqref{eq:omega_quantity}, the condition $(s+1)M^2\eta^2 \leq 1$ is guaranteed if we choose $\eta \leq \frac{1}{L\sqrt{2(s+1)(1+\mred{\hat{\omega}})}}$.
Since $\eta \geq \frac{\revise{8}(s-1)\rho_n}{\revise{s-2}}$, we finally get $\frac{\revise{8}(s-1)\rho_n}{\revise{s-2}} \leq \eta \leq \frac{\hat{\lambda}}{L}$ as in \eqref{eq:SAEG_para_update_v1} with $\hat{\lambda} := \frac{1}{\sqrt{2(s+1)(1+\mred{\hat{\omega}})}}$.

Moreover, utilizing \eqref{eq:SFOG_para_cond1} and Lemma~\ref{le:VFOG_key_estimate2} with $c_1 = \frac{s-2}{\revise{16}(s-1)}$, \revise{$c_2 = \frac{s-2}{8(s-1)}$}, $t_k = k+s+1$, and $\gamma = \frac{\revise{(s-2)}\eta}{\revise{32(s-1)}}$ as before,  we also have $\gamma_k := \frac{\gamma (t_k -1)}{t_k} = \frac{\revise{(s-2)}\eta(k+s)}{\revise{32(s-1)}(k+s+1)}$ and $\beta_k = \big[\frac{\revise{3}(s-2) \eta}{\revise{16}(s-1)} + 2\rho_n]\frac{k+1}{k+s+1} - \frac{\gamma_k}{k+s+1}$ as in \eqref{eq:SAEG_para_update_v1}.
Since  $\gamma = \frac{\revise{(s-2)}\eta}{\revise{32(s-1)}}$, we also have 
\begin{equation*}
\arraycolsep=0.2em
\begin{array}{lcl}
\phi_k & = & \big[ \frac{\revise{5}(s-2)}{\revise{8}}\eta - 4(s-1)\rho_n - 3\gamma \big](k + s) + \frac{(s-1)(\revise{3s+10})\eta}{\revise{16}} + 2(2s^2 - 3s + 1)\rho_n + \gamma \vspace{1ex}\\
& =  & \big[ \revise{\frac{(20s-23)(s-2)}{32(s-1)}}\eta - 4(s-1)\rho_n \big](k + s) + (s-1)\big[ \frac{(\revise{3s+10})\eta}{\revise{16}} + 2(2s-1)\rho_n\big] + \frac{\revise{(s-2)}\eta}{\revise{32(s-1)}}
\end{array}
\end{equation*}
as desired.
\Eproof

\beforesubsec
\subsection{\textbf{The proof of Corollary~\ref{co:complexity_of_VFOG}: Oracle complexity of \ref{eq:VFOG} using \eqref{eq:mini_batch_Gy}.}}
\label{apdx:co:complexity_of_VFOG}
\aftersubsec
Assume that $\widetilde{G}y^k$ satisfies the condition \eqref{eq:error_cond3} with $\kappa = 1$, $\Theta = 0$, and $\delta_k := \frac{\delta}{t_k^{1+\nu}}$ for a given $\delta > 0$ and $\nu > 0$.
Then, we have $\Sc_{\infty} := \sum_{k=0}^{\infty}\delta_k = \delta\sum_{k=0}^{\infty}\frac{1}{(k+s+1)^{1+\nu}} < +\infty$.
Moreover, by \eqref{eq:error_cond3} we also have $\Expsn{k}{\norms{e^k}^2} \leq \frac{\delta_k}{t_k^2} = \frac{\delta}{t_k^{3+\nu}}$.
However, since $\widetilde{G}y^k$ is constructed by \eqref{eq:mini_batch_Gy}, we have $\Expsn{k}{\norms{e^k}^2} \leq \frac{\sigma^2}{b_k}$, where $\sigma^2$ is the variance bound in Assumption~\ref{as:A0}.
Therefore, the condition $\Expsn{k}{\norms{e^k}^2} \leq  \frac{\delta}{t_k^{3+\nu}}$ holds if we impose $\frac{\sigma^2}{b_k} \leq \frac{\delta}{t_k^{3+\nu}}$.
Hence, we can choose $b_k := \big\lceil \frac{\sigma^2t_k^{3+\nu}}{\delta} \big\rceil = \big\lceil \frac{\sigma^2(k+s + 1)^{3+\nu}}{\delta} \big\rceil$ as stated.
In this case, the expected total number of oracle calls $\mbf{G}$ is at most 
\begin{equation*}
\arraycolsep=0.2em
\begin{array}{lcl}
\Expn{\Tc_K} & := & \sum_{k=0}^Kb_k \leq \frac{\sigma^2}{\delta} \sum_{k=0}^K (k+s+1)^{3+\nu} = \BigOs{\frac{\sigma^2K^{4 + \nu}}{\delta}}.
\end{array}
\end{equation*}
By \eqref{eq:SFOG_BigO_rates} of Theorem~\ref{th:VFOG1_convergence}, the number of iterations $K$ to reach an $\epsilon$-solution in expectation is $K = \BigO{\epsilon^{-1}}$.
Therefore, we obtain $\Expn{\Tc_K} = \BigOs{\frac{\sigma^2}{\delta\epsilon^{4 + \nu}}}$ as stated.
The expected total number of $J_{\eta T}$ evaluations is the same as $K$, which is $\BigO{\epsilon^{-1}}$.

Alternatively, if we choose $b_k := \big\lceil \frac{\sigma^2(k+s+1)^3\log(k+s+1)}{\delta} \big\rceil$, then $\Sc_K := \delta\sum_{k=0}^K \frac{1}{(k+s+1)\log(k+s+1)} \leq 1 + \log(\log(K+s+1))$.
Hence, from \eqref{eq:SFOG_BigO_rates} of Theorem~\ref{th:VFOG1_convergence}, to achieve $\Expn{\norms{Gx^k + v^k}^2} \leq \epsilon^2$, we require $\frac{C_0(\mcal{R}_0^2 + \Lambda_0 S_K)}{(K+s+1)^2} \leq \epsilon^2$.
Using the upper bound of $\Sc_K$, we can show that $K = \BigOs{\frac{1}{\epsilon\sqrt{\log(\epsilon^{-1})}}}$.
Using this estimation of $K$, we can estimate that $\Expn{\Tc_K} = \sum_{k=0}^Kb_k \leq \frac{\sigma^2}{\delta} \sum_{k=0}^K(k+s+1)^3\log(k+s+1) \leq \frac{\sigma^2(K+s+1)^4\log(K+s+1)}{4\delta} \sim \BigOs{\frac{1}{\epsilon^4 \log(\epsilon^{-1})}}$.
The expected total number of $J_{\eta T}$ evaluations is the same as $K$, which is $\widetilde{\mcal{O}}(\epsilon^{-1})$.
\Eproof

\revise{
\beforesubsec
\subsection{The proof of Proposition~\ref{pro:properties_of_G_and_T}.}\label{apdx:prop:h_function_smooth}
\aftersubsec
	We divide the proof of Proposition~\ref{pro:properties_of_G_and_T} into several steps.
	
	\vspace{0.5ex}
	\noindent\textit{\textbf{Step 1.}}
	First, we show that $h$ defined in \eqref{eq:portfolio_h_function} is $\gamma$-strongly convex and $L_h$-smooth with $L_h := \gamma + \nu \norms{C}^2$.
	Indeed, $\nabla h(y) = \gamma y + \nu C^\top \tanh(Cy)$ and $\nabla^2 h(y) = \gamma\Id + \nu C^\top D(y) C$, where $D(y) := \text{diag}\big( \operatorname{sech}^2(c_1^\top y), \cdots, \operatorname{sech}^2(c_r^\top y) \big)$.
	Since $0 \preceq D(y) \preceq \Id$, we have
	\begin{equation*}
		\gamma \Id \preceq \nabla^2 h(y) \preceq \gamma \Id + \nu C^\top C \preceq (\gamma + \nu \norms{C}^2)\Id.
	\end{equation*}
	Hence, $h$ is $\gamma$-strongly convex and $L_h$-smooth with $L_h := \gamma + \nu \norms{C}^2$.

	\vspace{0.5ex}
	\noindent\textit{\textbf{Step 2.}}
	Next, we show that $\nabla{\tilde{\Rc}}_{\lambda, \theta}(z)$ from \eqref{eq:portfolio_SCAD_complementary} is $\frac{1}{\theta-1}$-Lipschitz continuous.
	The derivative of $\tilde{R}_{\lambda, \theta}$ is
	$\tilde{R}'_{\lambda, \theta}(t) = -\lambda \boldsymbol{1}_{\sets{t < -\theta \lambda}} + \frac{t+\lambda}{\theta - 1} \boldsymbol{1}_{\sets{-\theta\lambda 
	\leq t < -\lambda}} + \frac{t - \lambda}{\theta - 1}\boldsymbol{1}_{\sets{\lambda < t \leq \theta \lambda}} + \lambda\boldsymbol{1}_{\sets{t > \theta \lambda}}$.
	It is continuous, monotone nondecreasing, and piece-wise affine, with slope either $0$ or $\frac{1}{\theta - 1}$.
	Therefore, $\tilde{R}_{\lambda, \theta}$ is convex and continuously differentiable, and $\nabla \tilde{\Rc}_{\lambda, \theta}(z)$ is $\frac{1}{\theta - 1}$-Lipschitz continuous.

	\vspace{0.5ex}
	\noindent\textit{\textbf{Step 3.}} 
	We now prove (i).
	Let $x = (z,y)$ and $x' = (z', y')$ be arbitrary.
	By Jensen's inequality and the $L_f$-average smoothness of $f$, $f$ is $L_f$-smooth.
	Thus, by \textit{\textbf{Steps 1--2}}, one can show that
	\begin{equation*}
	\begin{array}{llcl}
		& \norms{\nabla f(z) - \nabla f(z') + B(y - y') - \tau(\nabla \tilde{\Rc}_{\lambda, \theta}(z) - \nabla \tilde{\Rc}_{\lambda, \theta}(z'))}  & \leq &  \big(L_f + \frac{\tau}{\theta - 1}\big)\norms{z - z'} + \norms{B}\norms{y - y'}, \vspace{1ex}\\
		& \norms{\nabla h(y) - \nabla h(y') - B^\top(z - z')}  & \leq & \norms{B}\norms{z - z'} + L_h \norms{y - y'}.
	\end{array}
	\end{equation*}
	Therefore, with $\tilde{L}_f := L_f + \frac{\tau}{\theta - 1}$, we get
	\begin{equation*}
	\begin{array}{lcl}
		\begin{bmatrix} \norms{[Gx - Gx']_z} \\ \norms{[Gx - Gx']_y} \end{bmatrix} \leq  \begin{bmatrix} \tilde{L}_f & \norms{B} \\ \norms{B} & L_h \end{bmatrix}  \begin{bmatrix} \norms{z - z'} \\ \norms{y - y'} \end{bmatrix}.
	\end{array}
	\end{equation*}
	Hence,
	\begin{equation*}
	\begin{array}{lcl}
	\norms{Gx - Gx'} \leq \norm{\begin{bmatrix} \tilde{L}_f & \norms{B} \\ \norms{B} & L_h \end{bmatrix} }_2 \norms{x - x'}.
	\end{array}
	\end{equation*}
	The spectral norm of this symmetric matrix is its largest eigenvalue $L := \frac{\tilde{L}_f + L_h + \sqrt{(\tilde{L}_f - L_h)^2 + 4\norms{B}^2}}{2}$.

	\vspace{0.5ex}
	\noindent\textit{\textbf{Step 4.}} Next, we prove (ii).
	Let $x = (z,y)$ and $x' = (z', y')$, and take $w \in \Phi x$ and $w' \in \Phi x'$.
	By the definition of $T$, there exist $\xi \in \tau \lambda \partial \norms{z}_1 + \Nc_{\Zc}(z)$ and $\xi' \in \tau \lambda \partial \norms{z'}_1 + \Nc_{\Zc}(z')$ such that
	\begin{equation*}
		w =	\begin{bmatrix} \nabla f(z)	+ By - \tau\nabla \tilde{\Rc}_{\lambda, \theta}(z) + \xi \\ \nabla h(y) - B^\top z \end{bmatrix} \quad \text{and} \quad  
		w' = \begin{bmatrix} \nabla f(z') + By' - \tau\nabla \tilde{\Rc}_{\lambda, \theta}(z') + \xi' \\ \nabla h(y') - B^\top z' \end{bmatrix}.
	\end{equation*}
	Thus,
	\begin{equation}\label{eq:portfolio_cohypomonotone_proof0}
		w - w' = \begin{bmatrix} \big(\nabla f(z) - \nabla f(z')\big)	+ B(y - y') - \tau\big(\nabla \tilde{\Rc}_{\lambda, \theta}(z) - \nabla \tilde{\Rc}_{\lambda, \theta}(z')\big) + (\xi - \xi') \\ \big(\nabla h(y)-\nabla h(y')) - B^\top (z - z') 	\end{bmatrix}.
	\end{equation}
	Taking its inner product with $x - x'$ and using the convexity of $f$, the monotonicity of $\tau \lambda \partial \norms{\cdot}_1 + \partial \delta_{\Zc}$, and \textit{\textbf{Steps 1--2}}, we obtain
	\begin{equation}\label{eq:portfolio_cohypomonotone_proof1}
	\begin{array}{lcl}
		\iprods{w - w', x - x'} 
		&\geq& -\frac{\tau}{\theta - 1}\norms{z - z'}^2 + \gamma \norms{y - y'}^2.
	\end{array}
	\end{equation}
	We consider two cases.
	\begin{compactitem}
		\item\textit{Case 1.}
		If $\iprods{w - w', x - x'} \geq 0$, then the desired co-hypomonotonicity inequality holds trivially.
		
		\item\textit{Case 2.}
		If $\iprods{w - w', x - x'} < 0$, then \eqref{eq:portfolio_cohypomonotone_proof1} implies $\norms{y - y'}^2 \leq \frac{\tau}{\gamma(\theta - 1)} \norms{z - z'}^2$.
		Hence, using $\sigma := \sigma_{\min}(B^\top)$ and the $L_h$-smoothness of $h$,
		\begin{equation*}
			\begin{array}{lcl}
				\norms{w - w'} &\geq& \norms{\big(\nabla h(y)-\nabla h(y')) - B^\top (z - z')} \vspace{1ex}\\
				&\geq& \sigma\norms{z - z'} - L_h\norms{y - y'} \vspace{1ex}\\
				&\geq& \big(\sigma - L_h\big(\frac{\tau}{\gamma(\theta - 1)}\big)^{1/2}\big)\norms{z - z'} > 0.
			\end{array}
		\end{equation*}
		Combining this estimate with \eqref{eq:portfolio_cohypomonotone_proof1} gives
		\begin{equation*}
			\begin{array}{lcl}
				\iprods{w - w', x - x'} &\geq& -\frac{\tau}{\theta - 1}\norms{z - z'}^2 \geq -\frac{\tau}{(\theta - 1)\left(\sigma - L_h\sqrt{\tau/((\theta - 1)\gamma)}\right)^2}\norms{w - w'}^2.
			\end{array}
		\end{equation*}
		Thus, the desired inequality holds with $\rho := \frac{\tau}{(\theta - 1)\left(\sigma - L_h\sqrt{\tau/((\theta - 1)\gamma)}\right)^2}$.
	\end{compactitem}	
	Therefore, $\Phi := G+T$ is globally $\rho$-co-hypomonotone.

	\vspace{0.5ex}
	\noindent\textit{\textbf{Step 5.}} Finally, we prove (iii).
	By Minkowski's inequality, the $\frac{1}{\theta-1}$-Lipschitz continuity of $\nabla\tilde{\Rc}_{\lambda,\theta}$, and the average smoothness of $f$, for $\tilde{L}_f := L_f + \frac{\tau}{\theta - 1}$, we can show that
	\begin{equation*}
	\begin{array}{lcl}
		\left[\frac{1}{n}\sum_{i=1}^n \norms{\nabla f_i(z) - \nabla f_i(z') + B(y-y') - \tau(\nabla\tilde{\Rc}_{\lambda,\theta}(z) - \nabla\tilde{\Rc}_{\lambda,\theta}(z'))}^2\right]^{1/2} {\!\!\!}  \leq \tilde{L}_f \norms{z - z'} + \norms{B}\norms{y - y'}.
	\end{array}
	\end{equation*}
	Arguing as in \textbf{\textit{Step 3}} gives
	\begin{equation}\label{eq:portfolio_G_average_Lipschitz}
	\begin{array}{lcl}
		\frac{1}{n}\sum_{i=1}^n \norms{G_ix - G_ix'}^2 \leq L^2 \norms{x - x'}^2,
	\end{array}
	\end{equation}
	where $L = \frac{\tilde{L}_f + L_h + \sqrt{(\tilde{L}_f - L_h)^2 + 4\norms{B}^2}}{2}$.
	
	We next bound $\norms{x - x'}$ by $\norms{w - w'}$.
	It is easy to show that  
	$\nabla h(y) - \nabla h(y') = Q(y - y')$,
	where $Q := \int_0^1 \nabla^2 h(y' + t(y-y'))\mathrm{d}t$ satisfies $\gamma\Id \preceq Q \preceq L_h\Id$.
	Thus, from the second block of \eqref{eq:portfolio_cohypomonotone_proof0}, we get
	\begin{equation}\label{eq:Q_inv_bound}
	\begin{array}{lcl}
		y - y' = Q^{-1} \big([w - w']_y + B^\top (z - z')\big), \qquad \frac{1}{L_h}\Id \preceq Q^{-1} \preceq \frac{1}{\gamma}\Id.
	\end{array}
	\end{equation}
	Substituting this relation into the first block of \eqref{eq:portfolio_cohypomonotone_proof0}, we can derive that  
	\begin{equation*}
	\begin{array}{lcl}
		\iprods{[w - w']_z - B Q^{-1} [w - w']_y, z - z'} 
		&\geq& \frac{1}{L_h}\norms{B^\top(z - z')}^2 - \frac{\tau}{\theta - 1}\norms{z-z'}^2  \geq  \big(\frac{\sigma^2}{L_h} - \frac{\tau}{\theta - 1}\big)\norms{z-z'}^2.
	\end{array}
	\end{equation*}
	Since $\sigma^2(\theta - 1) > \tau L_h$, the Cauchy-Schwarz inequality and \eqref{eq:Q_inv_bound} yield
	\begin{equation*}
	\begin{array}{lcl}
		\norms{z-z'}  &\leq& \frac{L_h(\theta - 1)\sqrt{\gamma^2 + \norms{B}^2}}{\gamma\left[\sigma^2(\theta - 1) - \tau L_h\right]} \norms{w - w'}.
	\end{array}
	\end{equation*}
	Using \eqref{eq:Q_inv_bound} again, we also have
	\begin{equation*}
	\begin{array}{lcl}
		\norms{y - y'}   & \leq & \frac{\gamma\left[\sigma^2(\theta - 1) - \tau L_h\right] + L_h(\theta - 1)\norms{B}\sqrt{\gamma^2 + \norms{B}^2}}{\gamma^2\left[\sigma^2(\theta - 1) - \tau L_h\right]} \norms{w - w'}.
	\end{array}
	\end{equation*}
	Therefore, combining the last two estimates, we finally arrive at
	\begin{equation*}
	\begin{array}{lcl}
		\norms{x - x'}^2 &\leq& \frac{L_h^2(\theta - 1)^2(\gamma^2 + \norms{B}^2) + \left[\gamma\left(\sigma^2(\theta - 1) - \tau L_h\right) + L_h(\theta - 1)\norms{B}\sqrt{\gamma^2 + \norms{B}^2}\right]^2}{\gamma^2\left[\sigma^2(\theta - 1) - \tau L_h\right]^2} \norms{w - w'}^2. 
	\end{array}
	\end{equation*}
	Combining this estimate with \eqref{eq:portfolio_G_average_Lipschitz}, we obtain
	\begin{equation*}
	\begin{array}{lcl}
		\frac{1}{n}\sum_{i=1}^n \norms{G_ix - G_ix'}^2 \leq \tilde{L}^2\norms{w - w'}^2 \leq (1 + \tilde{L}^2)\norms{w - w'}^2,
	\end{array}
	\end{equation*}
	where $\hat{L}^2 := \frac{L^2\left\{\gamma^2 L_h^2(\theta - 1)^2\left(\gamma^2 + \norms{B}^2\right) + \left[\gamma\left(\sigma^2(\theta - 1) - \tau L_h\right) + L_h(\theta - 1)\norms{B}\sqrt{\gamma^2 + \norms{B}^2}\right]^2\right\}}{\gamma^4\left[\sigma^2(\theta - 1) - \tau L_h\right]^2}$.

	Finally, by (ii), for any $\kappa \geq 0$, it is obvious to show that
	\begin{equation*}
	\begin{array}{lcl}
		\iprods{w - w', x - x'}  &\geq& -(1+\kappa)\rho \norms{w-w'}^2 + \frac{\kappa\rho}{1+ \hat{L}^2}\cdot \frac{1}{n}\sum_{i=1}^n \norms{G_ix - G_ix'}^2.
	\end{array}
	\end{equation*}
	Let us denote $\rho_n := (1 + \kappa)\rho$ and $\rho_c := \frac{\kappa\rho}{1+ \hat{L}^2}$.
	Since $\rho_n - \rho_c = \rho\big(1 + \frac{\kappa\hat{L}^2}{1 + \hat{L}^2}\big) > 0$, $\Phi$ satisfies Assumption~\ref{as:A3} with $\rho_n > \rho_c \geq 0$.
\Eproof
}
\end{APPENDICES}

\bibliographystyle{informs2014} 

\end{document}